\documentclass[9pt]{article}
\usepackage{amssymb}
\usepackage{tikz}
\usepackage{xr}
\usetikzlibrary{arrows, calc, decorations.pathmorphing, backgrounds, positioning, fit, petri, automata}
\definecolor{grey1}{rgb}{0.5,0.5,0.5}
\usepackage[top=2.54cm, bottom=2.54cm, left=2.2cm, right=2.2cm]{geometry}
\usepackage{algorithmicx}
\usepackage{algcompatible}
\usepackage{algorithm}
\usepackage[shortlabels]{enumitem}
\usepackage{changebar}
\usepackage[authoryear]{natbib}

\usepackage{hyperref}
\hypersetup{ 
  colorlinks,
  citecolor=blue,
  linkcolor=red,
  urlcolor=magenta}

\usepackage{mathtools}
\usepackage{graphicx}
\usepackage{subfigure}
\usepackage{amsmath,amsthm,bm,bbm}
\usepackage{rotating}
\usepackage{mathrsfs}
\usepackage{chngcntr}
\usepackage{apptools}
\usepackage[titletoc,title]{appendix}
\usepackage{comment}

\usepackage{xcolor}         %for colors
\definecolor{grau}{rgb}{0.8,0.8,0.8}

\newcommand{\chen}[1]{\color{orange}}
\usepackage{setspace,lscape,longtable}
\usepackage{amsmath,amsthm,bm}
\usepackage{fullpage}  

\usepackage{titlesec}
\usepackage[english]{babel}
\usepackage{enumitem}
\setlist[enumerate]{itemsep=0mm}
\setlist[itemize]{itemsep=0mm}
\numberwithin{equation}{section}
\newtheorem{theorem}{Theorem}[section]

\newtheorem{lemma}[theorem]{Lemma}
\newtheorem{proposition}{Proposition}[section]
\theoremstyle{remark}

\newtheorem{assumption}{Assumption}
\newtheorem{example}{Example}

\newtheorem{result}{Result}[section]
\newtheorem*{excont}{Example \continuation}
\newcommand{\continuation}{??}
\newenvironment{continueexample}[1]
 {\renewcommand{\continuation}{\ref{#1}}\excont[continued]}
 {\endexcont}
\newtheorem{remark}{Remark}
\DeclareMathOperator*{\argmin}{arg\,min}

\DeclareMathOperator*{\arginf}{arg\,inf}

\newcommand{\prob}{{\mathbb{P}}}
\newcommand{\var}{{\mathrm{var}}}
\newcommand{\expect}{\mathbb{E}}

\newcommand{\transpose}{^{\mathrm{T}}}

\newcommand{\blambda}{{\bm{\lambda}}}

\newcommand{\calA}{{\mathcal{A}}}

\newcommand{\calD}{{\mathcal{D}}}
\newcommand{\calE}{{\mathcal{E}}}
\newcommand{\calF}{{\mathcal{F}}}
\newcommand{\calG}{{\mathcal{G}}}

\newcommand{\calL}{{\mathcal{L}}}

\newcommand{\calX}{{\mathcal{X}}}

\newcommand{\balpha}{{\boldsymbol{\alpha}}}

\newcommand{\bU}{{\mathbf{U}}}
\newcommand{\bP}{{\mathbf{P}}}
\newcommand{\bp}{{\mathbf{p}}}

\newcommand{\bv}{{\mathbf{v}}}
\newcommand{\bx}{{\mathbf{x}}}
\newcommand{\bX}{{\mathbf{X}}}
\newcommand{\bA}{{\mathbf{A}}}
\newcommand{\bB}{{\mathbf{B}}}

\newcommand{\bE}{{\mathbf{E}}}

\newcommand{\bG}{{\mathbf{G}}}
\newcommand{\bH}{{\mathbf{H}}}
\newcommand{\bM}{{\mathbf{M}}}

\newcommand{\bR}{{\mathbf{R}}}
\newcommand{\bQ}{{\mathbf{Q}}}
\newcommand{\bS}{{\mathbf{S}}}
\newcommand{\bW}{{\mathbf{W}}}
\newcommand{\bZ}{{\mathbf{Z}}}

\newcommand{\bg}{{\mathbf{g}}}

\newcommand{\br}{{\mathbf{r}}}

\newcommand{\bt}{{\mathbf{t}}}
\newcommand{\bu}{{\mathbf{u}}}

\newcommand{\by}{{\mathbf{y}}}
\newcommand{\bz}{{\mathbf{z}}}

\newcommand{\be}{{\mathbf{e}}}

\newcommand{\bV}{{\mathbf{V}}}

\newcommand{\bDelta}{{\bm{\Delta}}}
\newcommand{\bOmega}{{\bm{\Omega}}}

\newcommand{\bSigma}{{\bm{\Sigma}}}

\newcommand{\eye}{{\mathbf{I}}}
\newcommand{\one}{{\mathbf{1}}}

\newcommand{\bpsi}{{\bm{\psi}}}

\newcommand{\zero}{{\bm{0}}}
\newcommand{\eps}{\epsilon}

\author{Fangzheng Xie\footnotemark[1]\thanks{Department of Statistics, Indiana University} \and Dingbo Wu 
\footnotemark[1] \thanks{Correspondence should be addressed to Fangzheng Xie (fxie@iu.edu)}
}
\title{\bf Eigenvector-Assisted Statistical Inference for Signal-Plus-Noise Matrix Models}

\begin{document}

\allowdisplaybreaks

\maketitle

\begin{abstract}
In this paper, we develop a generalized Bayesian inference framework for a collection of signal-plus-noise matrix models arising in high-dimensional statistics and many applications. The framework is built upon an asymptotically unbiased estimating equation with the assistance of the leading eigenvectors of the data matrix. The solution to the estimating equation coincides with the maximizer of an appropriate statistical criterion function. The generalized posterior distribution is constructed by replacing the usual log-likelihood function in the Bayes formula with the criterion function. The proposed framework does not require the complete specification of the sampling distribution and is convenient for uncertainty quantification via a Markov Chain Monte Carlo sampler, circumventing the inconvenience of resampling the data matrix. Under mild regularity conditions, we establish the large sample properties of the estimating equation estimator and the generalized posterior distributions. In particular, the generalized posterior credible sets have the correct frequentist nominal coverage probability provided that the so-called generalized information equality holds. The validity and usefulness of the proposed framework are demonstrated through the analysis of synthetic datasets and the real-world ENZYMES network datasets. 
\end{abstract}

\noindent%
{\it Keywords:} Bernstein-von Mises theorem, eigenvector-assisted estimating equation, generalized Bayesian inference, Markov chain Monte Carlo, uncertainty analysis
% \vfill

\tableofcontents

\section{Introduction}
\label{sec:intro}

\subsection{Background}

In the era of data science, the emergence and the analysis of high-dimensional complex datasets have been a gigantic and rapidly developing field in recent decades. Low-rank matrix models, also known as \emph{signal-plus-noise matrix models}, have been broadly applied in numerous practical applications. Examples of such application domains include social network analysis \citep{HOLLAND1983109,young2007random}, signal processing and compressed sensing \citep{1614066,eldar2012compressed}, collaborative filtering and recommendation system \citep{bennett2007netflix,goldberg1992using}, neural science \citep{eichler2017complete}, camera sensor networks \citep{tron2009distributed}, and synchronization of wireless networks \citep{4177758}.
% The application and the analysis of 
% Low-rank matrix models are also of great interest to statisticians, machine learning scientists, applied mathematicians, computer scientists, bioinformaticians, or even physicists.

% Within the field of statistics and probability, 
Spectral methods are of fundamental interest in analyzing a broad range of signal-plus-noise matrix models. For example, the leading eigenvectors of the adjacency matrix of a stochastic block model encode the community structure of the vertices directly, leading to the renowned spectral clustering algorithm \citep{10.1214/19-AOS1854,lyzinski2014,rohe2011,sussman2012consistent}. Furthermore, spectral estimators can typically be directly applied to the subsequent inference tasks  \citep{ng2002spectral,868688,6565321,tang2013,doi:10.1080/10618600.2016.1193505,tang2017} or serve as ``warm-starts'' that initialize various optimization-based learning algorithms \citep{7029630,jain2013low,5466511}. On the theoretical side, the performance of spectral-based methods is backboned by the underlying matrix perturbation analysis \citep{10.1214/19-AOS1854,10.1214/17-AOS1541,cape2019signal,cape2017two,doi:10.1137/0707001,pmlr-v83-eldridge18a,fan2018eigenvector,doi:10.1080/01621459.2020.1751645,wedin1972perturbation,xie2021entrywise} and random matrix theory 
\citep{bai2010spectral,BENAYCHGEORGES2011494,PAUL20141,yao2015sample}. On the practical side, the computational cost of spectral estimators is typically low, which further popularizes them and their refinements in various contexts. 

% From the theoretical perspective, the success of spectral-based methods are backboned by the underlying matrix perturbation analysis \citep{doi:10.1137/0707001,wedin1972perturbation,10.1214/17-AOS1541} and random matrix theory 
% \citep{BENAYCHGEORGES2011494,PAUL20141,bai2010spectral,yao2015sample}. Motivated by the exact recovery in community detection for stochastic block models \citep{zhang2016minimax,7298436}, there has been a recent growing interest in investigating the entrywise behavior of the eigenvectors associated with signal-plus-noise matrix models \citep{cape2017two,cape2019signal,pmlr-v83-eldridge18a,10.1214/19-AOS1854,doi:10.1080/01621459.2020.1751645,xie2021entrywise,fan2018eigenvector}.

\subsection{Overview}

This paper proposes a general statistical inference framework for signal-plus-noise matrix models
 % by borrowing the idea in moment condition models \citep{10.2307/1912775,doi:10.1080/07350015.1996.10524656,10.1093/biomet/92.1.31,10.1214/009053606000001208}. 
% We derive an 
% quite flexible 
% The framework is 
based on a novel
eigenvector-assisted estimating equation.
 % that allows for heteroskedasticity through a user-defined weight function. 
The solution to the estimating equation can be alternatively viewed as the extremum of a general statistical criterion function. 
Examples of such a criterion function include the $M$-estimation objective function, the generalized method of moments objective function, and the exponentially tilted empirical likelihood. 
We propose to use the generalized posterior distribution to estimate the signal matrix, where the usual log-likelihood function in the Bayes formula is substituted by the aforementioned statistical criterion function of interest. Under mild regularity conditions, we establish the asymptotic normality of the eigenvector-assisted $Z$-estimator and the Bernstein-von Mises theorem of the generalized posterior distribution. 

Our proposed methodology enjoys several fascinating features:
\begin{enumerate}[(a), noitemsep, topsep = 0mm]
  \item The framework is likelihood-free and allows for various noise distributions.
  % robust to the potential misspecification of the sampling distribution. 

  \item The generalized posterior distribution can be computed via a standard Metropolis-Hastings algorithm, circumventing the inconvenience of nonconvex optimization problems. Furthermore, the Metropolis-Hastings algorithm can be implemented in parallel thanks to the separable structure of the criterion function (see Section \ref{sub:generalized_bayesian_estimation_with_moment_conditions} for details). 

  \item The generalized Bayesian method provides a convenient environment for uncertainty quantification through the Metropolis-Hastings algorithm. This advantage is in contrast to the frequentist approach for assessing the uncertainty via bootstrap because the resampling of signal-plus-noise matrices is not straightforward \citep{levin2019bootstrapping,10.1093/biomet/asaa006}. 

  \item The row-wise credible sets of the generalized posterior are well-calibrated. Namely, they have the correct frequentist coverage probability asymptotically, provided that the so-called generalized information equality holds (see Section \ref{sub:consequence_of_generalized_posterior} for details). 

  \item When the variance information of the noise is available, the practitioner can select the user-defined weight function in the estimating equation appropriately (see Section \ref{sub:moment_condition_models} for details), such that the resulting estimator has the minimum asymptotic covariance matrix in spectra among all eigenvector-assisted $Z$-estimators. 
\end{enumerate}

\subsection{Related work}

There are several recent papers addressing the theoretical properties of the eigenvectors of general signal-plus-noise matrix models. \cite{cape2019signal}
 % first used the term ``signal-plus-noise matrix model'' and 
explored the entrywise error bound and central limit theorem for the eigenvectors of signal-plus-noise matrices. \cite{10.1214/19-AOS1854} obtained sharper entrywise concentration bounds for the eigenvectors of symmetric random matrices with low expected rank. The asymptotic theory of the eigenvalues and linear functionals of the eigenvectors for the general random matrices with diverging leading eigenvalues was established by \cite{doi:10.1080/01621459.2020.1840990}. In the context of random graph inference, \cite{athreya2016limit}, \cite{tang2018}, and \cite{xie2021entrywise} studied the central limit theorems for the rows of the eigenvector matrix.
 % of the adjacency matrix 
\cite{xie2019efficient} and \cite{xie2021entrywise} proposed a one-step refinement for the eigenvectors and explored the corresponding entrywise limit theorem. \cite{agterberg2021entrywise} further extended the signal-plus-noise matrix framework to general rectangular matrices and allowed heteroskedasticity and dependence of the noise distributions. The asymptotic results obtained in the above work are with regard to frequentist estimators. While the uncertainty of a frequentist estimator can be assessed using bootstrap, the resampling of a signal-plus-noise matrix model is less straightforward than that of classical parametric models. This paper distinguishes itself from the aforementioned work as it provides a user-friendly environment for uncertainty quantification through the generalized Bayesian inference method. 

The idea of the generalized posterior distribution, which is obtained by replacing the usual log-likelihood function with a general statistical criterion function in the Bayes formula, is not entirely new in the literature. The convenience of the generalized posterior is that it does not require the full specification of the sampling distribution of the data. An early influential work is \cite{CHERNOZHUKOV2003293}, which established a systematic framework for studying the convergence of the generalized posteriors for a broad range of semiparametric econometrics models. There has also been some recent development on the Bernstein-von Mises theorem of the generalized posterior distributions \citep{kleijn2012bernstein,JMLR:v22:20-469,10.1093/biomet/asy054,syring2020gibbs}. These approaches, however, are not directly applicable to the signal-plus-noise matrix models. One contribution of the present paper is that we design appropriate statistical criterion functions for the signal-plus-noise matrix models by borrowing the idea of moment condition models with the assistance of the sample leading eigenvectors. In addition, the generalized posterior credible sets may not have the frequentist nominal coverage probability \citep{kleijn2012bernstein} and may require calibration \citep{10.1093/biomet/asy054} in general. In contrast, in our framework, the appropriate choice of the criterion function (e.g., the generalized method of moments criterion or the exponentially tilted empirical likelihood criterion) can provide the generalized posterior credible sets with the correct coverage probability.

Another line of the related literature is on the development of the moment condition models using the generalized method of moments \citep{10.2307/1912775}, the empirical likelihood \citep{owen1988empirical,10.1214/aos/1176347494}, the generalized empirical likelihood \citep{10.2307/2971718,10.2307/2171942,https://doi.org/10.1111/j.1468-0262.2004.00482.x}, and the exponentially tilted empirical likelihood \citep{doi:10.1080/01621459.2017.1358172,10.1093/biomet/92.1.31,10.1214/009053606000001208}. These papers tackle the higher-order properties of various point estimators for the low-dimensional parameters in general semiparametric moment condition
models that are popular in econometrics but do not apply directly to the high-dimensional signal-plus-noise matrix models. Our work fills this gap by developing a novel eigenvector-assisted estimation framework and the corresponding large sample properties. 
% by borrowing the idea of moment condition models. 
% , whereas our work puts the high-dimensional signal-plus-noise matrix models in the context of moment condition models. 

\subsection{Organization}

The rest of the paper is structured as follows. Section \ref{sec:signal_plus_noise_matrix_models} introduces the signal-plus-noise matrix model and presents several examples. Section \ref{sec:eigenvector_assisted_estimation} elaborates on the proposed eigenvector-assisted estimation framework. The main theoretical results of the proposed estimation procedure are established in Section \ref{sec:main_results}, including the large sample properties of the eigenvector-assisted $Z$-estimator and the generalized posterior distribution. Numerical examples are demonstrated in Section \ref{sec:numerical_examples}. We conclude the paper with a discussion in Section \ref{sec:discussion}. 

\vspace*{1ex}
\noindent\textbf{Notations:}
Given $n\in\mathbb{N}_+$, let $[n] = \{1,2,\ldots,n\}$. 
% For a vector-valued function $\bpsi$ taking values in a subset of $\mathbb{R}^d$, denote $[\bpsi]_k$ as the $k$th coordinate function of $\bpsi$. 
% For a vector valued differentiable function $\mathbf{f}(\bx):\mathbb{R}^d\to\mathbb{R}^d$, denote $D_\bx \mathbf{f}(\bx) = \partial\bf(\bx)/\partial\bx\transpose$ the Jacobian matrix of $\bf$ (with regard to $\bx$). 
For a scalar-valued $r$-times differentiable function $f(\bx):\mathbb{R}^d\to\mathbb{R}$ and a vector $\balpha = (\alpha_1,\ldots,\alpha_d)\in\mathbb{N}^d$ with $|\balpha|:=\sum_{k = 1}^d\alpha_k \leq r$, we use the notation $D^\balpha f(x_1,\ldots,x_d) = \partial^{|\balpha|}f(\bx)/\partial x_1^{\alpha_1}\ldots\partial x_d^{\alpha_d}$ to denote the corresponding $k$th-order mixed partial derivative associated with $\balpha$. 
For two non-negative sequences $(a_n)_{n = 1}^\infty,(b_n)_{n = 1}^\infty$, we write $a_n\lesssim b_n$,
% ($a_n\gtrsim b_n$, resp.), 
if $a_n\leq Cb_n$
% ($a_n\geq Cb_n$, resp.) 
for some constant $C > 0$. 
We use notations $C, c, C_1, C_2, \ldots$ to denote generic constants that may change from line to line but are independent of the asymptotic index $n$. 
With a slight abuse of notation, we say that a sequence of random variables $(X_n)_{n = 1}^\infty$ is upper bounded by a constant multiple of $\eps_n$ for a sequence $(\eps_n)_{n = 1}^\infty\subset\mathbb{R}$ with high probability, denoted by $X_n\lesssim \eps_n$ w.h.p. or $X_n = O(\eps_n)$ w.h.p., if for any $c > 0$, there exist constants $K_c > 0$ and $N_c\in\mathbb{N}_+$, such that $P(X_n\leq K_c\eps_n)\geq 1 - n^{-c}$ for all $n\geq N_c$. Similarly, a sequence of events $(\calE_n)_{n = 1}^\infty$ is said to occur with high probability (w.h.p.), if for all $c > 0$, there exists a constant $N_c\in\mathbb{N}_+$ depending on $c$, such that $\prob_0(\calE_n)\geq 1 - n^{-c}$ for all $n\geq N_c$. A sequence of events $(\calE_n)_{n = 1}^\infty$ is said to occur with probability approaching to one (w.p.a.1), if $\prob(\calE_n) \to 1$ as $n\to\infty$.
For $n,d\in\mathbb{N}_+$ with $n\geq d$, we denote $\eye_d$ the $d\times d$ identity matrix and $\mathbb{O}(n, d) = \{\bU\in\mathbb{R}^{n\times d}:\bU\transpose\bU = \eye_d\}$ the set of all orthonormal $d$-frames in $\mathbb{R}^n$, and we write $\mathbb{O}(d)$ when $n = d$. For a $n\times n$ symmetric matrix $\bA$, we denote $\lambda_k(\bA)$ its $k$th largest eigenvalue in magnitude, namely, $|\lambda_1(\bA)|\geq\ldots\geq|\lambda_n(\bA)|$. For a general rectangular $n\times d$ matrix $\bX$, we denote $\sigma_k(\bX)$ its $k$th largest singular value, such that $\sigma_1(\bX)\geq\ldots\geq\sigma_{\min(n, d)}(\bX)\geq 0$. For two positive semidefinite matrices $\bA$ and $\bB$, we denote $\bA\succeq\bB$ ($\bA\preceq \bB$, resp.) if $\bA - \bB$ is positive semidefinite (negative semidefinite, resp.). For a matrix $\bA = [A_{ij}]_{m\times n}$, we use $\|\bA\|_2$, $\|\bA\|_{\mathrm{F}}$, $\|\bA\|_{2\to\infty}$, and $\|\bA\|_\infty$ to denote the spectral norm, the Frobenius norm, the two-to-infinity norm defined by $\|\bA\|_{2\to\infty} = \max_{i\in [m]}(\sum_{j = 1}^nA_{ij}^2)^{1/2}$, and the matrix infinity norm defined by $\|\bA\|_\infty = \max_{i\in [m]}\sum_{j = 1}^n|A_{ij}|$, respectively. These norm notations also apply to (column) vectors in $\mathbb{R}^d$ for any $d\in\mathbb{N}_+$. For a (sub-Gaussian) random variable $A$, define the $\psi_2$-Orlicz norm of $A$ by
  $\|A\|_{\psi_2} = \sup_{p\geq 1}p^{-1/2}(\expect_0|A|^p)^{1/p}$ 
    (See, for example, \citealp{kosorok2008introduction} and \citealp{vershynin2010introduction}). 

\section{Signal-Plus-Noise Matrix Models} % (fold)
\label{sec:signal_plus_noise_matrix_models}

We first set the stage for the signal-plus-noise matrix model and review the basic properties of the spectral embedding in this section. Consider a symmetric positive semidefinite low-rank matrix $\bM\in\mathbb{R}^{n\times n}$ that can be written as $\bM = \rho_n\bX \bX \transpose$ for an $n\times d$ matrix $\bX$ and a scaling factor $\rho_n \in (0, 1]$, where $d\ll n$. The low-rank matrix $\bM$ represents the underlying signal matrix and is not accessible to the practitioners. Instead, only the noisy version $\bA$ of the signal matrix $\bM$ is observed. The signal-plus-noise matrix model specifies the following additive structure on $\bA$:
\begin{align}\label{eqn:signal_plus_noise}
\bA = \rho_n\bX\bX\transpose + \bE,
\end{align}
where $\bE = [E_{ij}]_{n\times n}$ is an $n\times n$ symmetric matrix of the noise and $(E_{ij}:1\leq i\leq j\leq n)$ are independent mean-zero random variables. The noise matrix $\bE$ is also referred to as the generalized Wigner matrix (see, for example, \citealp{yau2012universality}). The signal-plus-noise matrix model \eqref{eqn:signal_plus_noise} is flexible enough to include a broad range of popular statistical models, including 
% the stochastic block model \citep{HOLLAND1983109}, 
the random dot product graph \citep{young2007random} and the matrix completion problem \citep{candes2009exact}. We illustrate these special examples below in detail.

% \begin{example}[Stochastic block model]
% \label{example:SBM}
% Consider a network with $n$ vertices that are partitioned into $K$ communities for some $K\geq 1$, where $K$ is assumed to be much smaller than $n$. For any pair of vertices $(i, j)\in [n]\times [n]$, the probability that there exists an edge adjoining them is determined by the community memberships of vertex $i$ and vertex $j$. Formally, let $\sigma:[n]\to[K]$ be the community assignment function that assigns each vertex to a unique community label among the $K$ communities. Let $\bB = [B_{kl}]_{K\times K}\in(0, 1)^{K\times K}$ be a symmetric block probability matrix and $A_{ij}$ be the binary indicator of whether there is an edge linking vertices $i$ and $j$. Then the stochastic block model specifies that $A_{ij}\sim\mathrm{Bernoulli}(B_{\sigma(i)\sigma(j)})$ independently for $1\leq i\leq j\leq n$ and $A_{ji} = A_{ij}$. By converting the commmunity assignment $\sigma$ to a matrix $\bZ = [\mathbbm{1}\{\sigma(i) = k\}]_{n\times K}$, we see that the expected adjacency matrix $\bZ\bB\bZ\transpose$ is symmetric and low-rank. Furthermore, if $\bB$ is positive semidefinite with rank $d\leq K$ and can be factorized as $\bB = \bV\bV\transpose$ for a $K\times d$ matrix $\bV$, then $\bA$ can be represented using the signal-plus-noise matrix model \eqref{eqn:signal_plus_noise} with $\bE = [E_{ij}]_{n\times n}$, $E_{ij} = A_{ij} - B_{\sigma(i)\sigma(j)}$, $i,j\in [n]$, and $\rho_n^{1/2}\bX = \bZ\bV$. 
% \end{example}

\begin{example}[Random dot product graph]
\label{example:RDPG}
Consider a network with $n$ vertices labeled as $[n] = \{1,2,\ldots,n\}$. Each vertex $i\in [n]$ is assigned a $d$-dimensional Euclidean vector $\bx_i$, referred to as the latent position. The latent positions $\bx_1,\ldots,\bx_n$ are taken from the latent space $\calX\subset\mathbb{R}^d$ such that $\bx_i\transpose\bx_j\in [0, 1]$ for all $i, j\in [n]$. 
Let $\rho_n\in (0, 1]$ be the sparsity factor. 
Then the random dot product graph model generates a random adjacency matrix $\bA = [A_{ij}]_{n\times n}$ as follows: For each pair of vertices $(i, j)$, let 
% $A_{ij}$ be the binary indicator of the presence of an edge linking vertices $i$ and $j$. Then 
$A_{ij}\sim\mathrm{Bernoulli}(\rho_n\bx_i\transpose\bx_j)$ independently for $1\leq i\leq j\leq n$ and $A_{ji} = A_{ij}$. Clearly, 
with
% by denoting 
% $\rho_n^{1/2}\bx_{i} = \bx_i$ for all $i\in [n]$, 
$\bX = [\bx_{1},\ldots,\bx_{n}]\transpose$ and $\bE
 % = [E_{ij}]_{n\times n}
= [A_{ij} - \rho_n\bx_{i}\transpose\bx_{j}]_{n\times n}$, the random dot product graph model falls into the category of 
% the signal-plus-noise matrix 
model \eqref{eqn:signal_plus_noise}. 
\end{example}

\begin{example}[Symmetric noisy matrix completion]
\label{example:SNMC}
The general noisy matrix completion problem \citep{tight_oracle_inequalities,5466511} is described in the context of rectangular matrices, but the symmetric version of it also appears in certain applications, e.g., network cross-validation by edge sampling \citep{10.1093/biomet/asaa006}. Consider the signal-plus-noise matrix model \eqref{eqn:signal_plus_noise}, but the practitioners do not observe the complete matrix $\bA$. Instead, each entry $A_{ij}$ is observed with probability $p$ independently for $1\leq i\leq j\leq n$, and the missing entries of $\bA$ are replaced with zeros. Formally, let $z_{ij}\sim\mathrm{Bernoulli}(p)$ independently for all $1\leq i\leq j\leq n$, $z_{ji} = z_{ij}$, and denote $A_{ij}^{(\mathrm{obs})} = z_{ij}A_{ij}$. The matrix $\bA^{(\mathrm{obs})} = [A_{ij}^{(\mathrm{obs})}]_{n\times n}$ is biased for $\rho_n\bX\bX\transpose$, but $\bA^* = \bA^{(\mathrm{obs})} / p$ has the same expected value as $\rho_n\bX\bX\transpose$. Therefore, $\bA^*$ can be described by model \eqref{eqn:signal_plus_noise} as $\bA^* = \rho_n\bX\bX\transpose + \bE^*$, where $\bE^* = [z_{ij}A_{ij}/p - \expect A_{ij}]_{n\times n}$, $1\leq i\leq j\leq n$. 
\end{example}

In this work, we focus on estimating the signal matrix $\rho_n\bX\bX\transpose$ through the factor matrix $\bX\in\mathbb{R}^{n\times d}$.
 % that generates the underlying signal matrix $\bM$. 
Note that
 % without further restrictions, 
the signal-plus-noise matrix model \eqref{eqn:signal_plus_noise} is not identifiable in $\bX$. Firstly, 
% when $d$ is unknown, 
for any $d_1 > d$, there exists another matrix $\bX_1\in\mathbb{R}^{n\times d_1}$, such that $\bX\bX\transpose = \bX_1\bX_1\transpose$, and hence, they yield the same distribution on the observed random matrix $\bA$. This source of non-identifiability can be eliminated by requiring that $\sigma_d(\bX) > 0$. Secondly, 
% when $d$ is known, 
the factor matrix $\bX$ can only be identified up to an orthogonal matrix $\bW$ because $\bX\bX\transpose = (\bX\bW)(\bX\bW)\transpose$. The latter source of non-identifiability is inevitable without further constraints. Consequently, any estimator of $\bX$ can only recover $\bX$ up to an orthogonal transformation. 

% [Introduce the spectral embedding]
Perhaps the most straightforward estimator of $\bX$ is the spectral embedding estimator. It is formally defined as the solution to the least-squares problem 
\begin{align}\label{eqn:spectral_embedding}
\widetilde\bX = \arginf_{\bX\in\mathbb{R}^{n\times d}}\|\bA - \bX\bX\transpose\|_{\mathrm{F}}^2.
\end{align}
Conceptually, $\widetilde{\bX}\widetilde{\bX}\transpose$ is the projection of the noisy version of $\expect\bA$ to the space of all $n\times n$ rank-$d$ symmetric positive semidefinite matrices under the Frobenius norm metric. Practically, the spectral embedding $\widetilde{\bX}$ is simply the matrix concatenated by the top-$d$ scaled eigenvectors of $\bA$ \citep{Eckart1936}. Formally, let $\bA$ yield spectral decomposition
% \[
$\bA = \sum_{k = 1}^n\lambda_k(\bA)\widehat{\bu}_k\widehat{\bu}_k\transpose$,
% \]
where $\widehat{\bu}_k\transpose\widehat{\bu}_l = \mathbbm{1}(k = l)$ and $|\lambda_1(\bA)|\geq\ldots\geq|\lambda_n(\bA)|$. Then $\widetilde\bX$ can be taken as
% \[
$\widetilde\bX = \bU_\bA\bS_\bA^{1/2}$,
% \]
where $\bU_\bA = [\widehat{\bu}_1,\ldots,\widehat{\bu}_d]$ and $\bS_\bA = \mathrm{diag}\{\lambda_1(\bA),\ldots,\lambda_d(\bA)\}$. 

Although seemingly naive, the spectral embedding enjoys a collection of desirable features. In the context of stochastic block models, \cite{7298436}, \cite{10.1214/19-AOS1854}, \cite{lyzinski2014}, and \cite{sussman2012consistent} have shown that the spectral embedding can be applied to recover the community memberships of the underlying vertices.
% , a well-known procedure referred to as the spectral clustering. 
More generally, in the context of random dot product graphs, the asymptotic properties of the spectral embedding have been established, including the consistency \citep{6565321} and the central limit theorems \citep{athreya2016limit,tang2018,xie2021entrywise}. The eigenvector-based subsequent inference has also been studied, such as vertex classification \citep{tang2013} and hypothesis testing between graphs \citep{doi:10.1080/10618600.2016.1193505,tang2017}. For the generic signal-plus-noise matrix model \eqref{eqn:signal_plus_noise}, \cite{cape2019signal} has proved a sharp entrywise error bound for the unscaled eigenvectors
% of $\bA$ 
and a corresponding central limit theorem. 
% As will be seen in Section \ref{sec:main_results}, 
Their result is 
one of the building blocks
% of fundamental interest 
for developing the supporting theory of our proposed eigenvector-assisted estimation framework in Sections \ref{sec:eigenvector_assisted_estimation} and \ref{sec:main_results}. 

We close this subsection by constructing an appropriate orthogonal matrix $\bW$ to align the spectral embedding $\widetilde{\bX}$ with its estimand $\rho_n^{1/2}\bX$. 
% We denote $\bX_0$ the true value of $\bX$ giving rise to the distribution of $\bA$ through $\bA = \rho_n\bX_0\bX_0\transpose + \bE$.
 % $\rho_n^{1/2}\bX_0$. 
This alignment matrix is necessary for the theoretical analysis due to the orthogonal non-identifiability.
 % of the model \eqref{eqn:signal_plus_noise}. 
Nevertheless, the practitioners should be aware that it is not accessible because it requires the knowledge of the true value of $\bX$, which is not available in practice. 
To distinguish between the true value of $\bX$ and a generic $n\times d$ matrix $\bX$, we denote $\bX_0$ as the ground truth governing the distribution of the observed matrix $\bA$. 
Let $\rho_n^{1/2}\bX_0$ yield the singular value decomposition (SVD) $\rho_n^{1/2}\bX_0 = \bU_\bP\bS_\bP^{1/2}\bW_\bX$, where $\bU_\bP\in\mathbb{O}(n, d)$ and $\bW_\bX\in\mathbb{O}(d)$. 
Further let $\bU_\bP\transpose\bU_\bA$ have the SVD $\bU_\bP\transpose\bU_\bA = \bW_1\mathrm{diag}(\sigma_1,\ldots,\sigma_d)\bW_2\transpose$, where $\bW_1,\bW_2\in\mathbb{O}(d)$ and $\sigma_1 \geq\ldots\geq\sigma_d\geq 0$. Define the matrix sign \citep{10.1214/19-AOS1854,5714248} of $\bU_\bP\transpose\bU_\bA$ as $\bW^* = \mathrm{sgn}(\bU_\bP\transpose\bU_\bA) = \bW_1\bW_2\transpose$. Then the orthogonal alignment matrix between $\widetilde{\bX}$ and $\rho_n^{1/2}\bX_0$ is selected as $\bW = (\bW^*)\transpose\bW_\bX$. 
% As will be seen in Theorem \ref{thm:uniform_concentration_eigenvector}, 
\cite{tang2018} have shown that the choice of such an orthogonal alignment $\bW$ leads to the consistency result that $\|\widetilde{\bX}\bW - \rho_n^{1/2}\bX_0\|_{2\to\infty}\overset{\prob_0}{\to} 0$ under mild conditions.

% Let $\bP_0:=\rho_n\bX_0\bX_0\transpose$ yield spectral decomposition
% \[
% \bP_0 = \bU_\bP\bS_\bP\bU_\bP\transpose,
% \]
% where $\bU_\bP\transpose\bU_\bP = \eye_d$, $\bS_\bP = \mathrm{diag}\{\lambda_1(\bP_0),\ldots,\lambda_d(\bP_0)\}$. 

% subsection signal_plus_noise_matrix_models (end)

\section{Eigenvector-Assisted Estimation Framework} % (fold)
\label{sec:eigenvector_assisted_estimation}

% section eigenvector_assisted_estimation (end)

\subsection{Eigenvector-assisted estimating equation} % (fold)
\label{sub:moment_condition_models}

% [Background on moment condition models]
% [Motivate eigenvector-assisted estimating equations for moment conditions]

This subsection motivates the eigenvector-assisted estimating framework by 
constructing an asymptotically unbiased estimating equation. We first consider the problem of estimating a single row of $\bX_0$ when the remaining rows are available. Suppose we are interested in estimating the $i$th row of $\bX$ and assume that the remaining rows are readily available. Denote $\bX = [\bx_1,\ldots,\bx_n]\transpose$ and $\bX_0 = [\bx_{01},\ldots,\bx_{0n}]\transpose$. 
Namely, our goal is to estimate $\bx_{0i}$ given the information of $(\bx_{0j})_{j\neq i}$. Without loss of generality, we may consider estimating $\rho_n^{1/2}\bx_{0i}$ rather than $\bx_{0i}$ itself.
Then the signal-plus-noise matrix model \eqref{eqn:signal_plus_noise} implies that
 % $\expect(A_{ij} - \rho_n\bx_{i}\transpose\bx_{j}) = 0$, $i,j\in [n]$, where $\bX = [\bx_{1},\ldots,\bx_{n}]\transpose$. Note that when the parameter of interest is a single $\bx_i$, 
the data points $(\bx_j, A_{ij})_{j\neq i}$ come from the following linear regression model:
\begin{align*}
A_{ij} = \rho_n^{1/2}\bx_i\transpose\bx_{0j} + E_{ij},\quad j\in [n]\backslash \{i\},
\end{align*}
where $(\rho_n^{1/2}\bx_{0j})_{j\in [n]\backslash \{i\}}$ serve as the covariate vectors and $\bx_i$ is the unknown regression coefficient with the true value being $\rho_n^{1/2}\bx_{0i}$. 
% Unlike the classical linear model, 
% The model is potentially heteroskedastic because 
Note that the noise $(E_{ij})_{j\in [n]\backslash\{i\}}$ are independent but not necessarily identically distributed. To incorporate the potential heteroskedastic information, we consider a weight function $h_n(s, t):D\subset\mathbb{R}^2\to (0, +\infty)$ and the associated moment function
\begin{align*}
% $
\bg_{ij}(\bx_i) = (A_{ij} - \rho_n^{1/2}\bx_i\transpose\bx_{0j})h_n(\rho_n\bx_{0i}\transpose\bx_{0j}, \rho_n^{1/2}\bx_i\transpose\bx_{0j})\bx_{0j}.
% $
\end{align*}
Clearly, the moment conditions $\expect_0\{\bg_{ij}(\rho_n^{1/2}\bx_{0i})\} = \zero_d$ hold for $j\in [n]\backslash \{i\}$. 
Here, a canonical choice of the weight function $h_n$ is to require that $h_n(\rho_n\bx_{0i}\transpose\bx_{0j}, \rho_n\bx_{0i}\transpose\bx_{0j})^{-1} = \var_0(E_{ij})$ provided that the variance information of $E_{ij}$ is available and depends on $\bx_{0i}\transpose\bx_{0j}$. In general, the weight function $h_n(\cdot,\cdot)$ is quite flexible and can be designed according to the specific problem setup or the practitioners' expertise. 

The moment functions $\bg_{ij}(\bx_i)$, $j\in [n]\backslash\{i\}$ naturally lead to the unbiased generalized estimating equation (GEE)
\[
\frac{1}{n}\sum_{j \neq i}\bg_{ij}(\bx_i) = \zero_d.
\]
Solving the above GEE gives rise to a $Z$-estimator for $\rho_n^{1/2}\bx_{0i}$ provided that $(\bx_{0j})_{j\neq i}$ are accessible to the practitioners. However, obtaining the precise information of $(\bx_{0j})_{j\neq i}$ is non-trivial or even impossible for almost all real-world data problems. To this end, we introduce the eigenvector-assisted estimating equation
\begin{align}\label{eqn:EAEE}
\frac{1}{n}\sum_{j = 1}^n\widetilde{\bg}_{ij}(\bx_i) = \zero_d,\quad i \in [n],
\end{align}
where
% \[
$\widetilde{\bg}_{ij}(\bx_i) = (A_{ij} - \bx_i\transpose\widetilde{\bx}_j)h_n(\widetilde{\bx}_i\transpose\widetilde{\bx}_j, \bx_{i}\transpose\widetilde\bx_{j})\rho_n^{-1/2}\widetilde{\bx}_j$,
% \]
and $\widetilde{\bX} = [\widetilde{\bx}_1,\ldots,\widetilde{\bx}_n]\transpose$ is the spectral embedding defined in \eqref{eqn:spectral_embedding}. The eigenvector-assisted moment function $\widetilde{\bg}_{ij}(\bx_i)$ is obtained by replacing the unknown $(\bx_{0j})_{j\neq i}$ with their spectral embeddings $(\rho_n^{-1/2}\widetilde{\bx}_j)_{j\neq i}$. We refer to the solution $(\widehat{\bx}_i)_{i = 1}^n$ to the estimating equation \eqref{eqn:EAEE} as the eigenvector-assisted $Z$-estimator. 
% The flexibility of the choice of the weight function $h_n(\cdot,\cdot)$ allows some popular estimators developed in the literature can be formulated as the solutions to the equations of the form \eqref{eqn:EAEE}. 

\begin{remark}
An alternative strategy as opposed to replacing the unknown $(\bx_{0j})_{j\neq i}$ is to consider the following system of $n$ equations simultaneously:
\begin{align*}
% \begin{bmatrix}
% (1/n)
\frac{1}{n}
\sum_{j = 1}^n\bg_{1j}(\bx_1)
 = \zero_d
 ,\ldots,
% \\
% \vdots
% \\
% (1/n)
\frac{1}{n}\sum_{j = 1}^n\bg_{nj}(\bx_n)= \zero_d
% \end{bmatrix}
 % = \zero_{n\times d}
 . 
\end{align*}
However, the number of variables involved in this system is $n\times d$, and the computational cost of a solution may be expensive in general. In contrast, the eigenvector-assisted estimating equation \eqref{eqn:EAEE} can be solved for each $i\in [n]$ separately, where each sub-problem only contains $d$ variables. Consequently, the computation of the solutions to \eqref{eqn:EAEE} can be parallelized, which may further reduce the computational cost in practice. 
\end{remark}
Below, we provide two examples of the weight function $h_n(\cdot, \cdot)$ in the eigenvector-assisted estimating equation \eqref{eqn:EAEE}. These two choices of $h_n(\cdot,\cdot)$ lead to the spectral embedding defined in \eqref{eqn:spectral_embedding} and the one-step estimator for random dot product graphs \citep{xie2019efficient}. 
\begin{example}[Spectral embedding]
\label{example:Spectral_estimator}
The trivial choice that $h_n(s, t) = 1$ for all $(s, t)\in D$ results in 
% the $i$th row of 
the spectral embedding $\widetilde{\bX}$ as the corresponding eigenvector-assisted $Z$-estimator. To see this, denote $\be_i$ the $i$th standard basis vector in $\mathbb{R}^n$ whose coordinates are zeros except for the $i$th coordinate being one. Then the estimating equation \eqref{eqn:EAEE} implies
$\be_i\transpose (\bA \widetilde{\bX} -  \bX\widetilde{\bX}\transpose\widetilde{\bX}) = \zero_d\transpose$.
% reduces to
% \[
% \frac{1}{n}\sum_{j = 1}^n(A_{ij} - \bx_i\transpose\widetilde{\bx}_j)\widetilde{\bx}_j = 
% \frac{1}{n}\{\be_i\transpose (\bA \widetilde{\bX} -  \bX\widetilde{\bX}\transpose\widetilde{\bX})\}\transpose
%  = \zero_d. 
% \]
Note that $\bA\widetilde{\bX} = \bA\bU_\bA\bS_\bA^{1/2} = \bU_\bA\bS_\bA^{3/2} = \widetilde{\bX}(\widetilde{\bX}\transpose \widetilde{\bX})$ because $\widetilde{\bX} = \bU_\bA\bS_\bA^{1/2}$.
 % is the scaled eigenvector matrix corresponding to the top-$d$ largest eigenvalues of $\bA$ (in magnitude). 
The above estimating equation holds for all $i\in [n]$, implying that $\bA\bU_\bA\bS_\bA^{1/2} = \bX\bS_\bA$, and hence, $\bX = \bU_\bA\bS_\bA^{1/2}$ provided that $\bS_\bA$ is invertible. Therefore, the eigenvector-assisted $Z$-estimator coincides with the spectral embedding when $h_n(s, t) = 1$ for all $(s, t)\in D$. 
\end{example}

\begin{example}[One-step estimator for random dot product graphs]
\label{example:OSE_RDPG}
When $\bA$ is the adjacency matrix of a random graph, the spectral embedding $\widetilde{\bX}$ is also referred to as the adjacency spectral embedding (ASE) \citep{sussman2012consistent}. Although the ASE is practically useful because of the numerical stability and the ease of implementation, as pointed out by \cite{xie2019optimal} and \cite{xie2019efficient}, it is asymptotically sub-optimal because it does not
 % take the heteroskedasticity of the noise into account nor does it 
incorporate the information of the Bernoulli likelihood. Instead, \cite{xie2019efficient} proposed the following one-step estimator $\widehat{\bX}^{(\mathrm{OS})} = [\widehat{\bx}_1^{(\mathrm{OS})}, \ldots, \widehat{\bx}_n^{(\mathrm{OS})}]\transpose$ that improves upon the ASE:
\[
\widehat{\bx}_i^{(\mathrm{OS})} = \left\{\frac{1}{n}\sum_{j = 1}^n\frac{\widetilde\bx_j\widetilde\bx_j\transpose}{\widetilde\bx_i\transpose\widetilde\bx_j(1 - \widetilde\bx_i\transpose\widetilde\bx_j)}\right\}^{-1}\frac{1}{n}\sum_{j = 1}^n\frac{A_{ij}\widetilde\bx_j}{\widetilde\bx_i\transpose\widetilde\bx_j(1 - \widetilde\bx_i\transpose\widetilde\bx_j)},\quad i \in [n].
\]
It turns out that the one-step estimator coincides with the eigenvector-assisted $Z$-estimator when $h_n(s, t) = \{s(1 - s)\}^{-1}$. To see this, note that the estimating equation \eqref{eqn:EAEE} has the form
 % in this case can be equivalently written as
\[
\frac{1}{n}\sum_{j = 1}^n\frac{(A_{ij} - \bx_i\transpose\widetilde{\bx}_j)\widetilde\bx_j}{\widetilde\bx_i\transpose\widetilde\bx_j(1 - \widetilde\bx_i\transpose\widetilde\bx_j)} = \zero_d.
\]
Then a simple algebra shows that the solution to the above estimating equation coincides with the one-step estimator. The reason that the one-step estimator improves upon the ASE lies in the fact that $h_n(\rho_n\bx_{0i}\transpose\bx_{0j}, \rho_n\bx_{0i}\transpose\bx_{0j})^{-1} = \rho_n\bx_{0i}\transpose\bx_{0j}(1 - \rho_n\bx_{0i}\transpose\bx_{0j})$ is the same as the variance of $E_{ij} = A_{ij} - \rho_n\bx_{0i}\transpose\bx_{0j}$, $i,j \in [n]$. 
\end{example}

% \begin{example}[Weighted least squares for symmetric noisy matrix completion]
% \label{example:WLS_SNMC}
% \end{example}

\subsection{Generalized Bayesian estimation} % (fold)
\label{sub:generalized_bayesian_estimation_with_moment_conditions}

We now introduce the generalized Bayesian estimation method for the signal-plus-noise matrix model \eqref{eqn:signal_plus_noise} using the eigenvector-assisted estimating equation \eqref{eqn:EAEE}. Note that model \eqref{eqn:signal_plus_noise} does not specify a concrete likelihood function
 % unlike the classical fully Bayesian models, 
 due to its semiparametric nature. 
Therefore, we transform the zero-finding problem \eqref{eqn:EAEE} into a maximization problem and replace the usual log-likelihood function with the corresponding objective function. Specifically, let $\ell_{in}(\bx_i)$ be a criterion function whose maximizer is the solution to \eqref{eqn:EAEE} for each $i\in [n]$. Denote $\Theta$ the parameter space for $\bx_i$ and let $\pi(\bx_i)$ be the density of an absolutely continuous prior distribution on $\Theta$. Then we consider the following generalized posterior distribution associated with the criterion function $\ell_{in}(\bx_i)$:
\begin{align}\label{eqn:generalized_posterior}
\pi_{in}(\bx_i\mid\bA):=\frac{\exp\{\ell_{in}(\bx_i)\}\pi(\bx_i)}{\int_\Theta \exp\{\ell_{in}(\bx_i)\}\pi(\bx_i)\mathrm{d}\bx_i},\quad i \in [n].
\end{align}
Namely, the usual log-likelihood function for $\bx_i$ is substituted by the criterion function $\ell_{in}(\cdot)$ in the Bayes formula. 
Then the joint posterior distribution of $\bX
 % = [\bx_1,\ldots,\bx_n]\transpose
 $ is 
 % simply 
 obtained by taking the product: $\pi_n(\bX\mid\bA) = \prod_{i = 1}^n\pi_{in}(\bx_i\mid\bA)$. 
 In practice, the computation of the generalized posterior \eqref{eqn:generalized_posterior} can be implemented via a standard Metropolis-Hastings algorithm. The detailed algorithm is provided in the Supplementary Material. 
 
 Below, we consider three specific examples of the criterion function:
  % associated with \eqref{eqn:EAEE}: 
 the M-criterion function, the generalized method of moments (GMM) criterion function, and the exponentially tilted empirical likelihood (ETEL) criterion function. 

\vspace*{1ex}
\noindent\textbf{M-criterion.}
The most straightforward criterion function $\ell_{in}(\bx_i)$ is the indefinite integral of the estimating equation \eqref{eqn:EAEE} with respect to the argument $\bx_i$, leading to the following $M$-estimation criterion function:
\begin{align}
\label{eqn:M_estimation}
\ell_{in}(\bx_i) = \rho_n^{-1}\sum_{j = 1}^n\left\{\int_{t_0}^{\bx_i\transpose\widetilde{\bx}_j}(A_{ij} - t)h_n(\widetilde{\bx}_i\transpose\widetilde{\bx}_j, t)\mathrm{d}t\right\},
\end{align}
where $t_0\in\mathbb{R}$ is a fixed point such that $(\widetilde{\bx}_i\transpose\widetilde{\bx}_j, t)\in D$. The scaling factor $\rho_n^{-1}$ is added for technical considerations in Section \ref{sec:main_results} and does not change the maximizer of the criterion function. By the fundamental theorem of calculus, it is immediate to see that the gradient of the $M$-criterion function \eqref{eqn:M_estimation} coincides with $(1/n)\sum_{j = 1}^n\widetilde{\bg}_{ij}(\bx_i)$ up to a constant factor. 

\vspace*{1ex}
\noindent\textbf{GMM criterion.} The second choice of the criterion function that is maximized at the eigenvector-assisted $Z$-estimator is the generalized method of moments (GMM) criterion function:
\begin{align}\label{eqn:GMM}
\ell_{in}(\bx_i) = -\frac{n}{2}\left\{\frac{1}{n}\sum_{j = 1}^n\widetilde{\bg}_{ij}(\bx_i)\right\}\transpose\left\{\frac{1}{n}\sum_{j = 1}^n\widetilde{\bg}_{ij}(\widetilde{\bx}_i)\widetilde{\bg}_{ij}(\widetilde{\bx}_i)\transpose\right\}^{-1}\left\{\frac{1}{n}\sum_{j = 1}^n\widetilde{\bg}_{ij}(\bx_i)\right\}.
\end{align}
The GMM has been quite popular in econometrics \citep{10.2307/1912684,10.2307/2171802,10.2307/1912775,doi:10.1080/07350015.1996.10524656,10.2307/2971718}. It is clear that the maximizer of the GMM criterion \eqref{eqn:GMM} coincides with the zero to the estimating equation \eqref{eqn:EAEE} provided that $(1 / n)\sum_{j = 1}^n\widetilde{\bg}_{ij}(\widetilde\bx_i)\widetilde{\bg}_{ij}(\widetilde\bx_i)\transpose$ is positive definite. 

\vspace*{1ex}
\noindent\textbf{ETEL criterion.}
A popular Bayesian approach for moment condition models is the exponentially tilted empirical likelihood (ETEL) proposed in \cite{10.1093/biomet/92.1.31}. 
% The popularity is largely because 
In particular, \cite{10.1093/biomet/92.1.31} argued that
the ETEL could be interpreted as the limit of a nonparametric Bayesian procedure with a non-informative prior over the space of all distributions. For our purpose, we describe the ETEL criterion function in the context of the eigenvector-assisted estimating equation \eqref{eqn:EAEE}. Let $i\in [n]$ be a fixed row index. The ETEL is defined as the product of the empirical probabilities $\{p_{ij}(\bx_i)\}_{j = 1}^n$ for each observation $(A_{ij})_{j = 1}^n$:
% \begin{align*}
$L_{in}(\bx_i) = \prod_{i = 1}^np_{ij}(\bx_i)$.
% \end{align*}
Here, for each $i\in [n]$, $\{p_{ij}(\bx_i)\}_{j = 1}^n$ solve the constrained optimization problem
\begin{align}\label{eqn:ETEL_probabilities}
\begin{aligned}
\max_{[p_{i1},\ldots,p_{in}]\transpose\in \Psi_i(\bx_i)}&\sum_{j = 1}^n(-p_{ij}\log p_{ij}),
% \\
% \mbox{subject to }&\sum_{j = 1}^np_{ij} = 1,\quad\sum_{j = 1}^np_{ij}\widetilde{\bg}_{ij}(\bx_i) = \zero_d, \quad p_{ij}\geq 0,\quad j\in [n].
\end{aligned}
\end{align}
where $\Psi_i(\bx_i) = \{[p_{i1},\ldots,p_{in}]\transpose\in [0, 1]^n:\sum_{j = 1}^np_{ij} = 1, \sum_{j = 1}^np_{ij}\widetilde{\bg}_{ij}(\bx_i) = \zero_d\}$.
By the method of Lagrange multipliers, \cite{10.1093/biomet/asaa028} showed that $L_{in}(\bx_i) = \prod_{i = 1}^np_{ij}(\bx_i)$ is maximized at the eigenvector-assisted $Z$-estimator provided that the solution to the equation \eqref{eqn:EAEE} is well defined. Therefore, for each $i\in [n]$, the logarithmic ETEL 
% 
% $\ell_{in}(\bx_i) = \sum_{j = 1}^n\log\{p_{ij}(\bx_i)\}$ 
\begin{align}\label{eqn:ETEL}
\ell_{in}(\bx_i) = \sum_{i = 1}^n\log \{p_{ij}(\bx_i)\},
\end{align}
is also maximized at the solution to the equation \eqref{eqn:EAEE}. We refer to the criterion function \eqref{eqn:ETEL} as the ETEL criterion. In practice, for each $i\in [n]$ and any fixed $\bx_i\in\Theta$, the empirical probabilities $\{p_{ij}(\bx_i)\}_{j = 1}^n$ can be computed by solving the dual problem \citep{10.1214/009053606000001208}
\begin{align}\label{eqn:ETEL_probabilities_dual}
\begin{aligned}
p_{ij}(\bx_i) & = \frac{\exp\{\widehat\blambda_i(\bx_i)\transpose\widetilde{\bg}_{ij}(\bx_i)\}}{\sum_{l = 1}^n\exp\{\widehat\blambda_i(\bx_i)\transpose\widetilde{\bg}_{il}(\bx_i)\}},\quad j\in [n]\\
\widehat{\blambda}_i(\bx_i) &  = \argmin_{\blambda_i\in\mathbb{R}^d}\frac{1}{n}\sum_{j = 1}^n\exp\left\{\blambda_i(\bx_i)\transpose\widetilde{\bg}_{ij}(\bx_i)\right\}.
\end{aligned}
\end{align}

% [Why generalized Bayesian estimation? Natural uncertainty quantification; Avoid resampling or bootstrapping the data matrix. ]

% subsection generalized_bayesian_estimation_with_moment_conditions (end)

% section background_and_problem_formulation (end)

\section{Main Results} % (fold)
\label{sec:main_results}

\subsection{Large sample properties of the \texorpdfstring{$Z$}{Z}-estimator}
\label{sub:large_sample_properties_Z_estimator}

In this subsection, we establish the large sample properties of the eigenvector-assisted $Z$-estimator. 
% The rank $d$ of the mean matrix $\expect_0\bA$ is assumed to be known and does not change with the asymptotic index $n$ throughout. 
We first state the assumption for the signal-plus-noise matrix model \eqref{eqn:signal_plus_noise}.
\begin{assumption}[Sampling model]
\label{assumption:signal_plus_noise}
Model \eqref{eqn:signal_plus_noise} satisfies the following condition:
\begin{enumerate}[(i), noitemsep, topsep = 0mm]
  \item $\rho_n\in (0, 1]$, $\rho_n\to \rho\in [0, 1]$ exists, and $(\log n)^{4\xi}/(n\rho_n)\to 0$ for some constant $\xi > 1$.
%   there exists a constant $\xi > 1$ such that .

  \item There exist constants $C, c> 0$ such that $c(n\rho_n)\leq \lambda_d\{\expect_0(\bA)\}\leq \lambda_1\{\expect_0(\bA)\}\leq C(n\rho_n)$. 

  \item $\|\bE\|_2 = O\{(n\rho_n)^{1/2}\}$ w.h.p..
  
  \item There exist constants $C_E, \nu> 0$, such that for $\xi>1$ in (i) above, for all $n\geq N_0(C_E, \nu, \xi)$ depending on $C_E,\nu$, and $\xi$,
  \begin{align*}
  \prob_0\left[
  \bigcap_{m = 1}^{m(n) + 1}\bigcap_{k = 1}^d\left\{
  \left|\be_i\transpose\bE^m\bu_{0k}\right|\leq (C_En\rho_n)^{m/2}(\log n)^{m\xi}\|\bu_{0k}\|_\infty\right\}
  \right]\geq 1 - e^{-\nu(\log n)^\xi}.
  \end{align*}
  Here $m(n) = \lceil(\log n)/(\log n\rho_n)\rceil$ and $\bu_{0k}$ is the $k$th column vector of $\bU_\bP$. 

  \item The eigenvector matrix $\bU_\bP$ satisfies $\|\bU_\bP\|_{2\to\infty}\leq C_\mu\sqrt{d/n}$ for some constant $C_\mu \geq 1$. 

  \item $(A_{ij}:1\leq i\leq j\leq n)$ are independent; There exist constants $\sigma,\sigma_0 > 0$, such that $\expect_0|A_{ij} - \rho_n\bx_{0i}\transpose\bx_{0j}|\leq \sigma\rho_n$, $\var(A_{ij})\geq\sigma_0^2\rho_n$ for all $i,j\in [n]$, and either one of the following conditions holds:
  \begin{enumerate}[(a), noitemsep, topsep = 0mm]
    \item There exists a constant $M > 0$ such that $|A_{ij}|\leq M$ a.s., and $\var_0(A_{ij})\leq \sigma^2\rho_n$ for all $i,j\in [n]$. Without loss of generality we may assume that $M = 1$;
    \item $\max_{i,j\in [n]}\|A_{ij}\|_{\psi_2}\leq \sigma\rho_n^{1/2}$. 
  \end{enumerate}
\end{enumerate}
\end{assumption}
In Assumption \ref{assumption:signal_plus_noise} above, items (i) through (iv) have been adopted in \cite{cape2019signal} and are fundamental for the asymptotic normality of the rows of the unscaled eigenvector matrix $\bU_\bA$. Specifically, items (i) and (ii) introduce the scaling factor $\rho_n$ that governs the overall signal strength of $\expect_0\bA$. Item (iii) guarantees a concentration bound for the spectral norm of the noise matrix $\bE$, and item (iv) is a higher-order Bernstein-type concentration inequality for the row-wise behavior of $\bE$ and includes a broad class of generalized Wigner matrices \citep{cape2019signal,erdos2013,doi:10.1080/01621459.2020.1840990,doi:10.1080/01621459.2020.1751645}. In addition, item (v) is a delocalization condition for the  population unscaled eigenvector $\bU_\bP$ and appears in random graph inference \citep{JMLR:v18:17-448}, random matrix theory \citep{10.1215/00127094-3129809}, and matrix completion problems \citep{candes2009exact}. Item (vi) is a mild condition for the distribution of the noise matrix $\bE$. 

Next, Assumption \ref{assumption:regularity_condition} presents a standard regularity condition for the parameter space of $\bx_i$'s and the eigenvector-assisted estimating equation \eqref{eqn:EAEE}. 
\begin{assumption}[Regularity condition]
\label{assumption:regularity_condition}
Let $\Theta\subset\mathbb{R}^d$ be the parameter space for $\bx_1,\ldots,\bx_n$. 
\begin{enumerate}[(i), noitemsep, topsep = 0mm]
  \item $\Theta = \{\bx:\|\bx\|_2\leq r\}$ for some constant $r > 0$ and $\rho_n^{1/2}\bx_{0i}$ is inside the interior of $\Theta$. 

  \item The estimating equation \eqref{eqn:EAEE} has a unique solution $\widehat{\bx}_i$ inside the interior of $\Theta$ w.h.p.. 
\end{enumerate}
\end{assumption}
% In Assumption \ref{assumption:regularity_condition} (i) above, we select the parameter space $\Theta$ as a closed Euclidean ball in $\mathbb{R}^d$ with radius $r > 0$.
%  % because of the orthogonal non-identifiability of $\rho_n^{1/2}\bx_{0i}$. 
% Because $\rho_n^{1/2}\bx_{0i}$ is only identifiable up to an orthogonal alignment, we allow searching over a parameter space that is symmetric with regard to any orthogonal transformation. 
Assumption \ref{assumption:weight_functions} below is a Lipschitz condition for the weight function $h_n(\cdot, \cdot)$ in the estimating equation \eqref{eqn:EAEE} and can be satisfied, e.g., by the weight functions appearing in Examples \ref{example:Spectral_estimator} and \ref{example:OSE_RDPG}.
 % satisfy Assumption \ref{assumption:weight_functions}. 
\begin{assumption}[Weight functions]
\label{assumption:weight_functions}
There exist constants $c_1,c_2,K,\delta,r > 0$ such that for all $i,j\in [n]$, $(s, t)\in B(\rho_n\bx_{0i}\transpose\bx_{0j}, \rho_n\delta)\times [-r, r]$, the function $h_n$ is twice continuously differentiable, and
\begin{align*}
\begin{array}{lll}
c_1\leq h_n(s, t)\leq c_2,\quad&
|D^{(1, 0)}h_n(s, t)|\leq K\rho_n^{-1},\quad&
|D^{(0, 1)}h_n(s, t)|\leq K\rho_n,\\
|D^{(2, 0)}h_n(s, t)|\leq K\rho_n^{-2},\quad&
|D^{(1, 1)}h_n(s, t)|\leq K\rho_n^{-1},\quad&
|D^{(0, 2)}h_n(s, t)|\leq K\rho_n.
\end{array}
\end{align*}
\end{assumption}

We are now in a position to establish the large sample properties of the eigenvector-assisted $Z$-estimator.
% , including the consistency and the asymptotic normality. 
For notational convenience, denote
\begin{align}\label{eqn:matrix_formula}
\bG_{0in} = \frac{1}{n}\sum_{j = 1}^n\expect_0\left\{\frac{\partial\bg_{ij}}{\partial\bx_i\transpose}(\rho_n^{1/2}\bx_{0i})\right\}\quad\mbox{and}\quad
\bOmega_{0in} = \frac{1}{n}\sum_{j = 1}^n\expect_0\{\bg_{ij}(\rho_n^{1/2}\bx_{0i})\bg_{ij}(\rho_n^{1/2}\bx_{0i})\transpose\}.
\end{align}
\begin{theorem}
\label{thm:Large_sample_Z_estimator}
Suppose Assumptions \ref{assumption:signal_plus_noise}, \ref{assumption:regularity_condition}, and \ref{assumption:weight_functions} hold and
let $\widehat{\bx}_i$ be the solution to the estimating equation $(1/n)\sum_{j = 1}^n\widetilde{\bg}_{ij}(\bx_i) = \zero$ for each $i\in [n]$. 
% Denote
% \[
% \bG_{in}(\bx_i) = \frac{1}{n}\sum_{j = 1}^n\expect_0\left\{\frac{\partial\bg_{ij}}{\partial\bx_i\transpose}(\bx_i)\right\},\quad \bx_i\in\Theta.
% \]
Let $\bW\in\mathbb{O}(d)$ be the orthogonal alignment matrix between $\widetilde{\bX}$ and $\rho_n^{1/2}\bX_0$, where $\widetilde{\bX}$ is the spectral embedding in \eqref{eqn:spectral_embedding}. 
Then 
% for the $n$-dependent random matrix $\bW\in\mathbb{O}(d)$ in Theorem \ref{thm:uniform_concentration_eigenvector}, 
% \[
% \|\bW\transpose\widehat{\bx}_i - \rho_n^{1/2}\bx_{0i}\|_2 \lesssim \sqrt{\frac{(\log n)^{2\xi}}{(n\rho_n)}}\quad\mbox{w.h.p.}.
% \]
% Furthermore, 
\[
\sqrt{n}(\bW\transpose\widehat{\bx}_i - \rho_n^{1/2}\bx_{0i}) = -\frac{1}{\sqrt{n}}\sum_{j = 1}^n\bG_{0in}^{-1}\bg_{ij}(\rho_n^{1/2}\bx_{0i}) + O\left\{\frac{(\log n)^{2\xi}}{(n\rho_n)^{1/2}}\right\}\quad\mbox{w.h.p.},\quad i\in [n]
\]
and $\bOmega_{0in}^{-1/2}\bG_{0in}\sqrt{n}(\bW\transpose\widehat{\bx}_i - \rho_n^{1/2}\bx_{0i})\overset{\calL}{\to}\mathrm{N}_d(\zero_d, \eye_d)$. 
\end{theorem}
% Theorem \ref{thm:Large_sample_Z_estimator} suggests that
 % the following distributional approximation holds:
% \[
% $\sqrt{n}(\bW\transpose\widehat{\bx}_i - \rho_n^{1/2}\bx_{0i})\overset{\calL}{\approx}\mathrm{N}_d\left(\zero_d, \bG_{0in}^{-1}\bOmega_{0in}\bG_{0in}^{-1}\right)$
% for large $n$.
% \]
When the variance information of the noise $\bE$ is available and the weight function $h_n(s, t)$ satisfies $h_n(\rho_n\bx_{0i}\transpose\bx_{0j}, \rho_n\bx_{0i}\transpose\bx_{0j}) = \rho_n / \var_0(E_{ij})$ accordingly, Theorem \ref{thm:Large_sample_Z_estimator} further implies that
\[
\sqrt{n}(\bW\transpose\widehat{\bx}_i - \rho_n^{1/2}\bx_{0i})\overset{\calL}{\approx}\mathrm{N}_d\left(\zero_d, \left\{\frac{1}{n}\sum_{j = 1}^n\frac{\rho_n\bx_{0j}\bx_{0j}\transpose}{\var_0(E_{ij})}\right\}^{-1}\right).
\]
The following proposition shows that the asymptotic covariance matrix on the right-hand side of the above display is minimum in spectra among all eigenvector-assisted $Z$-estimators. 
\begin{proposition}\label{prop:optimal_weighting}
Suppose Assumptions \ref{assumption:signal_plus_noise} and \ref{assumption:regularity_condition} hold. Then for any weight function $h_n(s, t)$ satisfying Assumption \ref{assumption:weight_functions}, 
$\{(1/n)\sum_{j = 1}^n{\rho_n\bx_{0j}\bx_{0j}\transpose}/{\var_0(E_{ij})}\}^{-1}\preceq \bG_{0in}^{-1}\bOmega_{0in}\bG_{0in}^{-1}$.
% \[
% \left\{\frac{1}{n}\sum_{j = 1}^n\frac{\rho_n\bx_{0j}\bx_{0j}\transpose}{\var_0(E_{ij})}\right\}^{-1}\preceq \bG_{0in}^{-1}\bOmega_{0in}\bG_{0in}.
% \]
\end{proposition}

\begin{continueexample}{example:OSE_RDPG}
We now revisit Example \ref{example:OSE_RDPG} for illustration. In the context of random dot product graphs (Example \ref{example:RDPG}), with the weight function being $h_n(s, t) = \rho_n/\{s(1 - s)\}$, the eigenvector-assisted $Z$-estimator is the one-step estimator proposed in \cite{xie2019efficient}. Then 
% we have
% \begin{align*}
% \bG_{0in}  = -\frac{\rho_n^{1/2}}{n}\sum_{j = 1}^n\frac{\bx_{0j}\bx_{0j}\transpose}{\bx_{0i}\transpose\bx_{0j}(1 - \rho_n\bx_{0i}\transpose\bx_{0j})},\quad
% &\bOmega_{0in} 
% % \overset{\Delta}{=} \frac{1}{n}\sum_{j = 1}^n\expect_0\{\bg_{ij}(\rho_n^{1/2}\bx_{0i})\bg_{ij}(\rho_n^{1/2}\bx_{0i})\transpose\}
%  = \frac{\rho_n}{n}\sum_{j = 1}^n\frac{\bx_{0j}\bx_{0j}\transpose}{\bx_{0i}\transpose\bx_{0j}(1 - \rho_n\bx_{0i}\transpose\bx_{0j})}.
% \end{align*}
it follows immediately from Theorem \ref{thm:Large_sample_Z_estimator} that
% By the Lyapunov's central limit theorem and the Slutsky's theorem, we see that
\begin{align*}
% \left\{\frac{1}{n}\sum_{j = 1}^n\frac{\bx_{0j}\bx_{0j}\transpose}{\bx_{0i}\transpose\bx_{0j}(1 - \rho_n\bx_{0i}\transpose\bx_{0j})}\right\}^{1/2}
\sqrt{n}(\bW\transpose\widehat{\bx}_i - \rho_n^{1/2}\bx_{0i})\overset{\calL}{\approx}\mathrm{N}_d\left(\zero_d, \left\{\frac{1}{n}\sum_{j = 1}^n\frac{\bx_{0j}\bx_{0j}\transpose}{\bx_{0i}\transpose\bx_{0j}(1 - \rho_n\bx_{0i}\transpose\bx_{0j})}\right\}^{-1}\right).
\end{align*}
The above asymptotic normality coincides with Theorem 5 in \cite{xie2019efficient}. In addition, when the weight function is constantly one ($h_n(s, t) = 1$ for all $(s,t)$), the corresponding $Z$-estimator is the spectral embedding $\widetilde{\bX}$. Then Theorem \ref{thm:Large_sample_Z_estimator} implies that
% following the same reasoning, we have
% a simple algebra shows that $\bG_{in}(\rho_n^{1/2}\bx_{0i}) = (\rho_n^{1/2}/n)\bX_0\transpose\bX_0$, and therefore, 
\[
\sqrt{n}(\bW\transpose\widetilde{\bx}_i - \rho_n^{1/2}\bx_{0i})\overset{\calL}{\approx}\mathrm{N}_d\left(\zero_d, \bDelta_n^{-1}\left\{\frac{1}{n}\sum_{j = 1}^n{\bx_{0i}\transpose\bx_{0j}(1 - \rho_n\bx_{0i}\transpose\bx_{0j})}{\bx_{0j}\bx_{0j}\transpose}\right\}\bDelta_n^{-1}\right),
\]
where $\bDelta_n = (1/n)\bX_0\transpose\bX_0$. This recovers Theorem 1 in \cite{xie2019efficient}, which is rooted in \cite{athreya2016limit} and \cite{tang2018}. As shown in \cite{xie2019efficient}, the asymptotic covariance matrix of the spectral embedding is dominated by that of the one-step estimator in spectra because the weight function $h_n(s, t) = \rho_n/\{s(1 - s)\}$ adjusts for the heteroskedasticity of the noise matrix $\bE$. In contrast, the constant weight function $h_n(s, t) = 1$ ignores the variance information inherited from the Bernoulli likelihood. 
\end{continueexample}

\subsection{Convergence of the generalized posterior}
\label{sub:convergence_of_the_generalized_posterior}

We are now in a position to present the convergence properties of the generalized posterior  \eqref{eqn:generalized_posterior} with a generic criterion function $\ell_{in}(\bx_i)$.
% that is maximized at the eigenvector-assisted $Z$-estimator. 
Two necessary assumptions are in order. 
\begin{assumption}
\label{assumption:prior}
The prior density $\pi(\bx_i)$ is continuous over $\bx_i\in\Theta$ and there exist constants $c, C > 0$ such that $c\leq \pi(\bx_i)\leq C$ for all $\bx_i\in\Theta$. 
\end{assumption}
\begin{assumption}\label{assumption:criterion_function}
The criterion function $\ell_{in}(\bx_i)$ satisfies the following conditions: 
\begin{enumerate}[(i), noitemsep, topsep = 0mm]
  \item $\ell_{in}(\bx_i)$ is uniquely maximized at $\widehat{\bx}_i$ w.p.a.1, where $\widehat{\bx}_i$ solves equation \eqref{eqn:EAEE}.

  \item Let $\bW\in\mathbb{O}(d)$ be the orthogonal alignment matrix between the spectral embedding $\widetilde{\bX}$ and $\rho_n^{1/2}\bX_0$. There exist a positive definite matrix $\bSigma_{in}\in\mathbb{R}^{d\times d}$ whose eigenvalues are bounded away from $0$ and $+\infty$, and two positive sequences $(\eps_n)_{n = 1}^\infty$, $(\delta_n)_{n = 1}^\infty$, $\eps_n\leq \delta_n$ for all $n$, $n\eps_n^2\to +\infty$, $\max(\eps_n,\delta_n)\to 0$, such that for any row index $i\in [n]$ and $\alpha > 0$,
  \begin{align}
  \label{eqn:Hessian_A1}
  &\sup_{\bx_i\in B(\rho_n^{1/2}\bx_{0i}, \eps_n)}\left\|\frac{1}{n}\frac{\partial\ell_{in}}{\partial\bx_i\partial\bx_i\transpose}(\bW\bx_i) + \bW\bSigma_{in}\bW\transpose\right\|_2 = o\left(\frac{1}{n\eps_n^2}\right)\quad\mbox{w.p.a.1},\\
  \label{eqn:Hessian_A2}
  &\inf_{\bx_i\in B(\rho_n^{1/2}\bx_{0i}, 3\delta_n)}\lambda_{\min}\left\{-\frac{1}{n}\frac{\partial\ell_{in}}{\partial\bx_i\partial\bx_i\transpose}(\bW\bx_i)\right\}\gtrsim 1\quad\mbox{w.p.a.1},\\
  \label{eqn:identifiability}
  &\inf_{\bW\transpose\bx_i\notin B(\rho_n^{1/2}\bx_{0i}, \delta_n)}\{\ell_{in}(\widehat{\bx}_i) - \ell_{in}(\bx_i)\}\geq (1 + \alpha)d\log n\quad\mbox{w.p.a.1}.
  \end{align}
\end{enumerate}
\end{assumption}
Assumption \ref{assumption:prior}
% above, which is a mild condition, 
requires that the prior density is continuous and bounded away from $0$ and $+\infty$. Assumption \ref{assumption:criterion_function} is a requirement for the criterion function $\ell_{in}(\bx_i)$. As discussed in Section \ref{sub:generalized_bayesian_estimation_with_moment_conditions}, the maximizer of $\ell_{in}(\bx_i)$ needs to be the same as the solution to the estimating equation \eqref{eqn:EAEE}. Conditions \eqref{eqn:Hessian_A1} and \eqref{eqn:Hessian_A2} describe the local behavior of the Hessian of $\ell_{in}(\bx_i)$ in shrinking neighborhoods of the truth $\rho_n^{1/2}\bx_{0i}$. Specifically, in a shrinking neighborhood of $\rho_n^{1/2}\bx_{0i}$ with radius $\eps_n$, condition \eqref{eqn:Hessian_A1} requires that the negative Hessian of $(1/n)\ell_{in}$ is close to a deterministic $d\times d$ positive definite matrix $\bSigma_{in}$, and condition \eqref{eqn:Hessian_A2} guarantees that $\ell_{in}$ is strongly concave in a larger shrinking neighborhood of $\rho_n^{1/2}\bx_{0i}$ with radius $3\delta_n$. Finally, condition \eqref{eqn:identifiability} is an identifiability condition for the criterion function, which is standard in the literature on generalized Bayesian estimation (see, for example, \citealp{CHERNOZHUKOV2003293,doi:10.1080/01621459.2017.1358172,10.1093/biomet/asaa028,https://doi.org/10.1111/rssb.12342}). 

Below, Proposition \ref{prop:Criterion_satisfies_assumption} asserts that the M-criterion \eqref{eqn:M_estimation}, the GMM criterion \eqref{eqn:GMM}, and the ETEL criterion \eqref{eqn:ETEL} introduced in Section \ref{sub:generalized_bayesian_estimation_with_moment_conditions} satisfy Assumption \ref{assumption:criterion_function}.
\begin{proposition}\label{prop:Criterion_satisfies_assumption}
Suppose Assumptions \ref{assumption:signal_plus_noise}-\ref{assumption:weight_functions} hold. Then:
\begin{enumerate}[(a), noitemsep, topsep = 0mm]
  \item The M-criterion function \eqref{eqn:M_estimation} satisfies Assumption \ref{assumption:criterion_function} with $\bSigma_{in} = -\rho_n^{-1/2}\bG_{0in}$.
  % $\eps_n = (\log n)^{1/4}/\sqrt{n}$, $\delta_n = M_n\left\{(\log n)^{2\xi}/(n\rho_n)\right\}^{1/4}$, and $\bSigma_{in} = -\rho_n^{-1/2}\bG_{in}(\rho_n^{1/2}\bx_{0i})$, where $M_n = \log\log n$.

  \item The GMM criterion function \eqref{eqn:GMM} satisfies Assumption \ref{assumption:criterion_function} with
  $\bSigma_{in} = \bG_{0in}\transpose\bOmega_{0in}^{-1}\bG_{0in}$.
   % $\eps_n = (\log n)^{1/4}/\sqrt{n}$, $\delta_n = M_n\{(\log n)^{2\xi}/(n\rho_n)\}^{1/4}$, and  $\bSigma_{in} = \bG_{in}(\rho_n^{1/2}\bx_{0i})\transpose\bOmega_{in}(\rho_n^{1/2}\bx_{0i})\bG_{in}(\rho_n^{1/2}\bx_{0i})$, where $M_n = \log\log n$. 

  \item If further
%   Assumption \ref{assumption:signal_plus_noise}(vi) (a) and (b) is replaced by 
  Assumption \ref{assumption:signal_plus_noise} (vi) (b) holds (i.e., $\max_{i,j\in [n]}\|A_{ij}\|_{\psi_2}\leq \sigma\rho_n^{1/2}$), then the ETEL criterion function \eqref{eqn:ETEL} satisfies Assumption \ref{assumption:criterion_function} with 
  $\bSigma_{in} = \bG_{0in}\transpose\bOmega_{0in}^{-1}\bG_{0in}$.
  % $\eps_n = (\log n)^{\xi - 1}/\sqrt{M_nn}$, $\delta_n = M_n\sqrt{(\log n)^{2\xi + 1}/(n\rho_n)}$, and $\bSigma_{in} = \bG_{in}(\rho_n^{1/2}\bx_{0i})\transpose\bOmega_{in}(\rho_n^{1/2}\bx_{0i})^{-1}\bG_{in}(\rho_n^{1/2}\bx_{0i})$, where $M_n = \log\log n$. 
\end{enumerate} 
\end{proposition}

Theorem \ref{thm:BvM_generalized_posterior} below, which is the main result in this subsection, 
% guarantees that the generalized posterior \eqref{eqn:generalized_posterior} converges to a multivariate normal distribution 
% under Assumptions \ref{assumption:signal_plus_noise}-\ref{assumption:criterion_function}
establishes the large sample properties of the generalized posterior \eqref{eqn:generalized_posterior} under Assumptions \ref{assumption:signal_plus_noise}-\ref{assumption:criterion_function}.

\begin{theorem}[Convergence of the generalized posterior]
\label{thm:BvM_generalized_posterior}
Suppose Assumptions \ref{assumption:signal_plus_noise}-\ref{assumption:criterion_function} hold and
let $\widehat{\bx}_i$ be the solution to the estimating equation $(1/n)\sum_{j = 1}^n\widetilde{\bg}_{ij}(\bx_i) = \zero$ for each $i\in [n]$. 
% Consider the generalized posterior density \eqref{eqn:generalized_posterior}. 
% Let $\bW\in\mathbb{O}(d)$ be the $n$-dependent random matrix in Theorem \ref{thm:uniform_concentration_eigenvector},
Let $\bW\in\mathbb{O}(d)$ be the orthogonal alignment matrix between the spectral embedding $\widetilde{\bX}$ and $\rho_n^{1/2}\bX_0$,
 $\bt = \sqrt{n}\bW\transpose(\bx_i - \widehat{\bx}_i)$, and denote ${\pi}_{in}^*(\bt\mid\bA)$ the generalized posterior density of $\bt$ induced from $\pi_{in}(\bx_i\mid\bA)$ defined in \eqref{eqn:generalized_posterior}. Then for any $\alpha\geq 0$ and for each $i\in [n]$, 
\begin{align}\label{eqn:strong_convergence_generalized_posterior}
\int_{\mathbb{R}^d} (1 + \|\bt\|_2^\alpha)\left|{\pi}_{in}^*(\bt\mid\bA) - \frac{\exp(-\bt\transpose\bSigma_{in}\bt/2)}{\sqrt{\det(2\pi\bSigma_{in}^{-1})}}\right|\mathrm{d}\bt = o(1)\quad\mbox{w.p.a.1},
\end{align}
where $\bSigma_{in}$ is the $d\times d$ positive definite matrix in Assumption \ref{assumption:criterion_function}.
\end{theorem}
Theorem \ref{thm:BvM_generalized_posterior} implies that the total variation distance between the generalized posterior distribution of $\sqrt{n}\bW\transpose(\bx_i - \widehat{\bx}_i)$ and $\mathrm{N}_d(\zero_d, \bSigma_{in}^{-1})$ converges to $0$ in probability. 
% Namely, the shape of the generalized posterior distribution of $\bx_i$ can be approximated by a $d$-dimensional (random) multivariate normal distribution with mean $\widehat{\bx}_i$ and covariance matrix $(1/n)\bW\bSigma_{in}^{-1}\bW\transpose$. 
This result is also known as the Bernstein-von Mises theorem of the generalized posteriors \citep{CHERNOZHUKOV2003293,kleijn2012bernstein,JMLR:v22:20-469,10.1093/biomet/asy054,syring2020gibbs}.

\subsection{Generalized Bayesian inference}
\label{sub:consequence_of_generalized_posterior}

An important consequence of Theorem \ref{thm:BvM_generalized_posterior} is the asymptotic normality of the generalized posterior mean as a frequentist point estimator. Namely, the generalized posterior mean is asymptotically equivalent to the eigenvector-assisted $Z$-estimator up to the first order. This result is summarized in Theorem \ref{thm:GBE} below. 
\begin{theorem}[Generalized posterior mean]
\label{thm:GBE}
Assume the conditions of Theorem \ref{thm:BvM_generalized_posterior} hold.
Let $\bW\in\mathbb{O}(d)$ be the orthogonal alignment matrix between the spectral embedding $\widetilde{\bX}$ and $\rho_n^{1/2}\bX_0$.
% Let $\widehat{\bx}_i$ be the solution to the estimating equation $(1/n)\sum_{j = 1}^n\widetilde{\bg}_{ij}(\bx_i) = \zero$. 
% and 
For each $i\in [n]$, denote
% \[
$\bx_i^* = \int_\Theta\bx_i{\pi}_{in}(\bx_i\mid\bA)\mathrm{d}\bx_i$
% \]
the generalized posterior mean of $\pi_{in}(\bx_i\mid\bA)$ defined in \eqref{eqn:generalized_posterior}. Then for each $i\in [n]$,
  \begin{align}\label{eqn:asymptotic_normality_posterior_mean}
  \bOmega_{0in}^{-1/2}\bG_{0in}\sqrt{n}(\bW\transpose\bx_i^* - \rho_n^{1/2}\bx_{0i})\overset{\calL}{\to}\mathrm{N}_d(\zero_d, \eye_d).
  \end{align}
  % where
  % \begin{align*}
  % \bG_{0in} = \frac{1}{n}\sum_{j = 1}^n\expect_0\left\{\frac{\partial\bg_{ij}}{\partial\bx_i\transpose}(\rho_n^{1/2}\bx_{0i})\right\}\quad\mbox{and}\quad
  % \bOmega_{0in} = \frac{1}{n}\sum_{j = 1}^n\expect_0\left\{\bg_{ij}(\rho_n^{1/2}\bx_{0i})\bg_{ij}(\rho_n^{1/2}\bx_{0i})\transpose\right\}.
  % \end{align*}
\end{theorem}

% [Highlight the key contribution: Valid uncertainty quantification]
Another useful consequence of Theorem \ref{thm:BvM_generalized_posterior} is that the generalized posterior \eqref{eqn:generalized_posterior} provides a convenient approach for valid uncertainty quantification without bootstrapping the data matrix $\bA$, which is a fascinating feature of the eigenvector-assisted estimation framework. In order to produce a credible region with the correct coverage probability, we require that the following generalized information equality holds \citep{CHERNOZHUKOV2003293}:
\begin{align}\label{eqn:generalized_information_equality}
\lim_{n\to\infty}\bSigma_{in}^{-1}(\bG_{0in}\transpose\bOmega_{0in}^{-1}\bG_{0in}) = \eye_d. 
\end{align}
This equality guarantees that the asymptotic distribution of $\sqrt{n}(\bW\transpose\widehat{\bx}_i - \rho_n^{1/2}\bx_{0i})$ coincides with the Bernstein-von Mises limit of $\pi_{in}\{\sqrt{n}\bW\transpose(\bx_i - \widehat{\bx}_i)\in\cdot\mid\bA\}$. Proposition \ref{prop:Criterion_satisfies_assumption} shows that the GMM criterion \eqref{eqn:GMM} and the ETEL criterion \eqref{eqn:ETEL} satisfy the generalized information equality. For the $M$-criterion \eqref{eqn:M_estimation}, this equality holds provided that $-\rho_n^{1/2}\bG_{0in} = \bOmega_{0in}$. In particular, if the weight function $h_n(s, t)$ satisfies $h_n(\rho_n\bx_{0i}\transpose\bx_{0j}, \rho_n\bx_{0i}\transpose\bx_{0j}) = \rho_n/\var_0(E_{ij})$ for all $i,j\in [n]$, then equality \eqref{eqn:generalized_information_equality} holds for the $M$-criterion \eqref{eqn:M_estimation}.

Given a confidence level $\alpha\in (0, 1)$, we can construct a $(1 - \alpha)$ credible region for $\rho_n^{1/2}\bx_{0i}$ up to the orthogonal alignment $\bW\in\mathbb{O}(d)$ using the generalized posterior distribution \eqref{eqn:generalized_posterior}. Let $\widehat{\bV}_B$ be the covariance matrix of the generalized posterior $\eqref{eqn:generalized_posterior}$. In practice, $\widehat{\bV}_B$ can be estimated conveniently using the covariance matrix of the generalized posterior samples generated from the MCMC sampler. 
% Specifically, if $(\bx_i^{(t)})_{t = 1}^T$ are the output of the MCMC sampler of \eqref{eqn:generalized_posterior}, then we can estimate $\widehat{\bV}_B$ by
% \[
% \widehat{\bV}_B = \frac{1}{T}\sum_{t = 1}^T(\bx_i^{(t)} - \bar{\bx}_i)(\bx_i^{(t)} - \bar{\bx}_i)\transpose,\quad\mbox{where}\quad \bar{\bx}_i = \frac{1}{T}\sum_{t = 1}^T\bx_i^{(t)}.
% \]
Let $q_{(1 - \alpha)}$ be the $(1 - \alpha)$ quantile of the $\chi^2$ distribution with degree of freedom $d$. A large sample $(1 - \alpha)$ credible ellipse is then given by 
\begin{align}\label{eqn:GBI_credible_ellipse}
\calE_{in} = \left\{\bx_i:(\bx_i - \widehat{\bx}_i)\transpose\widehat{\bV}_B^{-1}(\bx_i - \widehat{\bx}_i)\leq q_{(1 - \alpha)}\right\}.
\end{align}
In what follows, Theorem \ref{thm:GBI} establishes that the credible ellipse \eqref{eqn:GBI_credible_ellipse} has an asymptotic valid $(1 - \alpha)$ coverage probability for $\rho_n^{1/2}\bx_{0i}$ up to an orthogonal transformation.
\begin{theorem}[Generalized posterior inference]
\label{thm:GBI}
Assume the conditions of Theorem \ref{thm:BvM_generalized_posterior} and the generalized information equality \eqref{eqn:generalized_information_equality} hold. 
% Let $\widehat{\bx}_i$ be the solution to the estimating equation $(1/n)\sum_{j = 1}^n\widetilde{\bg}_{ij}(\bx_i) = \zero$. 
Given $\alpha\in (0, 1)$, let $\calE_{in}$ be the $(1 - \alpha)$ credible ellipse defined in \eqref{eqn:GBI_credible_ellipse} and $\bW\in\mathbb{O}(d)$ be the orthogonal alignment matrix between the spectral embedding $\widetilde{\bX}$ and $\rho_n^{1/2}\bX_0$. Then 
% for the $n$-dependent random matrix $\bW\in\mathbb{O}(d)$ defined in \eqref{thm:uniform_concentration_eigenvector}, 
$\prob_0(\rho_n^{1/2}\bW\bx_{0i}\in\calE_{in})\to 1 - \alpha$ as $n\to\infty$. 
\end{theorem}

% section main_results (end)

\section{Numerical Examples} % (fold)
\label{sec:numerical_examples}

% \subsection{Random graph estimation} % (fold)
% \label{sub:random_graph_estimation}

\subsection{Synthetic examples} % (fold)
\label{sub:synthetic_examples}

We first illustrate the proposed eigenvector-assisted estimation framework using synthetic datasets. The MCMC sampler used here is the Metropolis-Hastings algorithm implemented in the \texttt{mcmc} \texttt{R} package \citep{mcmc} with parallelization over the row index $i\in [n]$. For each Markov chain, the first $1000$ iterations are discarded as the burn-in stage, and the subsequent $2000$ are collected as post-burn-in samples. The convergence diagnostics of the MCMC are provided in the Supplementary Material, and there are no signs of non-convergence. 

Below, we consider two simulation scenarios that fall into the category of the signal-plus-noise matrix model \eqref{eqn:signal_plus_noise}:
\begin{itemize}[noitemsep, topsep = 0mm]
  \item \textbf{Scenario I:} Random dot product graph model.
   % (see Example \ref{example:RDPG} for the definition). 
  The factor matrix $\bX_0$, also known as the latent position matrix, is generated from the curve $f(t) = 0.1 + 0.8\sin(\pi t)$, where $t\in [0, 1]$. Specifically, let $n = 800$, $\rho_n = 1$, $0 = t_1\leq t_2\leq\ldots\leq t_n = 1$ be equidistant points over $[0, 1]$, and the ground true $\bX_0$ be an $n\times 1$ matrix whose entries are $f(t_1),\ldots,f(t_n)$. 
  % Namely, the rank of the expected adjacency matrix is $1$. 
  Then for any $i,j\in [n]$, $i\leq j$, we generate the $(i, j)$th entry $A_{ij}$ of $\bA$ from $\mathrm{Bernoulli}(f(t_i)f(t_j))$ independently and we set $A_{ij} = A_{ji}$ for all $i > j$. The same example has also been considered in \cite{xie2019efficient}. 

  \item \textbf{Scenario II:} Symmetric noisy matrix completion.
   % (see Example \ref{example:SNMC}). 
  Let $\bX_0$ be a $n\times 1$ matrix defined in scenario I above, namely, $\bX_0 = [f(t_1),\ldots,f(t_n)]\transpose$, where $f(t) = 0.1 + 0.8\sin(\pi t)$ and $0 = t_1 \leq t_2\leq\ldots\leq t_n = 1$ are equidistant points over $[0, 1]$. 
  A symmetric random matrix $\bA^\star = [A_{ij}^\star]_{n\times n}$ is generated with $\bA^\star = \bX_0\bX_0\transpose + \bE$, where $\bE = [E_{ij}]_{n\times n}$, $(E_{ij}:1\leq i\leq j\leq n)$ are independent and identically distributed $\mathrm{N}(0, 1)$ random variables, and $E_{ij} = E_{ji}$ for all $i > j$. 
  % We set $E_{ij} = E_{ji}$ and $A_{ij}^\star = A_{ji}^\star$ for all $i > j$. 
  % The matrix $\bA^\star$ is observed with missing entries
  % , where the missing mechanism is \emph{missing completely at random} (MCAR) \citep{rubinmultiple} 
  Each $A_{ij}^\star$ is observed with probability $p$ independently for all $i\leq j$, $i,j\in [n]$. Formally, following the formulation in Example \ref{example:SNMC}, we let $z_{ij}\sim\mathrm{Bernoulli}(p)$ independently for $1\leq i\leq j\leq n$, $z_{ij} = z_{ji}$ for $i > j$, and $A_{ij} = A^\star_{ij}z_{ij}/p$. Namely, the $(i, j)$th entry of the matrix $\bA = [A_{ij}]_{n\times n}$ is $A_{ij}^\star/p$ if it is observed, and is $0$ if it is missing. Here we take $n = 400$ and $p = 0.6$. 
\end{itemize}

For each of the scenarios above, given a realization of the data matrix $\bA$, we consider the following approaches for estimating $\bX_0$: The spectral embedding (also known as the adjacency spectral embedding/ASE under scenario I), the eigenvector-assisted $Z$-estimate, and the three generalized Bayesian estimation methods associated with the $M$-criterion, the GMM criterion, and the ETEL criterion, respectively. For scenario I, we take $h_n(s, t) = 1/\{s(1 - t)\}$ as the weight function with the parameter space for $x_1,\ldots,x_n$ being $\Theta = [-1, 1]$. For scenario II, we let the weight function be $h_n(s, t) = p/\{(1 - p)t^2 + 1\}$ and the parameter space be $\Theta = [-1.2, 1.2]$. 
For the generalized posterior distributions, the posterior means are computed as the corresponding point estimates.
% with the three criterion functions. 
The same numerical experiment is repeated for $500$ independent Monte Carlo replicates for both scenario I and scenario II. 

We focus on the following inference objectives: The estimation accuracy of $\bX_0$ and the coverage probabilities of the (entrywise) generalized credible intervals for $\bX_0$. For the first objective, given one of the aforementioned estimates $\widehat{\bX}$ for $\bX_0$, we use the sum-of-squares error $\mathrm{SSE} = \|\widehat{\bX}\bW - \bX_0\|_{\mathrm{F}}^2$ 
as the evaluation metric, where $\bW\in\mathbb{O}(d)$ is the orthogonal alignment matrix between the spectral embedding and the ground truth. 
For the second objective, 
we compute the empirical coverage probabilities of the generalized posterior $95\%$ credible intervals for each $\bx_{0i}$ by taking the average number of credible intervals that cover the ground truth. 
% or the two aforementioned scenarios, the orthogonal alignment can be replaced by taking the absolute values of the entries of $\widehat{\bX}$. This is because we know that the latent positions are non-negative and one-dimensional orthogonal matrices are $\{\pm1\}$. 
% Therefore, we can use the absolute values of the entries to replace the orthogonal alignment. 
% The Metropolis-Hastings algorithm for the computation of the generalized posteriors is implemented with $1000$ burn-in iterations and $2000$ post-burn-in samples. The convergence diagnostics are provided in the Supplementary Material. 
\begin{figure}[t]
  \centerline{\includegraphics[width=1\textwidth]{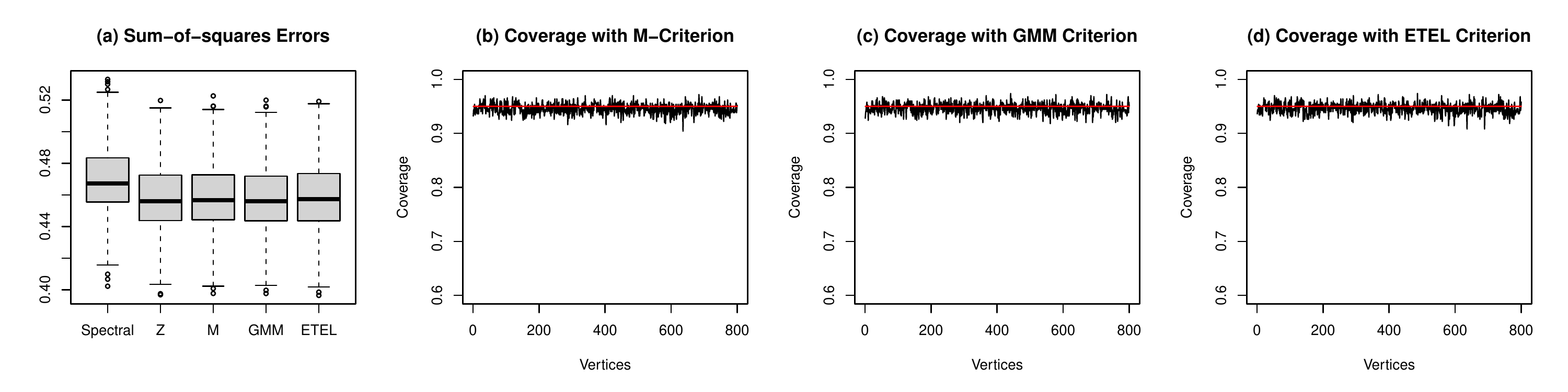}}
  \caption{Results for scenario I in Section \ref{sub:synthetic_examples}: Panel (a) presents boxplots of the sum-of-squares errors of the five point estimates involved across the $500$ Monte Carlo replicates; Panels (b), (c), and (d) display the empirical coverage probabilities of the (entrywise) $95\%$ credible intervals for $\bX_0$ obtained from the generalized posterior distributions with the $M$-criterion, the GMM criterion, and the ETEL criterion, respectively, where the red horizontal lines correspond to the nominal $95\%$ coverage probability.}
  \label{fig:graph_simulation}
\end{figure}
\begin{figure}[t]
  \centerline{\includegraphics[width=1\textwidth]{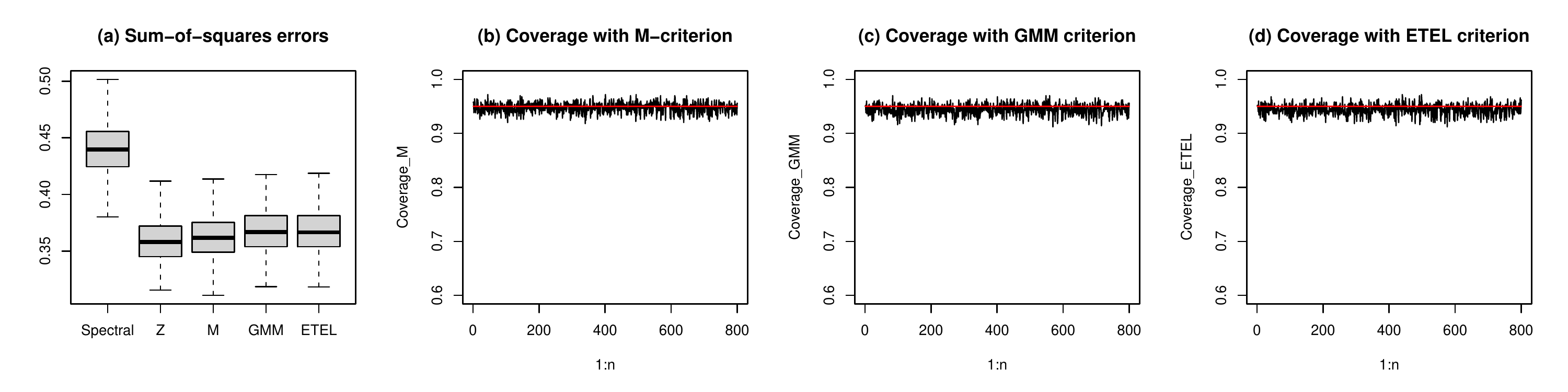}}
  \caption{Results for scenario II in Section \ref{sub:synthetic_examples}: Panel (a) presents boxplots of the sum-of-squares errors of the five point estimates involved across the $500$ Monte Carlo replicates; Panels (b), (c), and (d) display the empirical coverage probabilities of the (entrywise) $95\%$ credible intervals obtained from the generalized posterior distributions with the $M$-criterion, the GMM criterion, and the ETEL criterion, respectively, where the red horizontal lines correspond to the nominal $95\%$ coverage probability.}
  \label{fig:SNMC_simulation}
\end{figure}

Figure \ref{fig:graph_simulation} (a) and Figure \ref{fig:SNMC_simulation} (a) display the boxplots of the sum-of-squares errors of the aforementioned point estimates across the $500$ Monte Carlo replicates for scenarios I and II, respectively. The eigenvector-assisted $Z$-estimate and the generalized posterior means with the M-criterion, the GMM criterion, and the ETEL criterion have similar performance, and they all have smaller sum-of-squares errors than the spectral embedding. As discussed in Section \ref{sub:large_sample_properties_Z_estimator}, the improvement of the eigenvector-assisted estimates is because of the choice of the weight functions that encode the heteroskedastic variance information of the noise $\bE$, whereas the spectral embedding does not take it into account. The $p$-values of the two-sample $t$-tests among different sum-of-squares errors are tabulated in Table \ref{table:graph_ttest_pvalue}, which shows that the differences between the spectral embedding and the rest of the estimates are significant. 

\begin{table}[htbp]
  \caption{Results for Section \ref{sub:synthetic_examples}: the $p$-values of the two-sample $t$-tests among different sum-of-squares errors}
  \centering{%
  \begin{tabular}{c | c c c c }
    \hline\hline
    Comparison & Spectral vs $Z$ & Spectral vs $M$ & Spectral vs GMM & Spectral vs ETEL\\
    \hline
    Scenario I & $1.1\times 10^{-15}$ & $1.1\times 10^{-14}$ & $3.1\times 10^{-16}$ & $1.5\times 10^{-15}$ \\
    \hline
    Scenario II & $<2.2\times 10^{-16}$ & $<2.2\times 10^{-16}$ & $<2.2\times 10^{-16}$ & $<2.2\times 10^{-16}$ \\
    \hline\hline
  \end{tabular}%
  }
  \label{table:graph_ttest_pvalue}
  % \begin{tabnote}
% \end{tabnote}
\end{table}%
We also visualize the empirical coverage probabilities of the vertex-wise $95\%$ credible intervals using the generalized posteriors in Figures \ref{fig:graph_simulation} (b), \ref{fig:graph_simulation} (c), \ref{fig:graph_simulation} (d) under scenario I and Figures \ref{fig:SNMC_simulation} (b), \ref{fig:SNMC_simulation} (c), \ref{fig:SNMC_simulation} (d) under scenario II, respectively. Because the generalized information equality \eqref{eqn:generalized_information_equality} holds for both scenarios for the three criterion functions involved, the empirical coverage probabilities of the $95\%$ credible intervals obtained from the generalized posteriors are close to the nominal $95\%$ coverage probability. 
% We have also seen that the generalized information equality \eqref{eqn:generalized_information_equality} holds for the GMM criterion and the ETEL criterion. Therefore, the corresponding empirical coverage probabilities of their corresponding credible intervals are also close to the nominal $95\%$ coverage proability. 
These numerical findings validate the theoretical results in Section \ref{sec:main_results} empirically.

% subsection random_graph_estimation (end)

\subsection{Real-world network examples} % (fold)
\label{sub:real_data}

We now apply the proposed eigenvector-assisted estimation framework to real-world network examples. The datasets of interest are the ENZYMES networks taken from the BRENDA enzyme database \citep{schomburg2004brenda}. The networks are also publicly available at \url{https://networkrepository.com/index.php} \citep{nr}. 
These networks are graph representations of specific proteins.
% , and each graph is constructed as follows: 
The vertices represent the secondary structure elements that appear on certain amino acid sequences, and the existence of an edge linking two secondary structure elements means that the two elements appear as neighbors in the corresponding amino acid sequence or neighbors in the three-dimensional space \citep{10.1093/bioinformatics/bti1007}. In this study, we focus on the networks labeled ENZYMES 118, ENZYMES 123, ENZYMES 296, and ENZYMES 297. The summary statistics of these networks are provided in Table \ref{table:ENZYMES_summary} below. 
\begin{table}[htbp]
  \caption{Summary statistics of the ENZYMES networks}
  \centering{%
  \begin{tabular}{c | c c c c }
    \hline\hline
    Network label & ENZYMES 118 & ENZYMES 123 & ENZYMES 296 & ENZYMES 297\\
    \hline
    Number of vertices  & 95 & 90 & 125 & 121 \\
    \hline
    Number of edges & 121 & 127 & 141 & 149\\
    \hline
    Average degree & 5 & 9 & 5 & 7\\
    \hline\hline
  \end{tabular}%
  }
  \label{table:ENZYMES_summary}
  % \begin{tabnote}
% \end{tabnote}
\end{table}%

\begin{figure}[htbp]
  \centerline{\includegraphics[width=1\textwidth]{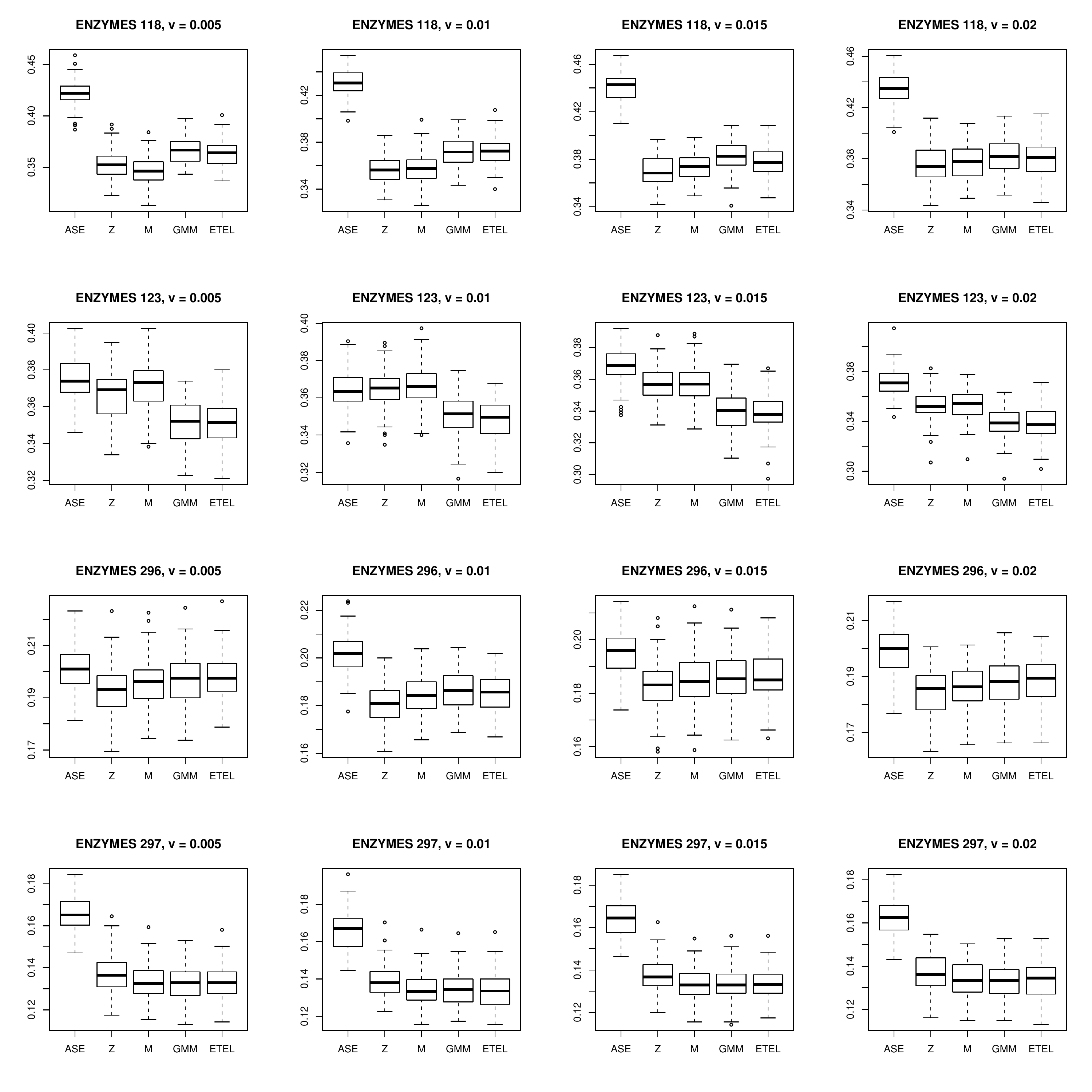}}
  \caption{The boxplots of the misclassification errors for the ENZYMES network data with different choice of $v$ across $50$ repeated experiments. }
  \label{fig:ENZYMES}
\end{figure}

We use the random dot product graph model as the working model for these ENZYMES networks. 
In addition to the observed adjacency matrices \emph{per se}, the class labels of the vertices are also available. Here, the inference goal of interest is the vertex classification when the observed network is contaminated by additional noise. The entire data analysis experiment consists of the following steps:
\begin{itemize}[noitemsep, topsep = 0mm]
  \item \textbf{Step 1: Noisy contamination of the data. }The adjacency matrix $\bA$ for each network is added with a symmetric noise matrix $\widetilde{\bE} = [\widetilde{E}_{ij}]_{n\times n}$ whose upper diagonal entries $(\widetilde{E}_{ij}:1\leq i\leq j\leq n)$ are independent and identically distributed $\mathrm{N}(0, v^2)$ random variables. The resulting data matrix $\widetilde{\bA} \overset{\Delta}{=} \bA + \widetilde{\bE}$ still falls into the category of the signal-plus-noise matrix model \eqref{eqn:signal_plus_noise} and has the same expected value as the original adjacency matrix $\bA$. 

  \item \textbf{Step 2: Dimensionality reduction. }Next, we estimate the latent position matrix using the following approaches: the adjacency spectral embedding (ASE), the eigenvector-assisted $Z$-estimate, the generalized posteriors with the $M$-criterion, the GMM criterion, and the ETEL criterion, respectively. Following the optimal weighting in Proposition \ref{prop:optimal_weighting}, we select the weight function as $h_n(s, t) = 1/\{s(1 - t) + v^2\}$ to match the reciprocal of the variance. We set the rank $d$ to be the same as the number of unique labels in each network. To compute the generalized posterior distributions, we implement the Metropolis-Hastings algorithm with $1000$ burn-in iterations, followed by another $2000$ post-burn-in MCMC samples. The convergence diagnostics are provided in the Supplementary Material, and they show no signs of non-convergence. We use the generalized posterior means as the point estimates. 

  \item \textbf{Step 3: Vertex classification. }The aforementioned five estimates are treated as the low-dimensional vertex features and fed into the $5$-nearest-neighbor classifier (5-NN) as the input variables for vertex classification. For each network, the 5-NN is implemented with approximately $75\%$ vertices as training data and the remaining vertices as testing data.
  % , where the partition is uniformly at random. 
  For each realization of the data matrix $\widetilde{\bA}$, the training-testing procedure is repeated independently for $100$ replicates, and the average misclassification errors on the testing data are reported. 
\end{itemize}
The range of the additional noise standard deviation $v$ is set to $\{0.005, 0.010, 0.015, 0.020\}$. For each fixed $v$, Steps 1-3 above are repeated for $50$ independent copies. The boxplots of the misclassification errors for the networks ENZYMES 118, ENZYMES 123, ENZYMES 296, and ENZYMES 297 with different choices of $v$ across $50$ repeated experiments are visualized in Figure \ref{fig:ENZYMES}. 
% It can be seen that 
Clearly,
for ENZYMES 118, ENZYMES 296, and ENZYMES 297, the proposed eigenvector-assisted estimates all outperform the baseline ASE significantly for different values of $v$. For ENZYMES 123, the generalized posterior means with the GMM criterion and the ETEL criterion have lower misclassification errors than those given by the baseline ASE, the eigenvector-assisted $Z$-estimate, and the generalized posterior mean with the $M$-criterion. Also, for ENZYMES 123, when $v$ increases, the eigenvector-assisted $Z$-estimate and the generalized posterior mean with the $M$-criterion outperform the ASE with lower misclassification errors. 
Overall, it is clear from the boxplots in Figure \ref{fig:ENZYMES} that the proposed methodology is more robust to the additional noisy contamination of the data matrix in terms of the vertex classification performance of the ENZYMES networks. 

% subsection vertex_classification_on_real_world_networks (end)

% section numerical_examples (end)

\section{Discussion}
\label{sec:discussion}

In this work, we propose a statistical inference framework for a broad range of signal-plus-noise matrix models using generalized posterior distributions based on a novel eigenvector-assisted estimating equation. 
The framework shares several fascinating properties. Firstly, it is quite flexible and allows the users to incorporate the heteroskedastic variance information of the noise. Secondly, it does not require the full specification of the noise distribution. Furthermore, from the computational perspective, the generalized posteriors can be computed via a Markov chain Monte Carlo sampler, which circumvents the potential challenging nonconvex optimization problems. In addition, the simulation-based inference algorithm also supplies the practitioners with a convenient environment for the uncertainty analysis and avoids the non-trivial resampling of the data matrix. Last but not least, our framework is backboned by solid theoretical support as we establish the large sample properties of the eigenvector-assisted $Z$-estimator and the generalized posterior distributions under mild regularity conditions. 

There are several potential extensions of the current framework. The large sample properties established in Section \ref{sec:main_results} may be applicable for certain subsequent inference tasks, such as testing whether two vertices in a stochastic block model are in the same community \citep{fan2019simple}. Our current signal-plus-noise matrix models are designed for symmetric random matrices with independent upper diagonal entries. There are, however, many high-dimensional statistical problems involving rectangular random matrices with low expected ranks, such as principal component analysis, high-dimensional clustering, compressed sensing, and collaborative filtering. It would be interesting to explore the singular-vector-assisted inference framework for general rectangular random matrices by taking advantage of the recent advance in the entrywise singular vector estimation \citep{agterberg2021entrywise,cape2017two}. On the practical side, our current computational strategy is a standard Metropolis-Hastings algorithm. The computational efficiency of such an algorithm will be hurt when the expected rank of the data matrix increases. This potential inconvenience leaves room for improving the practical performance of the algorithm if a more efficient MCMC sampler, such as a Hamiltonian Monte Carlo sampler, can be designed. We defer these interesting extensions to future research directions.

\bigskip
\begin{center}
{\large\bf SUPPLEMENTARY MATERIAL}
\end{center}

% \begin{description}

The supplementary material contains the proofs of the theoretical results in Section \ref{sec:main_results}, the detailed Metropolis-Hastings algorithm, and the MCMC convergence diagnostics for the numerical examples in Section \ref{sec:numerical_examples}. 

% \end{description}
\section*{Acknowledgements}
This research was supported in part by Lilly Endowment, Inc., through its support for the Indiana University Pervasive Technology Institute.

\clearpage

\begin{center}
  \begin{Large}
    \textbf{Supplementary Material for ``Eigenvector-Assisted Statistical Inference for Signal-Plus-Noise Matrix Models''}
  \end{Large}
\end{center}
\appendix

% \counterwithin{lemma}{section}
% \counterwithin{theorem}{section}
\allowdisplaybreaks

\begin{abstract}
This supplementary material file contains the proofs of the results in Section \ref{sec:main_results} of the manuscript and additional computational details, including the detailed Metropolis-Hastings algorithm and the convergence diagnostics of the numerical results in Section \ref{sec:numerical_examples} of the manuscript.
\end{abstract}

% \end{keyword}

% subsection notations (end)

% section introduction (end)

\section{Technical preparations}
\label{sec:technical_preparations}

\subsection{Large sample properties of the spectral emedding}
\label{sub:convergence_eigenvector}
% Under Assumption \ref{assumption:signal_plus_noise}, we are able to establish the large sample behavior of the rows of the spectral embedding in Theorem \ref{thm:uniform_concentration_eigenvector} below.
We begin by extending the results in \cite{cape2019signal} for the unscaled eigenvectors to the spectral embedding (scaled eigenvectors) in Theorem \ref{thm:uniform_concentration_eigenvector} below. 
\begin{theorem}
% [Uniform concentration bound for eigenvectors]
\label{thm:uniform_concentration_eigenvector}
Suppose Assumption \ref{assumption:signal_plus_noise} hold. 
% Suppose $\bW_1,\bW_2\in\mathbb{O}(d)$ are the left and right singular vectors of $\bU_\bP\transpose\bU_\bA$, and $\bV\in\mathbb{O}(d)$ is the right singular vector of $\bX_0$. Let $\bW^* = \bW_1\bW_2\transpose$ and $\bW = (\bW^*)\transpose\bV\transpose$. 
Let $\bW\in\mathbb{O}(d)$ be the orthogonal alignment matrix between $\widetilde{\bX}$ and $\rho_n^{1/2}\bX_0$. 
Then 
% there exists an $n$-dependent random matrix $\bW\in\mathbb{O}(d)$, 
% such that for any $c > 0$, there exists constants $K_c > 0, N_c\in\mathbb{N}_+$, such that 
% for all $n\geq N_c$, with probability at least $1 - n^{-c}$,
% \begin{align*}
$\|\widetilde\bX\bW - \rho_n^{1/2}\bX_0\|_{2\to\infty}\lesssim \sqrt{(\log n)^{2\xi}/n}$ w.h.p. and
\[
\widetilde\bX\bW - \rho_n^{1/2}\bX_0 = \rho_n^{-1/2}(\bA - \bP)\bX_0(\bX_0\transpose\bX_0)^{-1} + \bR_\bX,
\]
where $\|\bR_\bX\|_{2\to\infty}\lesssim (\log n)^{2\xi}/(n\rho_n^{1/2})$ w.h.p..

% \end{align*}
\end{theorem}

\begin{proof}[\bf Proof of Theorem \ref{thm:uniform_concentration_eigenvector}]
By Theorem 1 in \cite{cape2019signal}, we know that $\|\bU_\bA(\bW^*)\transpose - \bU_\bP\|_{2\to\infty}\lesssim (\log n)^\xi/(n\rho_n^{1/2})$ w.h.p.. Recall that $\bW=(\bW^*)\transpose\bW_\bX$. Now write
\begin{align*}
\|\widetilde{\bX}\bW - \rho_n^{1/2}\bX_0\|_{2\to\infty}
& = \|\bU_\bA\bS_\bA^{1/2}(\bW^*)\transpose  - \bU_\bP\bS_\bP^{1/2}\|_{2\to\infty}\\
&\leq \|\bU_\bA(\bS_\bA^{1/2}(\bW^*)\transpose - (\bW^*)\transpose\bS_\bP^{1/2})\|_{2\to\infty} + \|(\bU_\bA(\bW^*)\transpose - \bU_\bP)\bS_\bP^{1/2}\|_{2\to\infty}\\
&\leq (\|\bU_\bA(\bW^*)\transpose - \bU_\bP\|_{2\to\infty} + \|\bU_\bP\|_{2\to\infty})\|\bS_\bA^{1/2}(\bW^*)\transpose - (\bW^*)\transpose\bS_\bP^{1/2}\|_2\\
&\quad + \|\bU_\bA(\bW^*)\transpose - \bU_\bP\|_{2\to\infty}\|\bS_\bP^{1/2}\|_2.
\end{align*}
Since $\|\bU_\bA(\bW^*)\transpose - \bU_\bP\|_{2\to\infty}\lesssim (\log n)^\xi/(n\rho_n^{1/2})$ w.h.p., $\|\bU_\bP\|_{2\to\infty}\lesssim n^{-1/2}$, and $\|\bS_\bP^{1/2}\|_2\lesssim (n\rho_n)^{1/2}$ by Assumption \ref{assumption:signal_plus_noise}, then for the first assertion, it is sufficient to show that $\|\bS_\bA^{1/2}(\bW^*)\transpose - (\bW^*)\transpose\bS_\bP^{1/2}\|_2\lesssim (\log n)^{1/2}/(n\rho_n)^{1/2}$ w.h.p.. Following the derivation of equation (49) in \cite{JMLR:v18:17-448}, we have
\begin{align*}
\|\bW^*\bS_\bA - \bS_\bP\bW^*\|_2
&\leq \|\bW^* - \bU_\bP\transpose\bU_\bA\|_2(\|\bS_\bA\|_2 + \|\bS_\bP\|_2)\\
&\quad + \|\bU_\bP\transpose(\bA - \bP)\|_2\|\bU_\bA - \bU_\bP\bU_\bP\transpose\bU_\bA\|_2\\
&\quad + \|\bU_\bP\transpose(\bA - \bP)\bU_\bP\|_2\\
&\leq \|\sin\Theta(\bU_\bA,\bU_\bP)\|_2^2(\|\bS_\bA\|_2 + \|\bS_\bP\|_2)\\
&\quad + \|\bA - \bP\|_2\|\sin\Theta(\bU_\bA,\bU_\bP)\|_2\\
&\quad + \|\bU_\bP\transpose(\bA - \bP)\bU_\bP\|_2.
\end{align*}
Observe that the $(k, l)$th entry of $\bU_\bP\transpose(\bA - \bP)\bU_\bP$ can be written as $\sum_{i\leq j}\{A_{ij} - \expect_0(A_{ij})\}a_{ij}$, where the coefficients $a_{ij}$'s satisfy $\max_{i,j\in [n]}|a_{ij}|\lesssim 1/n$. 
Now we consider either one of the conditions hold in Assumption \ref{assumption:signal_plus_noise}(vi). If $|A_{ij}|$'s are bounded by $1$ almost surely, then Hoeffding's inequality and a union bound yield that $\|\bU_\bP\transpose(\bA - \bP)\bU_\bP\|_2\lesssim (\log n)^{1/2}$ w.h.p.. If $A_{ij}$'s are uniformly bounded in $\psi_2$-Orlicz norms, then by Proposition 5.10 in \cite{vershynin2010introduction}, we also have $\|\bU_\bP\transpose(\bA - \bP)\bU_\bP\|_2\lesssim (\log n)^{1/2}$ w.h.p.. Hence, we further obtain from Assumption \ref{assumption:signal_plus_noise} (ii), Assumption \ref{assumption:signal_plus_noise} (iii), Weyl's inequality, and Davis-Kahan theorem that
\begin{align*}
\|\bW^*\bS_\bA - \bS_\bP\bW^*\|_2
\lesssim \frac{n\rho_n}{n\rho_n} + \frac{(n\rho_n)^{1/2}}{n\rho_n}(n\rho_n)^{1/2} + (\log n)^{1/2}\lesssim (\log n)^{1/2}\quad\mbox{w.h.p..}
\end{align*}
To show the high probability bound for $\|\bS_\bA^{1/2}(\bW^*)\transpose - (\bW^*)\transpose\bS_\bP^{1/2}\|_2$, note that the the $(k, l)$th entry of the transpose of this matrix is
% \[
$w_{lk}^*\{\lambda_k(\bA) - \lambda_l(\bP)\}/\{\lambda_k(\bA)^{1/2} + \lambda_l(\bP)^{1/2}\}$,
% \]
where $w_{lk}^*$ is the $(l, k)$th entry of $\bW^*$. It follows directly that 
\[
\|\bS_\bA^{1/2}(\bW^*)\transpose - (\bW^*)\transpose\bS_\bP^{1/2}\|_{\mathrm{F}}
\leq \frac{1}{\lambda_d(\bP)^{1/2}}\|\bW^*\bS_\bA - \bS_\bP\bW^*\|_{\mathrm{F}}
\lesssim \frac{(\log n)^{1/2}}{(n\rho_n)^{1/2}}\quad\mbox{w.h.p..}
\]
This completes the proof of the first assertion. For the second assertion, by Theorem 2 in \cite{cape2019signal}, we have
% \[
$\bU_\bA - \bU_\bP\bW^* = \bE\bU_\bP\bS_\bP^{-1}\bW^* + \bR_\bU$,
% \]
where $\bE = \bA - \bP$ and $\|\bR_\bU\|_{2\to\infty}\lesssim (\log n)^{2\xi}/(n^{3/2}\rho_n)$ w.h.p.. Now write
\begin{align*}
\widetilde{\bX}\bW - \rho_n^{1/2}\bX_0 & = (\bU_\bA\bS_\bA^{1/2}(\bW^*)\transpose - \bU_\bP\bS_\bP^{1/2})\bW_\bX\\
& = \{\bU_\bA(\bS_\bA^{1/2}(\bW^*)\transpose - (\bW^*)\transpose\bS_\bP^{1/2}) + (\bU_\bA - \bU_\bP\bW^*)(\bW^*)\transpose\bS_\bP^{1/2}\}\bW_\bX\\
& = \bU_\bA(\bS_\bA^{1/2}(\bW^*)\transpose - (\bW^*)\transpose\bS_\bP^{1/2})\bW_\bX + (\bE\bU_\bP\bS_\bP^{-1/2}\bW_\bX + \bR_\bU(\bW^*)\transpose\bS_\bP^{1/2}\bW_\bX)\\
& = \rho_n^{-1/2}\bE\bX_0(\bX_0\transpose\bX_0)^{-1} + \bR_\bX,
\end{align*}
where 
\[
\bR_\bX = \bR_\bU(\bW^*)\transpose\bS_\bP^{1/2}\bW_\bX + \bU_\bA(\bS_\bA^{1/2}(\bW^*)\transpose - (\bW^*)\transpose\bS_\bP^{1/2})\bW_\bX.
\]
Note that $\rho_n^{-1/2}\bX_0(\bX_0\transpose\bX_0)^{-1} = \bU_\bP\bS_\bP^{-1/2}\bW_\bX$. We have already shown that 
\begin{align*}
\|\bU_\bA(\bS_\bA^{1/2}(\bW^*)\transpose - (\bW^*)\transpose\bS_\bP^{1/2})\bW_\bX\|_{2\to\infty}
&\leq (\|\bU_\bA - \bU_\bP\bW^*\|_{2\to\infty} + \|\bU_\bP\|_{2\to\infty})\|\bW^*\bS_\bA^{1/2} - \bS_\bP^{1/2}\bW^*\|_2
\\&
\lesssim \frac{(\log n)^{1/2}}{n\rho_n^{1/2}}\quad\mbox{w.h.p.}.
\end{align*}
It follows from the earlier derived high probability bounds that $\|\bR_\bX\|_{2\to\infty}\lesssim (\log n)^{2\xi}/(n\rho_n^{1/2})$ w.h.p.. This completes the proof of the second assertion. 
\end{proof}

\subsection{Some preliminary results}

\begin{result}
[Concentration of $\|\bA\|_\infty$]
\label{result:concentration_of_infinity_norm}
Under Assumption \ref{assumption:signal_plus_noise} and Assumption \ref{assumption:regularity_condition}, 
$\|\bA\|_\infty \lesssim n\rho_n$ w.h.p.. To see this, observe that
\[
\|\bA\|_\infty
 \leq \max_{i\in [n]}\sum_{j = 1}^n\{|E_{ij}| - \expect_0(|E_{ij}|)\} + \max_{i\in [n]}\sum_{j = 1}^n \expect_0(|E_{ij}|)+ \max_{i\in [n]}\sum_{j = 1}^n|\expect_0(A_{ij})|.
\] 
The third term is deterministic and is $O(n\rho_n)$ since $\Theta$ is compact for all $\bx_{0i}$. The second term is also deterministic and can be bounded by
\[
\max_{i\in [n]}\sum_{j = 1}^n \expect_0(|E_{ij}|) = \max_{i,j\in [n]}n\expect_0(|A_{ij} - \rho_n\bx_{0i}\transpose\bx_{0j}|)\lesssim n\rho_n
\]
under Assumption \ref{assumption:signal_plus_noise}(vi). For the first term, under Assumption \ref{assumption:signal_plus_noise}(vi)(a), we have
% \[
$\var_0(|E_{ij}|)\leq \expect_0(E_{ij}^2)\lesssim \rho_n$. It follows from Bernstein's inequality and a union bound over $i\in [n]$ that the first term is $O(n\rho_n)$ w.h.p.. Under Assumption \ref{assumption:signal_plus_noise} (vi)(b), we have $\||E_{ij}| - \expect_0(|E_{ij}|)\|_{\psi_2}\leq \|E_{ij}\|_{\psi_2} + \expect_0(|E_{ij}|)\lesssim \rho_n^{1/2}$. Therefore, by Proposition 5.10 in \cite{vershynin2010introduction} and a union bound over $i\in [n]$, the first term is also $O(n\rho_n)$ w.h.p.. Hence we conclude that $\|\bA\|_\infty \lesssim n\rho_n$ w.h.p..
% \]
\end{result}

\begin{result}[Uniform concentration of $\widetilde{\bX}$]
\label{result:pij_tilde_concentration}
Suppose Assumption \ref{assumption:signal_plus_noise}  and Assumption \ref{assumption:regularity_condition} holds. 
There exists a constant $\delta > 0$, such that for all $i,j\in [n]$, $\widetilde\bx_i\transpose\widetilde\bx_j\in B(\rho_n\bx_{0i}\transpose\bx_{0j}, \rho_n\delta)$, $\max_{j\in[n]}\|\widetilde{\bx}_j\|_2\lesssim \rho_n^{1/2}$, and $\max_{i,j\in [n]}|\widetilde\bx_i\transpose\widetilde\bx_j - \rho_n\bx_{0i}\transpose\bx_{0j}|\lesssim \rho_n^{1/2}(\log n)^{\xi}/n^{1/2}$ w.h.p.. To see these results, note that by Theorem \ref{thm:uniform_concentration_eigenvector}, Assumption \ref{assumption:signal_plus_noise} (i), and the fact that $\Theta$ is bounded, 
\[
\max_{j\in [n]}\|\widetilde{\bx}_j\|_2 \leq \|\widetilde\bX\bW - \rho_n^{1/2}\bX_0\|_{2\to\infty} + \rho_n^{1/2}\|\bX_0\|_{2\to\infty}\lesssim \frac{(\log n)^\xi}{\sqrt{n}} + \rho_n^{1/2} \lesssim \rho_n^{1/2}\quad\mbox{w.h.p..}
\]
Also, 
\begin{equation}
% \label{eqn:pij_tilde_concentration}
\begin{aligned}
\max_{i,j\in [n]}|\widetilde\bx_i\transpose\widetilde\bx_j - \rho_n\bx_{0i}\transpose\bx_{0j}|
&\leq \max_{i,j\in [n]}(\|\bW\transpose\widetilde\bx_i\|_2 + \rho_n^{1/2}\|\bx_{0j}\|_2)\|\widetilde{\bX}\bW - \rho_n^{1/2}\bX_{0}\|_{2\to\infty}\\
&\lesssim \rho_n^{1/2}\sqrt{\frac{(\log n)^{2\xi}}{n}}\quad\mbox{w.h.p..}
\end{aligned}
\nonumber
\end{equation}
Since $(\log n)^{4\xi}/(n\rho_n)\to 0$, it follows that $\widetilde\bx_i\transpose\widetilde\bx_j\in B(\rho_n\bx_{0i}\transpose\bx_{0j}, \rho_n\delta)$ for all $i,j\in [n]$ w.h.p..
\end{result}

\begin{result}[Bernstein-type concentration of $\bE^2$]
\label{result:Berstein_concentration_E_square}
Suppose Assumption \ref{assumption:signal_plus_noise}  and Assumption \ref{assumption:regularity_condition} holds. 
For each $n$, let $(\alpha_{nijl}:i,j,l\in [n])$ be a three-dimensional array of real numbers such that $\max_{i,j,l}|\alpha_{nijl}|\lesssim 1/n$. Then for any $p = 1,\ldots,\lceil(\log n)^\xi/2\rceil$, 
\[
\expect_0\left(\left|\sum_{j = 1}^n\sum_{l = 1}^n\alpha_{nijl}E_{ij}E_{jl}\right|^p\right)\leq (n\rho_n)^{p}(4\sigma p)^{2p}\max_{i,j,l\in [n]}|\alpha_{nijl}|^p.
\]
Furthermore, with $p = \lfloor(\log n)^{\xi}/(8\sigma)\rfloor$, we obtain by a higher-order Markov's inequality that
\begin{align*}
\prob_0\left\{\left|\sum_{j = 1}^n\sum_{l = 1}^n\alpha_{nijl}E_{ij}E_{jl}\right| > (n\rho_n)(\log n)^{2\xi}\max_{i,j,l\in [n]}|\alpha_{nijl}|\right\}\leq \left\{\frac{4\sigma p}{(\log n)^{\xi}}\right\}^{2p}\leq e^{-\nu(\log n)^\xi}
\end{align*}
for some constant $\nu > 0$. 
The proof is similar to those of Lemma 5.4 in \cite{doi:10.1080/01621459.2020.1751645}, Lemma 7.10 in \cite{erdos2013}, and Lemma B.1 in \cite{xie2019efficient}. Denote $H_{ij} = E_{ij}/(\sigma^2 n\rho_n)^{1/2}
% := (A_{ij} - \rho_n\bx_{0i}\transpose\bx_{0j})/(\sigma^2 n\rho_n)^{1/2}
$. To adapt the proofs there under Assumption \ref{assumption:signal_plus_noise} (vi), it is sufficient to show that $\expect_0(|H_{ij}^m|)\leq \sigma^2/n$ for all $2\leq m\leq \lceil(\log n)^{\xi}\rceil$. 
\begin{itemize}
  \item Under Assumption \ref{assumption:signal_plus_noise} (vi) (a), we have, $|H_{ij}|\leq |A_{ij} - \expect_0(A_{ij})|/(\sigma^2 n\rho_n)^{1/2}\leq 1$ because $n\rho_n\to\infty$. Therefore
  \[
  \expect_0(|H_{ij}|^m)\leq \expect_0(H_{ij}^2) = \frac{1}{n\rho_n\sigma^2}\var_0(A_{ij})\leq \frac{1}{n}.
  \]
  \item Under Assumption \ref{assumption:signal_plus_noise} (vi) (b), we have,
  \begin{align*}
  \expect_0(|H_{ij}|^m) = \frac{m^{m/2}}{(\sigma^2n\rho_n)^{m/2}}\left\{\frac{1}{\sqrt{m}}\expect_0(|E_{ij}|^m)^{1/m}\right\}^m\leq\frac{m^{m/2}}{(\sigma^2n\rho_n)^{m/2}}\|A_{ij}\|_{\psi_2}^m\leq \left(\frac{m}{n}\right)^{m/2}. 
  \end{align*}
  Because $m\leq n$, we obtain directly that $\expect_0(|H_{ij}|^m)\leq 1$ for all $2\leq m\leq \lceil(\log n)^\xi\rceil$. 
\end{itemize}
\end{result}

\begin{result}[Uniform concentration of $\widetilde{\bg}_{ij}$ and $\bg_{ij}$]
\label{result:uniform_concentration_g}
Suppose Assumptions \ref{assumption:signal_plus_noise}-\ref{assumption:weight_functions} hold. 
For each $i\in [n]$, 
\[
\max_{j\in [n]}\sup_{\bx_i\in\Theta}\{\|\widetilde{\bg}_{ij}(\bx_i)\|_2 + \|\bg_{ij}(\bx_i)\|_2\}\lesssim
\left\{
\begin{aligned}
& 1, \quad\mbox{w.h.p., if Assumption \ref{assumption:signal_plus_noise} (vi)(a) holds,}\\
& (\rho_n\log n)^{1/2}, \quad\mbox{w.h.p., if Assumption \ref{assumption:signal_plus_noise} (vi)(b) holds.}
\end{aligned}
\right.
\]
% for any sequence $(M_n)_{n = 1}^\infty$ with $M_n\to\infty$. 
To see why this holds, we first write 
\begin{align*}
&\widetilde{\bg}_{ij}(\bx_i) = A_{ij}h_n(\widetilde{\bx}_i\transpose\widetilde{\bx}_j, \bx_i\transpose\widetilde{\bx}_j)\rho_n^{-1/2}\widetilde{\bx}_j - \bx_i\transpose\widetilde{\bx}_jh_n(\widetilde{\bx}_i\transpose\widetilde{\bx}_j, \bx_i\transpose\widetilde{\bx}_j)\rho_n^{-1/2}\widetilde{\bx}_j,\\
&{\bg}_{ij}(\bx_i) = A_{ij}h_n(\rho_n\bx_{0i}\transpose\bx_{0j}, \rho_n^{1/2}\bx_i\transpose{\bx}_{0j})\bx_{0j} - \rho_n^{1/2}\bx_i\transpose\bx_{0j}h_n(\rho_n\bx_{0i}\transpose\bx_{0j}, \rho_n^{1/2}\bx_i\transpose{\bx}_{0j})\bx_{0j}
.
\end{align*}
By Result \ref{result:pij_tilde_concentration} and Assumption \ref{assumption:weight_functions}, 
\begin{align*}
&\max_{j\in [n]}\sup_{\bx_i\in\Theta}\|h_n(\rho_n\bx_{0i}\transpose\bx_{0j}, \rho_n^{1/2}\bx_i\transpose{\bx}_{0j})\bx_{0j}\|_2\lesssim 1,\\
&\max_{j\in [n]}\sup_{\bx_i\in\Theta}\|\rho_n^{1/2}\bx_i\transpose{\bx}_{0j}h_n(\rho_n\bx_{0i}\transpose\bx_{0j}, \rho_n^{1/2}\bx_i\transpose{\bx}_{0j})\bx_{0j}\|_2\lesssim \rho_n^{1/2},\\
&\max_{j\in [n]}\sup_{\bx_i\in\Theta}\|h_n(\widetilde{\bx}_i\transpose\widetilde{\bx}_j, \bx_i\transpose\widetilde{\bx}_j)\rho_n^{-1/2}\widetilde{\bx}_j\|_2\lesssim 1\quad\mbox{w.h.p.},\\
&\max_{j\in [n]}\sup_{\bx_i\in\Theta}\|\bx_i\transpose\widetilde{\bx}_jh_n(\widetilde{\bx}_i\transpose\widetilde{\bx}_j, \bx_i\transpose\widetilde{\bx}_j)\rho_n^{-1/2}\widetilde{\bx}_j\|_2\lesssim \rho_n^{1/2}\quad\mbox{w.h.p.}.
\end{align*}
By Assumption \ref{assumption:signal_plus_noise} (vi), Lemma 8.1 in \cite{kosorok2008introduction}, and a union bound over $j\in [n]$, we have
\[
\max_{j\in [n]}|A_{ij}|\lesssim
\left\{
\begin{aligned}
& 1, \quad\mbox{w.p.1., if Assumption \ref{assumption:signal_plus_noise} (vi)(a) holds,}\\
& (\rho_n\log n)^{1/2}, \quad\mbox{w.h.p., if Assumption \ref{assumption:signal_plus_noise} (vi)(b) holds.}
\end{aligned}
\right.
\]
Then the proof of Result \ref{result:uniform_concentration_g} is completed by combining the above high probability bounds. 
\end{result}

\begin{result}[Identifiability]
\label{result:identifiability}
Suppose Assumptions \ref{assumption:signal_plus_noise}, \ref{assumption:regularity_condition}, and \ref{assumption:weight_functions} holds. Then for each $i\in [n]$, $\bx_i = \rho_n^{1/2}\bx_{0i}$ is the unique solution to the population estimating equation
\[
\frac{1}{n}\sum_{j = 1}^n\expect_0\{\bg_{ij}(\bx_i)\} = \zero_d.
\]
Furthermore, there exists a constant $\delta > 0$, such that for any $\eps > 0$,
\[
\inf_{\|\bx_i - \rho_n^{1/2}\bx_{0i}\|_2 > \eps}\left\|\frac{1}{n}\sum_{j = 1}^n\expect_0\{\bg_{ij}(\bx_i)\}\right\|_2 > \rho_n^{1/2}\delta_0\eps. 
\]
Now we show this. Denote $h_{0nij}(\bx_i)=h_n(\rho_n\bx_{0i}\transpose\bx_{0j},\rho_n^{1/2}\bx_{i}\transpose\bx_{0j})$. Then
\begin{align*}
\frac{1}{n}\sum_{j=1}^n\expect_0\{\bg_{ij}(\bx_i)\} &= \frac{1}{n}\sum_{j=1}^n\expect_0\{(A_{ij}-\rho_n^{1/2}\bx_{i}\transpose\bx_{0j})h_{0nij}(\bx_i)\bx_{0j}\}\\
&= \frac{1}{n}\sum_{j=1}^n(\rho_n^{1/2}\bx_{0i}-\bx_i)\transpose\rho_n^{1/2}\bx_{0j}h_{0nij}(\bx_i)\bx_{0j}\\
&= \frac{1}{n} \rho_n^{1/2}\bX_0\transpose\mathrm{diag}\{h_{0ni1}(\bx_i),\ldots,h_{0nin}(\bx_i)\}\bX_0(\rho_n^{1/2}\bx_{0i}-\bx_i).
\end{align*}
By Assumption \ref{assumption:signal_plus_noise} (ii) and Assumption \ref{assumption:weight_functions},
\begin{align*}
\left\|\frac{1}{n}\sum_{j=1}^n\expect_0\{\bg_{ij}(\bx_i)\}\right\|_2
&= \left\|\frac{1}{n} \rho_n^{1/2}\bX_0\transpose\mathrm{diag}\{h_{0ni1}(\bx_i),\ldots,h_{0nin}(\bx_i)\}\bX_0(\rho_n^{1/2}\bx_{0i}-\bx_i)\right\|_2\\
&\geq \lambda_d\left\{\bX_0\transpose\mathrm{diag}\{h_{0ni1}(\bx_i),\ldots,h_{0nin}(\bx_i)\}\bX_0\right\}\frac{1}{n}\rho_n^{1/2}\|\rho_n^{1/2}\bx_{0i}-\bx_i\|_2\\
&\geq \rho_n^{1/2}\delta_0\|\rho_n^{1/2}\bx_{0i}-\bx_i\|_2
\end{align*}
for some constant $\delta_0>0$. So we have
$\inf_{\|\rho_n^{1/2}\bx_{0i}-\bx_i\|_2>\epsilon}\|{1}/{n}\sum_{j=1}^n\expect_0\{\bg_{ij}(\bx_i)\}\|_2 > \rho_n^{1/2}\delta_0\epsilon$, and that ${1}/{n}\sum_{j=1}^n\expect_0\{\bg_{ij}(\bx_i)\}=0$ implies $\bx_i=\rho_n^{1/2}\bx_{0i}$.
\end{result}

\begin{result}[Jacobian]
\label{result:Jacobian}
Suppose Assumptions \ref{assumption:signal_plus_noise}, \ref{assumption:regularity_condition}, and \ref{assumption:weight_functions} holds. Then there exists constants $\delta, c, C > 0$, such that for each $i\in [n]$, the matrix
\[
\bG_{in}(\bx_i) = \frac{1}{n}\sum_{j = 1}^n\expect_0\left\{\frac{\partial\bg_{ij}}{\partial\bx_i\transpose}(\bx_i)\right\}
\]
satisfies
\[
c\rho_n^{1/2}\leq \inf_{\|\bx_i - \rho_n^{1/2}\bx_{0i}\|_2\leq\delta}\lambda_d\{-\bG_{in}(\bx_i)\}\leq \sup_{\|\bx_i - \rho_n^{1/2}\bx_{0i}\|_2\leq\delta}\lambda_1\{-\bG_{in}(\bx_i)\}\leq C\rho_n^{1/2}.
\]
Now we show this. Denote $h_{0nij}(\bx_i)=h_n(\rho_n\bx_{0i}\transpose\bx_{0j},\rho_n^{1/2}\bx_{i}\transpose\bx_{0j})$. Then
\begin{align*}
\frac{1}{n}\sum_{j=1}^n\expect_0\left\{\frac{\partial\bg_{ij}}{\partial\bx_i}(\bx_i)\right\}
&= \frac{1}{n}\sum_{j=1}^n\expect_0\{(A_{ij}-\rho_n^{1/2}\bx_{i}\transpose\bx_{0j})D^{(0,1)}h_{0nij}(\bx_i)\rho_n^{1/2}\bx_{0j}\bx_{0j}\transpose - h_{0nij}(\bx_i)\rho_n^{1/2}\bx_{0j}\bx_{0j}\transpose\}\\
&= \frac{1}{n}\sum_{j=1}^n [(\rho_n^{1/2}\bx_{0i} - \bx_i)\transpose\bx_{0j}D^{(0,1)}h_{0nij}(\bx_i) - h_{0nij}(\bx_i)]\rho_n^{1/2}\bx_{0j}\bx_{0j}\transpose\\
&= -\frac{1}{n}\rho_n^{1/2}\bX_0\transpose\mathrm{diag}\left\{h_{0nij}(\bx_i) - (\rho_n^{1/2}\bx_{0i}-\bx_i)\transpose\bx_{0j}D^{(0,1)}h_{0nij}(\bx_i)\right\}_{j=1}^n\bX_0.
\end{align*}
By Assumptions \ref{assumption:regularity_condition} and \ref{assumption:weight_functions}, there exist constants $\delta, c', C' > 0$, such that
\begin{align*}
c'&\leq \inf_{\|\bx_i-\rho_n^{1/2}\bx_{0i}\|_2\leq\delta}\left\{h_{0nij}(\bx_i) - (\rho_n^{1/2}\bx_{0i}-\bx_i)\transpose\bx_{0j}D^{(0,1)}h_{0nij}(\bx_i)\right\}\\
&\leq \sup_{\|\bx_i-\rho_n^{1/2}\bx_{0i}\|_2\leq\delta}\left\{h_{0nij}(\bx_i) - (\rho_n^{1/2}\bx_{0i}-\bx_i)\transpose\bx_{0j}D^{(0,1)}h_{0nij}(\bx_i)\right\}
\leq C'.
\end{align*}
Then Result \ref{result:Jacobian} follows directly.
\end{result}

\begin{result}[Second-moment matrix]
\label{result:second_moment_matrix}
Suppose Assumptions \ref{assumption:signal_plus_noise}, \ref{assumption:regularity_condition}, and \ref{assumption:weight_functions} holds. Then there exists constants $c, C > 0$, such that for each $i\in [n]$, the matrix
\[
\bOmega_{in}(\bx_i) = \frac{1}{n}\sum_{j = 1}^n\expect_0\left\{ \bg_{ij}(\bx_i)\bg_{ij}(\bx_i)\transpose\right\}
\]
satisfies
\[
c\rho_n\leq \inf_{\bx_i\in\Theta}\lambda_d\{\bOmega_{in}(\bx_i)\}\leq \sup_{\bx_i\in\Theta}\lambda_1\{\bOmega_{in}(\bx_i)\}\leq C\rho_n.
\]
Now we show this. Denote $h_{0nij}(\bx_i)=h_n(\rho_n\bx_{0i}\transpose\bx_{0j},\rho_n^{1/2}\bx_{i}\transpose\bx_{0j})$. Then
\begin{align*}
\frac{1}{n}\sum_{j=1}^n\expect_0\left\{\bg_{ij}(\bx_i)\bg_{ij}(\bx_i)\transpose\right\}
&= \frac{1}{n}\sum_{j=1}^n\expect_0\left\{(A_{ij}-\rho_n^{1/2}\bx_i\transpose\bx_{0j})^2 h_{0nij}(\bx_i)^2\bx_{0j}\bx_{0j}\transpose\right\}\\
&= \frac{1}{n}\sum_{j=1}^n\left\{\var_0(A_{ij}) + \rho_n[(\rho_n^{1/2}\bx_{0i}-\bx_i)\transpose\bx_{0j}]^2\right\}h_{0nij}(\bx_i)^2\bx_{0j}\bx_{0j}\transpose.
\end{align*}
By Assumption \ref{assumption:signal_plus_noise} (vi), $\sigma_0^2\rho_n\leq \var_0(A_{ij}) \leq \sigma^2\rho_n$. By Assumption \ref{assumption:regularity_condition}, $0\leq \inf_{\bx_i\in\Theta}\rho_n[(\rho_n^{1/2}\bx_{0i}-\bx_i)\transpose\bx_{0j}]^2 \leq\sup_{\bx_i\in\Theta}\rho_n[(\rho_n^{1/2}\bx_{0i}-\bx_i)\transpose\bx_{0j}]^2 \leq C\rho_n$.
Then by Assumptions \ref{assumption:signal_plus_noise} (ii) and \ref{assumption:weight_functions}, there exist constants $c, C > 0$, such that
\begin{align*}
c\rho_n &\leq c_1\sigma^2\rho_n\lambda_d\left\{\frac{1}{n}\sum_{j=1}^n\bx_{0j}\bx_{0j}\transpose\right\} \leq \inf_{\bx_i\in\Theta}\lambda_d\left\{\bOmega_{in}(\bx_i)\right\} \leq \sup_{\bx_i\in\Theta}\lambda_1\left\{\bOmega_{in}(\bx_i)\right\}\\
&\leq c_2\sigma^2\rho_n\lambda_1\left\{\frac{1}{n}\sum_{j=1}^n\bx_{0j}\bx_{0j}\transpose\right\} \leq C\rho_n.
\end{align*}
\end{result}

\subsection{Law of Large Numbers}
\label{sub:LLN}

\begin{lemma}[Law of Large Numbers]
\label{lemma:LLN}
Suppose Assumptions \ref{assumption:signal_plus_noise}, \ref{assumption:regularity_condition}, and \ref{assumption:weight_functions} hold. Then for all $i\in [n]$,
\begin{align*}
&\left\|\frac{1}{n}\sum_{j = 1}^n\bW\transpose\frac{\partial\widetilde{\bg}_{ij}}{\partial\bx_i\transpose}(\rho_n^{1/2}\bW\bx_{0i})\bW - \bG_{in}(\rho_n^{1/2}\bx_{0i})\right\|_2 \lesssim \sqrt{\frac{(\log n)^{2\xi}}{n}}\quad\mbox{w.h.p.},\\
&\left\|\frac{1}{n}\sum_{j = 1}^n \bW\transpose\widetilde{\bg}_{ij}(\widetilde\bx_i)\widetilde{\bg}_{ij}(\widetilde\bx_i)\transpose\bW - \bOmega_{in}(\rho_n^{1/2}\bx_{0i})\right\|_2\lesssim \rho_n^{1/2}\sqrt{\frac{(\log n)^{2\xi}}{n}}\quad\mbox{w.h.p.},
\end{align*}
where $\bW\in\mathbb{O}(d)$ is the orthogonal alignment matrix between $\widetilde{\bX}$ and $\rho_n^{1/2}\bX_0$.
\end{lemma}
\begin{proof}[\bf Proof of Lemma \ref{lemma:LLN}]
$\blacksquare$ {\bf Proof of the first assertion.} First compute the Jacobian
\begin{align*}
\bW\transpose\frac{\partial\widetilde{\bg}_{ij}}{\partial\bx_i\transpose}(\bx_i)\bW
& = (A_{ij} - \bx_i\transpose\widetilde{\bx}_j)D^{(0, 1)}h_n(\widetilde{\bx}_i\transpose\widetilde{\bx}_j, \bx_i\transpose\widetilde{\bx}_j)\rho_n^{-1/2}\bW\transpose\widetilde{\bx}_j\widetilde{\bx}_j\transpose\bW\\
&\quad -h_n(\widetilde{\bx}_i\transpose\widetilde{\bx}_j, \bx_i\transpose\widetilde{\bx}_j)\rho_n^{-1/2}\bW\transpose\widetilde{\bx}_j\widetilde{\bx}_j\transpose\bW
 .
\end{align*}
Denote $E_{ij} = A_{ij} - \expect_0(A_{ij})$, $\widetilde{h}_{0nij} = h_n(\widetilde{\bx}_i\transpose\widetilde{\bx}_j, \rho_n^{1/2}\bx_{0i}\transpose\bW\transpose\widetilde{\bx}_j)$, and $h_{0nij} = h_n(\rho_n\bx_{0i}\transpose\bx_{0j}, \rho_n\bx_{0i}\transpose{\bx}_{0j})$. With a slight abuse of notations, we also denote $D^{(0, 1)}\widetilde{h}_{0nij} = D^{(0, 1)}h_n(\widetilde{\bx}_i\transpose\widetilde{\bx}_j, \rho_n^{1/2}\bx_{0i}\transpose\bW\transpose\widetilde{\bx}_j)$ and $D^{(0, 1)}{h}_{0nij} = D^{(0, 1)}h_n(\rho_n\bx_{0i}\transpose\bx_{0j}, \rho_n\bx_{0i}\transpose\bx_{0j})$. Then by triangle inequality and Cauchy-Schwarz inequality,
\begin{align*}
&\left\|\frac{1}{n}\sum_{j = 1}^n\bW\transpose\frac{\partial\widetilde{\bg}_{ij}}{\partial\bx_i\transpose}(\rho_n^{1/2}\bW\bx_{0i})\bW - \bG_{in}(\rho_n^{1/2}\bx_{0i})\right\|_2\\
&\quad\leq \left\|\frac{1}{n\rho_n^{1/2}}\sum_{j = 1}^nE_{ij} D^{(0, 1)}h_{0nij}\rho_n\bx_{0j}\bx_{0j}\transpose\right\|_2\\
&\quad\quad + \frac{1}{n\rho_n^{1/2}}\sum_{j = 1}^n|E_{ij}| \left\|D^{(0, 1)}\widetilde{h}_{0nij}\bW\transpose\widetilde{\bx}_j\widetilde{\bx}_j\transpose\bW - D^{(0, 1)}h_{0nij}\rho_n\bx_{0j}\bx_{0j}\transpose\right\|_2\\
&\quad\quad + \frac{1}{n\rho_n^{1/2}}\sum_{j = 1}^n\|\rho_n^{1/2}\bx_{0i}\|_2\|\bW\transpose\widetilde{\bx}_j - \rho_n^{1/2}\bx_{0j}\|_2 |D^{(0, 1)}\widetilde{h}_{0nij}| \|\widetilde{\bx}_j\|_2^2\\
&\quad\quad + \frac{1}{n}\sum_{j = 1}^n\left\|\widetilde{h}_{0nij}\rho_n^{-1/2}\bW\transpose\widetilde{\bx}_j\widetilde{\bx}_j\transpose\bW - h_{0nij}\rho_n^{1/2}\bx_{0j}\bx_{0j}\transpose\right\|_2.
\end{align*}
By Result \ref{result:pij_tilde_concentration}, Theorem \ref{thm:uniform_concentration_eigenvector}, and Assumption \ref{assumption:weight_functions}, the third term can be bounded as follows:
\begin{align*}
&\frac{1}{n\rho_n^{1/2}}\sum_{j = 1}^n\|\rho_n^{1/2}\bx_{0i}\|_2\|\bW\transpose\widetilde{\bx}_j - \rho_n^{1/2}\bx_{0j}\|_2|D^{(0, 1)}\widetilde{h}_{0nij}|\|\widetilde{\bx}_j\|_2^2\\
&\quad\lesssim \frac{1}{n\rho_n^{1/2}}\rho_n^{1/2}n\|\widetilde\bX\bW - \rho_n^{1/2}\bX_0\|_{2\to\infty}\rho_n^2\lesssim \rho_n^2\sqrt{\frac{(\log n)^{2\xi}}{n}}\quad\mbox{w.h.p..}
\end{align*}
Since $\widetilde\bx_i\transpose\widetilde\bx_j\in B(\rho_n\bx_{0i}\transpose\bx_{0j}, \rho_n\delta)$ for all $i,j\in [n]$ w.h.p. by Result \ref{result:pij_tilde_concentration}, it follows from Assumption \ref{assumption:weight_functions} (Lipschitz continuity of $h$) and Cauchy-Schwarz inequality that
\begin{align}\label{eqn:uniform_Lipschitz_h}
\begin{aligned}
\max_{i,j\in [n]}|\widetilde{h}_{0nij} - h_{0nij}|&\lesssim \rho_n^{-1}\max_{i,j\in [n]}|\widetilde\bx_i\transpose\widetilde\bx_j - \rho_n\bx_{0i}\transpose\bx_{0j}| + \rho_n\max_{i,j\in [n]}\rho_n^{1/2}\|\bx_{0i}\|_2\|\bW\transpose\widetilde{\bx}_j - \rho_n^{1/2}\bx_{0j}\|_2\\
&\lesssim \sqrt{\frac{(\log n)^{2\xi}}{n\rho_n}}\quad\mbox{w.h.p..}
\end{aligned}
\end{align}
Hence, the fourth term can be bounded using a similar approach:
\begin{align*}
&\frac{1}{n}\sum_{j = 1}^n\left\|\widetilde{h}_{0nij}\rho_n^{-1/2}\bW\transpose\widetilde{\bx}_j\widetilde{\bx}_j\transpose\bW - h_{0nij}\rho_n^{1/2}\bx_{0j}\bx_{0j}\transpose\right\|_2\\
&\quad\leq \frac{1}{n\rho_n^{1/2}}\sum_{j = 1}^n|\widetilde{h}_{0nij} - h_{0nij}|\|\widetilde{\bx}_j\|_2^2 + \frac{1}{n\rho_n^{1/2}}\sum_{j = 1}^n|h_{0nij}|\|\bW\transpose\widetilde{\bx}_j - \rho_n^{1/2}\bx_{0j}\|_2(\|\bW\transpose\widetilde{\bx}_j\|_2 + \rho_n^{1/2}\|\bx_{0j}\|_2)\\
&\quad\lesssim \frac{1}{n\rho_n^{1/2}}n\sqrt{\frac{(\log n)^{2\xi}}{n\rho_n}}\rho_n + \frac{1}{n\rho_n^{1/2}}n\sqrt{\frac{(\log n)^{2\xi}}{n}}\rho_n^{1/2} \lesssim \sqrt{\frac{(\log n)^{2\xi}}{n}}\quad\mbox{w.h.p.}.
\end{align*}
For the first term, by Assumption \ref{assumption:weight_functions}, we know that $|D^{(0, 1)}h_{0nij}| = O(\rho_n)$ and $\max_{j\in [n]}\|\bx_{0j}\|_2^2 = O(1)$. Then by either Bernstein's inequality under Assumption \ref{assumption:signal_plus_noise} (vi) (a) or Proposition 5.16 in \cite{vershynin2010introduction} under Assumption \ref{assumption:signal_plus_noise} (vi) (b), we obtain
\begin{align*}
\left\|\frac{1}{n\rho_n^{1/2}}\sum_{j = 1}^nE_{ij} D^{(0, 1)}h_{0nij}\rho_n\bx_{0j}\bx_{0j}\transpose\right\|_2 \lesssim \rho_n^2\sqrt{\frac{(\log n)^{2\xi}}{n}}\quad\mbox{w.h.p..}
\end{align*}
For the second term, we first observe that by Assumption \ref{assumption:weight_functions} and Result \ref{result:pij_tilde_concentration},
\[
\max_{i,j\in [n]}\left\|D^{(0, 1)}\widetilde{h}_{0nij}\bW\transpose\widetilde{\bx}_j\widetilde{\bx}_j\transpose\bW - D^{(0, 1)}h_{0nij}\rho_n\bx_{0j}\bx_{0j}\transpose\right\|_2 \lesssim \rho_n^{1/2}\sqrt{\frac{(\log n)^{2\xi}}{n}}\quad\mbox{w.h.p..}
\]
Then by Cauchy-Schwarz inequality and Assumption \ref{assumption:signal_plus_noise} (ii),
\begin{align*}
&\frac{1}{n\rho_n^{1/2}}\sum_{j = 1}^n|E_{ij}| \left\|D^{(0, 1)}\widetilde{h}_{0nij}\bW\transpose\widetilde{\bx}_j\widetilde{\bx}_j\transpose\bW - D^{(0, 1)}h_{0nij}\rho_n\bx_{0j}\bx_{0j}\transpose\right\|_2\\
&\quad\leq \frac{1}{n\rho_n^{1/2}}\|\bE\|_{\infty}\max_{i,j\in [n]}\left\|D^{(0, 1)}\widetilde{h}_{0nij}\bW\transpose\widetilde{\bx}_j\widetilde{\bx}_j\transpose\bW - D^{(0, 1)}h_{0nij}\rho_n\bx_{0j}\bx_{0j}\transpose\right\|_2\\
&\quad\lesssim \frac{1}{n\rho_n^{1/2}} n\rho_n \rho_n^{1/2}\sqrt{\frac{(\log n)^{2\xi}}{n}}\\
&\quad = \rho_n\sqrt{\frac{(\log n)^{2\xi}}{n}}\quad\mbox{w.h.p..}. 
\end{align*}
This completes the proof of the first assertion that
\[
\left\|\frac{1}{n}\sum_{j = 1}^n\bW\transpose\frac{\partial\widetilde{\bg}_{ij}}{\partial\bx_i\transpose}(\rho_n^{1/2}\bW\bx_i)\bW - \bG_{in}(\rho_n^{1/2}\bx_{0i})\right\|_2 \lesssim \sqrt{\frac{(\log n)^{2\xi}}{n}}\quad\mbox{w.h.p..}
\]
$\blacksquare$ {\bf Proof of the second assertion.}
By triangle inequality, with $\widetilde{h}_{nij} = h_n(\widetilde{\bx}_i\transpose\widetilde{\bx}_j,\widetilde{\bx}_i\transpose\widetilde{\bx}_j)$,
\begin{align*}
&\left\|\frac{1}{n}\sum_{j = 1}^n \bW\transpose\widetilde{\bg}_{ij}(\widetilde\bx_i)\widetilde{\bg}_{ij}(\widetilde\bx_i)\transpose\bW - \bOmega_{in}(\rho_n^{1/2}\bx_{0i})\right\|_2\\
&\quad\leq \frac{1}{n\rho_n}\sum_{j = 1}^n2|A_{ij}||\widetilde{\bx}_i\transpose\widetilde{\bx}_j - \rho_n\bx_{0i}\transpose\bx_{0j}|\widetilde{h}_{nij}^2\left\|\widetilde{\bx}_j\right\|_2^2\\
&\quad\quad + \frac{1}{n\rho_n}\sum_{j = 1}^n|(\widetilde{\bx}_i\transpose\widetilde{\bx}_j)^2 - (\rho_n\bx_{0i}\transpose\bx_{0j})^2|\widetilde{h}_{nij}^2\left\|\widetilde{\bx}_j\right\|_2^2\\
&\quad\quad + \frac{1}{n}\sum_{j = 1}^nE_{ij}^2 \left\|\rho_n^{-1}\widetilde{h}_{nij}^2\bW\transpose\widetilde{\bx}_j\widetilde{\bx}_j\transpose\bW - h_{0nij}^2\bx_{0j}\bx_{0j}\transpose\right\|_2\\
&\quad\quad + \left\|\frac{1}{n}\sum_{j = 1}^n\{(A_{ij} - \rho_n\bx_{0i}\transpose\bx_{0j})^2 - \var_0(A_{ij})\}h_{0nij}^2\bx_{0j}\bx_{0j}\transpose\right\|_2.
\end{align*}
For the first term, we apply Result \ref{result:concentration_of_infinity_norm}, Assumption \ref{assumption:weight_functions}, and Result \ref{result:pij_tilde_concentration} to obtain
\begin{align*}
\frac{1}{n\rho_n}\sum_{j = 1}^n2|A_{ij}||\widetilde{\bx}_i\transpose\widetilde{\bx}_j - \rho_n\bx_{0i}\transpose\bx_{0j}|\widetilde{h}_{nij}^2\left\|\widetilde{\bx}_j\right\|_2^2 \lesssim \rho_n^{3/2}\sqrt{\frac{(\log n)^{2\xi}}{n}}\quad\mbox{w.h.p.}.
\end{align*}
Following the same reasoning, the second term is $O(\rho_n^{3/2}(\log n)^\xi/\sqrt{n})$ w.h.p.. For the third term, we apply Assumption \ref{assumption:signal_plus_noise} (iii), Assumption \ref{assumption:weight_functions}, and Result \ref{result:pij_tilde_concentration} to obtain
\begin{align*}
&\frac{1}{n}\sum_{j = 1}^nE_{ij}^2 \left\|\rho_n^{-1}\widetilde{h}_{nij}^2\bW\transpose\widetilde{\bx}_j\widetilde{\bx}_j\transpose\bW - h_{0nij}^2\bx_{0j}\bx_{0j}\transpose\right\|_2\\
&\quad\leq \frac{1}{n}\|\bE\|_{2\to\infty}^2\max_{i,j}\left\|\rho_n^{-1}\widetilde{h}_{nij}^2\bW\transpose\widetilde{\bx}_j\widetilde{\bx}_j\transpose\bW - h_{0nij}^2\bx_{0j}\bx_{0j}\transpose\right\|_2\\
&\quad\lesssim \frac{1}{n}(n\rho_n)\sqrt{\frac{(\log n)^{2\xi}}{n\rho_n}} = \rho_n^{1/2}\sqrt{\frac{(\log n)^{2\xi}}{n}}\quad\mbox{w.h.p.}.
\end{align*}
For the fourth term, we consider two scenarios under Assumption \ref{assumption:signal_plus_noise} (vi). 
Note that the entries of $h_{0nij}^2\bx_{0j}\bx_{0j}\transpose$ are uniformly bounded. Under Assumption \ref{assumption:signal_plus_noise} (vi) (a), 
% \[
$\var_0\{(A_{ij} - \rho_n\bx_{0i}\transpose\bx_{0j})^2\}\leq \expect_0\{(A_{ij} - \rho_n\bx_{0i}\transpose\bx_{0j})^2\}\lesssim \rho_n$.
% \]
Then by Bernstein's inequality, the fourth term is $O((\rho_n\log n)^{1/2}/\sqrt{n})$ w.h.p.. Under Assumption \ref{assumption:signal_plus_noise} (vi) (b), 
$\|(A_{ij} - \rho_n\bx_{0i}\transpose\bx_{0j})^2 - \var_0(A_{ij})\|_{\psi_1}\leq \|A_{ij} - \rho_n\bx_{0i}\transpose\bx_{0j}\|_{\psi_2}^2 + \var_0(A_{ij}) \lesssim \rho_n$. Then by Proposition 5.16 in \cite{vershynin2010introduction}, 
the fourth term is $O(\rho_n(\log n)^{1/2}/\sqrt{n})$ w.h.p.. The proof of the second assertion is thus completed. 
\end{proof}

\subsection{Uniform Law of Large Numbers}
\label{sub:ULLN}

\begin{lemma}[Uniform Law of Large Numbers]
\label{lemma:ULLN}
Suppose Assumptions \ref{assumption:signal_plus_noise}, \ref{assumption:regularity_condition}, and \ref{assumption:weight_functions} hold. Then for all $i\in [n]$,
\begin{align*}
&\sup_{\bx_i\in\Theta}\left\|\frac{1}{n}\sum_{j = 1}^n[\bW\transpose\widetilde{\bg}_{ij}(\bW\bx_i) - \expect_0\{\bg_{ij}(\bx_i)\}]\right\|\lesssim \sqrt{\frac{(\log n)^{2\xi}}{n}}\quad\mbox{w.h.p.},\\
&\sup_{\bx_i\in\Theta}\left\|\frac{1}{n}\sum_{j = 1}^n\bW\transpose\frac{\partial\widetilde{\bg}_{ij}}{\partial\bx_i\transpose}(\bW\bx_i)\bW - \bG_{in}(\bx_i)\right\|\lesssim \sqrt{\frac{(\log n)^{2\xi}}{n}}\quad\mbox{w.h.p.}.
\end{align*}
\end{lemma}

\begin{proof}[\bf Proof of Lemma \ref{lemma:ULLN}]
Denote $\widetilde{h}_{nij}(\bx_i) = h_n(\widetilde{\bx}_i\transpose\widetilde{\bx}_j, \bx_i\transpose\widetilde{\bx}_j)$ and $h_{0nij}(\bx_i) = h_n(\rho_n\bx_{0i}\transpose\bx_{0j}, \rho_n^{1/2}\bx_i\transpose{\bx}_{0j})$ for notational simplicity. 

\vspace*{2ex}\noindent
$\blacksquare$ \textbf{Proof of the first assertion.} By triangle inequality and Cauchy-Schwarz inequality, with $E_{ij} := A_{ij} - \expect_0(A_{ij})$,
\begin{align*}
&\sup_{\bx_i\in\Theta}\left\|\frac{1}{n}\sum_{j = 1}^n\{\bW\transpose\widetilde{\bg}_{ij}(\bW\bx_i) - \bg_{ij}(\bx_i)\}\right\|_2\\
&\quad\leq \sup_{\bx_i\in\Theta}
\frac{1}{n}\sum_{j = 1}^n|E_{ij}|\left\| \{\widetilde{h}_{nij}(\bW\bx_i)\rho_n^{-1/2}\bW\transpose\widetilde{\bx}_j - h_{0nij}(\bx_i)\bx_{0j}\}\right\|_2\\
&\quad\quad + \sup_{\bx_i\in\Theta}
\frac{1}{n}\sum_{j = 1}^n|\bx_i\transpose\bW\transpose\widetilde{\bx}_j - \rho_n\bx_{0i}\transpose\bx_{0j}|\left\|\widetilde{h}_{nij}(\bW\bx_i)\rho_n^{-1/2}\bW\transpose\widetilde{\bx}_j - h_{0nij}(\bx_i)\bx_{0j}\right\|_2\\
&\quad\quad + \sup_{\bx_i\in\Theta}\|\bx_i\|_2\left\|\frac{1}{n}\sum_{j = 1}^n(\bW\transpose\widetilde\bx_j - \rho_n^{1/2}\bx_{0j})\bx_{0j}\transpose h_{0nij}(\bx_i)\right\|_2.
\end{align*}
By Assumption \ref{assumption:weight_functions}, Result \ref{result:pij_tilde_concentration}, and Cauchy-Schwarz inequality, 
\begin{align*}
\max_{i,j\in [n]}\sup_{\bx_i\in\Theta}|\widetilde{h}_{nij}(\bW\bx_i) - h_{0nij}(\bx_i)|
&\lesssim \rho_n^{-1}\max_{i,j\in [n]}|\rho_n\bx_{0i}\transpose\bx_{0j} - \widetilde{\bx}_i\transpose\widetilde{\bx}_j|
% \\&\quad
+ \rho_n\max_{i,j\in [n]}\sup_{\bx_i\in \Theta}\|\bx_i\|\|\bW\transpose\widetilde{\bx}_j - \rho_n^{1/2}\bx_{0j}\|_2\\
&\lesssim \sqrt{\frac{(\log n)^{2\xi}}{n\rho_n}}\quad\mbox{w.h.p.}.
\end{align*}
Therefore,
\begin{align*}
&\max_{i,j\in [n]}\sup_{\bx_i\in\Theta}\left\|\widetilde{h}_{nij}(\bW\bx_i)\rho_n^{-1/2}\bW\transpose\widetilde{\bx}_j - h_{0nij}(\bx_i)\bx_{0j}\right\|_2\\
&\quad\leq\max_{i,j\in [n]}\sup_{\bx_i\in\Theta}\left\{|\widetilde{h}_{nij}(\bW\bx_i) - h_{0nij}(\bx_i)|\|\rho_n^{-1/2}\widetilde{\bx}_j\| + h_{0nij}(\bx_i)\rho_n^{-1/2}\|\bW\transpose\widetilde{\bx}_j - \rho_n^{1/2}\bx_{0j}\|_2\right\}\\
&\quad\lesssim \sqrt{\frac{(\log n)^{2\xi}}{n\rho_n}}\quad\mbox{w.h.p..}
\end{align*}
For the first term, we apply Result \ref{result:concentration_of_infinity_norm} to obtain
\begin{align*}
&\sup_{\bx_i\in\Theta}
\left\|\frac{1}{n}\sum_{j = 1}^nE_{ij}\{\widetilde{h}_{nij}(\bW\bx_i)\rho_n^{-1/2}\bW\transpose\widetilde{\bx}_j - h_{0nij}(\bx_i)\bx_{0j}\}\right\|_2\\
&\quad\leq \frac{1}{n}\left(\sum_{j = 1}^n|E_{ij}|\right)\max_{i,j\in [n]}\sup_{\bx_i\in\Theta}\left\|\widetilde{h}_{nij}(\bW\bx_i)\rho_n^{-1/2}\bW\transpose\widetilde{\bx}_j - h_{0nij}(\bx_i)\bx_{0j}\right\|_2\\
&\quad\lesssim \frac{1}{n}(n\rho_n) \sqrt{\frac{(\log n)^{2\xi}}{n\rho_n}} = \rho_n^{1/2}\sqrt{\frac{(\log n)^{2\xi}}{n}}\quad\mbox{w.h.p..}
\end{align*}
Similarly, the second term can be bounded as follows:
\begin{align*}
&\sup_{\bx_i\in\Theta}
\frac{1}{n}\sum_{j = 1}^n|\bx_i\transpose\bW\transpose\widetilde{\bx}_j - \rho_n\bx_{0i}\transpose\bx_{0j}|\left\|\widetilde{h}_{nij}(\bW\bx_i)\rho_n^{-1/2}\bW\transpose\widetilde{\bx}_j - h_{0nij}(\bx_i)\bx_{0j}\right\|_2\\
&\quad\leq \max_{i,j\in [n]}\sup_{\bx_i\in\Theta}\left(\|\bx_i\|_2\|\widetilde{\bx}_j\| + \rho_n\|\bx_{0i}\|_2\|\bx_{0j}\|\right)\left\|\widetilde{h}_{nij}(\bW\bx_i)\rho_n^{-1/2}\bW\transpose\widetilde{\bx}_j - h_{0nij}(\bx_i)\bx_{0j}\right\|_2\\
&\quad\lesssim  \sqrt{\frac{(\log n)^{2\xi}}{n}}\quad\mbox{w.h.p..}
\end{align*}
The third term can be bounded by $\|\widetilde{\bX}\bW - \rho_n^{1/2}\bX_0\|_{2\to\infty} = O((\log n)^{\xi}/\sqrt{n})$ w.h.p.. To finish the proof, it is sufficient to show that
\[
\sup_{\bx_i\in\Theta}\left\|\frac{1}{n}\sum_{j = 1}^n[\bg_{ij}(\bx_i) - \expect_0\{\bg_{ij}(\bx_i)\}]\right\|_2\lesssim \rho_n^{1/2}\sqrt{\frac{\log n}{n}}\quad\mbox{w.h.p.}.
\]
For each $k\in [d]$, define a stochastic process $J_{kin}(\bx_i) = (1/n)\sum_{j = 1}^n[[\bg_{ij}(\bx_i)]_k - \expect_0\{[\bg_{ij}(\bx_i)]_k\}]$, where $[\cdot]_k$ denotes the $k$th coordinate of the vector. 
By definition, $\bg_{ij}(\bx_i) - \expect_0\{\bg_{ij}(\bx_i)\} = (A_{ij} - \rho_n\bx_{0i}\transpose\bx_{0j})h_n(\rho_n\bx_{0i}\transpose\bx_{0j}, \rho_n^{1/2}\bx_i\transpose\bx_{0j})\bx_{0j}$. By Assumption \ref{assumption:weight_functions}, we know that
\[
|h_n(\rho_n\bx_{0i}\transpose\bx_{0j}, \rho_n^{1/2}\bx_i\transpose\bx_{0j}) - h_n(\rho_n\bx_{0i}\transpose\bx_{0j}, \rho_n^{1/2}\by_i\transpose\bx_{0j})|\leq K\rho_n^{3/2}\sup_{\bx_j\in\Theta}\|\bx_j\|_2\|\bx_i - \by_i\|_2.
\]
Under Assumption \ref{assumption:signal_plus_noise} (vi), $A_{ij}$'s are uniformly bounded in $\psi_2$-Orlicz norms. Therefore, by Proposition 5.10 in \cite{vershynin2010introduction}, there exists a constant $C > 0$, such that for any $t > 0$ and $\bx_i,\by_i\in\Theta$,
\begin{align*}
\prob_0\left\{|J_{kin}(\bx_i) - J_{kin}(\by_i)| \geq t\right\}\leq e\exp\left( - \frac{Ct^2}{\rho_n^3\|\bx_i - \by_i\|_2^2/n}\right).
\end{align*}
Namely, there exists a constant $C_1 > 0$, such that $J_{kin}(\bx_i)$ is a sub-Gaussian process with regard to the metric $C_1(\rho_n^3/n)^{1/2}\|\cdot\|_2$. Since $\Theta$ is compact, then the packing entropy can also be bounded: There exists some constant $C_2$, such that
\begin{align*}
&\log\calD\left(\epsilon, \Theta, \frac{C_1\sqrt{\rho_n^3}\|\cdot\|_2}{\sqrt{n}}\right)\leq d\log\left(\frac{C_2\sqrt{\rho_n^3}}{\epsilon\sqrt{n}}\right),
\end{align*}
where, given a metric space $(T, \rho)$ and $\eps > 0$, the packing number $D(\eps, T, \rho)$ is the maximum number of disjoint balls with radius $\eps$ that are contained in $T$. 
% Note that by Lemma 5.9 in \cite{vershynin2010introduction}, 
% \[
% \|J_{kin}(\rho_n^{1/2}\bx_{0i})\|_{\psi_2}^2\lesssim \frac{1}{n^2}\sum_{j = 1}^n\|(A_{ij} - \rho_n\bx_{0i}\transpose\bx_{0j})h_{0nij}(\rho_n^{1/2}\bx_{0i})\bx_{0j}\|_{\psi_2}^2\lesssim \frac{1}{n}.
% \]
Since $\sup_{\bx_i,\by_i\in\Theta}\|\bx_i - \by_i\|_2 = C_3 < \infty$ for some constant $C_3 > 0$, we apply the maximal inequality for sub-Gaussian processes (Theorem 8.4 in \citealp{kosorok2008introduction}) to obtain
\begin{align*}
\left\|\sup_{\|\bx_i\|_2\leq 1}|J_{kin}(\bx_i) - J_{kin}(\rho_n^{1/2}\bx_{0i})|\right\|_{\psi_2} 
&\lesssim \int_0^{\frac{C_3\sqrt{\rho_n^3}}{\sqrt{n}}}\sqrt{\log\calD\left(\epsilon, \Theta, \frac{C_1\sqrt{\rho_n^3}\|\cdot\|_2}{\sqrt{n}}\right)}\mathrm{d}\eps\\
&\leq \int_0^{\frac{C_3\sqrt{\rho_n^3}}{\sqrt{n}}}\sqrt{
d\log\left(\frac{C_2\sqrt{\rho_n^3}}{\eps\sqrt{n}}\right)}\mathrm{d}\eps\lesssim \sqrt{\frac{\rho_n^3}{n}}.
\end{align*}
By Lemma 8.1 in \cite{kosorok2008introduction}, we obtain
\[
\sup_{\|\bx_i\|_2\leq 1}|J_{kin}(\bx_i) - J_{kin}(\rho_n^{1/2}\bx_{0i})| \lesssim \rho_n^{3/2}\sqrt{\frac{\log n}{n}}\quad\mbox{w.h.p.}.
\]
Now it is sufficient to consider $|J_{kin}(\rho_n^{1/2}\bx_{0i})|$ by triangle inequality. We consider the two scenarios under Assumption \ref{assumption:signal_plus_noise} (vi). If Assumption \ref{assumption:signal_plus_noise} (vi) (a) holds, then by Bernstein's inequality, we have, $|J_{kin}(\rho_n^{1/2}\bx_{0i})|\lesssim \rho_n^{1/2}\sqrt{(\log n)/n}$ w.h.p.. On the other hand, under Assumption \ref{assumption:signal_plus_noise} (vi) (b), we obtain from Proposition 5.16 in \cite{vershynin2010introduction} that $|J_{kin}(\rho_n^{1/2}\bx_{0i})|\lesssim \rho_n^{1/2}\sqrt{(\log n)/n}$ w.h.p. as well. Therefore, the proof is completed by the fact that
\begin{align*}
\sup_{\bx_i\in\Theta}\left\|\frac{1}{n}\sum_{j = 1}^n[\bg_{ij}(\bx_i) - \expect_0\{\bg_{ij}(\bx_i)\}]\right\|_2
&\leq \sum_{k=1}^d\left\{\sup_{\|\bx_i\|_2\leq 1}|J_{kin}(\bx_i) - J_{kin}(\rho_n^{1/2}\bx_{0i})| + |J_{kin}(\rho_n^{1/2}\bx_{0i})|\right\}\\
&\lesssim \rho_n^{1/2}\sqrt{\frac{\log n}{n}}\quad\mbox{w.h.p.}.
\end{align*}
\vspace*{2ex}\noindent
$\blacksquare$ \textbf{Proof of the second assertion. }
To begin with, we first compute the Jacobian
\begin{align*}
\bW\transpose\frac{\partial\widetilde{\bg}_{ij}}{\partial\bx_i\transpose}(\bW\bx_i)\bW
& = (A_{ij} - \bx_i\transpose\bW\transpose\widetilde{\bx}_j)D^{(0, 1)}h_n(\widetilde{\bx}_i\transpose\widetilde{\bx}_j, \bx_i\transpose\bW\transpose\widetilde{\bx}_j)\rho_n^{-1/2}\bW\transpose\widetilde{\bx}_j\widetilde{\bx}_j\transpose\bW\\
&\quad-h_n(\widetilde{\bx}_i\transpose\widetilde{\bx}_j, \bx_i\transpose\bW\transpose\widetilde{\bx}_j)\rho_n^{-1/2}\bW\transpose\widetilde{\bx}_j\widetilde{\bx}_j\transpose\bW.
\end{align*}
With a slight abuse of notations, we denote $D^{(0, 1)}\widetilde{h}_{nij}(\bx_i) = D^{(0, 1)}h_n(\widetilde{\bx}_i\transpose\widetilde{\bx}_j, \bx_i\transpose\widetilde{\bx}_j)$ and $D^{(0, 1)}h_{0nij}(\bx_i) = D^{(0, 1)}h_n(\rho_n\bx_{0i}\transpose\bx_{0j}, \rho_n^{1/2}\bx_i\transpose\bx_{0j})$. Then we have
\begin{align*}
&\bW\transpose\frac{\partial\widetilde{\bg}_{ij}}{\partial\bx_i\transpose}(\bW\bx_i)\bW - \frac{\partial\bg_{ij}}{\partial\bx_i\transpose}(\bx_i)\\
&\quad = (A_{ij} - \rho_n^{1/2}\bx_i\transpose\bx_{0j})\left\{D^{(0, 1)}\widetilde{h}_{nij}(\bW\bx_i)\rho_n^{-1/2}\bW\transpose\widetilde{\bx}_j\widetilde{\bx}_j\transpose\bW
- D^{(0, 1)}h_{0nij}(\bx_i)\rho_n^{1/2}\bx_{0j}\bx_{0j}\transpose
\right\}\\
&\quad\quad -\bx_i\transpose(\bW\transpose\widetilde{\bx}_j - \rho_n^{1/2}\bx_{0j}) D^{(0, 1)}\widetilde{h}_{nij}(\bW\bx_i)\rho_n^{-1/2}\bW\transpose\widetilde{\bx}_j\widetilde{\bx}_j\transpose\bW\\
&\quad\quad -\left\{\widetilde{h}_{nij}(\bW\bx_i)\rho_n^{-1/2}\bW\transpose\widetilde{\bx}_j\widetilde{\bx}_j\transpose\bW - h_{0nij}(\bx_i)\rho_n^{1/2}\bx_{0j}\bx_{0j}\transpose\right\}
 .
\end{align*}
Following the proof of the first assertion of Lemma \ref{lemma:ULLN}, we have
\begin{align*}
\max_{i,j\in [n]}\sup_{\bx_i\in\Theta}\left|\widetilde{h}_{nij}(\bW\bx_i) - h_{0nij}(\bx_i)\right|
&\lesssim \rho_n^{-1}\max_{i,j\in [n]}|\widetilde{\bx}_i\transpose\widetilde{\bx}_j - \rho_n\bx_{0i}\transpose\bx_{0j}| + \rho_n\sup_{\bx_i\in\Theta}\|\bx_i\|_2\|\widetilde\bX\bW - \rho_n^{1/2}\bX_0\|_{2\to\infty}\\
&\lesssim \sqrt{\frac{(\log n)^{2\xi}}{n\rho_n}}\quad\mbox{w.h.p.},
\end{align*}
and
\begin{align*}
&\max_{i,j\in [n]}\sup_{\bx_i\in\Theta}\left\|\widetilde{h}_{nij}(\bW\bx_i)\rho_n^{-1/2}\bW\transpose\widetilde{\bx}_j\widetilde{\bx}_j\transpose\bW - h_{0nij}(\bx_i)\rho_n^{1/2}\bx_{0j}\bx_{0j}\transpose\right\|\\
&\quad\lesssim \max_{i,j\in [n]}\left\{\rho_n^{-1/2}\|\widetilde{\bx}_j\|_2^2\sup_{\bx_i\in\Theta}\left|\widetilde{h}_{nij}(\bW\bx_i) - h_{0nij}(\bx_i)\right| + \rho_n^{-1/2}\|\bW\transpose\widetilde{\bx}_j - \rho_n^{1/2}\bx_{0j}\|(\|\bW\transpose\widetilde{\bx}_j\|_2 + \|\rho_n^{1/2}\bx_{0j}\|_2)\right\}\\
&\quad\lesssim \rho_n^{1/2}\sqrt{\frac{(\log n)^{2\xi}}{n\rho_n}} + \sqrt{\frac{(\log n)^{2\xi}}{n}} \asymp \sqrt{\frac{(\log n)^{2\xi}}{n}}\quad\mbox{w.h.p.}.
\end{align*}
Also, from the proof of the first assertion of Lemma \ref{lemma:ULLN} again, we have
\begin{align*}
&\max_{i,j\in [n]}\sup_{\bx_i\in\Theta}\left\|D^{(0, 1)}\widetilde{h}_{nij}(\bW\bx_i)\rho_n^{-1/2}\bW\transpose\widetilde{\bx}_j\widetilde{\bx}_j\transpose\bW - D^{(0, 1)}h_{0nij}(\bx_i)\rho_n^{1/2}\bx_{0j}\bx_{0j}\transpose\right\|_2\lesssim \sqrt{\frac{(\log n)^{2\xi}}{n}}\quad\mbox{w.h.p.}.
\end{align*}
It follows that
\begin{align*}
\sup_{\bx_i\in\Theta}\frac{1}{n}\sum_{j = 1}^n\left\|\bW\transpose\frac{\partial\widetilde{\bg}_{ij}}{\partial\bx_i\transpose}(\bW\bx_i)\bW - \frac{\partial{\bg}_{ij}}{\partial\bx_i\transpose}(\bx_i)\right\|_2
&\lesssim \sqrt{\frac{(\log n)^{2\xi}}{n}}\frac{1}{n}\sum_{j = 1}^n\left(|A_{ij}| + \rho_n^{1/2}\|\bx_{0j}\|_2\sup_{\bx_i\in\Theta}\|\bx_i\|_2\right)\\
&\quad + \sup_{\bx_i\in\Theta}\|\bx_i\|_2\|\widetilde{\bX}\bW - \rho_n^{1/2}\bX_0\|_{2\to\infty}\rho_n^{1/2}\|\widetilde{\bx}_j\|_2^2\\
&\quad + \sqrt{\frac{(\log n)^{2\xi}}{n}}\\
&\lesssim \rho_n^{1/2}\sqrt{\frac{(\log n)^{2\xi}}{n}} +  \rho_n^{3/2}\sqrt{\frac{(\log n)^{2\xi}}{n}} +  \sqrt{\frac{(\log n)^{2\xi}}{n}}\\
&\asymp \sqrt{\frac{(\log n)^{2\xi}}{n}}\quad\mbox{w.h.p.}.
\end{align*}
Following the proof of the first assertion above, by the maximal inequality for sub-Gaussian processes and Assumption \ref{assumption:weight_functions}, we have
\[
\sup_{\bx_i\in\Theta}\left\|\frac{1}{n}\sum_{j = 1}^n\left\{\frac{\partial{\bg}_{ij}}{\partial\bx_i\transpose}(\bx_i) - \expect_0\frac{\partial{\bg}_{ij}}{\partial\bx_i\transpose}(\bx_i)\right\}\right\|_2\lesssim \rho_n^{1/2}\sqrt{\frac{\log n}{n}}\quad\mbox{w.h.p.}.
\]
The proof is then completed by combining the two uniform concentration bounds. 
\end{proof}

\subsection{Central limit theorem}
\label{sub:CLT}
\begin{theorem}[Central Limit Theorem]
\label{thm:CLT}
Suppose Assumptions \ref{assumption:signal_plus_noise}, \ref{assumption:regularity_condition}, and \ref{assumption:weight_functions} hold. Then for all $i\in [n]$,
\[
\left\|\frac{1}{n}\sum_{j = 1}^n\bW\transpose\widetilde{\bg}_{ij}(\rho_n^{1/2}\bW\bx_{0i}) - \frac{1}{n}\sum_{j = 1}^n{\bg}_{ij}(\rho_n^{1/2}\bx_{0i})\right\|_2 \lesssim \frac{(\log n)^{2\xi}}{n}\quad\mbox{w.h.p.}.
\]
\end{theorem}
\begin{proof}[\bf Proof of Theorem \ref{thm:CLT}]
Denote $\widetilde{h}_{0nij} = h_n(\widetilde\bx_i\transpose\widetilde\bx_j, \rho_n^{1/2}\bx_{0i}\transpose\bW\transpose\widetilde{\bx}_j)$, $h_{0nij} = h_n(\rho_n\bx_{0i}\transpose\bx_{0j}, \rho_n\bx_{0i}\transpose\bx_{0j})$, and $\bpsi_n(\bx_i, \bu, \bv) = h_n(\rho_n\bu\transpose\bv, \rho_n^{1/2}\bx_i\transpose\bv)\bv$. 
For a vector $\bx$, let $[\bx]_k$ denote its $k$th coordinate. 
Simple calculation leads to 
\begin{align*}
&\frac{\partial \bpsi_n}{\partial\bx_i\transpose} = D^{(0, 1)}h_n\rho_n^{1/2}\bv\bv\transpose,
\quad
\frac{\partial \bpsi_n}{\partial\bu\transpose} = D^{(1, 0)}h_n\rho_n\bv\bv\transpose,\\
&\frac{\partial \bpsi_n}{\partial\bv\transpose}
 = D^{(1, 0)}h_n\rho_n\bv\bu\transpose + D^{(0, 1)}h_n\rho_n^{1/2}\bv\bx_i\transpose + h_n\eye_d,\\
&\frac{\partial^2 [\bpsi_n]_k}{\partial\bu\partial\bu\transpose} = (\be_k\transpose\bv)D^{(2, 0)}h_n\rho_n^2\bv\bv\transpose,\\
&\frac{\partial^2 [\bpsi_n]_k}{\partial\bu\partial\bv\transpose} = (\be_k\transpose\bv)\left\{D^{(2, 0)}h_n\rho_n^2\bv\bu\transpose + D^{(1, 0)}h_n\rho_n\eye_d + D^{(1, 1)}h_n\rho_n^{3/2}\bv\bx_i\transpose\right\} + D^{(1, 0)}h_n\rho_n\bv\be_k\transpose\\
&\frac{\partial^2 [\bpsi_n]_k}{\partial\bv\partial\bv\transpose} = 
\left\{
D^{(1, 0)}h_n\rho_n\be_k\bu\transpose + D^{(0, 1)}h_n\rho_n^{1/2}\be_k\bx_i\transpose
\right\}\\
&\quad\quad\quad\quad\quad + (\be_k\transpose\bv)
\left\{
D^{(2, 0)}h_n\rho_n^2\bu\bu\transpose + D^{(1, 1)}h_n\rho_n^{3/2}(\bu\bx_i\transpose + \bx_i\bu\transpose) + D^{(0, 2)}h_n\rho_n\bx_i\bx_i\transpose
\right\}\\
&\quad\quad\quad\quad\quad + \left\{D^{(1, 0)}h_n\rho_n\bu\be_k\transpose + D^{(0, 1)}h_n\rho_n^{1/2}\bx_i\be_k\transpose\right\},
\end{align*}
where we have suppressed the arguments $\bx_i,\bu,\bv$ for $h_n$ and $\bpsi$. Then by Assumption \ref{assumption:weight_functions}, there exists $\eps > 0$, such that for all $\bx_i\in \Theta$, $\bu\in B(\bx_{0i}, \eps)$, and $\bv\in B(\bx_{0j}, \eps)$, 
\begin{align*}
\left\|\frac{\partial^2 [\bpsi_n]_k}{\partial\bu\partial\bu\transpose}\right\|_2\lesssim 1,\quad
\left\|\frac{\partial^2 [\bpsi_n]_k}{\partial\bu\partial\bv\transpose}\right\|_2\lesssim 1,\quad
\left\|\frac{\partial^2 [\bpsi_n]_k}{\partial\bv\partial\bv\transpose}\right\|_2\lesssim 1.
\end{align*}
Denote
\[
\bB_{nij}^{(\bu)} = \frac{\partial \bpsi_n}{\partial\bu\transpose}(\rho_n^{1/2}\bx_{0i}, \bx_{0i}, \bx_{0j})\quad\mbox{and}\quad
\bB_{nij}^{(\bv)} = \frac{\partial \bpsi_n}{\partial\bv\transpose}(\rho_n^{1/2}\bx_{0i}, \bx_{0i}, \bx_{0j}).
\]
Clearly, the entries of $\bB_{nij}^{(\bu)}$ and $\bB_{nij}^{(\bv)}$ are uniformly bounded by a constant by Assumption \ref{assumption:weight_functions}. 
Therefore, by Theorem \ref{thm:uniform_concentration_eigenvector} and a Taylor expansion of $\bpsi_n$, we obtain
\begin{equation}
\label{eqn:Taylor_expansion_psi_function}
\begin{aligned}
\widetilde{h}_{0nij}\rho_n^{-1/2}\bW\transpose\widetilde{\bx}_j - h_{0nij}\bx_{0j}
& = \bpsi_n(\rho_n^{1/2}\bx_{0i}, \rho_n^{-1/2}\bW\transpose\widetilde{\bx}_i, \rho_n^{-1/2}\bW\transpose\widetilde{\bx}_j) - \bpsi_n(\rho_n^{1/2}\bx_{0i}, \bx_{0i}, \bx_{0j})\\
& = \bB_{nij}^{(\bu)}\rho_n^{-1/2}(\bW\transpose\widetilde{\bx}_i - \rho_n^{1/2}\bx_{0i}) + 
\bB_{nij}^{(\bv)}\rho_n^{-1/2}(\bW\transpose\widetilde{\bx}_j - \rho_n^{1/2}\bx_{0j})
% \\&\quad
 + \br_{nij}^{(\bpsi)},
\end{aligned}
\end{equation}
where $\max_{i,j\in [n]}\|\br_{nij}^{(\bpsi)}\|_2\lesssim \rho_n^{-1}\|\widetilde{\bX}\bW - \rho_n^{1/2}\bX_0\|_{2\to\infty}^2\lesssim (\log n)^{2\xi}/(n\rho_n)$ w.h.p.. Now write by triangle inequality and Cauchy-Schwarz inequality
\begin{align*}
&\left\|\frac{1}{n}\sum_{j = 1}^n\bW\transpose\widetilde{\bg}_{ij}(\rho_n^{1/2}\bW\bx_{0i})
 - \frac{1}{n}\sum_{j = 1}^n{\bg}_{ij}(\rho_n^{1/2}\bx_{0i})\right\|_2\\
&\quad\leq \frac{1}{n}\rho_n^{1/2}\|\bx_{0i}\|_2\left\|\sum_{j = 1}^n(\bW\transpose\widetilde{\bx}_j - \rho_n^{1/2}\bx_{0j})\bx_{0j}\transpose h_{0nij}\right\|_2\\
&\quad\quad + \frac{1}{n}\rho_n^{1/2}\|\bx_{0i}\|_2\left\|\sum_{j = 1}^n(\bW\transpose\widetilde{\bx}_j - \rho_n^{1/2}\bx_{0j})(\widetilde{h}_{0nij}\rho_n^{-1/2}\bW\transpose\widetilde{\bx}_j - h_{0nij}\bx_{0j})\transpose\right\|_2\\
&\quad\quad + \left\|\frac{1}{n}\sum_{j = 1}^nE_{ij}(\widetilde{h}_{0nij}\rho_n^{-1/2}\bW\transpose\widetilde{\bx}_j - h_{0nij}\bx_{0j})\right\|_2.
\end{align*}
For the first term, we apply Theorem \ref{thm:uniform_concentration_eigenvector} to write
\begin{align*}
\left\|\frac{1}{n}\sum_{j = 1}^n(\bW\transpose\widetilde{\bx}_j - \rho_n^{1/2}\bx_{0j})\bx_{0j}\transpose h_{0nij}\right\|_2
&\leq \left\|\frac{1}{n\rho_n^{1/2}}\sum_{j = 1}^n\sum_{l = 1}^nh_{0nij}(\bX_0\transpose\bX_0)^{-1}\bx_{0l}\bx_{0j}\transpose E_{jl}\right\|_2\\
&\quad + \frac{1}{n}\sum_{j = 1}^n|h_{0nij}|\|\bx_{0j}\|_2\|\bR_\bX\|_{2\to\infty}.
\end{align*}
Under Assumption \ref{assumption:signal_plus_noise} (vi) (a), by Bernstein's inequality, we have
\begin{align*}
\left\|\frac{1}{n\rho_n^{1/2}}\sum_{j = 1}^n\sum_{l = 1}^nh_{0nij}(\bX_0\transpose\bX_0)^{-1}\bx_{0l}\bx_{0j}\transpose E_{jl}\right\|_2
&\lesssim \frac{\log n}{n\rho_n^{1/2}}\max_{j,l\in [n]}\|(\bX_0\transpose\bX_0)^{-1}\|_2\|\bx_{0j}\|_2\|\bx_{0l}\|_2\\
&\quad + \frac{(\rho_n\log n)^{1/2}}{n\rho_n^{1/2}}\left(\sum_{j = 1}^n\sum_{l = 1}^n\|(\bX_0\transpose\bX_0)^{-1}\|_2^2\|\bx_{0j}\|_2^2\|\bx_{0l}\|_2^2\right)^{1/2}\\
&\lesssim \frac{\log n}{n^2\rho_n^{1/2}} + \frac{(\rho_n\log n)^{1/2}}{n\rho_n^{1/2}}\lesssim \frac{(\log n)^{1/2}}{n}\quad\mbox{w.h.p.}.
\end{align*}
Under Assumption \ref{assumption:signal_plus_noise} (vi) (b), by Proposition 5.10 in \cite{vershynin2010introduction}, we have
\[
\left\|\frac{1}{n\rho_n^{1/2}}\sum_{j = 1}^n\sum_{l = 1}^nh_{0nij}(\bX_0\transpose\bX_0)^{-1}\bx_{0l}\bx_{0j}\transpose E_{jl}\right\|_2\lesssim \frac{(\log n)^{1/2}}{n}\quad\mbox{w.h.p.}.
\]
Applying the second assertion of Theorem \ref{thm:uniform_concentration_eigenvector} yields
\[
\frac{1}{n}\sum_{j = 1}^n|h_{0nij}|\|\bx_{0j}\|_2\|\bR_\bX\|_{2\to\infty}\lesssim \frac{(\log n)^{2\xi}}{n\rho_n^{1/2}}\quad\mbox{w.h.p.}.
\]
Hence the first term is $O((\log n)^{2\xi}/n)$ w.h.p..
Also, by \eqref{eqn:uniform_Lipschitz_h}, we know that
\[
\max_{i,j\in [n]}\|\widetilde{h}_{0nij}\rho_n^{-1/2}\bW\transpose\widetilde{\bx}_j - h_{0nij}\bx_{0j}\|_2\lesssim \sqrt{\frac{(\log n)^{2\xi}}{n\rho_n}}.
\]
Therefore, the second term is also $O((\log n)^{2\xi}/{n})$ w.h.p. by the same reasoning. It suffices to show that the third term is $O((\log n)^{2\xi}/n)$ w.h.p.. By \eqref{eqn:Taylor_expansion_psi_function}, we have
\begin{align*}
\left\|\frac{1}{n}\sum_{j = 1}^nE_{ij}(\widetilde{h}_{0nij}\rho_n^{-1/2}\bW\transpose\widetilde{\bx}_j - h_{0nij}\bx_{0j})\right\|_2
&\leq \left\|\frac{1}{n\rho_n^{1/2}}\sum_{j = 1}^n\bB_{nij}^{(\bu)}E_{ij}\right\|_2\|\widetilde{\bX}\bW - \rho_n^{1/2}\bX_0\|_{2\to\infty}\\
&\quad + \left\|\frac{1}{n\rho_n^{1/2}}\sum_{j = 1}^nE_{ij}\bB_{nij}^{(\bv)}(\bW\transpose\widetilde{\bx}_j - \rho_n^{1/2}\bx_{0j})\right\|_2
\\&\quad
+ \frac{1}{n}\sum_{j = 1}^n|E_{ij}|\|\br_{nij}^{(\bpsi)}\|_2.
\end{align*}
Since $\max_{i,j\in [n]}\|\br_{nij}^{(\bpsi)}\|_2\lesssim (\log n)^{2\xi}/(n\rho_n)$ w.h.p., then 
\[
\frac{1}{n}\sum_{j = 1}^n|E_{ij}|\|\br_{nij}^{(\bpsi)}\|_2\leq \frac{1}{n}\|\bE\|_\infty\|\max_{i,j\in [n]}\|\br_{nij}^{(\bpsi)}\|_2\lesssim \frac{(\log n)^{2\xi}}{n}\quad\mbox{w.h.p.}.
\]
Also, by Theorem \ref{thm:uniform_concentration_eigenvector} and either Bernstein's inequality under Assumption \ref{assumption:signal_plus_noise} (vi) (i) or Proposition 5.10 in \cite{vershynin2010introduction} under Assumption \ref{assumption:signal_plus_noise} (vi) (ii), 
\begin{align*}
\left\|\frac{1}{n\rho_n^{1/2}}\sum_{j = 1}^n\bB_{nij}^{(\bu)}E_{ij}\right\|_2\|\widetilde{\bX}\bW - \rho_n^{1/2}\bX_0\|_{2\to\infty}
\lesssim \frac{(n\rho_n\log n)^{1/2}}{n\rho_n^{1/2}}\sqrt{\frac{(\log n)^{2\xi}}{n}}\leq \frac{(\log n)^{2\xi}}{n}\quad\mbox{w.h.p.}.
\end{align*}
Now we focus on the remaining term. First write by Theorem \ref{thm:uniform_concentration_eigenvector} and Result \ref{result:concentration_of_infinity_norm} that
\begin{align*}
&\left\|\frac{1}{n\rho_n^{1/2}}\sum_{j = 1}^nE_{ij}\bB_{nij}^{(\bv)}(\bW\transpose\widetilde{\bx}_j - \rho_n^{1/2}\bx_{0j})\right\|_2\\
&\quad\leq \left\|\frac{1}{n\rho_n}\sum_{j = 1}^n\sum_{l = 1}^nE_{ij}E_{jl}\bB_{nij}^{(\bv)}(\bX_0\transpose\bX_0)^{-1}\bx_{0l}\right\|_2
% \\&\quad\quad
+ \frac{1}{n\rho_n^{1/2}}\sum_{j = 1}^n|E_{ij}|\|\bB_{nij}^{(\bv)}\|_2\|\bR_\bX\|_{2\to\infty}\\
&\quad\lesssim \left\|\frac{1}{n^2\rho_n}\sum_{j = 1}^n\sum_{l = 1}^nE_{ij}E_{jl}\bB_{nij}^{(\bv)}\left(\frac{1}{n}\bX_0\transpose\bX_0\right)^{-1}\bx_{0j}\right\|_2
% \\&\quad
 + \frac{1}{n\rho_n^{1/2}}(n\rho_n)\frac{(\log n)^{2\xi}}{n\rho_n^{1/2}}\\
&\quad = \left\|\frac{1}{n^2\rho_n}\sum_{j = 1}^n\sum_{l = 1}^nE_{ij}E_{jl}\bB_{nij}^{(\bv)}\left(\frac{1}{n}\bX_0\transpose\bX_0\right)^{-1}\bx_{0j}\right\|_2 + \frac{(\log n)^{2\xi}}{n}\quad\mbox{w.h.p.}.
\end{align*}
Since the entries of $\bB_{nij}^{(\bv)}$ and $(\frac{1}{n}\bX_0\transpose\bX_0)^{-1}$ are uniformly bounded, applying Result \ref{result:Berstein_concentration_E_square} completes the proof.
\end{proof}

\section{Proofs of The Main Results}
\label{sec:proofs_main_results}

\subsection{Proof of Theorem \ref{thm:Large_sample_Z_estimator}}
\label{sub:proof_Z_estimator}
We first show the following weaker consistency result
\[
\|\bW\transpose\widehat{\bx}_i - \rho_n^{1/2}\bx_{0i}\|_2 \lesssim \sqrt{\frac{(\log n)^{2\xi}}{(n\rho_n)}}\quad\mbox{w.h.p.}
\]
and then establish the asymptotic normality based on this convergence rate result. 
\begin{proof}[\bf Proof of consistency]
Let $M_{in}(\bx) = \|(1/n)\sum_{j = 1}^n\expect_0\bg_{ij}(\bx_i)\}\|_2$. By Result \ref{result:identifiability}, $M_{in}$ is uniquely minimized at $\rho_n^{1/2}\bx_{0i}$ and for all $\eps > 0$,
\begin{align}\label{eqn:Mfunction_isolated_maximizer}
\sup_{\|\bx_i - \rho_n^{1/2}\bx_{0i}\|_2 > \eps}\{-M_{in}(\bx_i)\} + \rho_n^{1/2} \delta_0\eps \leq -M_{in}(\rho_n^{1/2}\bx_{0i}) = 0.
\end{align}
Now denote $\widetilde{M}_{in}(\bx_i) = \|(1/n)\sum_{j = 1}^n\bW\transpose\widetilde{\bg}_{ij}(\bW\bx_i)\}\|_2$. By Assumption \ref{assumption:regularity_condition} (ii), $\bW\transpose\widehat{\bx}_i$ is the unique minimizer of $\widetilde{M}_{in}$ inside the interior of $\Theta$ w.h.p.. In addition, by Lemma \ref{lemma:ULLN}, 
% \begin{align}\label{eqn:uniform_convergence_surrogate_likelihood}
% \sup_{\|\bx_i\|_2\in \Theta}\left| \widetilde{M}_{in}(\bx_i) - M_{in}(\bx_i)\right| \lesssim\sqrt{\frac{\log n}{n}}\quad\mbox{w.h.p.}.
% \end{align}
% Namely,
% \[
${M}_{in}(\rho_n^{1/2}\bx_{0i}) - \widetilde{M}_{in}(\rho_n^{1/2}\bx_{0i}) = r_{in}^{(M)}$, where $|r_{in}^{(M)}|\lesssim  \sqrt{{(\log n)^{2\xi}}/{n}}$ w.h.p..
% \] 
Since $\bW\transpose\widehat{\bx}_i$ is the minimizer of $\widetilde{M}_{in}(\bx_i)$, it follows again by Lemma \ref{lemma:ULLN} that
\begin{align*}
M_{in}(\bW\transpose\widehat{\bx}_i) - M_{in}(\rho_n^{1/2}\bx_{0i})
& = M_{in}(\bW\transpose\widehat{\bx}_i) - \widetilde{M}_{in}(\rho_n^{1/2}\bx_{0i}) - r_{in}^{(M)} \\
&\leq M_{in}(\bW\transpose\widehat{\bx}_i) - \widetilde{M}_{in}(\bW\transpose\widehat{\bx}_i) + |r_{in}^{(M)}|\\
% & = \widetilde{M}_{in}(\bW(\bW\transpose\widehat{\bx}_i)) - M_{in}(\bW\transpose\widehat{\bx}_i) + r_{in}^{(M)}\\
&\leq \sup_{\bx_i\in\Theta}|\widetilde{M}_{in}(\bx_i) - M_{in}(\bx_i)| + |r_{in}^{(M)}|.
\end{align*}
This implies that $M_{in}(\bW\transpose\widehat{\bx}_i) - M_{in}(\rho_n^{1/2}\bx_{0i})\lesssim \sqrt{(\log n)^{2\xi}/n}$ w.h.p. by Lemma \ref{lemma:ULLN}.
By \eqref{eqn:Mfunction_isolated_maximizer}, for all $\eps > 0$ and for any $\by$ with $\|\by - \rho_n^{1/2}\bx_{0i}\|_2 >  \eps$, we have $M_{in}(\by) - M_{in}(\rho_n^{1/2}\bx_{0i})\geq \rho_n^{1/2}\delta_0\eps$. The proof is thus completed by taking $\eps = C\sqrt{(\log n)^{2\xi}/(n\rho_n)}$ for an appropriate constant $C > 0$.
% and taking the contrapositive of the previous sentence.
\end{proof}

\begin{proof}[\bf Proof of asymptotic normality]
By the consistency result in the aforementioned proof, we know that $\|\bW\transpose\widehat{\bx}_i - \rho_n^{1/2}\bx_{0i}\|_2 = o_{\prob_0}(1)$. 
Let $[\bx]_k$ denote the $k$th coordinate of a vector $\bx$. 
By Assumption \ref{assumption:regularity_condition} (ii) and Taylor's theorem,
\begin{align*}
\zero_d = \frac{1}{n}\sum_{j = 1}^n\bW\transpose\widetilde{\bg}_{ij}(\widehat{\bx}_i)& = 
\frac{1}{n}\sum_{j = 1}^n\bW\transpose\widetilde{\bg}_{ij}(\rho_n^{1/2}\bW\bx_{0i}) + \frac{1}{n}\sum_{j = 1}^n\bW\transpose\frac{\partial\widetilde{\bg}_{ij}}{\partial\bx_i\transpose}(\rho_n^{1/2}\bW\bx_{0i})\bW(\bW\transpose\widehat{\bx}_i - \rho_n^{1/2}\bx_{0i})\\
&\quad + \bW\transpose\br_n^{(\bg)},
\end{align*}
where for each $k\in [d]$, there exists $\theta_{ik}\in [0, 1]$, such that $\bar{\bx}_i^{(k)} = (1 - \theta_{ik})\rho_n^{1/2}\bW\bx_{0i} + \theta_{ik}\widehat{\bx}_{i}$, and
\[
[\br_n^{(\bg)}]_k = (\widehat{\bx}_i - \rho_n^{1/2}\bW\bx_{0i})\transpose\frac{1}{n}\sum_{j = 1}^n\frac{\partial^2[\widetilde{\bg}_{ij}]_k}{\partial\bx_i\partial\bx_i\transpose}(\bar{\bx}_i^{(k)})(\widehat{\bx}_i - \rho_n^{1/2}\bW\bx_{0i}).
\]
Denote $\widetilde{h}_{nij}(\bx_i) = h_n(\widetilde{\bx}_i\transpose\widetilde{\bx}_j, \bx_i\transpose\widetilde{\bx}_j)$ and $h_{0nij}(\bx_i) = h_n(\rho_n\bx_{0i}\transpose\bx_{0j}, \rho_n^{1/2}\bx_i\transpose\bx_{0j})$. With a slight abuse of notations, we also denote $D^{(0, 1)}\widetilde{h}_{nij}(\bx_i) = D^{(0, 1)}h_n(\widetilde{\bx}_i\transpose\widetilde{\bx}_j, \bx_i\transpose\widetilde{\bx}_j)$ and $D^{(0, 2)}\widetilde{h}_{nij}(\bx_i) = D^{(0, 2)}h_n(\widetilde{\bx}_i\transpose\widetilde{\bx}_j, \bx_i\transpose\widetilde{\bx}_j)$. 
We then have 
\begin{align*}
\frac{1}{n}\sum_{j = 1}^n\frac{\partial^2[\widetilde{\bg}_{ij}]_k}{\partial\bx_i\partial\bx_i\transpose}(\bar{\bx}_i^{(k)})
& = -2\frac{1}{n\rho_n^{1/2}}\sum_{j = 1}^n[\widetilde{\bx}_j]_k D^{(0, 1)}\widetilde{h}_{nij}(\bar{\bx}_i^{(k)})\widetilde{\bx}_j\widetilde{\bx}_j\transpose\\
&\quad + \frac{1}{n\rho_n^{1/2}}\sum_{j = 1}^n(\rho_n\bx_{0i}\transpose\bx_{0j} - \bar{\bx}_i^{(k)}\widetilde{\bx}_{j}) [\widetilde{\bx}_j]_kD^{(0, 2)}\widetilde{h}_{nij}(\bar{\bx}_i^{(k)})\widetilde{\bx}_j\widetilde{\bx}_j\transpose\\
&\quad + \frac{1}{n\rho_n^{1/2}}\sum_{j = 1}^nE_{ij}[\widetilde{\bx}_j]_kD^{(0, 2)}\widetilde{h}_{nij}(\bar{\bx}_i^{(k)})\widetilde{\bx}_j\widetilde{\bx}_j\transpose.
\end{align*}
Note that $\|\bW\transpose\bar{\bx}_i^{(k)} - \rho_n^{1/2}\bx_{0i}\|_2\leq \theta_{ik}\|\bW\transpose\widehat{\bx}_{ik} - \rho_n^{1/2}\bx_{0i}\|_2 \lesssim \sqrt{(\log n)^{2\xi}/(n\rho_n)}$ w.h.p. by the previously proved consistency result. This implies that $(\bW\transpose\bar{\bx}_i^{(k)})\transpose\bW\transpose\widetilde{\bx}_j\in [-r,r]$ w.h.p.. Then by Assumption \ref{assumption:weight_functions} and Result \ref{result:pij_tilde_concentration}, 
\begin{align*}
\max_{i,j}|D^{(0, 1)}\widetilde{h}_{nij}(\bar{\bx}_i^{(k)})|\lesssim \rho_n,\quad\max_{i,j}|D^{(0, 2)}\widetilde{h}_{nij}(\bar{\bx}_i^{(k)})|\lesssim \rho_n\quad\mbox{w.h.p.}.
\end{align*}
It follows that
\begin{align*}
&\left\|\frac{1}{n\rho_n^{1/2}}\sum_{j = 1}^n[\widetilde{\bx}_j]_k D^{(0, 1)}\widetilde{h}_{nij}(\bar{\bx}_i^{(k)})\widetilde{\bx}_j\widetilde{\bx}_j\transpose\right\|_2\lesssim \frac{1}{n\rho_n^{1/2}}n\rho_n^{5/2} = \rho_n^2\quad\mbox{w.h.p.},\\
&\left\|\frac{1}{n\rho_n^{1/2}}\sum_{j = 1}^n(\rho_n\bx_{0i}\transpose\bx_{0j} - \bar{\bx}_i^{(k)}\widetilde{\bx}_{j}) [\widetilde{\bx}_j]_kD^{(0, 2)}\widetilde{h}_{nij}(\bar{\bx}_i^{(k)})\widetilde{\bx}_j\widetilde{\bx}_j\transpose\right\|_2\lesssim \frac{1}{n\rho_n^{1/2}}n\rho_n^{3} = \rho_n^{5/2}\quad\mbox{w.h.p.}.
\end{align*}
We also obtain from Result \ref{result:concentration_of_infinity_norm}, Result \ref{result:pij_tilde_concentration}, and Assumption \ref{assumption:weight_functions} that
\begin{align*}
&\frac{1}{n\rho_n^{1/2}}\sum_{j = 1}^n|E_{ij}||[\widetilde{\bx}_j]_k||D^{(0, 2)}\widetilde{h}_{nij}(\bar{\bx}_i^{(k)})|\|\widetilde{\bx}_j\|_2^2\\
&\quad\leq \frac{1}{n\rho_n^{1/2}}\|\bE\|_\infty \max_{i,j\in [n]}\|\widetilde{\bx}_j\|_2^3|D^{(0, 2)}\widetilde{h}_{nij}(\bar{\bx}_i^{(k)})|\lesssim \frac{1}{n\rho_n^{1/2}}n\rho_n^{5/2} = \rho_n^2\quad\mbox{w.h.p.}.
\end{align*}
Therefore, $\left\|(1/n)\sum_{j = 1}^n\partial^2[\widetilde{\bg}_{ij}]_k(\bar{\bx}_i^{(k)})/\partial\bx_i\partial\bx_i\transpose\right\|_2 = O(\rho_n^2)$ w.h.p.. Namely, 
\[
\bW\transpose\br_n^{(\bg)} = O\left(\rho_n^{3/2}\sqrt{\frac{(\log n)^{2\xi}}{n}}\right)(\bW\transpose\widehat{\bx}_i - \rho_n^{1/2}\bx_{0i})\quad\mbox{w.h.p.}.
\]
Hence, by Theorem \ref{thm:CLT} and Lemma \ref{lemma:LLN}, with $\bg_{0ij}:=\bg_{ij}(\rho_n^{1/2}\bx_{0i})$ and $\bG_{0in} = \bG_{in}(\rho_n^{1/2}\bx_{0i})$, we have
\begin{align*}
-\frac{1}{n}\sum_{j = 1}^n\bg_{0ij} + O\left(\frac{(\log n)^{2\xi}}{n}\right)
& = \left\{\bG_{0in} + O\left(\sqrt{\frac{(\log n)^{2\xi}}{n}}\right)\right\}(\bW\transpose\widehat{\bx}_i - \rho_n^{1/2}\bx_{0i})\\
&\quad + O\left(\rho_n^{3/2}\sqrt{\frac{(\log n)^{2\xi}}{n}}\right)(\bW\transpose\widehat{\bx}_i - \rho_n^{1/2}\bx_{0i})\\
& = \left\{\bG_{0in} + O\left(\sqrt{\frac{(\log n)^{2\xi}}{n}}\right)\right\}(\bW\transpose\widehat{\bx}_i - \rho_n^{1/2}\bx_{0i})\quad\mbox{w.h.p.}.
\end{align*}
By Woodbury matrix identity, matrix series expansion of $(\eye_d + \bB)^{-1}$ for $\|\bB\|_2 < 1$, and Result \ref{result:Jacobian}, 
\[
\left\{\bG_{0in} + O\left(\sqrt{\frac{(\log n)^{2\xi}}{n}}\right)\right\}^{-1} = \bG_{0in}^{-1} + O\left(\sqrt{\frac{(\log n)^{2\xi}}{n\rho_n^2}}\right). 
\]
It follows after simple simplification that
\[
\bW\transpose\widehat{\bx}_i - \rho_n^{1/2}\bx_{0i} = -\frac{1}{n}\sum_{j = 1}^n\bG_{0in}^{-1}\bg_{0ij} + O\left(\frac{(\log n)^{2\xi}}{n\rho_n^{1/2}}\right)\quad\mbox{w.h.p.}.
\]
The proof is completed by multiplying $\sqrt{n}$ on both sides of the above equation.
\end{proof}

\subsection{Proof of Theorem \ref{thm:BvM_generalized_posterior}}
\label{sub:proof_of_BvM}

\begin{proof}[\bf Proof of Theorem \ref{thm:BvM_generalized_posterior}]
Since $\bt = \sqrt{n}\bW\transpose(\bx_i - \widehat{\bx}_i)$, then $\bx_i = \widehat{\bx}_i + \bW\bt/\sqrt{n}$. Denote
\begin{align*}
d_{in} & = \int_{\mathbb{R}^d} \exp\left\{\ell_{in}\left(\widehat{\bx}_i + \frac{\bW\bt}{\sqrt{n}}\right) - \ell_{in}(\widehat{\bx}_i)\right\}\pi\left(\widehat{\bx}_i+\frac{\bW\bt}{\sqrt{n}}\right)\mathbbm{1}\left(\bt\in\widehat{\Theta}_i \right)\mathrm{d}\bt,
\end{align*}
where $\widehat{\Theta}_i = \{\bt:\bW\transpose\widehat{\bx}_i + n^{-1/2}\bt\in\Theta\}$. Note that $\sup_{\bt\in\widehat{\Theta}_i}\|\bt\|_2\lesssim \sqrt{n}$. 
Clearly, by definition, we have
\begin{align*}
\pi_{in}^*(\bt\mid\bA) & = \frac{1}{d_{in}}\exp\left\{\ell_{in}\left(\widehat{\bx}_i + \frac{\bW\bt}{\sqrt{n}}\right) - \ell_{in}(\widehat{\bx}_i)\right\}\pi\left(\widehat{\bx}_i+\frac{\bW\bt}{\sqrt{n}}\right)\mathbbm{1}\left(\bt\in\widehat{\Theta}_i\right).
\end{align*}
It is sufficient to show that
\begin{align}
% \begin{aligned}
&\int_{\mathbb{R}^d}(1 + \|\bt\|_2^\alpha)\left|\exp\left\{\ell_{in}\left(\widehat{\bx}_i + \frac{\bW\bt}{\sqrt{n}}\right) - \ell_{in}(\widehat{\bx}_i)\right\}\pi\left(\widehat{\bx}_i+\frac{\bW\bt}{\sqrt{n}}\right)\mathbbm{1}\left(\bt\in\widehat{\Theta}_i\right) - e^{-\bt\transpose\bSigma_{in}\bt/2}\pi(\rho_n^{1/2}\bW\bx_{0i})\right|\mathrm{d}\bt\nonumber\\
\label{eqn:BvM_sufficient_condition}
&\quad = o(1)\quad\mbox{w.p.a.1}.
% \end{aligned}
\end{align}
To see this, observe that the left-hand side of  \eqref{eqn:strong_convergence_generalized_posterior} can be written as
\begin{align*}
&\frac{1}{d_{in}}\int_{\mathbb{R}^d} (1 + \|\bt\|_2^\alpha)
\left|
\exp\left\{\ell_{in}\left(\widehat{\bx}_i + \frac{\bW\bt}{\sqrt{n}}\right) - \ell_{in}(\widehat{\bx}_i)\right\}\pi\left(\widehat{\bx}_i+\frac{\bW\bt}{\sqrt{n}}\right)\mathbbm{1}\left( \bt\in\widehat{\Theta}_i\right)
- \frac{d_{in}e^{-\bt\transpose\bSigma_{in}\bt/2}}{\det(2\pi\bSigma_{in}^{-1})^{1/2}}
\right|\mathrm{d}\bt\\
&\quad\leq \frac{1}{d_{in}}\int_{\mathbb{R}^d} (1 + \|\bt\|_2^\alpha)
\left|
\exp\left\{\ell_{in}\left(\widehat{\bx}_i + \frac{\bW\bt}{\sqrt{n}}\right) - \ell_{in}(\widehat{\bx}_i)\right\}\pi\left(\widehat{\bx}_i+\frac{\bW\bt}{\sqrt{n}}\right)\mathbbm{1}\left(\bt\in\widehat{\Theta}_i\right)
- e^{-\bt\transpose\bSigma_{in}\bt/2}\pi(\rho_n^{1/2}\bW\bx_{0i})
\right|\mathrm{d}\bt\\
&\quad\quad + \frac{1}{d_{in}}\left|\pi(\rho_n^{1/2}\bW\bx_{0i})
- \frac{d_{in}}{\det(2\pi\bSigma_{in}^{-1})^{1/2}}
\right|\int_{\mathbb{R}^d} (1 + \|\bt\|_2^\alpha)
e^{-\bt\transpose\bSigma_{in}\bt/2}\mathrm{d}\bt.
\end{align*}
Since \eqref{eqn:BvM_sufficient_condition} implies that $d_{in}=\det(2\pi\bSigma_{in}^{-1})^{1/2}\pi(\rho_n^{1/2}\bW\bx_{0i})+o(1)$ w.p.a.1. (by taking $\alpha=0$), it can be seen that \eqref{eqn:BvM_sufficient_condition} implies that the two terms on the right hand side of the previous display are $o(1)$ w.p.a.1.. Hence, we are left with establishing \eqref{eqn:BvM_sufficient_condition}. 

\vspace*{2ex}\noindent
Let $\eps_n,\delta_n$ be the sequences given by Assumption \ref{assumption:criterion_function} and consider the following partition:
\begin{align*}
\calA_1 &= \{\bt\in\widehat{\Theta}_i:\|\bt\|_2 \leq \sqrt{n}\eps_n/2\},\quad\calA_2 = \{\bt\in\widehat{\Theta}_i:\sqrt{n}\eps_n/2\leq \|\bt\|_2 <  2\sqrt{n}\delta_n\},\\
\calA_3 &= \{\bt\in\widehat{\Theta}_i:\|\bt\|_2> 2\sqrt{n}\delta_n\}.
\end{align*}
Let $\calA_4 = \mathbb{R}^d\backslash\widehat{\Theta}_i$. 
We first consider the integral of \eqref{eqn:BvM_sufficient_condition} over $\calA_4$. 
By Assumption \ref{assumption:regularity_condition}, there exists some $\eta > 0$, such that $B(\rho_n^{1/2}\bx_{0i}, \eta)\subset\Theta$. 
By Theorem \ref{thm:Large_sample_Z_estimator}, $\|\bW\transpose\widehat{\bx}_i - \rho_n^{1/2}\bx_{0i}\|_2 = o(1)$ w.h.p.. Then $\bt\notin \widehat{\Theta}_i$ implies that
\[
\eta\leq \left\|\bW\transpose\widehat{\bx}_i + \frac{\bt}{\sqrt{n}} - \rho_n^{1/2}\bx_{0i}\right\|_2\leq \|\bW\transpose\widehat{\bx}_i - \rho_n^{1/2}\bx_{0i}\|_2 + \frac{1}{\sqrt{n}}\|\bt\|_2. 
\]
Therefore $\|\bt\|_2\geq \sqrt{n}\eta/2$ w.h.p. when $\bt\notin\widehat{\Theta}_i$.
In this case, we have
\begin{align}
% \begin{aligned}
&\int_{\calA_4} (1 + \|\bt\|_2^\alpha)\left|\exp\left\{\ell_{in}\left(\widehat{\bx}_i + \frac{\bW\bt}{\sqrt{n}}\right) - \ell_{in}(\widehat{\bx}_i)\right\}\pi\left(\widehat{\bx}_i+\frac{\bW\bt}{\sqrt{n}}\right)\mathbbm{1}\left(\bt\in\widehat{\Theta}_i\right) - e^{-\bt\transpose\bSigma_{in}\bt/2}\pi(\rho_n^{1/2}\bW\bx_{0i})\right|\mathrm{d}\bt\nonumber\\
&\quad = \int_{\calA_4} (1 + \|\bt\|_2^\alpha)  e^{-\bt\transpose\bSigma_{in}\bt/2}\pi(\rho_n^{1/2}\bW\bx_{0i})\mathrm{d}\bt\nonumber
\\&\quad
\label{eqn:BvM_sufficient_condition_IV}
\leq \pi(\rho_n^{1/2}\bW\bx_{0i})\int_{\{\bt:\|\bt\|_2\geq\sqrt{n}\eta/2\}} (1 + \|\bt\|_2^\alpha) e^{-\lambda_d(\bSigma_{in})\|\bt\|_2^2/2}\mathrm{d}\bt = o(1)\quad\mbox{w.h.p.}.
% \end{aligned}
\end{align}
Now we turn to the integral of \eqref{eqn:BvM_sufficient_condition} over $\calA_3$. Recall from Theorem \ref{thm:Large_sample_Z_estimator} that $\|\bW\transpose\widehat{\bx}_i - \rho_n^{1/2}\bx_{0i}\|_2\leq \delta_n$ w.p.a.1. because $\sqrt{n}\delta_n\to + \infty$. Then $\bt\in\calA_3$ implies that
\[
\|\bW\transpose\bx_i - \rho_n^{1/2}\bx_{0i}\|\geq \frac{\|\bt\|_2}{\sqrt{n}} - \|\bW\transpose\widehat{\bx}_i - \rho_n^{1/2}\bx_{0i}\|_2\geq \delta_n\quad\mbox{w.p.a.1.}
\]
% Namely, $\bW\transpose\bx_i\notin B(\rho_n^{1/2}\bx_{0i}, \delta_n)$. 
By \eqref{eqn:identifiability} in Assumption \ref{assumption:criterion_function}, \ref{assumption:regularity_condition}, \ref{assumption:prior},
\begin{align}
% \begin{aligned}
&\int_{\calA_3} (1 + \|\bt\|_2^\alpha)\left|\exp\left\{\ell_{in}\left(\widehat{\bx}_i + \frac{\bW\bt}{\sqrt{n}}\right) - \ell_{in}(\widehat{\bx}_i)\right\}\pi\left(\widehat{\bx}_i+\frac{\bW\bt}{\sqrt{n}}\right)\mathbbm{1}\left(\bt\in\widehat{\Theta}_i\right) - e^{-\bt\transpose\bSigma_{in}\bt/2}\pi(\rho_n^{1/2}\bW\bx_{0i})\right|\mathrm{d}\bt\nonumber\\
&\quad\leq \int_{\calA_3}(1 + \|\bt\|_2^\alpha)\exp\left[\sup_{\bW\transpose\bx_i\notin B(\rho_n^{1/2}\bx_{0i}, \delta_n)}\left\{\ell_{in}\left(\bx_i\right) - \ell_{in}(\widehat{\bx}_i)\right\}\right]\pi\left(\widehat{\bx}_i+\frac{\bW\bt}{\sqrt{n}}\right)\mathrm{d}\bt\nonumber\\
&\quad\quad+ \int_{\calA_3}(1 + \|\bt\|_2^\alpha) e^{-\bt\transpose\bSigma_{in}\bt/2}\pi(\rho_n^{1/2}\bW\bx_{0i})\mathrm{d}\bt\nonumber\\
&\quad\lesssim 
% \left\{\sup_{\bx_i\in\Theta}\pi(\bx_i)\right\}
\int_{\{\bt:\|\bt\|_2\leq c\sqrt{n}\}}(1 + \|\bt\|_2^\alpha)\left(\frac{1}{n}\right)^{(1 + \alpha)d}\mathrm{d}\bt
% \nonumber\\&\quad\quad
+ 
% \pi(\rho_n^{1/2}\bW\bx_{0i})
\int_{\{\bt:\|\bt\|_2\geq 2\sqrt{n}\delta_n\}}(1 + \|\bt\|_2^\alpha) e^{-\bt\transpose\bSigma_{in}\bt/2}\mathrm{d}\bt\nonumber\\
% &\quad\leq n^{(1+\alpha)(1/2-d)}\int_{\Theta}\|\bx_i-\widehat{\bx}_i\|_2^\alpha\pi(\bx_i)\mathrm{d}\bx_i + \pi(\rho_n^{1/2}\bW\bx_{0i})\int_{\{\bt:\|\bt\|_2>2\sqrt{n}\delta_n\}}\|\bt\|_2^\alpha e^{-\bt\transpose\bSigma_{in}\bt/2}\mathrm{d}\bt\nonumber\\
\label{eqn:BvM_sufficient_condition_III}
&\quad = o(1)\quad\mbox{w.p.a.1.}
% \end{aligned}
\end{align}
for some constant $c > 0$. 

\vspace*{2ex}\noindent
It is now sufficient to consider the integral of \eqref{eqn:BvM_sufficient_condition} over $\calA_2$ and $\calA_1$. 
Because $\widehat{\bx}_i$ is the maximizer of $\ell_{in}(\bx_i)$ and is inside the interior of $\Theta$ with probability going to one, then by Taylor's theorem,
\begin{align*}
\ell_{in}\left(\widehat{\bx}_i + \frac{\bW\bt}{\sqrt{n}}\right) - \ell_{in}(\widehat{\bx}_i)
& = \frac{\partial \ell_{in}}{\partial\bx_i}(\widehat{\bx}_i)\transpose\frac{\bW\bt}{\sqrt{n}}
    + \frac{1}{2}\bt\transpose\bW\transpose\frac{1}{n}\frac{\partial^2\ell_{in}}{\partial\bx_i\partial\bx_i\transpose}(\bar{\bx}_i)\bW\bt = \frac{1}{2}\bt\transpose\bW\transpose\frac{1}{n}\frac{\partial^2\ell_{in}}{\partial\bx_i\partial\bx_i\transpose}(\bar{\bx}_i)\bW\bt,
\end{align*}
where there exists some $\theta_i(\bt)\in [0, 1]$ for each $i$ and $\bt$, such that $\bar{\bx}_i = \widehat{\bx}_i + \theta_i(\bt)\bW\bt/\sqrt{n}$. 
% because $\Theta$ has non-empty interior, implying that $\inf_{\bt\notin\widehat{\Theta}_i}\|\bt\|_2\to +\infty$. 
We next focus on the integral of \eqref{eqn:BvM_sufficient_condition} over $\calA_2$. By Theorem \ref{thm:Large_sample_Z_estimator}, for all $\bt\in\calA_2$, we have
\[
\|\bW\transpose\bar{\bx}_i - \rho_n^{1/2}\bx_{0i}\|_2\leq \|\bW\transpose\widehat{\bx}_i - \rho_n^{1/2}\bx_{0i}\|_2 + \frac{\|\bt\|_2}{\sqrt{n}}\leq 3\delta_n\quad\mbox{w.p.a.1.}
\]
Then by \eqref{eqn:Hessian_A2} in Assumption \ref{assumption:criterion_function}, for all $\bt\in\calA_2$,
\begin{align*}
\ell_{in}\left(\widehat{\bx}_i + \frac{\bW\bt}{\sqrt{n}}\right) - \ell_{in}(\widehat{\bx}_i)
&\leq-\frac{1}{2}\|\bt\|_2^2\lambda_{\min}\left\{-\frac{1}{n}\frac{\partial^2\ell_{in}}{\partial\bx_i\partial\bx_i\transpose}(\bar{\bx}_i)\right\}\\
&\leq-\frac{1}{2}\|\bt\|_2^2\inf_{\bW\transpose\bx_i\in B(\rho_n^{1/2}\bx_{0i}, 3\delta_n)}\lambda_{\min}\left\{
-\frac{1}{n}\frac{\partial^2\ell_{in}}{\partial\bx_i\partial\bx_i\transpose}({\bx}_i)
\right\}\leq - C\|\bt\|_2^2\quad\mbox{w.p.a.1,}
\end{align*}
where $C$ is some constant independent of $n$. It follows that
\begin{align}
% \begin{aligned}
&\int_{\calA_2} (1 + \|\bt\|_2^\alpha)\left|\exp\left\{\ell_{in}\left(\widehat{\bx}_i + \frac{\bW\bt}{\sqrt{n}}\right) - \ell_{in}(\widehat{\bx}_i)\right\}\pi\left(\widehat{\bx}_i+\frac{\bW\bt}{\sqrt{n}}\right)\mathbbm{1}\left(\bt\in\widehat{\Theta}_i\right) - e^{-\bt\transpose\bSigma_{in}\bt/2}\pi(\rho_n^{1/2}\bW\bx_{0i})\right|\mathrm{d}\bt\nonumber\\
\label{eqn:BvM_sufficient_condition_II}
&\quad\lesssim \int_{\{\|\bt\|_2 > \sqrt{n}\eps_n/2\}}(1 + \|\bt\|_2^\alpha) e^{-C\|\bt\|_2^2}\mathrm{d}\bt
% \nonumber\\
% &\quad\quad+ 
% \pi(\rho_n^{1/2}\bW\bx_{0i})
+ \int_{\{\|\bt\|_2 > \sqrt{n}\eps_n/2\}}(1 + \|\bt\|_2^\alpha) e^{-\bt\transpose\bSigma_{in}\bt/2}\mathrm{d}\bt
 = o(1)\quad\mbox{w.p.a.1.}
% \end{aligned}
\end{align}
We finally consider the integral of \eqref{eqn:BvM_sufficient_condition} over $\calA_1$. For all $\bt\in \calA_1$, by Theorem \ref{thm:Large_sample_Z_estimator}, we have
\[
\|\bW\transpose\by_i - \rho_n^{1/2}\bx_{0i}\|_2\leq \|\bW\transpose\widehat{\bx}_i - \rho_n^{1/2}\bx_{0i}\|_2 + \frac{\|\bt\|_2}{\sqrt{n}}\leq \eps_n\quad\mbox{w.p.a.1.} 
\]
for both $\by_i=\bar{\bx}_i$ and $\by_i=\bx_i\equiv\widehat\bx_i+\bW\bt/\sqrt{n}$. Denote
\[
\bR_{in}^{(\ell)}(\bx_i) = \frac{1}{n}\bW\transpose\frac{\partial^2\ell_{in}}{\partial\bx_i\partial\bx_i\transpose}(\bx_i)\bW + \bSigma_{in}.
\]
Then by \eqref{eqn:Hessian_A1} in Assumption \ref{assumption:criterion_function},
\begin{align*}
\sup_{\bt\in \calA_1}\|\bR_{in}^{(\ell)}(\bar{\bx}_i)\|_2\leq \sup_{\bW\transpose\bx_i\in B(\rho_n^{1/2}\bx_{0i}, \eps_n)}\|\bR_{in}^{(\ell)}({\bx}_i)\|_2 \leq \frac{\gamma_n}{n\eps_n^2}\quad\mbox{w.p.a.1.},
\end{align*}
where $(\gamma_n)_{n = 1}^\infty$ is a positive sequence converging to $0$. It follows that
\begin{align*}
&\int_{\calA_1}(1 + \|\bt\|_2^\alpha)\left|\exp\left\{\ell_{in}\left(\widehat{\bx}_i + \frac{\bW\bt}{\sqrt{n}}\right) - \ell_{in}(\widehat{\bx}_i)\right\}\pi\left(\widehat{\bx}_i+\frac{\bW\bt}{\sqrt{n}}\right)\mathbbm{1}(\bt\in\widehat\Theta_i) - e^{-\bt\transpose\bSigma_{in}\bt/2}\pi(\rho_n^{1/2}\bW\bx_{0i})\right|\mathrm{d}\bt\\
& \quad = \int_{\calA_1}(1 + \|\bt\|_2^\alpha)\left|\exp\left\{ -\frac{1}{2}\bt\transpose\bSigma_{in}\bt + \frac{1}{2}\bt\transpose\bR_{in}^{(\ell)}(\bar{\bx}_i)\bt\right\}\pi\left(\widehat{\bx}_i+\frac{\bW\bt}{\sqrt{n}}\right) - e^{-\bt\transpose\bSigma_{in}\bt/2}\pi(\rho_n^{1/2}\bW\bx_{0i})\right|\mathrm{d}\bt\\
& \quad = \int_{\calA_1}(1 + \|\bt\|_2^\alpha)\left|\exp\left\{\frac{1}{2}\bt\transpose\bR_{in}^{(\ell)}(\bar{\bx}_i)\bt\right\} - \frac{\pi(\rho_n^{1/2}\bW\bx_{0i})}{\pi(\widehat{\bx}_i+{\bW\bt}/{\sqrt{n}})}\right|e^{-\bt\transpose\bSigma_{in}\bt/2}\pi\left(\widehat{\bx}_i+\frac{\bW\bt}{\sqrt{n}}\right)\mathrm{d}\bt\\
& \quad \leq \left\{\left|\exp\left\{\frac{1}{2}\sup_{\|\bt\|_2\leq \sqrt{n}\eps_n}\|\bt\|_2^2\sup_{\bt\in\calA_1}\left\|\bR_{in}^{(\ell)}\left(\bar{\bx}_i\right)\right\|_2\right\} - 1\right| + \sup_{\bW\transpose\bx_i\in B(\rho_n^{1/2}\bx_{0i},\epsilon_n)}\left|1-\frac{\pi(\rho_n^{1/2}\bW\bx_{0i})}{\pi(\widehat{\bx}_i+{\bW\bt}/{\sqrt{n}})}\right|\right\}\\
&\quad\quad \times\int_{\calA_1} (1 + \|\bt\|_2^\alpha) e^{-\bt\transpose\bSigma_{in}\bt/2}\pi\left(\widehat{\bx}_i+\frac{\bW\bt}{\sqrt{n}}\right)\mathrm{d}\bt\\
& \quad \leq\left\{\left|e^{\gamma_n/2} - 1\right| + \sup_{\bW\transpose\bx_i\in B(\rho_n^{1/2}\bx_{0i},\epsilon_n)}C_1\left|\pi\left(\widehat{\bx}_i+\frac{\bW\bt}{\sqrt{n}}\right)-\pi(\rho_n^{1/2}\bW\bx_{0i})\right|\right\}C_2\int (1 + \|\bt\|_2^\alpha) e^{-\bt\transpose\bSigma_{in}\bt/2}\mathrm{d}\bt\\
& \quad = o(1)\quad\mbox{w.p.a.1.}
\end{align*}
for some constants $C_1, C_2 > 0$. 
This shows that 
\begin{equation}\label{eqn:BvM_sufficient_condition_I}
\begin{aligned}
&\int_{\calA_1}(1 + \|\bt\|_2^\alpha)\left|\exp\left\{\ell_{in}\left(\widehat{\bx}_i + \frac{\bW\bt}{\sqrt{n}}\right) - \ell_{in}(\widehat{\bx}_i)\right\}\pi\left(\widehat{\bx}_i+\frac{\bW\bt}{\sqrt{n}}\right) - e^{-\bt\transpose\bSigma_{in}\bt/2}\pi(\rho_n^{1/2}\bW\bx_{0i})\right|\mathrm{d}\bt\\
&\quad = o(1)\quad\mbox{w.p.a.1.}.
\end{aligned}
\end{equation}
The proof of \eqref{eqn:BvM_sufficient_condition} is thus completed by combining \eqref{eqn:BvM_sufficient_condition_IV}, \eqref{eqn:BvM_sufficient_condition_III}, \eqref{eqn:BvM_sufficient_condition_II}, and \eqref{eqn:BvM_sufficient_condition_I}. 
\end{proof}

\subsection{Proof of Theorem \ref{thm:GBE}}
\label{sub:proof_of_GBE}

\begin{proof}[\bf Proof of Theorem \ref{thm:GBE}]
We observe that with $\alpha = 1$, Theorem \ref{thm:BvM_generalized_posterior} implies
\begin{align*}
\|\sqrt{n}(\bx_i^* - \widehat{\bx}_i)\|_2
& = \left\|
\int_\Theta\sqrt{n}\bW\transpose(\bx_i - \widehat{\bx}_i)\pi_{in}(\bx_i\mid\bA)\mathrm{d}\bx_i
\right\|_2\\
& = \left\|
\int\bt{\pi}_{in}^*(\bt\mid\bA)\mathrm{d}\bt - \int\bt \det(2\pi\bSigma_{in}^{-1})^{-1/2}e^{-\bt\transpose\bSigma_{in}\bt/2}\mathrm{d}\bt
\right\|_2\\
&\leq \int \|\bt\|_2\left|\pi_{in}^*(\bt\mid\bA) - \det(2\pi\bSigma_{in}^{-1})^{-1/2}e^{-\bt\transpose\bSigma_{in}\bt/2}\right|\mathrm{d}\bt\\
&\leq \int (1 + \|\bt\|_2)\left|{\pi}_{in}^*(\bt\mid\bA) - \det(2\pi\bSigma_{in}^{-1})^{-1/2}e^{-\bt\transpose\bSigma_{in}\bt/2}\right|\mathrm{d}\bt\overset{\prob_0}{\to}0.
\end{align*}
Namely, $\sqrt{n}(\bx_i^* - \widehat{\bx}_i) = o_{\prob_0}(1)$. 
Then the asymptotic normality of $\bOmega_{0in}^{-1/2}\bG_{0in}\sqrt{n}(\bW\transpose\bx_i^* - \rho_n^{1/2}\bx_{0i})$ directly follows from Theorem \ref{thm:Large_sample_Z_estimator}, the result that $\sqrt{n}(\bx_i^* - \widehat{\bx}_i) = o_{\prob_0}(1)$, Result \ref{result:Jacobian}, and Result \ref{result:second_moment_matrix}. 
\end{proof}

\subsection{Proof of Theorem \ref{thm:GBI}}
\label{sub:proof_of_GBI}
\begin{proof}[\bf Proof of Theorem \ref{thm:GBI}]
We first show that $n\bW\transpose \widehat{\bV}_B\bW = \bSigma_{in}^{-1} + o_{\prob_0}(1)$. Denote $\expect_\bt$ the expected value with regard to the posterior distribution of $\bt = \sqrt{n}\bW\transpose(\bx_i - \widehat{\bx}_i)$. Let $\bx_i^*$ denote the posterior mean of $\bx_i$. From the proof of Theorem \ref{thm:GBE}, we know that $\expect_\bt\bt = o_{\prob_0}(1)$. 
By definition of $\widehat{\bV}_B$, we have
\begin{align*}
n\bW\transpose \widehat{\bV}_B\bW
& = \bW\transpose\int n\bx_i\bx_i\transpose\pi_{in}(\bx_i\mid \bA)\mathrm{d}\bx_i\bW - \bW\transpose(\sqrt{n}\bx_i^*)(\sqrt{n}\bx_i^*)\transpose\bW\\
& = \bW\transpose\left[\expect_\bt\{(\sqrt{n}\widehat{\bx}_i + \bW\bt)(\sqrt{n}\widehat{\bx}_i + \bW\bt)\transpose\}
- \{\sqrt{n}\widehat{\bx}_i + \bW(\expect_\bt\bt)\}\{\sqrt{n}\widehat{\bx}_i + \bW(\expect_\bt\bt)\}\transpose\right]\bW\\
& = \bW\transpose\left[
n\widehat{\bx}_i\widehat{\bx}_i\transpose + \bW(\expect_\bt\bt)\sqrt{n}\widehat{\bx}_i\transpose + \sqrt{n}\widehat{\bx}_i(\expect_\bt\bt)\transpose\bW\transpose + \bW\expect_\bt(\bt\bt\transpose)\bW\transpose\right]\bW\\
&\quad
- \bW\transpose\{n\widehat{\bx}_i\widehat{\bx_i}\transpose + \bW(\expect_\bt\bt)\sqrt{n}\widehat{\bx}_i\transpose + \sqrt{n}\widehat{\bx}_i(\expect_\bt\bt)\transpose\bW\transpose + \bW(\expect_\bt\bt)(\expect_\bt\bt)\transpose\bW\transpose
\}\bW\\
& = \expect_\bt(\bt\bt\transpose) + o_{\prob_0}(1)\\
& = \left\{\int \bt\bt\transpose\pi_{in}^*(\bt\mid\bA)\mathrm{d}\bt - \bSigma_{in}^{-1}\right\} + \bSigma_{in}^{-1} + o_{\prob_0}(1)\\
& = \left\{\int \bt\bt\transpose\pi_{in}^*(\bt\mid\bA)\mathrm{d}\bt - \int \bt\bt\transpose\frac{\exp(-\bt\transpose\bSigma_{in}\bt/2)}{\det(2\pi\bSigma_{in}^{-1})}\mathrm{d}\bt\right\} + \bSigma_{in}^{-1} + o_{\prob_0}(1)\\
& = \int \bt\bt\transpose\left\{\pi_{in}^*(\bt\mid\bA) - \frac{\exp(-\bt\transpose\bSigma_{in}\bt/2)}{\det(2\pi\bSigma_{in}^{-1})}\right\}\mathrm{d}\bt + \bSigma_{in}^{-1} + o_{\prob_0}(1).
\end{align*}
By Theorem \ref{thm:BvM_generalized_posterior}, we have
\begin{align*}
\left\|\int \bt\bt\transpose\left\{\pi_{in}^*(\bt\mid\bA) - \frac{\exp(-\bt\transpose\bSigma_{in}\bt/2)}{\det(2\pi\bSigma_{in}^{-1})}\right\}\mathrm{d}\bt\right\|_2
&\leq \int \|\bt\|_2^2\left\{\pi_{in}^*(\bt\mid\bA) - \frac{\exp(-\bt\transpose\bSigma_{in}\bt/2)}{\det(2\pi\bSigma_{in}^{-1})}\right\}\mathrm{d}\bt = o_{\prob_0}(1).
\end{align*}
Therefore, we conclude that $n\bW\transpose \widehat{\bV}_B\bW = \bSigma_{in}^{-1} + o_{\prob_0}(1)$, and hence, $(\bW\transpose \widehat{\bV}_B\bW)^{-1} = n\bSigma_{in} + o_{\prob_0}(n)$. Now let $\bZ_{in}\sim N(\rho_n^{1/2}\bx_{0i}, (n\bSigma_{in})^{-1})$. Then $\bW\transpose\widehat\bx_i = \bZ_{in}+o_{\prob_0}(n^{-1/2})$, and $\bZ_{in}=O_{\prob_0}(n^{-1/2})$. Therefore,
\begin{align*}
&\mathbb{P}_0\left\{\rho_n^{1/2}\bW\bx_{0i}\in\calE_{in}\right\} \\
&\quad= \mathbb{P}_0\left\{(\rho_n^{1/2}\bW\bx_{0i}-\widehat\bx_i)\transpose \widehat\bV_B^{-1} (\rho_n^{1/2}\bW\bx_{0i}-\widehat\bx_i) \leq q_{(1-\alpha)}\right\} \\
&\quad= \mathbb{P}_0\left\{(\bW\transpose\widehat\bx_i - \rho_n^{1/2}\bx_{0i})\transpose (\bW\transpose\widehat\bV_B\bW)^{-1} (\bW\transpose\widehat\bx_i - \rho_n^{1/2}\bx_{0i}) \leq q_{(1-\alpha)}\right\} \\
&\quad= \mathbb{P}_0\left\{(\bZ_{in} - \rho_n^{1/2}\bx_{0i} + o_{\prob_0}(n^{-1/2}))\transpose \{n\bSigma_{in}+o_{\prob_0}(n)\} (\bZ_{in} - \rho_n^{1/2}\bx_{0i} + o_{\prob_0}(n^{-1/2})) \leq q_{(1-\alpha)}\right\} \\
&\quad= \mathbb{P}_0\left[\{\sqrt{n}(\bZ_{in} - \rho_n^{1/2}\bx_{0i}) + o_{\prob_0}(1)\}\transpose \bSigma_{in}\{\sqrt{n}(\bZ_{in} - \rho_n^{1/2}\bx_{0i}) + o_{\prob_0}(1)\} + o_{\prob_0}(1)\leq q_{(1-\alpha)}\right\} 
% \\
% &\quad= F_{\chi^2_d + o_p(1)}(q_{(1-\alpha)}) 
% \\&\quad
\to 
% F_{\chi^2_d}(q_{(1-\alpha)}) 
% \\&\quad
= 1-\alpha.
\end{align*}
The proof is thus completed.
% where $F_{\chi^2_d + o_p(1)}(q_{(1-\alpha)})\to F_{\chi^2_d}(q_{(1-\alpha)})$ follows from Polya's Theorem.
\end{proof}

\subsection{Proof of Proposition \ref{prop:optimal_weighting}}
\label{sub:proof_of_prop_optimal_weighting}
\begin{proof}[\bf Proof of Proposition \ref{prop:optimal_weighting}]
Denote $h_{0nij}=h_n(\rho_n\bx_{0i}\transpose\bx_{0j},\rho_n\bx_{0i}\transpose\bx_{0j})$, $\bH_{in}=\mathrm{diag}\{h_{0ni1},\ldots,h_{0nin}\}$,, and $\bV_{in}=\mathrm{diag}\{\var_0(E_{i1}), \ldots,\var_0(E_{in)}\}$. By simple calculation, we have
\[
\bG_{0in}=-\frac{\rho_n^{1/2}}{n}\bX_0\transpose\bH_{in}\bX_0
\quad\mbox{and}\quad
\bOmega_{0in}=\frac{1}{n}\bX_0\transpose\bH_{in}\bV_{in}\bH_{in}\bX_0.
\]
Now
\[
\bG_{0in}^{-1}\bOmega_{0in}\bG_{0in}^{-1} = n\rho_n^{-1}(\bX_0\transpose\bH_{in}\bX_0)^{-1}\bX_0\transpose\bH_{in}\bV_{in}\bH_{in}\bX_0(\bX_0\transpose\bH_{in}\bX_0)^{-1},
\]
and
\[
\left\{\frac{1}{n}\sum_{j=1}^n\frac{\rho_n\bx_{0j}\bx_{0j}\transpose}{\var_0(E_{ij})}\right\}^{-1} = n\rho_n^{-1}(\bX_0\transpose\bV_{in}^{-1}\bX_0)^{-1}.
\]
Then
\begin{align*}
&\bG_{0in}^{-1}\bOmega_{0in}\bG_{0in}^{-1} - \left\{\frac{1}{n}\sum_{j=1}^n\frac{\rho_n\bx_{0j}\bx_{0j}\transpose}{\var_0(E_{ij})}\right\}^{-1} \\
&\quad= n\rho_n^{-1}(\bX_0\transpose\bH_{in}\bX_0)^{-1}\bX_0\transpose\bH_{in}\bV_{in}\bH_{in}\bX_0(\bX_0\transpose\bH_{in}\bX_0)^{-1} - n\rho_n^{-1}(\bX_0\transpose\bV_{in}^{-1}\bX_0)^{-1} \\
&\quad= n\rho_n^{-1} (\bX_0\transpose\bH_{in}\bX_0)^{-1}\left(\bX_0\transpose\bH_{in}\bV_{in}\bH_{in}\bX_0 - \bX_0\transpose\bH_{in}\bX_0 (\bX_0\transpose\bV_{in}^{-1}\bX_0)^{-1} \bX_0\transpose\bH_{in}\bX_0\right)(\bX_0\transpose\bH_{in}\bX_0)^{-1}\\
&\quad= n\rho_n^{-1} (\bX_0\transpose\bH_{in}\bX_0)^{-1}\bX_0\transpose\bH_{in}\bV_{in}^{1/2}\left(\eye_n - \bV_{in}^{-1/2}\bX_0(\bX_0\transpose\bV_{in}^{-1}\bX_0)^{-1}\bX_0\transpose\bV_{in}^{-1/2}\right) \bV_{in}^{1/2}\bH_{in}\bX_0(\bX_0\transpose\bH_{in}\bX_0)^{-1}\\
&\quad= n\rho_n^{-1} (\bX_0\transpose\bH_{in}\bX_0)^{-1}\bX_0\transpose\bH_{in}\bV_{in}^{1/2}\left(\eye_n - \bV_{in}^{-1/2}\bX_0[(\bV_{in}^{-1/2}\bX_0)\transpose\bV_{in}^{-1/2}\bX_0]^{-1}(\bV_{in}^{-1/2}\bX_0)\transpose\right) \\
&\quad\quad\times \bV_{in}^{1/2}\bH_{in}\bX_0(\bX_0\transpose\bH_{in}\bX_0)^{-1}\\
&\quad= n\rho_n^{-1} \bS\left(\eye_n - \bP_{\bV_{in}^{-1/2}\bX_0}\right)\bS\transpose
% \\&\quad= n\rho_n^{-1} \bS\left(\eye_n - \bP_{\bV_{in}^{-1/2}\bX_0}\right) \left[\bS\left(\eye_n - \bP_{\bV_{in}^{-1/2}\bX_0}\right)\right]\transpose \\
% &\quad
\succeq \bf{0}_{d\times d},
\end{align*}
where $\bS = (\bX_0\transpose\bH_{in}\bX_0)^{-1}\bX_0\transpose\bH_{in}\bV_{in}^{1/2}$, and $\bP_{\bV_{in}^{-1/2}\bX_0} = \bV_{in}^{-1/2}\bX_0[(\bV_{in}^{-1/2}\bX_0)\transpose\bV_{in}^{-1/2}\bX_0]^{-1}(\bV_{in}^{-1/2}\bX_0)\transpose$ is the projection matrix onto the subspace spanned by the columns of $\bV_{in}^{-1/2}\bX_0$.
\end{proof}

\section{Proof of Proposition \ref{prop:Criterion_satisfies_assumption}}
\label{sub:proof_of_criterion_proposition}

The proof of Proposition \ref{prop:Criterion_satisfies_assumption} is lengthy and quite technical. We breakdown the proof for the $M$-criterion, the GMM criterion, and the ETEL criterion into Subsection \ref{sub:proof_of_prop_M_criterion}, Subsection \ref{sub:proof_of_prop_GMM_criterion}, and Subsection \ref{sub:proof_of_prop_ETEL_criterion}, respectively. 

\subsection{Proof of Proposition \ref{prop:Criterion_satisfies_assumption} (a)}
\label{sub:proof_of_prop_M_criterion}

\begin{proof}[\bf Proof of Proposition \ref{prop:Criterion_satisfies_assumption} (a)]
Let $\eps_n = (\log n)^{1/4}/\sqrt{n}$, $M_n = \log\log n$, and $\delta_n = M_n\{(\log n)^{2\xi}/(n\rho_n)\}^{1/4}$. Denote $\bG_{0in} = \bG_{in}(\rho_n^{1/2}\bx_{0i})$. By construction, we have 
\begin{align*}
\bW\transpose\frac{1}{n}\frac{\partial^2\ell_{in}}{\partial\bx_i\partial\bx_i\transpose}(\bW\bx_i)\bW
& = \frac{1}{n\rho_n^{1/2}}\sum_{j = 1}^n\bW\transpose\frac{\partial\widetilde{\bg}_{ij}}{\partial\bx_i}(\bW\bx_i)\bW\\
& = \rho_n^{-1/2}\bG_{0in} + \rho_n^{-1/2}\{\bG_{in}(\bx_i) - \bG_{0in}\}
\\&\quad
+ \rho_n^{-1/2}\left\{\frac{1}{n}\sum_{j = 1}^n\bW\transpose\frac{\partial\widetilde{\bg}_{ij}}{\partial\bx_i}(\bW\bx_i)\bW - \bG_{in}(\bx_i)\right\}.
\end{align*}
Denote $h_{0nij}(\bx_i) = h_n(\rho_n\bx_{0i}\transpose\bx_{0j}, \rho_n^{1/2}\bx_{i}\transpose\bx_{0j})$ and $D^{(0, 1)}h_{0nij}(\bx_i) = D^{(0, 1)} h_n(\rho_n\bx_{0i}\transpose\bx_{0j}, \rho_n^{1/2}\bx_{i}\transpose\bx_{0j})$ with a slight abuse of notations. 
By definition of $\bG_{in}(\bx_i)$, 
\begin{align*}
&\sup_{\bx_i\in B(\rho_n^{1/2}\bx_{0i}, 3\delta_n)}\rho_n^{-1/2}\|\bG_{in}(\bx_i) - \bG_{0in}\|_2\\
&\quad\leq 
\sup_{\bx_i\in B(\rho_n^{1/2}\bx_{0i}, 3\delta_n)}\frac{1}{n}\sum_{j = 1}^n|h_{0nij}(\bx_i) - h_{0nij}(\rho_n^{1/2}\bx_{0i})|\|\bx_{0j}\|_2^2\\
&\quad\quad + \sup_{\bx_i\in B(\rho_n^{1/2}\bx_{0i}, 3\delta_n)}\rho_n^{1/2}\frac{1}{n}\sum_{j = 1}^n\|\rho_n^{1/2}\bx_{0i} - \bx_{i}\|_2\|\bx_{0j}\|_2|D^{(0, 1)}h_{0nij}(\bx_i)|\|\bx_{0j}\|_2^2\\
&\quad\lesssim \rho_n^{3/2}\delta_n \leq \rho_n^{1/2}M_n\rho_n^{1/4}\left\{\frac{(\log n)^{2\xi}}{n}\right\}^{1/4}.
\end{align*}
We then obtain from Lemma \ref{lemma:ULLN} that
\begin{align*}
\sup_{\bx_i\in B(\rho_n^{1/2}\bx_{0i}, 3\delta_n)}\left\|\bW\transpose\frac{1}{n}\frac{\partial^2\ell_{in}}{\partial\bx_i\partial\bx_i\transpose}(\bW\bx_i)\bW - \rho_n^{-1/2}\bG_{0in}\right\|_2\lesssim M_n\left\{\frac{(\log n)^{2\xi}}{n\rho_n}\right\}^{1/4} = o\left(\frac{1}{n\eps_n^2}\right)\quad\mbox{w.h.p.}.
\end{align*}
Since $\eps_n\leq \delta_n$ and the eigenvalues of $-\rho_n^{-1/2}\bG_{0in}$ are bounded away from $0$ and $\infty$, this completes the proof of \eqref{eqn:Hessian_A1} and \eqref{eqn:Hessian_A2} in Assumption \ref{assumption:criterion_function} simulatenously. We now focus on the verification of \eqref{eqn:identifiability}. Without loss of generality, we can take $t_0 = 0$. Denote
\[
M_{in}(\bx_i) = \frac{1}{n\rho_n}\sum_{j = 1}^n\int_0^{\rho_n^{1/2}\bx_i\transpose\bx_{0j}}(\rho_n\bx_{0i}\transpose\bx_{0j} - t)h_{n}(\rho_n\bx_{0i}\transpose\bx_{0j}, t)\mathrm{d}t. 
\]
By triangle inequality, we have
\begin{align*}
\left|\frac{1}{n}\ell_{in}(\bW\bx_i) - M_{in}(\bx_i)\right|
& \leq\left|\frac{1}{n\rho_n}\sum_{j= 1}^n\int_{\rho_n^{1/2}\bx_i\transpose\bx_{0j}}^{\bx_i\transpose\bW\transpose\widetilde{\bx}_{j}}(A_{ij} - t)h_n(\widetilde{\bx}_i\transpose\widetilde{\bx}_j, t)\mathrm{d}t\right|\\
&\quad + \left|\frac{1}{n\rho_n}\sum_{j= 1}^n\int_0^{\rho_n^{1/2}\bx_i\transpose\bx_{0j}}(A_{ij} - t)\{h_n(\widetilde{\bx}_i\transpose\widetilde{\bx}_j, t) - h_{n}(\rho_n\bx_{0i}\transpose\bx_{0j}, t)\}\mathrm{d}t\right|\\
&\quad + \left|\frac{1}{n\rho_n}\sum_{j= 1}^n(A_{ij} - \rho_n\bx_{0i}\transpose\bx_{0j})\int_0^{\rho_n^{1/2}\bx_i\transpose\bx_{0j}}h_{n}(\rho_n\bx_{0i}\transpose\bx_{0j}, t)\mathrm{d}t\right|.
% \bx_i\transpose(\bW\transpose\widetilde{\bx}_j - \rho_n^{1/2}\bx_{0j})
\end{align*}
For the first term, by the mean-value theorem, for each $j\in [n]$, there exists some $t_{ij}(\bx_i)\in [-r,r]$ adjoining $\rho_n^{1/2}\bx_i\transpose\bx_{0j}$ and $\bx_i\transpose\bW\transpose\widetilde{\bx}_j$, such that
\begin{align*}
&\sup_{\bx_i\in \Theta}\left|\frac{1}{n\rho_n}\sum_{j= 1}^n\int_{\rho_n^{1/2}\bx_i\transpose\bx_{0j}}^{\bx_i\transpose\bW\transpose\widetilde{\bx}_{j}}(A_{ij} - t)h_n(\widetilde{\bx}_i\transpose\widetilde{\bx}_j, t)\mathrm{d}t\right|\\
&\quad\leq\sup_{\bx_i\in \Theta}\frac{1}{n\rho_n}\sum_{j = 1}^n|\bx_i\transpose(\bW\transpose\widetilde{\bx}_j - \rho_n^{1/2}\bx_{0j})||A_{ij}- t_{ij}(\bx_i)|h_n(\widetilde{\bx}_i\transpose\widetilde{\bx}_j, t_{ij}(\bx_i))\\
&\quad\lesssim \sup_{\bx_i\in \Theta}\frac{1}{n\rho_n }\|\bx_i\|_2 \max_{j\in [n]}\|\bW\transpose\widetilde{\bx}_j - \rho_n^{1/2}\bx_{0j}\|_2 \left\{\|\bA\|_\infty + \sum_{j = 1}^n\|\bx_i\|_2(\|\rho_n^{1/2}\bx_{0j}\|_2 + \|\bW\transpose\widetilde{\bx}_j\|_2)\right\}\\
&\quad\lesssim \sqrt{
\frac{(\log n)^{2\xi}}{n\rho_n}
}\quad\mbox{w.h.p.}.
\end{align*}
For the second term, by the mean-value theorem, Assumption \ref{assumption:weight_functions}, and Result \ref{result:pij_tilde_concentration}, for each $j\in [n]$, there exists some $t_{ij}(\bx_i)$ adjoining $\rho_n^{1/2}\bx_i\transpose\bx_{0j}$ and $0$, such that
\begin{align*}
&\sup_{\bx_i\in \Theta}\left|\frac{1}{n\rho_n}\sum_{j = 1}^n\int_0^{\rho_n^{1/2}\bx_i\transpose\bx_{0j}}(A_{ij} - t)\{h_n(\widetilde{\bx}_i\transpose\widetilde{\bx}_j, t) - h_{n}(\rho_n\bx_{0i}\transpose\bx_{0j}, t)\}\mathrm{d}t\right|\\
&\quad\leq \frac{1}{n\rho_n}\sum_{j = 1}^n\left(|A_{ij}| + \rho_n^{1/2}\sup_{\bx_i\in \Theta}\|\bx_i\|_2\|\bx_{0j}\|_2\right)\sup_{\bx_i\in \Theta}\int_0^{\rho_n^{1/2}\bx_i\transpose\bx_{0j}}|h_n(\widetilde{\bx}_i\transpose\widetilde{\bx}_j, t) - h_{n}(\rho_n\bx_{0i}\transpose\bx_{0j}, t)|\mathrm{d}t\\
&\quad = \frac{1}{n\rho_n}\sum_{j = 1}^n\left(|A_{ij}| + \rho_n^{1/2}\sup_{\bx_i\in \Theta}\|\bx_i\|_2\|\bx_{0j}\|_2\right)\sup_{\bx_i\in \Theta}|\rho_n^{1/2}\bx_i\transpose\bx_{0j}||h_n(\widetilde{\bx}_i\transpose\widetilde{\bx}_j, t_{ij}(\bx_i)) - h_{n}(\rho_n\bx_{0i}\transpose\bx_{0j}, t_{ij}(\bx_i))|\\
&\quad\lesssim \frac{1}{n\rho_n}(n\rho_n + n\rho_n^{1/2})\rho_n^{-1/2}\max_{j\in [n]}|\widetilde{\bx}_i\transpose\widetilde{\bx}_j - \rho_n\bx_{0i}\transpose\bx_{0j}|\lesssim \sqrt{
\frac{(\log n)^{2\xi}}{n\rho_n}
}\quad\mbox{w.h.p.}.
\end{align*}
For the third term, we apply the maximal inequality for sub-Gaussian processes. Define the function $\kappa_{nij}(\bx_i) = \int_0^{\rho_n^{1/2}\bx_i\transpose\bx_{0j}}h_n(\rho_n\bx_{0i}\transpose\bx_{0j}, t)\mathrm{d}t$ and the stochastic process $J_{in}(\bx_i) = \sum_{j = 1}^n(A_{ij} - \rho_n\bx_{0i}\transpose\bx_{0j})\kappa_{nij}(\bx_i)$. 
By Assumption \ref{assumption:regularity_condition}, $|\kappa_{nij}(\bx_i) - \kappa_{nij}(\bx_i')|\lesssim \rho_n^{1/2}\|\bx_i - \bx_i'\|_2$. Observe that $(A_{ij} - \expect_0A_{ij})_{j = 1}^n$ are uniformly bounded in sub-Gaussian norms. 
Then by Proposition 5.10 in \cite{vershynin2010introduction}, for any $t > 0$
\begin{align*}
\prob_0\left\{|J_{in}(\bx_i) - J_{in}(\bx_i')| > t\right\}\leq e\exp\left\{-\frac{t^2}{Cn\rho_n\|\bx_i - \bx_i'\|_2^2}\right\}.
\end{align*}
Namely, $\{(n\rho_n)^{-1/2}J_{in}(\bx_i):\bx_i\in\Theta\}$ is a sub-Gaussian process with respect to the distance $C\|\cdot\|_2$ for some constant $C > 0$. By Theorem 8.4 in \cite{kosorok2008introduction}, 
\[
\left\|\sup_{\bx_i\in\Theta}\frac{1}{n\rho_n}\sum_{j = 1}^n(A_{ij} - \rho_n\bx_{0i}\transpose\bx_{0j})\kappa_{nij}(\bx_i)\right\|_{\psi_2}\lesssim \frac{1}{\sqrt{n\rho_n}}.
\]
By Lemma 8.1 in \cite{kosorok2008introduction}, 
\[
\left|\sup_{\bx_i\in\Theta}\frac{1}{n\rho_n}\sum_{j = 1}^n(A_{ij} - \rho_n\bx_{0i}\transpose\bx_{0j})\kappa_{nij}(\bx_i)\right|\lesssim \sqrt{\frac{\log n}{n\rho_n}}\quad\mbox{w.h.p.}.
\]
We conclude from the three pieces of the concentration bounds obtained earlier that
\[
\sup_{\bx_i\in\Theta}\left|\frac{1}{n}\ell_{in}(\bW\bx_i) - M_{in}(\bx_i)\right|\lesssim \sqrt{\frac{(\log n)^{2\xi}}{n\rho_n}}\quad\mbox{w.h.p.}.
\]
Note that by Assumption \ref{assumption:weight_functions} and Assumption \ref{assumption:signal_plus_noise} (ii), for any $\bx_i\notin B(\rho_n^{1/2}\bx_{0i}, \delta_n)$,
\begin{align*}
&\inf_{\bx_i\notin B(\rho_n^{1/2}\bx_{0i}, \delta_n)}\{M_{in}(\rho_n^{1/2}\bx_{0i}) - M_{in}(\bx_i)\}\\
&\quad = \inf_{\bx_i\notin B(\rho_n^{1/2}\bx_{0i}, \delta_n)}\frac{1}{n\rho_n}\sum_{j = 1}^n\int_{\rho_n\bx_{0i}\transpose\bx_{0j}}^{\rho_n^{1/2}\bx_i\transpose\bx_{0j}}(t - \rho_n\bx_{0i}\transpose\bx_{0j})h_n(\rho_n\bx_{0i}\transpose\bx_{0j}, t)\mathrm{d}t\\
&\quad\gtrsim \inf_{\bx_i\notin B(\rho_n^{1/2}\bx_{0i}, \delta_n)}\frac{1}{n\rho_n}\sum_{j = 1}^n\int_{\rho_n\bx_{0i}\transpose\bx_{0j}}^{\rho_n^{1/2}\bx_i\transpose\bx_{0j}}(t - \rho_n\bx_{0i}\transpose\bx_{0j})\mathrm{d}t\\
&\quad = \inf_{\bx_i\notin B(\rho_n^{1/2}\bx_{0i}, \delta_n)}\frac{1}{2n\rho_n}\sum_{j = 1}^n(\rho_n^{1/2}\bx_i\transpose\bx_{0j} - \rho_n\bx_{0i}\transpose\bx_{0j})^2\\
&\quad\geq \inf_{\bx_i\notin B(\rho_n^{1/2}\bx_{0i}, \delta_n)}\frac{1}{2}\|\bx_i - \rho_n^{1/2}\bx_{0i}\|_2^2\lambda_{\min}\left(\frac{1}{n}\sum_{j = 1}^n\bx_{0j}\bx_{0j}\transpose\right)\gtrsim \delta_n^2.
\end{align*}
By Theorem \ref{thm:Large_sample_Z_estimator}, for $\gamma_n = \sqrt{(\log n)/(n\rho_n)}$, $\bW\transpose\widehat{\bx}_i\in B(\rho_n^{1/2}\bx_{0i}, \gamma_n)$ w.p.a.1. Therefore, following a similar reasoning,
\begin{align*}
|M_{in}(\rho_n^{1/2}\bx_{0i}) - M_{in}(\bW\transpose\widehat{\bx}_i)|
&\lesssim \frac{1}{n\rho_n}\sum_{j = 1}^n\int_{\rho_n\bx_{0i}\transpose\bx_{0j}}^{\rho_n^{1/2}\widehat{\bx}_i\transpose\bW\bx_{0j}}(t - \rho_n\bx_{0i}\transpose\bx_{0j})\mathrm{d}t\lesssim \|\bW\transpose\widehat{\bx}_i - \rho_n^{1/2}\bx_{0i}\|_2^2\lesssim \gamma_n\quad\mbox{w.p.a.1}.
\end{align*}
Hence, we obtain
\begin{align*}
&\inf_{\bx_i\notin B(\rho_n^{1/2}\bx_{0i}, \delta_n)}\frac{1}{n}\{\ell_{in}(\widehat{\bx}_i) - \ell_{in}(\bW\bx_i)\}\\
&\quad\geq \left\{\frac{1}{n}\ell_{in}(\widehat{\bx}_i) - M_{in}(\bW\transpose\widehat{\bx}_i)\right\} + \{M_{in}(\bW\transpose\widehat{\bx}_i) - M_{in}(\rho_n^{1/2}\bx_{0i})\}\\
&\quad\quad + \inf_{\bx_i\notin B(\rho_n^{1/2}\bx_{0i}, \delta_n)}\{M_{in}(\rho_n^{1/2}\bx_{0i}) - M_{in}(\bx_i)\}
% \\&\quad
 + \inf_{\bx_i\notin B(\rho_n^{1/2}\bx_{0i}, \delta_n)}\left\{M_{in}(\bx_i) - \frac{1}{n}\ell_{in}(\bW\bx_i)\right\}\\
&\quad\geq \inf_{\bx_i\notin B(\rho_n^{1/2}\bx_{0i}, \delta_n)}\{M_{in}(\rho_n^{1/2}\bx_{0i}) - M_{in}(\bx_i)\} - 2\sup_{\bx_i\in\Theta}\left|M_{in}(\bx_i) - \frac{1}{n}\ell_{in}(\bW\bx_i)\right|\\
&\quad\quad- |M_{in}(\rho_n^{1/2}\bx_{0i}) - M_{in}(\bW\transpose\widehat{\bx}_i)|\\
&\quad\geq C_1\delta_n^2 - C_2\sqrt{\frac{(\log n)^{2\xi}}{n\rho_n}} - C_3\gamma_n\geq \frac{(1 + \alpha)d\log n}{n}\quad\mbox{w.p.a.1},
\end{align*}
where $C_1,C_2,C_3 > 0$ are constants. The proof is thus completed. 
\end{proof}

\subsection{Proof of Proposition \ref{prop:Criterion_satisfies_assumption} (b)}
\label{sub:proof_of_prop_GMM_criterion}

\begin{lemma}\label{lemma:Sample_moments_g}
Suppose Assumptions \ref{assumption:signal_plus_noise}, \ref{assumption:regularity_condition}, \ref{assumption:weight_functions} hold. Let $[\bx]_k$ denote the $k$th coordinate of a vector $\bx$. Then
\begin{align*}
&\sup_{\bx_i\in\Theta}\frac{1}{n}\sum_{j = 1}^n\|\widetilde{\bg}_{ij}(\bx_i)\|_2\leq \left\{
\sup_{\bx_i\in\Theta}\frac{1}{n}\sum_{j = 1}^n\|\widetilde{\bg}_{ij}(\bx_i)\|_2^2
\right\}^{1/2}\lesssim \rho_n^{1/2}\quad\mbox{w.h.p.},\\
&\sup_{\bx_i\in\Theta}\frac{1}{n}\sum_{j = 1}^n\left\|\frac{\partial\widetilde{\bg}_{ij}}{\partial\bx_i\transpose}(\bx_i)\right\|_2\leq \left\{
\sup_{\bx_i\in\Theta}\frac{1}{n}\sum_{j = 1}^n\left\|\frac{\partial\widetilde{\bg}_{ij}}{\partial\bx_i\transpose}(\bx_i)\right\|_2^2
\right\}^{1/2}\lesssim \rho_n^{1/2}\quad\mbox{w.h.p.},\\
&\sup_{\bx_i\in\Theta}\frac{1}{n}\sum_{j = 1}^n\left\|\frac{\partial^2[\widetilde{\bg}_{ij}]_k}{\partial\bx_i\partial\bx_i\transpose}(\bx_i)\right\|_2\lesssim \rho_n^2\quad\mbox{w.h.p., }k\in [d], \\
&\sup_{\bx_i\in\Theta}\frac{1}{n}\sum_{j = 1}^n\|\widetilde{\bg}_{ij}(\bx_i)\|_2^3\lesssim\rho_n^{3/2}\quad\mbox{w.p.a.1.}\quad\mbox{under Assumption \ref{assumption:signal_plus_noise} (vi) (b)}.
\end{align*}
\end{lemma}

\begin{proof}[\bf Proof of Lemma \ref{lemma:Sample_moments_g}]
By Cauchy-Schwarz inequality (or, equivalently, Jensen's inequality), for the first and second inequalities, it is sufficient to prove the latter upper bounds. For the first inequality, by Assumption \ref{assumption:signal_plus_noise} (iii) and Assumption \ref{assumption:signal_plus_noise} (v), we have
\begin{align*}
\sup_{\bx_i\in B(\rho_n^{1/2}\bx_{0i}, \delta_n)}
\frac{1}{n}\sum_{j = 1}^n\|\widetilde{\bg}_{ij}(\bW\bx_i)\|_2^2
&\lesssim
\sup_{\bx_i\in\Theta}\frac{1}{n}\sum_{j = 1}^n\left\{(A_{ij}^2 + \|\bx_i\|_2^2\|\widetilde{\bx}_j\|_2^2)\widetilde{h}_{nij}(\bW\bx_i)^2\right\}\rho_n^{-1}\|\widetilde{\bx}_j\|_2^2\\
&\lesssim \frac{1}{n}(\|\bE\|_{2\to\infty}^2 + \|\bU_\bP\|_{2\to\infty}^2\|\bS_\bP\|_2^2) + \rho_n\lesssim \rho_n\quad\mbox{w.h.p.}.
\end{align*}
For the second inequality, since 
\[
\left\|\frac{\partial\widetilde{\bg}_{ij}}{\partial\bx_i\transpose}(\bx_i)\right\|_2
\leq (|A_{ij}| + \|\bx_i\|_2\|\widetilde{\bx}_j\|_2)|D^{(0, 1)}\widetilde{h}_{nij}(\bx_i)|\rho_n^{-1/2}\|\widetilde{\bx}_j\|_2^2 + \widetilde{h}_{nij}(\bx_i)\rho_n^{-1/2}\|\widetilde{\bx}_j\|_2^2,
\]
then by Result \ref{result:concentration_of_infinity_norm} that $\|\bA\|_\infty\lesssim n\rho_n$ w.h.p., Result \ref{result:pij_tilde_concentration}, and Assumption \ref{assumption:weight_functions}, we have
\begin{align*}
\sup_{\bx_i\in \Theta}\frac{1}{n}\sum_{j = 1}^n\left\|\frac{\partial\widetilde{\bg}_{ij}}{\partial\bx_i\transpose}(\bx_i)\right\|_2^2
&\leq \frac{3}{n}\sum_{j = 1}^n\left(A_{ij}^2 + \sup_{\bx_i\in \Theta}\|\bx_i\|_2^2\|\widetilde{\bx}_j\|_2^2\right)\sup_{\bx_i\in \Theta}D^{(0, 1)}\widetilde{h}_{nij}(\bW\bx_i)^2\rho_n^{-1}\|\widetilde{\bx}_j\|_2^4\\
&\quad + \frac{3}{n}\sum_{j = 1}^n\sup_{\bx_i\in \Theta}\widetilde{h}_{nij}(\bW\bx_i)^2\rho_n^{-1}\|\widetilde{\bx}_j\|_2^4\\
&\lesssim \frac{1}{n}(\|\bE\|_{2\to\infty}^2 + \|\bU_\bP\|_{2\to\infty}^2\|\bS_\bP\|_2^2 + n\rho_n)\rho_n^3 + \rho_n\lesssim \rho_n\quad\mbox{w.h.p.}.
\end{align*}
For the third inequality, write
\begin{align*}
\frac{\partial^2[\widetilde{\bg}_{ij}]_k}{\partial\bx_i\partial\bx_i\transpose}(\bx_i)
& = \{(A_{ij} - \bx_i\transpose\widetilde{\bx}_j)D^{(0, 2)}\widetilde{h}_{nij}(\bx_i) - 2D^{(0, 1)}\widetilde{h}_{nij}(\bx_i)\}\rho_n^{-1/2}\widetilde{x}_{jk}\widetilde{\bx}_j\widetilde{\bx}_j\transpose.
\end{align*}
By Result \ref{result:concentration_of_infinity_norm}, $\sum_{j= 1}^n|A_{ij}| = \|\bA\|_\infty\lesssim n\rho_n$ w.h.p.. It follows from Assumption \ref{assumption:weight_functions} and Result \ref{result:pij_tilde_concentration} that, 
\begin{align*}
&\sup_{\bx_i\in\Theta}\frac{1}{n}\sum_{j = 1}^n\left\|\frac{\partial^2[\widetilde{\bg}_{ij}]_k}{\partial\bx_i\partial\bx_i\transpose}(\bx_i)
\right\|_2\\
&\quad
\lesssim \sup_{\bx_i\in\Theta}\frac{1}{n\rho_n^{1/2}}\sum_{j = 1}^n\left\{(|A_{ij}| + \|\bx_i\|_2\|\widetilde{\bx}_j\|_2)|D^{(0, 2)}\widetilde{h}_{nij}(\bx_i)| + 2|D^{(0, 1)}\widetilde{h}_{nij}(\bx_i)|\right\}\|\widetilde{\bx}_j\|_2^3\\
&\quad\lesssim \frac{1}{n\rho_n^{1/2}}\left(\rho_n\|\bA\|_\infty + n\rho_n^{3/2} + n\rho_n\right)\rho_n^{3/2}\lesssim \rho_n^2 \quad\mbox{w.h.p., }k\in[d].
\end{align*}
For the last inequality, we have
\begin{align*}
\sup_{\bx_i\in\Theta}\frac{1}{n}\sum_{j = 1}^n\left\|\widetilde{\bg}_{ij}(\bx_i)\right\|_2^3 &=
\sup_{\bx_i\in\Theta}\frac{1}{n}\sum_{j = 1}^n\left\|(A_{ij}-\bx_i\transpose\widetilde\bx_j)\widetilde{h}_{nij}(\bx_i)\rho_n^{-1/2}\widetilde\bx_j\right\|_2^3 \\
&\leq \sup_{\bx_i\in\Theta}\frac{1}{n}\sum_{j = 1}^n \left|A_{ij}-\bx_i\transpose\widetilde\bx_j\right|^3 \max_{i,j\in[n]}\sup_{\bx_i\in\Theta}\left|\widetilde{h}_{nij}(\bx_i)\right|^3 \max_{j\in[n]}\|\rho_n^{-1/2}\widetilde\bx_j\|_2^3 \\
&\lesssim \sup_{\bx_i\in\Theta}\frac{1}{n}\sum_{j = 1}^n\left\{|A_{ij}|^3 + 3|A_{ij}|^2\|\bx_i\|_2\|\widetilde\bx_j\|_2 + 3|A_{ij}|\|\bx_i\|_2^2\|\widetilde\bx_j\|_2^2 + \|\bx_i\|_2^3\|\widetilde\bx_j\|_2^3\right\}\\
&\lesssim \frac{1}{n}\sum_{j = 1}^n|A_{ij}|^3 + \frac{1}{n}\sum_{j = 1}^n|A_{ij}|^2\rho_n^{1/2} + \frac{1}{n}\left\|\bA\right\|_{\infty}\rho_n + \frac{1}{n}n\rho_n^{3/2}.
\end{align*}
For $(1/n)\sum_{j = 1}^n|A_{ij}|^3$, we have 
\begin{align*}
\prob_0\left\{\frac{1}{n}\left|\sum_{j=1}^n(|A_{ij}|^3-\expect_0|A_{ij}|^3)\right|\geq\rho_n^{3/2}\right\}
&\leq \frac{1}{n^2\rho_n^3}\sum_{j = 1}^n\var_0(A_{ij}^3)\leq \frac{1}{n^2\rho_n^3}\sum_{j=1}^n\expect_0(A_{ij}^6)\lesssim \frac{1}{n\rho_n^3}\max_{j\in [n]}\|A_{ij}\|_{\psi_2}^6 \lesssim \frac{1}{n}
\end{align*}
by Chebyshev's Inequality and 
\[
\max_{j\in [n]}\expect_0|A_{ij}|^3\lesssim \max_{j\in [n]}\|A_{ij}\|_{\psi_2}^3\lesssim \rho_n^{3/2}
\]
under Assumption \ref{assumption:signal_plus_noise} (vi) (b), so $(1/n)\sum_{j = 1}^n|A_{ij}|^3\lesssim\rho_n^{3/2}$ w.p.a.1.\ by triangle inequality.
For $({1}/{n})\sum_{j = 1}^n|A_{ij}|^2$, we have 
$(1/n)\sum_{j = 1}^n|A_{ij}|^2\leq (1/n)(2\|\bE\|_{2\to\infty}^2 + 2\|\bU_\bP\|_{2\to\infty}^2\|\bS_\bP\|_2^2)\lesssim n\rho_n$ w.h.p..
% $\|A_{ij}^2\|_{\psi_1\leq}\|A_{ij}\|_{\psi_2}^2$ by Lemma 5.14 in \cite{vershynin2010introduction}, $\frac{1}{n}\left|\sum_{j=1}^n(|A_{ij}|^3-\expect_0|A_{ij}|^2)\right|\lesssim\rho_n$ w.h.p.\ by Proposition 5.16 in \cite{vershynin2010introduction} and $\expect_0|A_{ij}|^2\lesssim\rho_n$ under Assumption \ref{assumption:signal_plus_noise} (vi) (b), so $\frac{1}{n}\sum_{j = 1}^n|A_{ij}|^2\lesssim\rho_n$ w.h.p.\ by triangle inequality.
By Result \ref{result:concentration_of_infinity_norm}, $\|\bA\|_{\infty}\lesssim n\rho_n$ w.h.p.. 
Therefore, we conclude that $\sup_{\bx_i\in\Theta}({1}/{n})\sum_{j = 1}^n\left\|\widetilde{\bg}_{ij}(\bx_i)\right\|_2^3\lesssim\rho_n^{3/2}$ w.p.a.1..
The proof is thus completed.
\end{proof}

\begin{proof}[\bf Proof of Proposition \ref{prop:Criterion_satisfies_assumption} (b)]
Let $\eps_n = (\log n)^{1/4}/\sqrt{n}$, $M_n = \log\log n$, and $\delta_n = M_n\{(\log n)^{2\xi}/(n\rho_n)\}^{1/4}$. 
Denote $\widetilde{\bV}_{in} = \{1/n\sum_{j = 1}^n\bW\transpose\widetilde{\bg}_{ij}(\widetilde{\bx}_i)\widetilde{\bg}_{ij}(\widetilde{\bx}_i)\transpose\bW\}^{-1}$. 
By definition of the GMM criterion function \eqref{eqn:GMM}, we have
\begin{align*}
\frac{\partial\ell_{in}}{\partial\bx_i\transpose}(\bx_i)
& = -n\left\{\frac{1}{n}\sum_{j = 1}^n\bW\transpose\widetilde{\bg}_{ij}(\bx_i)\right\}\transpose\widetilde{\bV}_{in}\left\{\frac{1}{n}\sum_{j = 1}^n\bW\transpose\frac{\partial\widetilde{\bg}_{ij}}{\partial\bx_i\transpose}(\bx_i)\right\},\\
\frac{1}{n}\bW\transpose\frac{\partial^2\ell_{in}}{\partial\bx_i\partial\bx_i\transpose}(\bx_i)\bW
& = -\left\{\frac{1}{n}\sum_{j = 1}^n\bW\transpose\frac{\partial\widetilde{\bg}_{ij}}{\partial\bx_i\transpose}(\bx_i)\bW\right\}\transpose\widetilde{\bV}_{in}\left\{\frac{1}{n}\sum_{j = 1}^n\bW\transpose\frac{\partial\widetilde{\bg}_{ij}}{\partial\bx_i\transpose}(\bx_i)\bW\right\}\\
&\quad - \bW\transpose\sum_{k,l = 1}^d\left\{\frac{1}{n}\sum_{j = 1}^n[\widetilde{\bg}_{ij}]_k(\bx_i)\right\}\transpose[\bW\widetilde{\bV}_{in}\bW\transpose]_{kl}\left\{\frac{1}{n}\sum_{j = 1}^n\frac{\partial^2[\widetilde{\bg}_{ij}]_l}{\partial\bx_i\partial\bx_i\transpose}(\bx_i)\right\}\bW.
\end{align*}
Write
\begin{align*}
&\frac{1}{n}\sum_{j = 1}^n\bW\transpose\widetilde{\bg}_{ij}(\bx_i) = \frac{1}{n}\sum_{j = 1}^n\expect_0\{\bg_{ij}(\rho_n^{1/2}\bx_{0i})\} + \bR_{in}^{(\bg)}(\bx_i),\\
&\frac{1}{n}\sum_{j = 1}^n\bW\transpose\frac{\partial\widetilde{\bg}_{ij}}{\partial\bx_i\transpose}(\bx_i)\bW = \bG_{0in} + \bR_{in}^{(\bG)}(\bx_i),\\
&\frac{1}{n}\sum_{j = 1}^n\bW\transpose\widetilde{\bg}_{ij}(\widetilde{\bx}_i)\widetilde{\bg}_{ij}(\widetilde{\bx}_i)\transpose\bW = \bOmega_{0in} + \bR_{in}^{(\bOmega)},\\
&\widetilde{\bV}_{in} = \bOmega_{0in}^{-1} + \bR_{in}^{(\bV)},\\
&\frac{1}{n}\bW\transpose\frac{\partial^2\ell_{in}}{\partial\bx_i\partial\bx_i\transpose}(\bx_i)\bW = -\bG_{0in}\transpose\bOmega_{0in}^{-1}\bG_{0in} + \bR_{in}^{(\ell)}(\bx_i),
\end{align*}
where $\bG_{0in}\overset{\Delta}{=}\bG_{in}(\rho_n^{1/2}\bx_{0i})$ and $\bOmega_{0in}\overset{\Delta}{=}\bOmega_{in}(\rho_n^{1/2}\bx_{0i})$. By Lemma \ref{lemma:ULLN} and Assumption \ref{assumption:weight_functions},
\begin{align*}
\sup_{\bx_i\in B(\rho_n^{1/2}\bx_{0i}, 3\delta_n)}\|\bR_{in}^{(\bg)}(\bW\bx_i)\|_2
& \lesssim \sqrt{\frac{(\log n)^{2\xi}}{n}} + \sup_{\bx_i\in B(\rho_n^{1/2}\bx_{0i}, 3\delta_n)}\left\|\frac{1}{n}\sum_{j = 1}^n\expect_0\{\bg_{ij}(\bx_i) - \bg_{ij}(\rho_n^{1/2}\bx_{0i})\}\right\|_2\\
& \lesssim M_n\rho_n^{1/4}\left\{\frac{(\log n)^{2\xi}}{n}\right\}^{1/4}\quad\mbox{w.h.p.},\\
\sup_{\bx_i\in B(\rho_n^{1/2}\bx_{0i}, 3\delta_n)}\|\bR_{in}^{(\bG)}(\bW\bx_i)\|_2
& \lesssim \sqrt{\frac{(\log n)^{2\xi}}{n}} + \sup_{\bx_i\in B(\rho_n^{1/2}\bx_{0i}, 3\delta_n)}\left\|\bG_{in}(\bx_i) - \bG_{0in}\right\|_2\\
& \lesssim \sqrt{\frac{(\log n)^{2\xi}}{n}} + \sup_{\bx_i\in B(\rho_n^{1/2}\bx_{0i}, 3\delta_n)}\rho_n^{1/2}\frac{1}{n}\sum_{j = 1}^n|h_{0nij}(\bx_i) - h_{0nij}(\rho_n^{1/2}\bx_{0i})|\|\bx_{0j}\|_2^2\\
&\quad + \sup_{\bx_i\in B(\rho_n^{1/2}\bx_{0i}, 3\delta_n)}\rho_n\frac{1}{n}\sum_{j = 1}^n\|\rho_n^{1/2}\bx_{0i} - \bx_{i}\|_2\|\bx_{0j}\|_2|D^{(0, 1)}h_{0nij}(\bx_i)|\|\bx_{0j}\|_2^2\\
&\lesssim M_n\rho_n^{1/4}\left\{\frac{(\log n)^{2\xi}}{n}\right\}^{1/4}\quad\mbox{w.h.p.}
\end{align*}
By Lemma \ref{lemma:LLN}, $\|\bR_{in}^{(\bOmega)}\|_2\lesssim \rho_n^{1/2}\sqrt{(\log n)^{2\xi}/n}$ w.h.p.. Then by Assumption \ref{assumption:regularity_condition} and Result \ref{result:second_moment_matrix}, $\|\widetilde{\bV}_{in}\|_2\lesssim \rho_n^{-1}$ and $\|\widetilde{\bV}_{in}^{-1}\|_2\lesssim \rho_n$ w.h.p.. Also, we have $\|\bOmega_{0in}\|_2\lesssim \rho_n$ and $\|\bOmega_{0in}^{-1}\|_2\lesssim \rho_n^{-1}$ deterministically. Therefore,
\[
\|\bR_{in}^{(\bV)}\|_2 \leq \|\widetilde{\bV}_{in}\|_2\|\widetilde{\bV}_{in}^{-1} - \bOmega_{0in}\|_2\|\bOmega_{0in}^{-1}\|_2\lesssim \rho_n^{-1}\sqrt{\frac{(\log n)^{2\xi}}{n\rho_n}}\quad\mbox{w.h.p.}.
\]
Observe that $(1/n)\sum_{j = 1}^n\expect_0\{\bg_{ij}(\rho_n^{1/2}\bx_{0i})\} = \zero_d$. Then by Lemma \ref{lemma:Sample_moments_g}, for all $k,l\in [d]$, we have
\begin{align*}
&\sup_{\bx_i\in B(\rho_n^{1/2}\bx_{0i}, 3\delta_n)}\left\|\left\{\frac{1}{n}\sum_{j = 1}^n[\widetilde{\bg}_{ij}]_k(\bx_i)\right\}\transpose[\bW\widetilde{\bV}_{in}\bW\transpose]_{kl}\left\{\frac{1}{n}\sum_{j = 1}^n\frac{\partial^2[\widetilde{\bg}_{ij}]_l}{\partial\bx_i\partial\bx_i\transpose}(\bx_i)\right\}\right\|_2
\\
&\quad\leq \|\widetilde{\bV}_{in}\|_{\mathrm{F}}\sup_{\bx_i\in B(\rho_n^{1/2}\bx_{0i}, 3\delta_n)}\left\|\frac{1}{n}\sum_{j = 1}^n\bW\transpose\widetilde{\bg}_{ij}(\bx_i)\right\|_2\sup_{\bx_i\in\Theta}\frac{1}{n}\sum_{j = 1}^n\left\|\frac{\partial^2[\widetilde{\bg}_{ij}]_l}{\partial\bx_i\partial\bx_i\transpose}(\bx_i)\right\|_2\\
&\quad\lesssim M_n\rho_n^{5/4}\left\{\frac{(\log n)^{2\xi}}{n}\right\}^{1/4}\quad\mbox{w.h.p.}.
\end{align*}
Now we focus on the Hessian of the GMM criterion function $\ell_{in}(\bx_i)$. By the previous computation, we have
\begin{align*}
&\left\{\frac{1}{n}\sum_{j = 1}^n\bW\transpose\frac{\partial\widetilde{\bg}_{ij}}{\partial\bx_i\transpose}(\bx_i)\bW\right\}\transpose\widetilde{\bV}_{in}\left\{\frac{1}{n}\sum_{j = 1}^n\bW\transpose\frac{\partial\widetilde{\bg}_{ij}}{\partial\bx_i\transpose}(\bx_i)\bW\right\}\\
&\quad = \{\bG_{0in} + \bR_{in}^{(\bG)}(\bx_i)\}\transpose\{\bOmega_{0in}^{-1} + \bR_{in}^{(\bV)}\}\{\bG_{0in} + \bR_{in}^{(\bG)}(\bx_i)\}\\
&\quad = \bG_{0in}\transpose\bOmega_{0in}^{-1}\bG_{0in} + \bG_{0in}\transpose\bR_{in}^{(\bV)}\bG_{0in} + \{\bR_{in}^{(\bG)}(\bx_i)\}\transpose\widetilde{\bV}_{in}\bG_{0in}\\
&\quad\quad + \bG_{0in}\transpose\widetilde{\bV}_{in}\bR_{in}^{(\bG)}(\bx_i) + \{\bR_{in}^{(\bG)}(\bx_i)\}\transpose\widetilde{\bV}_{in}\bR_{in}^{(\bG)}(\bx_i).
\end{align*}
It follows directly from the previous results and Assumption \ref{assumption:regularity_condition} that
\begin{align*}
&\sup_{\bx_i\in B(\rho_n^{1/2}\bx_{0i}, 3\delta_n)}\left\|\frac{1}{n}\bW\transpose\frac{\partial^2\ell_{in}}{\partial\bx_i\partial\bx_i\transpose}(\bW\bx_i)\bW + \bG_{0in}\transpose\bOmega_{0in}^{-1}\bG_{0in}\right\|_2\\
&\quad\leq \|\bG_{0in}\|_2^2\|\bR_{in}^{(\bV)}\|_2 + \sup_{\bx_i\in B(\rho_n^{1/2}\bx_{0i}, 3\delta_n)}\left\{2\|\bR_{in}^{(\bG)}(\bW\bx_i)\|_2\|\widetilde{\bV}_{in}\|_2\|\bG_{0in}\| + 
\|\bR_{in}^{(\bG)}(\bW\bx_i)\|_2^2\|\widetilde{\bV}_{in}\|_2\right\}\\
&\quad\quad + \sum_{k,l\in [d]}\sup_{\bx_i\in B(\rho_n^{1/2}\bx_{0i}, 3\delta_n)}\left\|\left\{\frac{1}{n}\sum_{j = 1}^n[\widetilde{\bg}_{ij}]_k(\bx_i)\right\}\transpose[\bW\widetilde{\bV}_{in}\bW\transpose]_{kl}\left\{\frac{1}{n}\sum_{j = 1}^n\frac{\partial^2[\widetilde{\bg}_{ij}]_l}{\partial\bx_i\partial\bx_i\transpose}(\bx_i)\right\}\right\|_2\\
&\quad\lesssim M_n\left\{\frac{(\log n)^{2\xi}}{n\rho_n}\right\}^{1/4} = M_n\left\{\frac{(\log n)^{2\xi + 2}}{n\rho_n}\right\}^{1/4}\frac{1}{n\eps_n^2} = o\left(\frac{1}{n\eps_n^2}\right)\quad\mbox{w.h.p.}.
\end{align*}
Since $\eps_n\leq\delta_n$ and the eigenvalues of $\bG_{0in}\transpose\bOmega_{0in}^{-1}\bG_{0in}$ are bounded away from $0$ and $\infty$ by Result \ref{result:Jacobian} and Result \ref{result:second_moment_matrix}, the above concentration bound completes the proofs of \eqref{eqn:Hessian_A1} and \eqref{eqn:Hessian_A2} simultaneously. It is now sufficient to establish \eqref{eqn:identifiability}. Define the function $M_{in}(\bx_i)\overset{\Delta}{=}\|\bOmega_{0in}^{-1/2}(1/n)\sum_{j = 1}^n\expect_0\{\bg_{ij}(\bx_i)\}\|_2^2$. A simple algebra leads to
\begin{align*}
&\left|\frac{2}{n}\ell_{in}(\bW\bx_i) + M_{in}(\bx_i)\right|\\
&\quad = \left\{
\left\|\widetilde{\bV}_{in}^{1/2}\frac{1}{n}\sum_{j = 1}^n\bW\transpose\widetilde{\bg}_{ij}(\bW\bx_i)\right\|_2 + \left\|\bOmega_{0in}^{-1/2}\frac{1}{n}\sum_{j = 1}^n\expect_0\{\bg_{ij}(\bx_i)\}\right\|_2
\right\}\\
&\quad\quad\times \left|\left\|\widetilde{\bV}_{in}^{1/2}\frac{1}{n}\sum_{j = 1}^n\bW\transpose\widetilde{\bg}_{ij}(\bW\bx_i)\right\|_2 - \left\|\bOmega_{0in}^{-1/2}\frac{1}{n}\sum_{j = 1}^n\expect_0\{\bg_{ij}(\bx_i)\}\right\|_2\right|\\
&\quad\leq \left\{\|\widetilde{\bV}_{in}\|_2^{1/2}\left\|\frac{1}{n}\sum_{j = 1}^n\bW\transpose\widetilde{\bg}_{ij}(\bW\bx_i)\right\|_2 + \|\bOmega_{0in}^{-1}\|_2^{1/2}\left\|\frac{1}{n}\sum_{j = 1}^n\expect_0\{\bg_{ij}(\bx_i)\}\right\|_2\right\}\\
&\quad\quad\times \left\{\|\widetilde{\bV}_{in}^{1/2} - \bOmega_{0in}^{-1/2}\|_2
\left\|\frac{1}{n}\sum_{j = 1}^n\bW\transpose\widetilde{\bg}_{ij}(\bW\bx_i)\right\|_2 + \|\bOmega_{0in}^{-1}\|_2^{1/2}\left\|\frac{1}{n}\sum_{j = 1}^n[\bW\transpose\widetilde{\bg}_{ij}(\bW\bx_i) - \expect_0\{\bg_{ij}(\bx_i)\}]\right\|_2
\right\}.
\end{align*}
To bound $\|\widetilde{\bV}_{in}^{1/2} - \bOmega_{0in}^{-1/2}\|_2$, consider an eigenvector $\by$ of $\widetilde{\bV}_{in}^{1/2} - \bOmega_{0in}^{-1/2}$ associated with an eigenvalue $\mu$ with $\|\by\|_2 = 1$. Clearly,
\begin{align*}
\by\transpose(\widetilde{\bV}_{in} - \bOmega_{0in}^{-1})\by& = \by\transpose(\widetilde{\bV}_{in}^{1/2} - \bOmega_{0in}^{-1/2})\widetilde{\bV}_{in}^{1/2}\by + \by\transpose\bOmega_{0in}^{-1/2}(\widetilde{\bV}_{in}^{1/2} - \bOmega_{0in}^{-1/2})\by\\
& = \mu\by\transpose(\widetilde{\bV}_{in}^{1/2} + \bOmega_{0in}^{-1/2})\by.
\end{align*}
It follows that
\[
\lambda_{\min}(\bOmega_{0in}^{-1/2})|\mu|\leq|\mu||\by\transpose(\widetilde{\bV}_{in}^{1/2} + \bOmega_{0in}^{-1/2})\by| = |\by\transpose(\widetilde{\bV}_{in} - \bOmega_{0in}^{-1})\by|\leq \|\widetilde{\bV}_{in} - \bOmega_{0in}^{-1}\|_2.
\]
In particular, $\mu$ can be selected such that $|\mu| = \|\widetilde{\bV}_{in}^{1/2} - \bOmega_{0in}^{-1/2}\|_2$. Therefore, we obtain that
\[
\|\widetilde{\bV}_{in}^{1/2} - \bOmega_{0in}^{-1/2}\|_2\leq \frac{1}{\lambda_{\min}(\bOmega_{0in}^{-1/2})}\|\widetilde{\bV}_{in} - \bOmega_{0in}^{-1}\|_2\lesssim \rho_n^{-1/2}\sqrt{\frac{(\log n)^{2\xi}}{n\rho_n}}\quad\mbox{w.h.p.}.
\]
Also, by Lemma \ref{lemma:ULLN},
\begin{align*}
\sup_{\bx_i\in\Theta}\left\|\frac{1}{n}\sum_{j = 1}^n\bW\transpose\widetilde{\bg}_{ij}(\bW\bx_i)\right\|_2
&\leq \sup_{\bx_i\in\Theta}\left\|\frac{1}{n}\sum_{j = 1}^n[\bW\transpose\widetilde{\bg}_{ij}(\bW\bx_i) - \expect_0\{\bg_{ij}(\bx_i)\}]\right\|_2 + \sup_{\bx_i\in\Theta}\left\|\frac{1}{n}\sum_{j = 1}^n\expect_0\{\bg_{ij}(\bx_i)\}\right\|_2\\
&\lesssim \sqrt{\frac{(\log n)^{2\xi}}{n}} + \rho_n^{1/2}\lesssim \rho_n^{1/2}\quad\mbox{w.h.p.}.
\end{align*}
Hence, we conclude from the previous concentration bounds that
\begin{align*}
\sup_{\bx_i\in\Theta}\left|\frac{2}{n}\ell_{in}(\bW\bx_i) + M_{in}(\bx_i)\right|
&\lesssim \sqrt{\frac{(\log n)^{2\xi}}{n\rho_n}}\quad\mbox{w.h.p.}.
\end{align*}
By Assumption \ref{assumption:regularity_condition} and the fact that $M_{in}(\rho_n^{1/2}\bx_{0i}) = 0$, for $\gamma_n = \sqrt{(\log n)/(n\rho_n)}$ and for sufficiently large $n$,
\begin{align*}
&\inf_{\bx_i\notin B(\rho_n^{1/2}\bx_{0i}, \delta_n)}\{M_{in}(\bx_i) - M_{in}(\rho_n^{1/2}\bx_{0i})\}\geq M_n\sqrt{\frac{(\log n)^{2\xi}}{n\rho_n}},\\
&\sup_{\bx_i\in B(\rho_n^{1/2}\bx_{0i}, \gamma_n)}\{M_{in}(\bx_i) - M_{in}(\rho_n^{1/2}\bx_{0i})\}\leq M_n\sqrt{\frac{(\log n)}{n\rho_n}}.
\end{align*}
Observe that by Theorem \ref{thm:Large_sample_Z_estimator}, $\bW\transpose\widehat{\bx}_i \in B(\rho_n^{1/2}\bx_{0i}, \gamma_n)$ w.p.a.1. 
It follows that
\begin{align*}
\inf_{\bW\transpose\bx_i\notin B(\rho_n^{1/2}\bx_{0i}, \delta_n)}\{\ell_{in}(\widehat{\bx}_i) - \ell_{in}(\bx_i)\}
&\geq \left\{\ell_{in}(\widehat{\bx}_i) + \frac{n}{2} M_{in}(\bW\transpose\widehat{\bx}_i)\right\}\\
&\quad + \inf_{\bW\transpose\bx_i\notin B(\rho_n^{1/2}\bx_{0i}, \delta_n)}\left\{-\frac{n}{2}M_{in}(\bW\transpose\widehat{\bx}_i) +\frac{n}{2} M_{in}(\bW\transpose\bx_i)\right\}\\
&\quad  + \inf_{\bW\transpose\bx_i\notin B(\rho_n^{1/2}\bx_{0i}, \delta_n)}\left\{ - \frac{n}{2}M_{in}(\bW\transpose\bx_i) - \ell_{in}(\bx_i)\right\}\\
&\geq \frac{n}{2}\inf_{\bW\transpose\bx_i\notin B(\rho_n^{1/2}\bx_{0i}, \delta_n)}\{M_{in}(\bW\transpose\bx_i) - M_{in}(\bW\transpose\widehat{\bx}_i)\}\\
&\quad - n\sup_{\bz\in\Theta}\left|\frac{2}{n}\ell_{in}(\bW\bz) + M_{in}(\bz)\right|\\
&\geq \frac{n}{2}\inf_{\bz\notin B(\rho_n^{1/2}\bx_{0i}, \delta_n)}\{M_{in}(\bz) - M_{in}(\rho_n^{1/2}\bx_{0i})\}\\
&\quad - \frac{n}{2}\left|M_{in}(\bW\transpose\widehat{\bx}_i) - M_{in}(\rho_n^{1/2}\bx_{0i})\right|
% \\&\quad
 - n\sup_{\bz\in\Theta}\left|\frac{2}{n}\ell_{in}(\bW\bz) + M_{in}(\bz)\right|\\
 &\geq \frac{n}{2}\left\{M_n\sqrt{\frac{(\log n)^{2\xi}}{n\rho_n}} - M_n\sqrt{\frac{(\log n)}{n\rho_n}} - 2C\sqrt{\frac{(\log n)^{2\xi}}{n\rho_n}}\right\}\\
&\geq \frac{n}{4}M_n\sqrt{\frac{(\log n)^{2\xi}}{n\rho_n}}\geq (1 + \alpha)d\log n\quad\mbox{w.p.a.1}
\end{align*}
for any $\alpha > 0$, where $C > 0$ is a constant. The proof is thus completed.
\end{proof}

\subsection{Auxiliary results for ETEL}
\label{sub:auxiliar_results_ETEL}
The most technical part of the proof of Proposition \ref{prop:Criterion_satisfies_assumption} is the analysis of the ETEL criterion. In preparation for doing so, we provide a collection of auxiliary results for the ETEL in this subsection. 

\begin{lemma}\label{lemma:lambda_g_ETEL}
Suppose Assumptions \ref{assumption:signal_plus_noise}, \ref{assumption:regularity_condition}, and \ref{assumption:weight_functions} hold. Let $\Lambda_n = \{\blambda\in\mathbb{R}^d:\|\blambda\|_2\leq \zeta_n\}$, where $(\zeta_n)_{n = 1}^\infty$ is a positive sequence. Then
\[
\max_{i,j\in [n]}\sup_{(\blambda_i,\bx_i)\in\Lambda_n\times\Theta}|\blambda_i\transpose\widetilde{\bg}_{ij}(\bx_i)|\lesssim 
\left\{\begin{aligned}
&\zeta_n &\quad&\mbox{w.h.p., if Assumption \ref{assumption:signal_plus_noise} (vi) (a) holds,}\\
&(\rho_n\log n)^{1/2}\zeta_n &\quad&\mbox{w.h.p., if Assumption \ref{assumption:signal_plus_noise} (vi) (b) holds}.
\end{aligned}
\right.
\]
\end{lemma}
\begin{proof}[\bf Proof of Lemma \ref{lemma:lambda_g_ETEL}]
% By triangle inequality and Cauchy-Schwarz inequality, 
% \[
% \|\widetilde{\bg}_{ij}(\bx_i)\|_2\leq |A_{ij}|\rho_n^{-1/2}\widetilde{h}_{nij}(\bx_i)\|\widehat{\bx}_j\|_2 + \rho_n^{-1/2}\widetilde{h}_{nij}(\bx_i)\|\bx_i\|_2\|\widetilde{\bx}_j\|_2^2,
% \]
% where $\widetilde{h}_{nij}(\bx_i):=h_n(\widetilde{\bx}_i\transpose\widetilde{\bx}_j, \bx_i\transpose\widetilde{\bx}_j)$. By Assumption \ref{assumption:weight_functions} and Result \ref{result:pij_tilde_concentration}, $\sup_{\bx_i\in\Theta}|\widetilde{h}_{nij}(\bx_i)|\lesssim 1$ and $\|\widetilde{\bx}_j\|_2^2\lesssim \rho_n$ w.h.p.. By Assumption \ref{assumption:signal_plus_noise} (vi) and Lemma 8.1 in \cite{kosorok2008introduction}, 
% \[
% \max_{i,j\in [n]}|A_{ij}|\lesssim
% \left\{\begin{aligned}
% &1&\quad&\mbox{w.p.1., if Assumption \ref{assumption:signal_plus_noise} (vi) (a) holds,}\\
% &(\rho_n\log n)^{1/2}&\quad&\mbox{w.h.p., if Assumption \ref{assumption:signal_plus_noise} (vi) (b) holds}.
% \end{aligned}
% \right.
% \]
% By Result \ref{result:pij_tilde_concentration}, this implies that
By Result \ref{result:uniform_concentration_g}, we have
\[
\max_{i,j\in [n]}\sup_{\bx_i\in\Theta}\|\widetilde{\bg}_{ij}(\bx_i)\|_2\lesssim \left\{\begin{aligned}
&1&\quad&\mbox{w.h.p., if Assumption \ref{assumption:signal_plus_noise} (vi) (a) holds,}\\
&(\rho_n\log n)^{1/2} &\quad&\mbox{w.h.p., if Assumption \ref{assumption:signal_plus_noise} (vi) (b) holds}.
\end{aligned}
\right.
\] It follows from Cauchy-Schwarz inequality that
\begin{align*}
\max_{i,j\in [n]}\sup_{(\blambda_i,\bx_i)\in\Lambda_n\times\Theta}|\blambda_i\transpose\widetilde{\bg}_{ij}(\bx_i)|
&\leq \sup_{\blambda_i\in\Lambda_n}\|\blambda_i\|_2\max_{i,j\in[n]}\sup_{\bx_i\in\Theta}\|\widetilde{\bg}_{ij}(\bx_i)\|_2\\
&\lesssim \left\{\begin{aligned}
&\zeta_n &\quad&\mbox{w.h.p., if Assumption \ref{assumption:signal_plus_noise} (vi) (a) holds,}\\
&(\rho_n\log n)^{1/2}\zeta_n &\quad&\mbox{w.h.p., if Assumption \ref{assumption:signal_plus_noise} (vi) (b) holds}.
\end{aligned}
\right.
\end{align*}
\end{proof}

\begin{lemma}\label{lemma:covariance_bound_ETEL}
Suppose Assumptions \ref{assumption:signal_plus_noise}, \ref{assumption:regularity_condition}, and \ref{assumption:weight_functions} hold. Let $(\delta_n)_{n = 1}^\infty$ be a positive sequence converging to $0$. Then
\begin{align*}
\sup_{\bx_i\in B(\rho_n^{1/2}\bx_{0i}, \delta_n)}\left\|\frac{1}{n}\sum_{j = 1}^n\bW\transpose\widetilde{\bg}_{ij}(\bW\bx_i)\widetilde{\bg}_{ij}(\bW\bx_i)\transpose\bW - \bOmega_{in}(\bx_i)\right\|_2\lesssim \rho_n\left(\delta_n + \sqrt{\frac{(\log n)^{2\xi}}{n\rho_n}}\right)\quad\mbox{w.h.p.}.
\end{align*}
Consequently, by Result \ref{result:second_moment_matrix}, 
\[
\inf_{\bx_i\in B(\rho_n^{1/2}\bx_{0i}, \delta_n)}\lambda_{\min}\left\{\frac{1}{n}\sum_{j = 1}^n\widetilde{\bg}_{ij}(\bW\bx_i)\widetilde{\bg}_{ij}(\bW\bx_i)\transpose\right\}\gtrsim \rho_n\quad\mbox{w.h.p.}.
\]
\end{lemma}
\begin{proof}[\bf Proof of Lemma \ref{lemma:covariance_bound_ETEL}]
With $\widetilde{h}_{nij}(\bx_i) = h_n(\widetilde{\bx}_i\transpose\widetilde{\bx}_j, \bx_i\transpose\widetilde{\bx}_j)$ and ${h}_{0nij}(\bx_i) = h_n(\rho_n\bx_{0i}\transpose\bx_{0j}, \rho_n^{1/2}\bx_i\transpose\bx_{0j})$, we have
\begin{align*}
&\frac{1}{n}\sum_{j = 1}^n\bW\transpose\widetilde{\bg}_{ij}(\bW\bx_i)\widetilde{\bg}_{ij}(\bW\bx_i)\transpose\bW\\
&\quad = \frac{1}{n\rho_n}\sum_{j = 1}^n(A_{ij} - \bx_i\transpose\bW\transpose\widetilde{\bx}_j)^2\widetilde{h}_{nij}^2(\bW\bx_i)\bW
\transpose\widetilde{\bx}_j\widetilde{\bx}_j\transpose\bW\\
&\quad = \frac{1}{n}\sum_{j = 1}^n(A_{ij} - \rho_n^{1/2}\bx_{i}\transpose\bx_{0j})^2h_{0nij}^2(\bx_{i})\bx_{0j}\bx_{0j}\transpose\\
&\quad\quad + \frac{1}{n}\sum_{j = 1}^n(A_{ij} - \rho_n^{1/2}\bx_i\transpose\bx_{0j})^2\left\{
\widetilde{h}_{nij}^2(\bW\bx_{i})\rho_n^{-1}\bW\transpose\widetilde{\bx}_{j}\widetilde{\bx}_{j}\transpose\bW - 
h_{0nij}^2(\bx_{i})\bx_{0j}\bx_{0j}\transpose\right\}\\
&\quad\quad + \frac{1}{n}\sum_{j = 1}^n\left\{(A_{ij} - \bx_{i}\transpose\bW\transpose\widetilde{\bx}_{j})^2 - (A_{ij} - \rho_n^{1/2}\bx_{i}\transpose\bx_{0j})^2\right\}
\widetilde{h}_{nij}^2(\bW\bx_{i})\rho_n^{-1}\bW\transpose\widetilde{\bx}_{j}\widetilde{\bx}_{j}\transpose\bW.
\end{align*}
Following the same reasoning for Result \ref{result:pij_tilde_concentration} and \eqref{eqn:uniform_Lipschitz_h}, 
\begin{align*}
&\max_{i,j\in [n]}\sup_{\bx_i\in \Theta}|\bx_{i}\transpose\bW\transpose\widetilde{\bx}_{j} - \rho_n^{1/2}\bx_{i}\transpose\bx_{0j}|\lesssim \sqrt{\frac{(\log n)^{2\xi}}{n}}\quad\mbox{w.h.p.},\\
&\max_{i,j\in [n]}\sup_{\bx_i\in \Theta}|(\bx_{i}\transpose\bW\transpose\widetilde{\bx}_{j})^2 - (\rho_n^{1/2}\bx_{i}\transpose\bx_{0j})^2|\lesssim \rho_n^{1/2}\sqrt{\frac{(\log n)^{2\xi}}{n}}\quad\mbox{w.h.p.},\\
&\max_{i,j\in [n]}\sup_{\bx_i\in \Theta}|\widetilde{h}_{nij}(\bW\bx_i) - h_{0nij}(\bx_{i})|\lesssim \sqrt{\frac{(\log n)^{2\xi}}{n\rho_n}}\quad\mbox{w.h.p.},\\
&\max_{i,j\in [n]}\sup_{\bx_i\in \Theta}|\widetilde{h}_{nij}^2(\bW\bx_i) - h_{0nij}^2(\bx_i)|\lesssim \sqrt{\frac{(\log n)^{2\xi}}{n\rho_n}}\quad\mbox{w.h.p.},\\
&\max_{i,j\in [n]}\sup_{\bx_i\in \Theta}\|\widetilde{h}_{nij}^2(\bW\bx_i)\rho_n^{-1}\bW\transpose\widetilde{\bx}_j\widetilde{\bx}_j\transpose\bW - h_{0nij}^2(\bx_{i})\bx_{0j}\bx_{0j}\transpose\|_2\lesssim \sqrt{\frac{(\log n)^{2\xi}}{n\rho_n}}\quad\mbox{w.h.p.}.
\end{align*}
Therefore, the second term can be bounded as follows:
\begin{align*}
&\sup_{\bx_i\in \Theta}\left\|\frac{1}{n}\sum_{j = 1}^n(A_{ij} - \rho_n^{1/2}\bx_{i}\transpose\bx_{0j})^2\left\{
\widetilde{h}_{nij}^2(\bW\bx_{i})\rho_n^{-1}\bW\transpose\widetilde{\bx}_{j}\widetilde{\bx}_{j}\transpose\bW - 
h_{0nij}^2(\bx_{i})\bx_{0j}\bx_{0j}\transpose\right\}\right\|_2\\
&\quad
\lesssim \frac{1}{n}\left(\|\bE\|_{2\to\infty}^2 + \rho_n^{1/2}\|\bE\|_\infty + n\rho_n\right)\max_{i,j\in [n]}\sup_{\bx_i\in \Theta}\|\widetilde{h}_{nij}^2(\bW\bx_i)\rho_n^{-1}\bW\transpose\widetilde{\bx}_j\widetilde{\bx}_j\transpose\bW - h_{0nij}^2(\bx_{i})\bx_{0j}\bx_{0j}\transpose\|_2\\
&\quad\lesssim \frac{1}{n}(n\rho_n)\sqrt{\frac{(\log n)^{2\xi}}{n\rho_n}} = \rho_n^{1/2}\sqrt{\frac{(\log n)^{2\xi}}{n}}\quad\mbox{w.h.p.}.
\end{align*}
Also, observe that 
\[
|(A_{ij} - \bx_{i}\transpose\bW\transpose\widetilde{\bx}_{j})^2 - (A_{ij} - \rho_n^{1/2}\bx_{i}\transpose\bx_{0j})^2|
\leq 2|A_{ij}||\bx_{i}\transpose\bW\transpose\widetilde{\bx}_{j} - \rho_n^{1/2}\bx_{i}\transpose\bx_{0j}| + |(\bx_{i}\transpose\bW\transpose\widetilde{\bx}_{j})^2 - (\rho_n^{1/2}\bx_{i}\transpose\bx_{0j})^2|.
\]
Then by Result \ref{result:concentration_of_infinity_norm}, the third term can be bounded as follows:
\begin{align*}
&\sup_{\bx_i\in \Theta}\left\|\frac{1}{n}\sum_{j = 1}^n\left\{(A_{ij} - \bx_{i}\transpose\bW\transpose\widetilde{\bx}_{j})^2 - (A_{ij} - \rho_n^{1/2}\bx_{i}\transpose\bx_{0j})^2\right\}
\widetilde{h}_{nij}^2(\bW\bx_{i})\rho_n^{-1}\bW\transpose\widetilde{\bx}_{j}\widetilde{\bx}_{j}\transpose\bW\right\|_2\\
&\quad\lesssim  \frac{1}{n}\max_{i,j\in [n]}\sup_{\bx_i\in \Theta}\left\{|\bx_{i}\transpose\bW\transpose\widetilde{\bx}_{j} - \rho_n^{1/2}\bx_{i}\transpose\bx_{0j}|\|\bA\|_\infty + \max_{i,j\in [n]}n|(\bx_{i}\transpose\bW\transpose\widetilde{\bx}_{j})^2 - (\rho_n^{1/2}\bx_{i}\transpose\bx_{0j})^2|\right\}\\
&\quad\lesssim \frac{1}{n}\left\{(n\rho_n)\sqrt{\frac{(\log n)^{2\xi}}{n}} + n\rho_n^{1/2}\sqrt{\frac{(\log n)^{2\xi}}{n}}\right\} =  \rho_n^{1/2}\sqrt{\frac{(\log n)^{2\xi}}{n}}\quad\mbox{w.h.p.}.
\end{align*}
We now focus on the first term. Denote $\omega_{ij}^{(kl)}(\bx_i) = (A_{ij} - \rho_n^{1/2}\bx_i\transpose\bx_{0j})^2h_{0nij}^2(\bx_i)x_{0jk}x_{0jl}$. Write
\begin{align*}
\frac{\partial}{\partial\bx_i}\omega_{ij}^{(kl)}(\bx_i)& = \frac{\partial}{\partial\bx_i}(A_{ij} - \rho_n^{1/2}\bx_i\transpose\bx_{0j})^2h_{0nij}^2(\bx_i)x_{0jk}x_{0jl}\\
& = 2(A_{ij} - \rho_n^{1/2}\bx_i\transpose\bx_{0j})^2h_{0nij}(\bx_i)D^{(0, 1)}h_{0nij}(\bx_i)\bx_{0j}\rho_n^{1/2}x_{0jk}x_{0jl}\\
&\quad - 2(A_{ij} - \rho_n^{1/2}\bx_{i}\transpose\bx_{0j})h_{0nij}^2(\bx_i)\rho_n^{1/2}\bx_{0j}x_{0jk}x_{0jl}.
\end{align*}
Observe that by Assumptions \ref{assumption:signal_plus_noise} (iii), \ref{assumption:regularity_condition}, and \ref{assumption:weight_functions}, together with Result \ref{result:concentration_of_infinity_norm}, we have
\begin{align*}
\frac{1}{n}\sum_{j = 1}^n\sup_{\bx_i\in\Theta}\left\|\frac{\partial}{\partial\bx_i}\omega_{ij}^{(kl)}(\bx_i)\right\|_2
&\lesssim \frac{\rho_n^{3/2}}{n}\sum_{j = 1}^n\sup_{\bx_i\in\Theta}(A_{ij} - \rho_n^{1/2}\bx_i\transpose\bx_{0j})^2 + \frac{\rho_n^{1/2}}{n}\sum_{j = 1}^n\sup_{\bx_i\in\Theta}(|A_{ij}| + |\rho_n^{1/2}\bx_i\transpose\bx_{0j}|)\\
&\lesssim \frac{\rho_n^{3/2}}{n}(\|\bE\|_{2\to\infty}^2 + n\rho_n) + \frac{\rho_n^{1/2}}{n}\|\bA\|_\infty + \rho_n\lesssim \rho_n\quad\mbox{w.h.p.}.
\end{align*}
Denote
\begin{align*}
\Delta_{in}^{(kl)}
& = \frac{1}{n}\sum_{j = 1}^n\left\{\sup_{\bx_i\in B(\rho_n^{1/2}\bx_{0i}, \delta_n)}\omega_{ij}^{(kl)}(\bx_i) - \inf_{\bx_i\in B(\rho_n^{1/2}\bx_{0i}, \delta_n)}\omega_{ij}^{(kl)}(\bx_i)\right\}.
\end{align*}
By mean-value theorem and the previous result, we obtain
\begin{align*}
\Delta_{in}^{(kl)}
&\leq \frac{1}{n}\sum_{j = 1}^n2\sup_{\bx_i\in\Theta}\left\|\frac{\partial}{\partial\bx_i}\omega_{ij}^{(kl)}(\bx_i)\right\|_2\delta_n\lesssim \rho_n\delta_n\quad\mbox{w.h.p.}.
\end{align*}
Observe that
\[
\left|\inf_{\bx_i\in B(\rho_n^{1/2}\bx_{0i}, \delta_n)}\omega_{ij}^{(kl)}(\bx_i)\right|\lesssim A_{ij}^2 + \rho_n\quad\mbox{and}\quad
\left|\sup_{\bx_i\in B(\rho_n^{1/2}\bx_{0i}, \delta_n)}\omega_{ij}^{(kl)}(\bx_i)\right|\lesssim A_{ij}^2 + \rho_n
\]
regardless of the sign of $x_{0jk}x_{0jl}$. Then
under Assumption \ref{assumption:signal_plus_noise} (vi) (a), we have
\begin{align*}
\max_{j\in [n]}\expect_0\left\{\left[\sup_{\bx_i\in B(\rho_n^{1/2}\bx_{0i}, \delta_n)}\omega_{ij}^{(kl)}(\bx_i)\right]^2\right\}\lesssim \max_{j\in [n]}\expect_0(A_{ij})^4 + \rho_n^2\lesssim \rho_n.
\end{align*}
Under Assumption \ref{assumption:signal_plus_noise} (vi) (b), we have
\begin{align*}
&\max_{j\in [n]}\left\|\sup_{\bx_i\in B(\rho_n^{1/2}\bx_{0i}, \delta_n)}\omega_{ij}^{(kl)}(\bx_i) 
- \expect_0\left\{\sup_{\bx_i\in B(\rho_n^{1/2}\bx_{0i}, \delta_n)}\omega_{ij}^{(kl)}(\bx_i)\right\}\right\|_{\psi_1}\\
&\quad\leq \max_{j\in [n]}\left\|\sup_{\bx_i\in B(\rho_n^{1/2}\bx_{0i}, \delta_n)}\omega_{ij}^{(kl)}(\bx_i) \right\|_{\psi_1} + \left\|
\expect_0\left\{\sup_{\bx_i\in B(\rho_n^{1/2}\bx_{0i}, \delta_n)}\omega_{ij}^{(kl)}(\bx_i)\right\}\right\|_{\psi_1}\\
&\quad\leq \max_{j\in [n]}\left\|\left|\sup_{\bx_i\in B(\rho_n^{1/2}\bx_{0i}, \delta_n)}\omega_{ij}^{(kl)}(\bx_i)\right| \right\|_{\psi_1} + 
\expect_0\left|\sup_{\bx_i\in B(\rho_n^{1/2}\bx_{0i}, \delta_n)}\omega_{ij}^{(kl)}(\bx_i)\right|\\
&\quad
\leq 2\max_{j\in [n]}\left\|\left|\sup_{\bx_i\in B(\rho_n^{1/2}\bx_{0i}, \delta_n)}\omega_{ij}^{(kl)}(\bx_i)\right| \right\|_{\psi_1}
\lesssim \max_{j\in [n]}\|A_{ij}\|_{\psi_2}^2 + \rho_n\lesssim \rho_n.
\end{align*}
It follows either from Bernstein's inequality under Assumption \ref{assumption:signal_plus_noise} (vi) (a) or from Proposition 5.16 in \cite{vershynin2010introduction} under Assumption \ref{assumption:signal_plus_noise} (vi) (b) that 
\begin{align*}
\frac{1}{n}\sum_{j = 1}^n\sup_{\bx_i\in B(\rho_n^{1/2}\bx_{0i}, \delta_n)}\omega_{ij}^{(kl)}(\bx_i) = \frac{1}{n}\sum_{j = 1}^n\expect_0\left\{\sup_{\bx_i\in B(\rho_n^{1/2}\bx_{0i}, \delta_n)}\omega_{ij}^{(kl)}(\bx_i)\right\} + O\left(\rho_n^{1/2}\sqrt{\frac{\log n}{n}}\right)\quad\mbox{w.h.p.}.
\end{align*}
Similarly, we also have
\begin{align*}
\frac{1}{n}\sum_{j = 1}^n\inf_{\bx_i\in B(\rho_n^{1/2}\bx_{0i}, \delta_n)}\omega_{ij}^{(kl)}(\bx_i) = \frac{1}{n}\sum_{j = 1}^n\expect_0\left\{\inf_{\bx_i\in B(\rho_n^{1/2}\bx_{0i}, \delta_n)}\omega_{ij}^{(kl)}(\bx_i)\right\} + O\left(\rho_n^{1/2}\sqrt{\frac{\log n}{n}}\right)\quad\mbox{w.h.p.}.
\end{align*}
Hence, we obtain the following result:
\begin{align*}
&\sup_{\bx_i\in B(\rho_n^{1/2}\bx_{0i},\delta_n)}\left[\frac{1}{n}\sum_{j = 1}^n\{\omega_{ij}^{(kl)}(\bx_i) - \expect_0\omega_{ij}^{(kl)}(\bx_i)\}
\right]\\
&\quad\leq \frac{1}{n}\sum_{j = 1}^n\sup_{\bx_i\in B(\rho_n^{1/2}\bx_{0i},\delta_n)}\omega_{ij}^{(kl)}(\bx_i) - \frac{1}{n}\sum_{j = 1}^n\expect_0\left\{\inf_{\bx_i\in B(\rho_n^{1/2}\bx_{0i},\delta_n)}\omega_{ij}^{(kl)}(\bx_i)\right\}\\
&\quad \leq \frac{1}{n}\sum_{j = 1}^n\sup_{\bx_i\in B(\rho_n^{1/2}\bx_{0i},\delta_n)}\omega_{ij}^{(kl)}(\bx_i) - \frac{1}{n}\sum_{j = 1}^n\inf_{\bx_i\in B(\rho_n^{1/2}\bx_{0i},\delta_n)}\omega_{ij}^{(kl)}(\bx_i) + \left|O\left(\rho_n^{1/2}\sqrt{\frac{\log n}{n}}\right)\right|\\
&\quad = \Delta_{in}^{(kl)} + \left|O\left(\rho_n^{1/2}\sqrt{\frac{\log n}{n}}\right)\right| \lesssim \rho_n\left(\delta_n+ \sqrt{\frac{\log n}{n\rho_n}}\right)\quad\mbox{w.h.p.}.
\end{align*}
By an analogous argument, we also have
\[
\inf_{\bx_i\in B(\rho_n^{1/2}\bx_{0i},\delta_n)}\left[\frac{1}{n}\sum_{j = 1}^n\{\omega_{ij}^{(kl)}(\bx_i) - \expect_0\omega_{ij}^{(kl)}(\bx_i)\}
\right] \gtrsim -\rho_n\left(\delta_n+ \sqrt{\frac{\log n}{n\rho_n}}\right)\quad\mbox{w.h.p.}.
\]
This implies that
\[
\sup_{\bx_i\in B(\rho_n^{1/2}\bx_{0i}, \delta_n)}\left\|\frac{1}{n}\sum_{j = 1}^n(A_{ij} - \rho_n^{1/2}\bx_{i}\transpose\bx_{0j})^2h_{0nij}^2(\bx_{i})\bx_{0j}\bx_{0j}\transpose - \bOmega_{in}(\bx_i)\right\|_2\lesssim \rho_n\left(\delta_n + \sqrt{\frac{\log n}{n\rho_n}}\right)\quad\mbox{w.h.p.}.
\]
Combining the concentration bounds for the three terms, we obtain
\[
\sup_{\bx_i\in B(\rho_n^{1/2}\bx_{0i}, \delta_n)}\left\|
\frac{1}{n}\sum_{j = 1}^n\bW\transpose\widetilde{\bg}_{ij}(\bW\bx_i)\widetilde{\bg}_{ij}(\bW\bx_i)\transpose\bW - \bOmega_{in}(\bx_i)
\right\|_2\lesssim \rho_n\left(\delta_n + \sqrt{\frac{(\log n)^{2\xi}}{n\rho_n}}\right) = o(\rho_n)\quad\mbox{w.h.p.}.
\]
This completes the proof of the first assertion. For the second assertion, note that $\inf_{\bx_i\in\Theta}\lambda_d\{\bOmega_{in}(\bx_{i})\}\gtrsim \rho_n$ by Result \ref{result:second_moment_matrix}. We therefore conclude that
\begin{align*}
&\inf_{\bx_i\in B(\rho_n^{1/2}\bx_{0i}, \delta_n)}
\lambda_{\min}\left\{\frac{1}{n}\sum_{j = 1}^n\widetilde{\bg}_{ij}(\bW\bx_i)\widetilde{\bg}_{ij}(\bW\bx_i)\transpose
\right\}\\
&\quad = 
\inf_{\bx_i\in B(\rho_n^{1/2}\bx_{0i}, \delta_n)}\lambda_{\min}\left\{
\frac{1}{n}\sum_{j = 1}^n\bW\transpose\widetilde{\bg}_{ij}(\bW\bx_i)\widetilde{\bg}_{ij}(\bW\bx_i)\transpose\bW
\right\}\\
&\quad\geq \inf_{\bx_i\in\Theta}\lambda_{\min}\left\{\bOmega_{in}(\bx_{i})\right\} - \left|o\left(\rho_n\right)\right|\gtrsim \rho_n\quad\mbox{w.h.p.}.
\end{align*}
The proof of the second assertion is thus completed. 
\end{proof}

\begin{lemma}
\label{lemma:Local_ULLN}
Suppose Assumptions \ref{assumption:signal_plus_noise}, \ref{assumption:regularity_condition}, and \ref{assumption:weight_functions} hold. Let $(\delta_n)_{n = 1}^\infty$ be a positive sequence converging to $0$. Then for all $i\in [n]$,
\begin{align*}
\sup_{\bx_i\in B(\rho_n^{1/2}\bx_{0i}, \delta_n)}\left\|\frac{1}{n}\sum_{j = 1}^n[\bW\transpose\widetilde{\bg}_{ij}(\bW\bx_i) - \expect_0\{\bg_{ij}(\bx_i)\}]\right\|\lesssim (\delta_n + \rho_n^{1/2})\sqrt{\frac{(\log n)^{2\xi}}{n}}\quad\mbox{w.h.p.}.
\end{align*}
\end{lemma}

\begin{proof}[\bf Proof of Lemma \ref{lemma:Local_ULLN}]
Denote $\widetilde{h}_{nij}(\bx_i) = h_n(\widetilde{\bx}_i\transpose\widetilde{\bx}_j, \bx_i\transpose\widetilde{\bx}_j)$ and $\widetilde{h}_{0nij}(\bx_i) = h_n(\rho_n\bx_{0i}\transpose\bx_{0j}, \rho_n^{1/2}\bx_i\transpose{\bx}_{0j})$. The proof is almost the same as that of Lemma \ref{lemma:ULLN} except for some small modifications.
% \vspace*{2ex}\noindent
% $\blacksquare$ {\bf Proof of the first assertion.}
By triangle inequality and Cauchy-Schwarz inequality, 
\begin{align*}
&\sup_{\bx_i\in  B(\rho_n^{1/2}\bx_{0i}, \delta_n)}\left\|\frac{1}{n}\sum_{j = 1}^n[\bW\transpose\widetilde{\bg}_{ij}(\bW\bx_i) - \expect_0\{\bg_{ij}(\bx_i)\}]\right\|_2\\
&\quad\leq \sup_{\bx_i\in \Theta}
\frac{1}{n}\sum_{j = 1}^n|E_{ij}|\left\|\widetilde{h}_{nij}(\bW\bx_i)\rho_n^{-1/2}\bW\transpose\widetilde{\bx}_j - h_{0nij}(\bx_i)\bx_{0j}\right\|_2\\
&\quad\quad + \sup_{\bx_i\in  B(\rho_n^{1/2}\bx_{0i}, \delta_n)}
\frac{1}{n}\sum_{j = 1}^n|\bx_i\transpose\bW\transpose\widetilde{\bx}_j - \rho_n\bx_{0i}\transpose\bx_{0j}|\left\|\widetilde{h}_{nij}(\bW\bx_i)\rho_n^{-1/2}\bW\transpose\widetilde{\bx}_j - h_{0nij}(\bx_i)\bx_{0j}\right\|_2\\
&\quad\quad + \sup_{\bx_i\in  B(\rho_n^{1/2}\bx_{0i}, \delta_n)}\|\bx_i\|_2\left\|\frac{1}{n}\sum_{j = 1}^n(\bW\transpose\widetilde\bx_j - \rho_n^{1/2}\bx_{0j})\bx_{0j}\transpose h_{0nij}(\bx_i)\right\|_2\\
&\quad\quad + \sup_{\bx_i\in\Theta}\left\|\frac{1}{n}\sum_{j = 1}^n[\bg_{ij}(\bx_i) - \expect_0\{\bg_{ij}(\bx_i)\}]\right\|_2.
\end{align*}
By the proof of Lemma \ref{lemma:ULLN}, the first and fourth terms are $O(\rho_n^{1/2}\sqrt{(\log n)^{2\xi}/n})$ w.h.p.. The third term is bounded by a constant multiple of $(\delta_n + \rho_n^{1/2})\|\widetilde{\bX}\bW - \rho_n^{1/2}\bX_0\|_{2\to\infty} = O\{(\delta_n + \rho_n^{1/2})\sqrt{(\log n)^{2\xi}/n}\}$ w.h.p.. For the second term, we first recall in the proof of Lemma \ref{lemma:ULLN} that
\begin{align*}
\max_{i,j\in [n]}\sup_{\bx_i\in\Theta}\left\|\widetilde{h}_{nij}(\bW\bx_i)\rho_n^{-1/2}\bW\transpose\widetilde{\bx}_j - h_{0nij}(\bx_i)\bx_{0j}\right\|_2\lesssim \sqrt{\frac{(\log n)^{2\xi}}{n\rho_n}}\quad\mbox{w.h.p..}
\end{align*}
Then the second term can be bounded as follows:
\begin{align*}
&\sup_{\bx_i\in B(\rho_n^{1/2}\bx_{0i}, \delta_n)}
\frac{1}{n}\sum_{j = 1}^n|\bx_i\transpose\bW\transpose\widetilde{\bx}_j - \rho_n\bx_{0i}\transpose\bx_{0j}|\left\|\widetilde{h}_{nij}(\bW\bx_i)\rho_n^{-1/2}\bW\transpose\widetilde{\bx}_j - h_{0nij}(\bx_i)\bx_{0j}\right\|_2\\
&\quad\lesssim \max_{i,j\in [n]}\sup_{\bx_i\in B(\rho_n^{1/2}\bx_{0i}, \delta_n)}\left(\|\bx_i - \rho_n^{1/2}\bx_{0i}\|_2\|\bW\transpose\widetilde{\bx}_j\|_2 + \rho_n^{1/2}\|\bx_{0i}\|_2\|\bW\transpose\widetilde{\bx}_j - \rho_n^{1/2}\bx_{0j}\|_2\right)\sqrt{\frac{(\log n)^{2\xi}}{n\rho_n}}\\
&\quad\lesssim  (\delta_n + \rho_n^{1/2})\sqrt{\frac{(\log n)^{2\xi}}{n}}\quad\mbox{w.h.p..}
\end{align*}
Combining the high probability bounds for the four terms above completes the proof. 
\end{proof}

\begin{lemma}\label{lemma:convergence_lambda_ETEL}
Suppose Assumptions \ref{assumption:signal_plus_noise}, \ref{assumption:regularity_condition}, and \ref{assumption:weight_functions} hold. Further assume that Assumption \ref{assumption:signal_plus_noise} (vi) is strengthened to Assumption \ref{assumption:signal_plus_noise} (vi) (b). Let $\widehat{\blambda}_i(\bx_i)$ be the Lagrange multiplier given by \eqref{eqn:ETEL_probabilities_dual}
% , $\gamma_n = \sqrt{(\log n)/n}$, 
and $\delta_n = M_n\sqrt{{(\log n)^{2\xi + 1}}/{(n\rho_n)}}$, where $M_n = \log\log n$. 
% Then
% If furthermore under Assumption \ref{assumption:regularity_condition}(e)(i), $\sqrt{{(\log n)^{3\xi}}/{(n\rho_n^2)}} = O(1)$ holds, then
% \begin{itemize}
Then
%   \item[(a)] Under Assumption \ref{assumption:regularity_condition} (e)(i), if $\rho_n^{-1/2}\delta_n = O(1)$, then
%   \[
%   \sup_{\bW\transpose\bx_i\in B(\rho_n^{1/2}\bx_{0i}, \delta_n)}\|\widehat{\blambda}_i(\bx_i)\|_2\lesssim \sqrt{\frac{(\log n)^{3\xi}}{n\rho_n^2}}\quad\mbox{w.h.p.}.
%   \]
%   \item[(b)] Under Assumption \ref{assumption:regularity_condition} (e)(ii), 
  \begin{align*}
    % &\sup_{\bW\transpose\bx_i\in B(\rho_n^{1/2}\bx_{0i}, \eps_n)}\|\widehat{\blambda}_i(\bx_i)\|_2\lesssim \sqrt{\frac{(\log n)^{2\xi}}{n\rho_n}}\quad\mbox{w.h.p.},\\
  &\sup_{\bW\transpose\bx_i\in B(\rho_n^{1/2}\bx_{0i}, \delta_n)}\|\widehat{\blambda}_i(\bx_i)\|_2\lesssim M_n\sqrt{\frac{(\log n)^{2\xi + 1}}{n\rho_n^2}}\quad\mbox{w.h.p.}.
  \end{align*}
\end{lemma}

\begin{proof}[\bf Proof of Lemma \ref{lemma:convergence_lambda_ETEL}]
Denote 
\begin{align*}
P_{in}(\blambda_i, \bx_i) & = \frac{1}{n}\sum_{j = 1}^n\exp\{\blambda_i\transpose\widetilde{\bg}_{ij}(\bx_i)\},\\
\bQ_{in}(\blambda_i, \bx_i) & = \frac{\partial P_{in}}{\partial\blambda_i}(\blambda_i, \bx_i) = \frac{1}{n}\sum_{j = 1}^n\exp\{\blambda_i\transpose\widetilde{\bg}_{ij}(\bx_i)\}\widetilde{\bg}_{ij}(\bx_i).
\end{align*}
Let $\bar{\blambda}_i(\bx_i) = \argmin_{\blambda_i\in\Lambda_n}P_{in}(\blambda_i, \bx_i)$, where $\Lambda_n:= \{\blambda_i\in\mathbb{R}^d:\|\blambda_i\|_2\leq \sqrt{(\log n)^{4\xi - 1}/(n\rho_n^2)}\}$. We claim that the event $\{\bar{\blambda}_i(\bx_i) = \widehat{\blambda}_i(\bx_i)\mbox{ for all }\bW\transpose\bx_i\in B(\rho_n^{1/2}\bx_{0i}, \delta_n)\}$ w.h.p.. In fact, since $P_{in}$ is convex in $\blambda_i$, it is sufficient to show that the event
% \[
$\{\bar{\lambda}_i(\bx_i)\mbox{ is in the interior of }\Lambda_n\mbox{ for all }\bW\transpose\bx_i\in B(\rho_n^{1/2}\bx_{0i}, \delta_n)\}$
% \]
occurs w.h.p.. By the minimization property of $\bar{\blambda}_i$ and Taylor's theorem, we have
\begin{align*}
P_{in}(\zero_d, \bx_i)\geq P_{in}(\bar{\blambda}_i(\bx_i), \bx_i)
= P_{in}(\zero_d, \bx_i) + \bQ_{in}(\zero_d, \bx_i)\transpose\bar{\blambda}_i(\bx_i) + \frac{1}{2}\bar{\blambda}_i(\bx_i)\transpose\frac{\partial^2 P_{in}}{\partial\bx_i\partial\bx_i\transpose}(\theta_i(\bx_i)\bar{\blambda}_i(\bx_i), \bx_i)\bar{\blambda}_i(\bx_i),
\end{align*}
where $\theta_i(\bx_i)\in [0, 1]$. It follows from Cauchy-Schwarz inequality and a simple algebra that
\begin{align*}
\frac{1}{2}\bar{\blambda}_i(\bx_i)\transpose\left[\frac{1}{n}\sum_{j = 1}^n\exp\{\theta_i\bar{\blambda}_i(\bx_i)\transpose\widetilde{\bg}_{ij}(\bx_i)\}\widetilde{\bg}_{ij}(\bx_i)\widetilde{\bg}_{ij}(\bx_i)\transpose\right]\bar{\blambda}_i(\bx_i)
\leq \|\bar{\blambda}_i(\bx_i)\|_2\left\|\frac{1}{n}\sum_{j = 1}^n\widetilde{\bg}_{ij}(\bx_i)\right\|_2.
\end{align*}
By Lemma \ref{lemma:lambda_g_ETEL} and Assumption \ref{assumption:signal_plus_noise} (vi) (b),
$
\max_{j\in [n]}\sup_{\blambda_i\in\Lambda_n, \bx_i\in\Theta}|\blambda_i\transpose\widetilde{\bg}_{ij}(\bx_i)| = o(1)
$ 
w.h.p..
% \quad\mbox{w.h.p.},
% \]
Now let
\[
\calE_n = \left\{\min_{i,j\in [n]}\inf_{\bx_i\in\Theta}\exp\left\{\theta_i(\bx_i)\bar{\lambda}_i(\bx_i)\widetilde{\bg}_{ij}(\bx_i)\right\}\geq \frac{1}{2}\right\}.
\]
Then for any $c > 0$, there exists some $N_c\in\mathbb{N}_+$, such that $\prob_0(\calE_n) \geq 1 - n^{-c}$ for all $n\geq N_c$.
% which implies that $\min_{i,j\in [n]}\inf\{e^{\theta_i\bar{\blambda}_i\transpose\widetilde{\bg}_{ij}(\bx_i)}:\bx_i\in \Theta\}\gtrsim1$ w.h.p.. Namely, for all $\bx_i\in \Theta$, 
Observe that over $\calE_n$, for all $\bx_i\in\Theta$,
the left-hand side of the above inequality can be lower bounded by 
\[
\frac{1}{4}\|\bar{\blambda}_i(\bx_i)\|_2^2\lambda_{\min}\left\{\frac{1}{n}\sum_{j = 1}^n\widetilde{\bg}_{ij}(\bx_i)\widetilde{\bg}_{ij}(\bx_i)\transpose\right\}.
\]
It follows that over $\calE_n$, for all $\bx_i\in\Theta$, we have
\[
\|\bar{\blambda}_i(\bx_i)\|_2\leq 4\left\|\left\{\frac{1}{n}\sum_{j = 1}^n\widetilde{\bg}_{ij}(\bx_i)\widetilde{\bg}_{ij}(\bx_i)\transpose\right\}^{-1}\right\|_2\left\|\frac{1}{n}\sum_{j = 1}^n\widetilde{\bg}_{ij}(\bx_i)\right\|_2.
\]
By Lemma \ref{lemma:covariance_bound_ETEL}, we know that for any fixed $c > 0$, there exists a constant $K_c > 0$, such that the event
\begin{align*}
% &\sup_{\bW\transpose\bx_i\in B(\rho_n^{1/2}\bx_{0i}, \eps_n)}\left\|\left\{\frac{1}{n}\sum_{j = 1}^n\widetilde{\bg}_{ij}(\bx_i)\widetilde{\bg}_{ij}(\bx_i)\transpose\right\}^{-1}\right\|_2\\
% &\quad\leq 
\calF_n(K_c) = 
\left\{\sup_{\bW\transpose\bx_i\in B(\rho_n^{1/2}\bx_{0i}, \delta_n)}\left\|\left\{\frac{1}{n}\sum_{j = 1}^n\widetilde{\bg}_{ij}(\bx_i)\widetilde{\bg}_{ij}(\bx_i)\transpose\right\}^{-1}\right\|_2\leq K_c \rho_n^{-1}\right\}
% \quad\mbox{w.h.p.}.
\end{align*}
with probability at least $1 - n^{-c}$ for sufficiently large $n$. 
By Lemma \ref{lemma:Local_ULLN}, we know that
\begin{align*}
\sup_{\bW\transpose\bx_i\in B(\rho_n^{1/2}\bx_{0i}, \delta_n)}\left\|\frac{1}{n}\sum_{j = 1}^n\widetilde{\bg}_{ij}(\bx_i)\right\|_2
& = \sup_{\bx_i\in B(\rho_n^{1/2}\bx_{0i}, \delta_n)}\left\|\frac{1}{n}\sum_{j = 1}^n\bW\transpose\widetilde{\bg}_{ij}(\bW\bx_i)\right\|_2\\
&\leq \sup_{\bx_i\in B(\rho_n^{1/2}\bx_{0i}, \delta_n)}\left\|\frac{1}{n}\sum_{j = 1}^n[\bW\transpose\widetilde{\bg}_{ij}(\bW\bx_i) - \expect_0\bg_{ij}(\bx_i)]\right\|_2\\
&\quad + \sup_{\bx_i\in B(\rho_n^{1/2}\bx_{0i}, \delta_n)}\left\|\frac{1}{n}\sum_{j = 1}^n\rho_n^{1/2}(\bx_i - \rho_n^{1/2}\bx_{0i})\transpose\bx_{0j}h_{0nij}(\bx_i)\bx_{0j}\transpose\right\|_2\\
&\lesssim (\delta_n + \rho_n^{1/2})\sqrt{\frac{(\log n)^{2\xi}}{n}} + \rho_n^{1/2}\delta_n\asymp \rho_n^{1/2}\delta_n\quad\mbox{w.h.p.}
% ,
.
\end{align*}
This implies that for any fixed $c > 0$, there exists a constant $C_c > 0$, such that the event
\[
\calG_n(C_c) = \left\{
\sup_{\bW\transpose\bx_i\in B(\rho_n^{1/2}\bx_{0i}, \delta_n)}\left\|\frac{1}{n}\sum_{j = 1}^n\widetilde{\bg}_{ij}(\bx_i)\right\|_2\leq C_c \rho_n^{1/2}\delta_n
\right\}
\]
occurs with probability at least $1 - n^{-c}$ for sufficiently large $n$. 
% and similarly, we also have
% \begin{align*}
% \sup_{\bW\transpose\bx_i\in B(\rho_n^{1/2}\bx_{0i}, \eps_n)}\left\|\frac{1}{n}\sum_{j = 1}^n\widetilde{\bg}_{ij}(\bx_i)\right\|_2
% % & = \sup_{\bx_i\in B(\rho_n^{1/2}\bx_{0i}, \eps_n)}\left\|\frac{1}{n}\sum_{j = 1}^n\bW\transpose\widetilde{\bg}_{ij}(\bW\bx_i)\right\|_2\\
% % &\leq \sup_{\bx_i\in B(\rho_n^{1/2}\bx_{0i}, \eps_n)}\left\|\frac{1}{n}\sum_{j = 1}^n[\bW\transpose\widetilde{\bg}_{ij}(\bW\bx_i) - \expect_0\bg_{ij}(\bx_i)]\right\|_2\\
% % &\quad + \sup_{\bx_i\in B(\rho_n^{1/2}\bx_{0i}, \eps_n)}\left\|\frac{1}{n}\sum_{j = 1}^n\rho_n^{1/2}(\bx_i - \rho_n^{1/2}\bx_{0i})\transpose\bx_{0j}h_{0nij}(\bx_i)\bx_{0j}\bx_{0j}\transpose\right\|_2\\
% % &\lesssim (\delta_n + \rho_n^{1/2})\sqrt{\frac{(\log n)^{2\xi}}{n}} + \rho_n^{1/2}\delta_n\asymp \rho_n^{1/2}\delta_n
% \lesssim \rho_n^{1/2}\sqrt{\frac{(\log n)^{2\xi}}{n}}
% \quad\mbox{w.h.p.},
% \end{align*}
Note that over the event $\calE_n\cap \calF_n(K_c) \cap \calG_n(C_c)$, we have
\begin{align*}
% &\sup_{\bW\transpose\bx_i\in B(\rho_n^{1/2}\bx_{0i}, \eps_n)}\|\bar{\blambda}_i(\bx_i)\|_2\lesssim \sqrt{\frac{(\log n)^{2\xi}}{n\rho_n}}\quad\mbox{w.h.p.},\\
\sup_{\bW\transpose\bx_i\in B(\rho_n^{1/2}\bx_{0i}, \delta_n)}\|\bar{\blambda}_i(\bx_i)\|_2\leq 4K_cC_cM_n\sqrt{\frac{(\log n)^{2\xi + 1}}{n\rho_n^2}}.
\end{align*}
The event $\calE_n\cap \calF_n(K_c) \cap \calG_n(C_c)$ occurs with probability at least $1 - 3n^{-c}$.
This shows that the event 
$\{\bar{\blambda}_i(\bx_i)\mbox{ is in the interior of }\Lambda_n\mbox{ for all }\bW\transpose\bx_i\in B(\rho_n^{1/2}\bx_{0i}, \delta_n)\}$ occurs w.h.p.. and that
\[
\sup_{\bW\transpose\bx_i\in B(\rho_n^{1/2}\bx_{0i}, \delta_n)}\|\bar{\blambda}_i(\bx_i)\|_2\lesssim M_n\sqrt{\frac{(\log n)^{2\xi + 1}}{n\rho_n^2}}\quad\mbox{w.h.p.}.
\]
Replacing $\bar{\blambda}_i(\bx_i)$ by $\widehat{\blambda}_i(\bx_i)$ in the above concentration bound completes the proof.
\end{proof}

\begin{lemma}\label{lemma:EP_convergence_ETEL}
Suppose Assumptions \ref{assumption:signal_plus_noise}, \ref{assumption:regularity_condition}, and \ref{assumption:weight_functions} hold. Further assume that Assumption \ref{assumption:signal_plus_noise} (vi) is strengthened to Assumption \ref{assumption:signal_plus_noise} (vi) (b). Let $\widehat{\blambda}_i(\bx_i)$ be the Lagrange multiplier given by \eqref{eqn:ETEL_probabilities_dual}, $M_n = \log\log n$, and $\delta_n = M_n\sqrt{(\log n)^{2\xi + 1}/(n\rho_n)}$. Then
\begin{align*}
% &\max_{j\in [n]}\sup_{\bW\transpose\bx_i\in B(\rho_n^{1/2}\bx_{0i}, \eps_n)}|1 - np_{ij}(\bx_i)|\lesssim   \sqrt{\frac{(\log n)^{2\xi + 1}}{n}},\\
&\max_{j\in [n]}\sup_{\bW\transpose\bx_i\in B(\rho_n^{1/2}\bx_{0i}, \delta_n)}|1 - np_{ij}(\bx_i)|\lesssim M_n
\sqrt{\frac{(\log n)^{2\xi + 2}}{n\rho_n}},
\end{align*}
where $p_{ij}(\bx_i)$, $j\in [n]$ are the empirical probabilities given by \eqref{eqn:ETEL_probabilities}. 
\end{lemma}
\begin{proof}[\bf Proof of Lemma \ref{lemma:EP_convergence_ETEL}]
% We first focus on the concentration bound for 
% \[
% \max_{j\in [n]}\sup_{\bx_i\in B(\rho_n^{1/2}\bx_{0i}, \eps_n)}|\exp\{\widehat{\blambda}_i(\bx_i)\transpose\widetilde{\bg}_{ij}(\bx_i)\} - 1|\quad\mbox{and}\quad
% \max_{j\in [n]}\sup_{\bx_i\in B(\rho_n^{1/2}\bx_{0i}, \delta_n)}|\exp\{\widehat{\blambda}_i(\bx_i)\transpose\widetilde{\bg}_{ij}(\bx_i)\} - 1|.
% \]
By mean-value theorem, there exists some $s_{ij}(\bx_i)$ adjoining $0$ and $\widehat{\blambda}_i(\bx_i)\transpose\widetilde{\bg}_{ij}(\bx_i)$, such that
\begin{align*}
% \max_{j\in [n]}\sup_{\bx_i\in B(\rho_n^{1/2}\bx_{0i}, \eps_n)}|\exp\{\widehat{\blambda}_i(\bx_i)\transpose\widetilde{\bg}_{ij}(\bx_i)\} - 1|
% &\leq \max_{j\in [n]}\sup_{\bx_i\in B(\rho_n^{1/2}\bx_{0i}, \eps_n)}e^{s_{ij}(\bx_i)}|\widehat{\blambda}_i(\bx_i)\transpose\widetilde{\bg}_{ij}(\bx_i)|,\\
\max_{j\in [n]}\sup_{\bW\transpose\bx_i\in B(\rho_n^{1/2}\bx_{0i}, \delta_n)}|\exp\{\widehat{\blambda}_i(\bx_i)\transpose\widetilde{\bg}_{ij}(\bx_i)\} - 1|
&\leq \max_{j\in [n]}\sup_{\bW\transpose\bx_i\in B(\rho_n^{1/2}\bx_{0i}, \delta_n)}e^{s_{ij}(\bx_i)}|\widehat{\blambda}_i(\bx_i)\transpose\widetilde{\bg}_{ij}(\bx_i)|.
\end{align*}
By Lemma \ref{lemma:convergence_lambda_ETEL} and Result \ref{result:uniform_concentration_g}, $\sup_{\bx_i\in B(\rho_n^{1/2}\bx_{0i}, \delta_n)}e^{s_{ij}(\bx_i)}\leq 2$ w.h.p.. Then again, by Lemma \ref{lemma:convergence_lambda_ETEL}, Cauchy-Schwarz inequality, and Result \ref{result:uniform_concentration_g}, we have
\begin{align*}
\max_{j\in [n]}\sup_{\bW\transpose\bx_i\in B(\rho_n^{1/2}\bx_{0i}, \delta_n)}|\exp\{\widehat{\blambda}_i(\bx_i)\transpose\widetilde{\bg}_{ij}(\bx_i)\} - 1|
&\lesssim \max_{j\in [n]}\sup_{\bW\transpose\bx_i\in B(\rho_n^{1/2}\bx_{0i}, \delta_n)}\|\widehat{\blambda}_i(\bx_i)\|_2\|\widetilde{\bg}_{ij}(\bx_i)\|_2\\
&\lesssim M_n\sqrt{\frac{(\log n)^{2\xi + 2}}{n\rho_n}}\quad\mbox{w.h.p.},
\end{align*}
% and similarly,
% \begin{align*}
% \max_{j\in [n]}\sup_{\bx_i\in B(\rho_n^{1/2}\bx_{0i}, \delta_n)}|\exp\{\widehat{\blambda}_i(\bx_i)\transpose\widetilde{\bg}_{ij}(\bx_i)\} - 1|&\lesssim M_n\sqrt{\frac{(\log n)^{2\xi + 1}}{n\rho_n}}\quad\mbox{w.h.p.}.
% \end{align*}
Observe that
\begin{align*}
|1 - np_{ij}|
& = \left|\frac{(1/n)\sum_{m = 1}^n\{(e^{\widehat{\blambda}_i\transpose\widetilde{\bg}_{im}} - 1) + (1 - e^{\widehat{\blambda}_i\transpose\widetilde{\bg}_{ij}})\} }{(1/n)\sum_{m = 1}^n\{(e^{\widehat{\blambda}_i\transpose\widetilde{\bg}_{im}} - 1) + 1\}}\right|
% \\&
\leq \frac{2\max_{m\in [n]}|e^{\widehat{\blambda}_i\transpose \widetilde{\bg}_{im} } - 1|}{1 - 2\max_{m\in [n]}|e^{\widehat{\blambda}_i\transpose \widetilde{\bg}_{im} } - 1|},
\end{align*}
where we have suppressed the argument $\bx_i$. It follows from the previous concentration bounds that
\begin{align*}
% \max_{j\in [n]}\sup_{\bW\transpose\bx_i\in B(\rho_n^{1/2}\bx_{0i}, \eps_n)}|1 - np_{ij}(\bx_i)|
% &\lesssim 
% \sqrt{\frac{(\log n)^{2\xi + 1}}{n}}\quad\mbox{w.h.p.},\\
\max_{j\in [n]}\sup_{\bW\transpose\bx_i\in B(\rho_n^{1/2}\bx_{0i}, \delta_n)}|1 - np_{ij}(\bx_i)|
&\lesssim M_n
\sqrt{\frac{(\log n)^{2\xi + 2}}{n\rho_n}}\quad\mbox{w.h.p.}.
\end{align*}
This completes the proof of the lemma.
\end{proof}

\begin{lemma}\label{lemma:lambda_gradient_convergence_ETEL}
Suppose Assumptions \ref{assumption:signal_plus_noise}, \ref{assumption:regularity_condition}, and \ref{assumption:weight_functions} hold. Further assume that Assumption \ref{assumption:signal_plus_noise} (vi) is strengthened to Assumption \ref{assumption:signal_plus_noise} (vi) (b). Let $\widehat{\blambda}_i(\bx_i)$ be the Lagrange multiplier given by \eqref{eqn:ETEL_probabilities_dual}
% , $\eps_n = \sqrt{(\log n)/n}$, 
and $\delta_n = M_n\sqrt{(\log n)^{2\xi + 1}/(n\rho_n)}$, where $M_n = \log\log n$. Denote
\begin{align*}
\widetilde{\bG}_{in}(\bx_i) & = \frac{1}{n}\sum_{j = 1}^ne^{\widehat{\blambda}_i(\bx_i)\transpose\widetilde{\bg}_{ij}(\bx_i)}\{\eye_d + \bW\transpose\widetilde{\bg}_{ij}(\bx_i)\widehat{\blambda}_i(\bx_i)\transpose\bW\}\bW\transpose\frac{\partial\widetilde{\bg}_{ij}}{\partial\bx_i\transpose}(\bx_i)\bW,\\
\widetilde{\bOmega}_{in}(\bx_i) & = \frac{1}{n}\sum_{j = 1}^ne^{\widehat{\blambda}_i(\bx_i)\transpose\widetilde{\bg}_{ij}(\bx_i)}\bW\transpose\widetilde{\bg}_{ij}(\bx_i)\widetilde{\bg}_{ij}(\bx_i)\transpose\bW.
\end{align*}
Then
\begin{align*}
&\sup_{\bx_i\in B(\rho_n^{1/2}\bx_{0i}, \delta_n)}\|\widetilde{\bG}_{in}(\bW\bx_i) - \bG_{in}(\rho_n^{1/2}\bx_{0i})\|_2
\lesssim M_n\sqrt{\frac{(\log n)^{2\xi + 2}}{n}}\quad\mbox{w.h.p.},\\
&\sup_{\bx_i\in B(\rho_n^{1/2}\bx_{0i}, \delta_n)}\|\widetilde{\bOmega}_{in}(\bW\bx_i) - \bOmega_{in}(\rho_n^{1/2}\bx_{0i})\|_2\lesssim \rho_n^{1/2}M_n\sqrt{\frac{(\log n)^{2\xi + 2}}{n}} \quad\mbox{w.h.p.},\\
&\sup_{\bx_i\in B(\rho_n^{1/2}\bx_{0i}, \delta_n)}\left\|\bW\transpose\frac{\partial\widehat{\blambda}_i}{\partial\bx_i\transpose}(\bW\bx_i)\bW + \bOmega_{in}(\rho_n^{1/2}\bx_{0i})^{-1}\bG_{in}(\rho_n^{1/2}\bx_{0i})\right\|_2
\lesssim \rho_n^{-1}M_n\sqrt{\frac{(\log n)^{2\xi + 2}}{n}}\quad\mbox{w.h.p.},\\
&\sup_{\bx_i\in B(\rho_n^{1/2}\bx_{0i}, \delta_n)}\left\|\frac{\partial\widetilde{\bG}_{in}}{\partial x_{ik}}(\bW\bx_i)\right\|_2\lesssim \rho_n^{1/2}\quad\mbox{w.h.p.},\quad k\in [d],\\
&\sup_{\bx_i\in B(\rho_n^{1/2}\bx_{0i}, \delta_n)}\left\|\frac{\partial\widetilde{\bOmega}_{in}}{\partial x_{ik}}(\bW\bx_i)\right\|_2\lesssim \rho_n \quad\mbox{w.p.a.1.},\quad k\in [d],\\
&\sup_{\bx_i\in B(\rho_n^{1/2}\bx_{0i}, \delta_n)}\left\|\frac{\partial}{\partial x_{ik}}\left\{\frac{\partial\widehat{\blambda}_{i}}{\partial \bx_{i}\transpose}(\bW\bx_i)\right\}\right\|_2\lesssim \rho_n^{-1/2} \quad\mbox{w.p.a.1.},\quad k\in [d],\\
&\sup_{\bx_i\in B(\rho_n^{1/2}\bx_{0i}, \delta_n)}\frac{1}{n}\sum_{j = 1}^n\left\|\frac{\partial(\widehat{\blambda}_i\transpose\widetilde{\bg}_{ij})}{\partial\bx_i}(\bW \bx_i)\right\|_2^2\lesssim 1\quad\mbox{w.h.p.},\\
&\sup_{\bx_i\in B(\rho_n^{1/2}\bx_{0i}, \delta_n)}\frac{1}{n}\sum_{j = 1}^n\left\|\frac{\partial^2(\widehat{\blambda}_{i}\transpose\widetilde{\bg}_{ij})}{\partial \bx_i\partial\bx_i\transpose}(\bW\bx_i)\right\|_2\lesssim 1 \quad\mbox{w.p.a.1.}
\end{align*}
\end{lemma}

\begin{proof}[\bf Proof of Lemma \ref{lemma:lambda_gradient_convergence_ETEL}]
{$\blacksquare$ \bf Proof of the first assertion.}
By triangle inequality and Cauchy-Schwarz inequality, we write
\begin{align*}
\|\widetilde{\bG}_{in}(\bW\bx_i) - \bG_{in}(\rho_n^{1/2}\bx_{0i})\|_2
&\leq\frac{1}{n}\sum_{j = 1}^n\left\{\left|e^{\widehat{\blambda}_i(\bW\bx_i)\transpose\widetilde{\bg}_{ij}(\bW\bx_i)} - 1\right|
 + \|\widehat{\blambda}_i(\bW\bx_i)\|_2\|\widetilde{\bg}_{ij}(\bW\bx_i)\|_2
 \right\}
\left\|\frac{\partial\widetilde{\bg}_{ij}}{\partial\bx_i\transpose}(\bW\bx_i)\right\|_2\\
&\quad + \left\|\frac{1}{n}\sum_{j = 1}^n\bW\transpose\frac{\partial\widetilde{\bg}_{ij}}{\partial\bx_i\transpose}(\bW\bx_i)\bW - \bG_{in}(\bx_i)\right\|_2 + \|\bG_{in}(\bx_i) - \bG_{in}(\rho_n^{1/2}\bx_{0i})\|_2.
\end{align*}
By the second assertion of Lemma \ref{lemma:ULLN}, the second term is $O(\sqrt{(\log n)^{2\xi}/n})$ w.h.p. uniformly in $\bx_i\in\Theta$. For any $\delta > 0$, the third term can be bounded by
\begin{align*}
&\sup_{\bx_i\in B(\rho_n^{1/2}\bx_{0i}, \delta)}\|\bG_{in}(\bx_i) - \bG_{in}(\rho_n^{1/2}\bx_{0i})\|_2\\
&\quad\leq \sup_{\bx_i\in B(\rho_n^{1/2}\bx_{0i}, \delta)}\frac{1}{n}\sum_{j = 1}^n\rho_n\|\rho_n^{1/2}\bx_{0i} - \bx_{i}\|_2\|\bx_{0j}\|_2|D^{(0, 1)}h_{0nij}(\bx_i)|\|\bx_{0j}\|_2^2\\
&\quad\quad + \sup_{\bx_i\in B(\rho_n^{1/2}\bx_{0i}, \delta)}\frac{1}{n}\sum_{j = 1}^n\rho_n^{1/2}|h_{0nij}(\bx_i) - h_{0nij}(\rho_n^{1/2}\bx_{0i})|\|\bx_{0j}\|_2^2\lesssim {\color{red}\rho_n^2}\delta
\end{align*}
by Assumption \ref{assumption:weight_functions}. For the first term, by Lemma \ref{lemma:convergence_lambda_ETEL} and Result \ref{result:uniform_concentration_g}, under Assumption \ref{assumption:signal_plus_noise} (vi) (b),
\[
\max_{j\in [n]}\sup_{\bx_i\in B(\rho_n^{1/2}\bx_{0i}, \delta_n)}|\widehat{\blambda}_i(\bW\bx_i)\transpose\widetilde{\bg}_{ij}(\bW\bx_i)|\leq \max_{j\in [n]}\sup_{\bW\transpose\bx_i\in B(\rho_n^{1/2}\bx_{0i}, \delta_n)}\|\widehat{\blambda}_i(\bx_i)\|_2\sup_{\bx_i\in\Theta}\|\widetilde{\bg}_{ij}(\bx_i)\|_2\lesssim 1\quad\mbox{w.h.p.},
\]
and by mean-value theorem and the proof of Lemma \ref{lemma:EP_convergence_ETEL},
\begin{align*}
% &\max_{j\in [n]}\sup_{\bx_i\in B(\rho_n^{1/2}\bx_{0i}, \eps_n)}
% \left\{\left|e^{\widehat{\blambda}_i(\bW\bx_i)\transpose\widetilde{\bg}_{ij}(\bW\bx_i)} - 1\right|
%  + \|\widehat{\blambda}_i(\bW\bx_i)\|_2\|\widetilde{\bg}_{ij}(\bW\bx_i)\|_2
%  \right\}\lesssim\sqrt{\frac{(\log n)^{2\xi + 1}}{n}}\quad\mbox{w.h.p.},\\
&\max_{j\in [n]}\sup_{\bx_i\in B(\rho_n^{1/2}\bx_{0i}, \delta_n)}
\left\{\left|e^{\widehat{\blambda}_i(\bW\bx_i)\transpose\widetilde{\bg}_{ij}(\bW\bx_i)} - 1\right|
 + \|\widehat{\blambda}_i(\bW\bx_i)\|_2\|\widetilde{\bg}_{ij}(\bW\bx_i)\|_2
 \right\}
\lesssim M_n\sqrt{\frac{(\log n)^{2\xi + 2}}{n\rho_n}}\quad\mbox{w.h.p.},
\end{align*}
Also, by Lemma \ref{lemma:Sample_moments_g}, 
\begin{align*}
\sup_{\bx_i\in \Theta}\frac{1}{n}\sum_{j = 1}^n\left\|\frac{\partial\widetilde{\bg}_{ij}}{\partial\bx_i\transpose}(\bx_i)\right\|_2
% &\leq 
% \sup_{\bx_i\in \Theta}\max_{j\in [n]}\left\|\frac{\partial\widetilde{\bg}_{ij}}{\partial\bx_i\transpose}(\bx_i)\right\|_2
% \\&
\lesssim \rho_n^{1/2}\quad\mbox{w.h.p.}.
\end{align*}
We then obtain from combining the above concentration bounds that
\begin{align*}
&\sup_{\bx_i\in B(\rho_n^{1/2}\bx_{0i}, \delta_n)}\|\widetilde{\bG}_{in}(\bW\bx_i) - \bG_{in}(\rho_n^{1/2}\bx_{0i})\|_2
\lesssim M_n\sqrt{\frac{(\log n)^{2\xi + 2}}{n}}\quad\mbox{w.h.p.}
.
% ,\\
% &\sup_{\bx_i\in B(\rho_n^{1/2}\bx_{0i}, \delta_n)}\|\widetilde{\bG}_{in}(\bW\bx_i) - \bG_{in}(\rho_n^{1/2}\bx_{0i})\|_2
% \lesssim \max(1,\rho_n\log n)\sqrt{\frac{(\log n)^{3\xi}}{n\rho_n}}\quad\mbox{w.h.p.}.
\end{align*}
\vspace*{2ex}
{$\blacksquare$ \bf Proof of the second assertion.} By triangle inequality, we have
\begin{align*}
\|\widetilde{\bOmega}_{in}(\bW\bx_i) - \bOmega_{in}(\rho_n^{1/2}\bx_{0i})\|_2
&\leq \frac{1}{n}\sum_{j = 1}^n|\exp\{\widehat{\blambda}_i(\bW\bx_i)\transpose\widetilde{\bg}_{ij}(\bW\bx_i)\} - 1|\|\widetilde{\bg}_{ij}(\bW\bx_i)\|_2^2\\
&\quad + \left\|\frac{1}{n}\sum_{j = 1}^n\bW\transpose\widetilde{\bg}(\bW\bx_i)\widetilde{\bg}(\bW\bx_i)\transpose\bW - \bOmega_{in}(\bx_i)\right\|_2 + \|\bOmega_{in}(\bx_i) - \bOmega_{in}(\rho_n^{1/2}\bx_{0i})\|_2
\end{align*}
By Lemma \ref{lemma:covariance_bound_ETEL}, we see that the second term is $O\{\rho_n^{1/2}M_n\sqrt{(\log n)^{2\xi + 1}/n}\}$ w.h.p. uniformly over $\bx_i\in B(\rho_n^{1/2}\bx_{0i}, \delta_n)$. For the third term, we have
\begin{align*}
&\sup_{\bx_i\in B(\rho_n^{1/2}\bx_{0i}, \delta_n)}\|\bOmega_{in}(\bx_i) - \bOmega_{in}(\rho_n^{1/2}\bx_{0i})\|_2\\
&\quad\leq \sup_{\bx_i\in B(\rho_n^{1/2}\bx_{0i}, \delta_n)}\frac{1}{n}\sum_{j = 1}^n\rho_n\|\bx_i - \rho_n^{1/2}\bx_{0i}\|_2^2\|\bx_{0j}\|_2^4h_{0nij}^2(\bx_i)\\
&\quad\quad + \sup_{\bx_i\in B(\rho_n^{1/2}\bx_{0i}, \delta_n)}\frac{1}{n}\sum_{j = 1}^n\var_0(A_{ij})\|\bx_{0j}\|_2^2|h_{0nij}^2(\bx_i) - h_{0nij}^2(\rho_n^{1/2}\bx_{0i})|\\
&\quad\lesssim \rho_n\delta_n = \rho_n^{1/2}M_n\sqrt{\frac{(\log n)^{2\xi + 1}}{n}}.
\end{align*}
% By Lemma \ref{lemma:covariance_bound_ETEL} and its proof, the second term is $O(\rho_n^{1/2}\sqrt{(\log n)^{2\xi}/n})$ w.h.p. uniformly in $\bx_i\in B(\rho_n^{1/2}\bx_{0i}, \delta_n)$. 
For the first term, by Lemma \ref{lemma:Sample_moments_g},
\begin{align*}
\sup_{\bx_i\in B(\rho_n^{1/2}\bx_{0i}, \delta_n)}
\frac{1}{n}\sum_{j = 1}^n\|\widetilde{\bg}_{ij}(\bW\bx_i)\|_2^2
\lesssim \rho_n\quad\mbox{w.h.p.}.
\end{align*}
Then it follows from the proof of the first assertion that
\[
\sup_{\bx_i\in B(\rho_n^{1/2}\bx_{0i}, \delta_n)}\frac{1}{n}\sum_{j = 1}^n|\exp\{\widehat{\blambda}_i(\bW\bx_i)\transpose\widetilde{\bg}_{ij}(\bW\bx_i)\} - 1|\|\widetilde{\bg}_{ij}(\bW\bx_i)\|_2^2\lesssim \rho_n^{1/2} M_n\sqrt{\frac{(\log n)^{2\xi + 2}}{n}}\quad\mbox{w.h.p.}.
\]
Therefore, we conclude that 
\[
\sup_{\bx_i\in B(\rho_n^{1/2}\bx_{0i}, \delta_n)}\|\widetilde{\bOmega}_{in}(\bW\bx_i) - \bOmega_{in}(\rho_n^{1/2}\bx_{0i})\|_2\lesssim \rho_n^{1/2}M_n\sqrt{\frac{(\log n)^{2\xi + 2}}{n}}\quad\mbox{w.h.p.}.
\]

\vspace*{2ex}
\noindent
{$\blacksquare$ \bf Proof of the third assertion.} Since $\widehat{\blambda}_i(\bx_i)$ is the Lagrange multiplier defined by \eqref{eqn:ETEL_probabilities_dual}, then it satisfies the equation
\[
\frac{1}{n}\sum_{j = 1}^n\exp\{\widehat{\blambda}_i(\bx_i)\transpose\widetilde{\bg}_{ij}(\bx_i)\}\widetilde{\bg}_{ij}(\bx_i) = 0.
\]
By the implicit function theorem, 
\[
\bW\transpose\frac{\partial\widehat{\blambda}_i}{\partial\bx_i\transpose}(\bx_i)\bW = -\widetilde{\bOmega}_{in}^{-1}(\bx_i)\widetilde{\bG}_{in}(\bx_i).
\]
Denote $\bOmega_{0in} = \bOmega_{in}(\rho_n^{1/2}\bx_{0i})$ and $\bG_{0in}(\rho_n^{1/2}\bx_{0i})$. By Cauchy-Schwarz inequality, we write
\begin{align*}
&\sup_{\bx_i\in B(\rho_n^{1/2}\bx_{0i}, \delta_n)}\left\|\bW\transpose\frac{\partial\widehat{\blambda}_i}{\partial\bx_i\transpose}(\bW\bx_i)\bW + \bOmega_{0in}^{-1}\bG_{0in}\right\|_2\\
&\quad\leq \sup_{\bx_i\in B(\rho_n^{1/2}\bx_{0i}, \delta_n)}\|\widetilde{\bOmega}_{in}(\bW\bx_i)^{-1}\|_2\|\widetilde{\bG}_{in}(\bW\bx_i) - \bG_{0in}\|_2\\
&\quad\quad + \sup_{\bx_i\in B(\rho_n^{1/2}\bx_{0i}, \delta_n)}\|\widetilde{\bOmega}_{in}(\bW\bx_i)^{-1}\|_2\|\widetilde{\bOmega}_{in}(\bW\bx_i) - \bOmega_{0in}\|_2\|\bOmega_{0in}^{-1}\|_2\|\bG_{0in}\|_2.
\end{align*}
By Assumption \ref{assumption:regularity_condition} and the second assertion, $\widetilde{\bOmega}_{in}(\bW\bx_i) = \bOmega_{0in}(\bx_i) + o(\rho_n)$ uniformly in $\bx_i\in B(\rho_n^{1/2}\bx_{0i}, \delta_n)$ w.h.p. and $\bOmega_{0in}$ has eigenvalues bounded from below and above by constant multiples of $\rho_n$. It follows that $\|\widetilde{\bOmega}_{in}(\bW\bx_i)\|_2 = O(\rho_n)$ and $\|\widetilde{\bOmega}_{in}(\bW\bx_i)^{-1}\|_2 = O(\rho_n^{-1})$ w.h.p. uniformly in $\bx_i \in B(\rho_n^{1/2}\bx_{0i}, \delta_n)$. Hence, from the conclusions of the first and second assertions, we have
\[
\sup_{\bx_i\in B(\rho_n^{1/2}\bx_{0i}, \delta_n)}\left\|\bW\transpose\frac{\partial\widehat{\blambda}_i}{\partial\bx_i\transpose}(\bW\bx_i)\bW + \bOmega_{0in}^{-1}\bG_{0in}\right\|_2
\lesssim \rho_n^{-1}M_n\sqrt{\frac{(\log n)^{2\xi + 2}}{n}}\quad\mbox{w.h.p.}.
\]
\vspace*{2ex}
\noindent
{$\blacksquare$ \bf Proof of the fourth assertion.} 
We suppress the argument $\bW\bx_i$ for notational simplicity and compute by Cauchy-Schwarz inequality:
\begin{align*}
\left\|\frac{\partial\widetilde{\bG}_{in}}{\partial x_{ik}}\right\|_2
& \leq \left\|\frac{1}{n}\sum_{j = 1}^ne^{\widehat{\blambda}_i\transpose\widetilde{\bg}_{ij}}\left(\frac{\partial\widehat{\blambda}_i\transpose}{\partial x_{ik}}\widetilde{\bg}_{ij} + \widehat{\blambda}_i\transpose\frac{\partial\widetilde{\bg}_{ij}}{\partial x_{ik}}\right)(\eye_d +\bW\transpose \widetilde{\bg}_{ij}\widehat{\blambda}_i\transpose\bW)\bW\transpose\frac{\partial\widetilde{\bg}_{ij}}{\partial\bx_i\transpose}\right\|_2\\
&\quad + \left\|\frac{1}{n}\sum_{j = 1}^ne^{\widehat{\blambda}_i\transpose\widetilde{\bg}_{ij}}\left(\frac{\partial\widetilde{\bg}_{ij}}{\partial x_{ik}}\widehat{\blambda}_i\transpose + \widetilde{\bg}_{ij}\frac{\partial\widehat{\blambda}_i\transpose}{\partial x_{ik}}\right)\frac{\partial\widetilde{\bg}_{ij}}{\partial\bx_i\transpose}\right\|_2\\
&\quad + \left\|\frac{1}{n}\sum_{j = 1}^ne^{\widehat{\blambda}_i\transpose\widetilde{\bg}_{ij}}(\eye_d + \bW\transpose\widetilde{\bg}_{ij}\widehat{\blambda}_i\transpose\bW)\bW\transpose\frac{\partial}{\partial x_{ik}}\left(\frac{\partial\widetilde{\bg}_{ij}}{\partial\bx_i\transpose}\right)\right\|_2\\
&\leq \max_{j\in [n]}e^{\widehat{\blambda}_i\transpose\widetilde{\bg}_{ij}}(2 + \|\widehat{\blambda}_i\|_2\|\widetilde{\bg}_{ij}\|_2)\\
&\quad\times\left\{\left\|\frac{\partial\widehat{\blambda}_i\transpose}{\partial \bx_{i}\transpose}\right\|_2\left(\frac{1}{n}\sum_{j = 1}^n\|\widetilde{\bg}_{ij}\|_2^2\right)^{1/2}\left(\frac{1}{n}\sum_{j = 1}^n\left\|\frac{\partial\widetilde{\bg}_{ij}}{\partial\bx_i\transpose}\right\|_2^2\right)^{1/2} + \|\widehat{\blambda}_i\|\frac{1}{n}\sum_{j = 1}^n\left\|\frac{\partial\widetilde{\bg}_{ij}}{\partial\bx_i\transpose}\right\|_2^2\right\}\\
&\quad + \max_{j\in [n]}e^{\widehat{\blambda}_i\transpose\widetilde{\bg}_{ij}}(1 + \|\widehat{\blambda}_i\|_2\|\widetilde{\bg}_{ij}\|_2)\frac{1}{n}\sum_{j = 1}^n\left\|\frac{\partial}{\partial x_{ik}}\left(\frac{\partial\widetilde{\bg}_{ij}}{\partial\bx_i\transpose}\right)\right\|_2.
\end{align*}
By Lemma \ref{lemma:Sample_moments_g},
\begin{align*}
\sup_{\bx_i\in\Theta}\frac{1}{n}\sum_{k = 1}^d\sum_{j = 1}^n\left\|\frac{\partial^2[\widetilde{\bg}_{ij}]_k}{\partial\bx_i\partial\bx_i\transpose}(\bx_i)
\right\|_2\lesssim \rho_n^2 \quad\mbox{w.h.p.},
\end{align*}
where $[\cdot]_k$ denotes the $k$th coordinate of the vector. 
From the proof of the first assertion, we have
\[
\sup_{\bx_i\in B(\rho_n^{1/2}\bx_{0i}, \delta_n)}\max_{j\in [n]}\left\{e^{\widehat{\blambda}_i(\bW\bx_i)\transpose\widetilde{\bg}_{ij}(\bW\bx_i)} + \|\widehat{\blambda}_i(\bW\bx_i)\|_2\|\widetilde{\bg}_{ij}(\bW\bx_i)\|_2\right\}\lesssim 1\quad\mbox{w.h.p.}.
\]
This shows that the second term is $O(\rho_n^2)$ uniformly in $\bx_i\in B(\rho_n^{1/2}\bx_{0i}, \delta_n)$ w.h.p.. Also, by the third assertion, Result \ref{result:Jacobian}, and Result \ref{result:second_moment_matrix}, 
\[
\sup_{\bx_i\in B(\rho_n^{1/2}\bx_{0i}, \delta_n)}\left\|\frac{\partial\widehat{\blambda}_i}{\partial\bx_i\transpose}(\bW\bx_i)\right\|_2\lesssim \rho_n^{-1/2}M_n\sqrt{\frac{(\log n)^{2\xi + 2}}{n\rho_n}} +  \|\bOmega_{0in}^{-1}\bG_{0in}\|_2\lesssim \rho_n^{-1/2}\quad\mbox{w.h.p..}
\]
By Lemma \ref{lemma:Sample_moments_g},
\[
\sup_{\bx_i\in\Theta}\frac{1}{n}\sum_{j = 1}^n\left\|\widetilde{\bg}_{ij}(\bx_i)\right\|_2^2\lesssim \rho_n\quad\mbox{w.h.p.}
\quad\mbox{and}\quad
% % \sup_{\bx_i\in\Theta}\frac{1}{n}\sum_{j = 1}^n\left\|\frac{\partial\widetilde{\bg}_{ij}}{\partial\bx_i\transpose}(\bx_i)\right\|_2^2\lesssim \rho_n \quad\mbox{w.h.p.}
% .
% \]
% Following the similar reasoning, we have
% \[
\sup_{\bx_i\in\Theta}\frac{1}{n}\sum_{j = 1}^n\left\|\frac{\partial\widetilde{\bg}_{ij}}{\partial\bx_i\transpose}(\bx_i)\right\|_2^2\lesssim \rho_n\quad\mbox{w.h.p.}
\]
Also, recall from Lemma \ref{lemma:convergence_lambda_ETEL} that 
\[
\sup_{\bx_i\in B(\rho_n^{1/2}\bx_{0i}, \delta_n)}\|\widehat{\blambda}_i(\bW\bx_i)\|_2\lesssim M_n\sqrt{\frac{(\log n)^{2\xi + 1}}{n\rho_n^2}}\quad\mbox{w.h.p.}.
\]
It follows that the first term is $O(\rho_n^{1/2})$ w.h.p.. Combining the above concentration bounds yields that
\[
\max_{k\in [d]}\sup_{\bx_i\in B(\rho_n^{1/2}\bx_{0i}, \delta_n)}\left\|\frac{\partial\widetilde{\bG}_{in}}{\partial x_{ik}}(\bW\bx_i)\right\|_2\lesssim \rho_n^{1/2}\quad\mbox{w.h.p.}.
\]

\vspace*{2ex}
\noindent
{$\blacksquare$ \bf Proof of the fifth assertion.} Compute the derivative:
\begin{align*}
\frac{\partial\widetilde{\bOmega}_{in}}{\partial x_{ik}}(\bx_i) &=
\frac{1}{n}\sum_{j = 1}^ne^{\widehat{\blambda}_i(\bx_i)\transpose\widetilde{\bg}_{ij}(\bx_i)}\left\{\frac{\partial\widehat\blambda_i}{\partial x_{ik}}(\bx_i)\transpose\widetilde\bg_{ij}(\bx_i) + \widehat\blambda_i(\bx_i)\transpose\frac{\partial\widetilde\bg_{ij}}{\partial x_{ik}}(\bx_i)\right\}\bW\transpose\widetilde\bg_{ij}(\bx_i)\widetilde\bg_{ij}(\bx_i)\transpose\bW \\
&\quad+ \frac{1}{n}\sum_{j = 1}^ne^{\widehat{\blambda}_i(\bx_i)\transpose\widetilde{\bg}_{ij}(\bx_i)}\left\{\bW\transpose\frac{\partial\widetilde\bg_{ij}}{\partial x_{ik}}(\bx_i)\widetilde\bg_{ij}(\bx_i)\transpose\bW + \bW\transpose\widetilde\bg_{ij}(\bx_i)\frac{\partial\widetilde\bg_{ij}}{\partial x_{ik}}(\bx_i)\transpose\bW\right\}
\end{align*}
Then by Cauchy-Schwarz inequality, 
\begin{align*}
&\sup_{\bx_i\in B(\rho_n^{1/2}\bx_{0i}, \delta_n)}\left\|\frac{\partial\widetilde{\bOmega}_{in}}{\partial x_{ik}}(\bW\bx_i)\right\|_2\\
&\quad \leq \sup_{\bx_i\in B(\rho_n^{1/2}\bx_{0i}, \delta_n)}\frac{1}{n}\sum_{j = 1}^ne^{\widehat{\blambda}_i(\bW\bx_i)\transpose\widetilde{\bg}_{ij}(\bW\bx_i)}\left\|\frac{\partial\widehat{\blambda}_i}{\partial x_{ik}}(\bW\bx_i)\right\|_2\|\widetilde{\bg}_{ij}(\bW\bx_i)\|_2^3\\
&\quad\quad + \sup_{\bx_i\in B(\rho_n^{1/2}\bx_{0i}, \delta_n)}\frac{1}{n}\sum_{j = 1}^n e^{\widehat{\blambda}_i(\bW\bx_i)\transpose\widetilde{\bg}_{ij}(\bW\bx_i)}\left\|\widehat{\blambda}_i(\bW\bx_i)\right\|_2\left\|\frac{\partial\widetilde{\bg}_{ij}}{\partial x_{ik}}(\bW\bx_i)\right\|_2 \left\|\widetilde\bg_{ij}(\bW\bx_i)\right\|_2^2\\
&\quad\quad + \sup_{\bx_i\in B(\rho_n^{1/2}\bx_{0i}, \delta_n)}\frac{1}{n}\sum_{j = 1}^n 2e^{\widehat{\blambda}_i(\bW\bx_i)\transpose\widetilde{\bg}_{ij}(\bW\bx_i)}\left\|\frac{\partial\widetilde{\bg}_{ij}}{\partial x_{ik}}(\bW\bx_i)\right\|_2 \left\|\widetilde\bg_{ij}(\bW\bx_i)\right\|_2\\
&\quad\leq \sup_{\bx_i\in B(\rho_n^{1/2}\bx_{0i}, \delta_n)}\max_{j\in [n]}e^{\widehat{\blambda}_i(\bW\bx_i)\transpose\widetilde{\bg}_{ij}(\bW\bx_i)}\left\|\frac{\partial\widehat{\blambda}_i}{\partial x_{ik}}(\bW\bx_i)\right\|_2\times\sup_{\bx_i\in\Theta}\left\{\frac{1}{n}\sum_{j = 1}^n\|\widetilde{\bg}_{ij}(\bx_i)\|_2^3\right\}\\
&\quad\quad + \sup_{\bx_i\in B(\rho_n^{1/2}\bx_{0i}, \delta_n)}\max_{j\in [n]}e^{\widehat{\blambda}_i(\bW\bx_i)\transpose\widetilde{\bg}_{ij}(\bW\bx_i)}\left\|\widetilde\bg_{ij}(\bW\bx_i)\right\|_2\left\|\widehat{\blambda}_i(\bW\bx_i)\right\|_2\\
&\quad\quad\quad\times \sup_{\bx_i\in\Theta}\left(\frac{1}{n}\sum_{j = 1}^n\left\|\frac{\partial\widetilde{\bg}_{ij}}{\partial x_{ik}}(\bx_i)\right\|_2 ^2\right)^{1/2}\left(\frac{1}{n}\sum_{j = 1}^n \left\|\widetilde\bg_{ij}(\bx_i)\right\|_2^2\right)^{1/2}\\
&\quad\quad + \sup_{\bx_i\in B(\rho_n^{1/2}\bx_{0i}, \delta_n)}\max_{j\in [n]}2e^{\widehat{\blambda}_i(\bW\bx_i)\transpose\widetilde{\bg}_{ij}(\bW\bx_i)}\sup_{\bx_i\in\Theta}\left(\frac{1}{n}\sum_{j = 1}^n \left\|\frac{\partial\widetilde{\bg}_{ij}}{\partial x_{ik}}(\bx_i)\right\|_2^2\right)^{1/2}\left(\frac{1}{n}\sum_{j = 1}^n \left\|\widetilde\bg_{ij}(\bx_i)\right\|_2^2\right)^{1/2}.
\end{align*}
By the proof of first assertion, we have $\sup_{\bx_i\in B(\rho_n^{1/2}\bx_{0i}, \delta_n)}\max_{j\in[n]}e^{\widehat{\blambda}_i(\bW\bx_i)\transpose\widetilde{\bg}_{ij}(\bW\bx_i)}\lesssim1$ w.h.p.. By the proof of fourth assertion, we have $\sup_{\bx_i\in B(\rho_n^{1/2}\bx_{0i}, \delta_n)}\left\|\frac{\partial\widehat{\blambda}_i}{\partial x_{ik}}(\bW\bx_i)\right\|_2\lesssim\rho_n^{-1/2}$ w.h.p..
By Lemma \ref{lemma:convergence_lambda_ETEL}, $\sup_{\bx_i\in B(\rho_n^{1/2}\bx_{0i}, \delta_n)}\left\|\widehat{\blambda}_i(\bW\bx_i)\right\|_2\lesssim\rho_n^{-1/2}\delta_n$ w.h.p..
By Result \ref{result:uniform_concentration_g}, we have 
\[
\sup_{\bx_i\in B(\rho_n^{1/2}\bx_{0i}, \delta_n)}\max_{j\in[n]}\left\|\widetilde\bg_{ij}(\bW\bx_i)\right\|_2\lesssim(\rho_n\log{n})^{1/2}\quad\mbox{w.h.p.}.
\]
By Lemma \ref{lemma:Sample_moments_g}, we have
\begin{align*}
&\sup_{\bx_i\in\Theta}\frac{1}{n}\sum_{j = 1}^n\|\widetilde{\bg}_{ij}(\bx_i)\|_2^3\lesssim\rho_n^{3/2}\quad\mbox{w.p.a.1.},\\
&\sup_{\bx_i\in\Theta}\left\{
\frac{1}{n}\sum_{j = 1}^n\|\widetilde{\bg}_{ij}(\bx_i)\|_2^2
\right\}^{1/2}\lesssim\rho_n^{1/2},\quad\mbox{w.h.p.},\\
&\sup_{\bx_i\in\Theta}\left\{
\frac{1}{n}\sum_{j = 1}^n\left\|\frac{\partial\widetilde{\bg}_{ij}}{\partial\bx_i\transpose}(\bx_i)\right\|_2^2\right\}^{1/2}\lesssim \rho_n^{1/2}\quad\mbox{w.h.p.}.
\end{align*}
Therefore, we obtain
\begin{align*}
\sup_{\bx_i\in B(\rho_n^{1/2}\bx_{0i}, \delta_n)}\left\|\frac{\partial\widetilde{\bOmega}_{in}}{\partial x_{ik}}(\bW\bx_i)\right\|_2
&\lesssim \rho_n^{-1/2}\rho_n^{3/2} + \rho_n^{-1/2}\delta_n\rho_n^{1/2}(\log{n})^{1/2}\rho_n^{1/2}\rho_n^{1/2} + \rho_n^{1/2}\rho_n^{1/2}\\
&\lesssim \rho_n\quad\mbox{w.p.a.1.}.
\end{align*}

\vspace*{2ex}
\noindent
{$\blacksquare$ \bf Proof of the sixth assertion.} By definition,
\begin{align*}
\frac{\partial}{\partial x_{ik}}\left\{\frac{\partial\widehat{\blambda}_i}{\partial\bx_i}(\bx_i)\right\}
& = \widetilde{\bOmega}_{in}(\bx_i)^{-1}\frac{\partial\widetilde{\bOmega}_{in}}{\partial x_{ik}}(\bx_i)\widetilde{\bOmega}_{in}(\bx_i)^{-1}\widetilde{\bG}_{in}(\bx_i) - \widetilde{\bOmega}_{in}(\bx_i)^{-1}\frac{\partial\widetilde{\bG}_{in}}{\partial x_{ik}}(\bx_i).
\end{align*}
By Assumption \ref{assumption:regularity_condition} and the previous assertions, we have
\begin{align*}
&\sup_{\bx_i\in B(\rho_n^{1/2}\bx_{0i}, \delta_n)}\|\widetilde{\bG}_{in}(\bW\bx_i)\|_2\lesssim \rho_n^{1/2}\quad\mbox{w.h.p.},\\
&\sup_{\bx_i\in B(\rho_n^{1/2}\bx_{0i}, \delta_n)}\|\widetilde{\bOmega}_{in}(\bW\bx_i)^{-1}\|_2\lesssim \rho_n^{-1}\quad\mbox{w.h.p.},\\
&\sup_{\bx_i\in B(\rho_n^{1/2}\bx_{0i}, \delta_n)}\left\|\frac{\partial\widetilde{\bG}_{in}}{\partial x_{ik}}(\bW\bx_i)\right\|_2\lesssim \rho_n ^{1/2}\quad\mbox{w.h.p.},\\
&\sup_{\bx_i\in B(\rho_n^{1/2}\bx_{0i}, \delta_n)}\left\|\frac{\partial\widetilde{\bOmega}_{in}}{\partial x_{ik}}(\bW\bx_i)\right\|_2\lesssim \rho_n\quad\mbox{w.p.a.1}.
\end{align*}
It follows directly that
\begin{align*}
\sup_{\bx_i\in B(\rho_n^{1/2}\bx_{0i}, \delta_n)}\left\|\frac{\partial}{\partial x_{ik}}\left\{\frac{\partial\widehat{\blambda}_i}{\partial\bx_i}(\bW\bx_i)\right\}\right\|_2
&\lesssim \rho_n^{-1/2}\quad\mbox{w.p.a.1}. 
\end{align*}

\vspace*{2ex}
\noindent
{$\blacksquare$ \bf Proof of the seventh and eighth assertion.} By definition, we have
\begin{align*}
\frac{\partial}{\partial\bx_i}\widehat{\blambda}_i(\bx_i)\transpose\widetilde{\bg}_{ij}(\bx_i)
& = \left\{\frac{\partial\widehat{\blambda}_i\transpose}{\partial\bx_i}(\bx_i)\right\}\widetilde{\bg}_{ij}(\bx_i) + \left\{\frac{\partial\widetilde{\bg}_{ij}\transpose}{\partial\bx_i}(\bx_i)\right\}\widehat{\blambda}_i(\bx_i),\\
\frac{\partial}{\partial x_{ik}}\left\{\frac{\partial}{\partial\bx_i}\widehat{\blambda}_i(\bx_i)\transpose\widetilde{\bg}_{ij}(\bx_i)\right\}
& = \frac{\partial}{\partial x_{ik}}\left\{\frac{\partial\widehat{\blambda}_i\transpose}{\partial\bx_i}(\bx_i)\right\}\widetilde{\bg}_{ij}(\bx_i) + \left\{\frac{\partial\widehat{\blambda}_i\transpose}{\partial\bx_i}(\bx_i)\right\}\frac{\partial\widetilde{\bg}_{ij}}{\partial x_{ik}}(\bx_i)\\
&\quad + \frac{\partial}{\partial x_{ik}}\left\{\frac{\partial\widetilde{\bg}_{ij}\transpose}{\partial\bx_i}(\bx_i)\right\}\widehat{\blambda}_i(\bx_i) + \left\{\frac{\partial\widetilde{\bg}_{ij}\transpose}{\partial\bx_i}(\bx_i)\right\}\frac{\partial\widehat{\blambda}_i}{\partial x_{ik}}(\bx_i)
\end{align*}
% By Lemma \ref{lemma:Sample_moments_g}, we know that
% \begin{align*}
% &\sup_{\bx_i\in \Theta}\frac{1}{n}\sum_{j = 1}^n\left\|\widetilde{\bg}_{ij}(\bx_i)\right\|_2
% \lesssim \rho_n^{1/2}\quad\mbox{w.h.p.},\quad
% \sup_{\bx_i\in \Theta}\frac{1}{n}\sum_{j = 1}^n\left\|\frac{\partial\widetilde{\bg}_{ij}}{\partial\bx_i\transpose}(\bx_i)\right\|_2
% \lesssim \rho_n^{1/2}\quad\mbox{w.h.p.},\\
% &\sup_{\bx_i\in\Theta}\frac{1}{n}\sum_{j = 1}^n\left\|\frac{\partial}{\partial x_{ik}}\left\{\frac{\partial\widetilde{\bg}_{ij}}{\partial\bx_i\transpose}(\bx_i)\right\}\right\|_2\lesssim \rho_n^2\quad\mbox{w.h.p.}.
% \end{align*}
It then follows directly from the third and sixth assertion, together with Lemma \ref{lemma:Sample_moments_g}, that 
\begin{align*}
&\sup_{\bx_i\in B(\rho_n^{1/2}\bx_{0i}, \delta_n)}\frac{1}{n}\sum_{j = 1}^n\left\|\frac{\partial}{\partial\bx_i}(\widehat{\blambda}_i\transpose\widetilde{\bg}_{ij})(\bW\bx_i)\right\|_2^2\lesssim 1\quad\mbox{w.h.p.},\\
&\sup_{\bx_i\in B(\rho_n^{1/2}\bx_{0i}, \delta_n)}\frac{1}{n}\sum_{j = 1}^n\left\|\frac{\partial}{\partial x_{ik}}\left\{\frac{\partial}{\partial\bx_i}(\widehat{\blambda}_i\transpose\widetilde{\bg}_{ij})\right\}(\bW\bx_i)\right\|_2\lesssim 1\quad\mbox{w.p.a.1},\quad k\in [d].
\end{align*}
The proof is completed by applying a union bound over $k\in [d]$.
\end{proof}

\subsection{Proof of Proposition \ref{prop:Criterion_satisfies_assumption} (c)}
\label{sub:proof_of_prop_ETEL_criterion}

\begin{lemma}\label{lemma:ETEL_Hessian_term_III}
Suppose Assumptions \ref{assumption:signal_plus_noise}, \ref{assumption:regularity_condition}, and \ref{assumption:weight_functions} hold. Further assume that Assumption \ref{assumption:signal_plus_noise} (vi) is strengthened to Assumption \ref{assumption:signal_plus_noise} (vi) (b). Let $\widehat{\blambda}_i(\bx_i)$ be the Lagrange multiplier given by \eqref{eqn:ETEL_probabilities_dual}
% , $\eps_n = \sqrt{(\log n)/n}$, 
and $\delta_n = M_n\sqrt{(\log n)^{2\xi + 1}/(n\rho_n)}$, where $M_n = \log\log n$. Then
\begin{align*}
% &\sup_{\bx_i\in B(\rho_n^{1/2}\bx_{0i}, \eps_n)}
% \left\|\frac{1}{n}\sum_{j = 1}^nnp_{ij}(\bW\bx_i)\left\{\frac{\partial}{\partial\bx_i}\widehat{\blambda}_i(\bx_i)\transpose\widetilde{\bg}_{ij}(\bW\bx_i)\right\}\right\|_2\lesssim \sqrt{\frac{(\log n)^{2\xi + 1}}{n}}
%  \quad\mbox{w.h.p.},\\
&\sup_{\bx_i\in B(\rho_n^{1/2}\bx_{0i}, \delta_n)}
\left\|\frac{1}{n}\sum_{j = 1}^nnp_{ij}(\bW\bx_i)\left\{\frac{\partial(\widehat{\blambda}_i \transpose\widetilde{\bg}_{ij})}{\partial\bx_i} (\bW\bx_i)\right\}\right\|_2\lesssim M_n\sqrt{\frac{(\log n)^{2\xi + 2}}{n\rho_n}}
 \quad\mbox{w.h.p.},
\end{align*}
where $p_{ij}(\bx_i)$, $j\in [n]$ are the empirical probabilities given by \eqref{eqn:ETEL_probabilities}. 
\end{lemma}

\begin{proof}[\bf Proof of Lemma \ref{lemma:ETEL_Hessian_term_III}]
By the computation of the gradient of $\widehat{\blambda}_i(\bx_i)\transpose\widetilde{\bg}_{ij}(\bx_i)$, triangle inequality, and Cauchy-Schwarz inequality, we have
\begin{align*}
&\sup_{\bx_i\in B(\rho_n^{1/2}\bx_{0i}, \delta_n)}\left\|\frac{1}{n}\sum_{j = 1}^nnp_{ij}(\bW\bx_i)\left\{\frac{\partial(\widehat{\blambda}_i\transpose\widetilde{\bg}_{ij})}{\partial\bx_i}(\bW\bx_i)\right\}\right\|_2\\
&\quad\leq \sup_{\bx_i\in B(\rho_n^{1/2}\bx_{0i}, \delta_n)}\max_{j\in [n]}|1 - np_{ij}(\bW\bx_i)|\times \sup_{\bx_i\in B(\rho_n^{1/2}\bx_{0i}, \delta_n)}\left\{\frac{1}{n}\sum_{j = 1}^n\left\|\frac{\partial(\widehat{\blambda}_i\transpose\widetilde{\bg}_{ij})}{\partial\bx_i}(\bW\bx_i)\right\|_2^2\right\}^{1/2}\\
&\quad\quad + \sup_{\bx_i\in B(\rho_n^{1/2}\bx_{0i}, \delta_n)}\left\|\widehat{\blambda}_i(\bW\bx_i)\right\|_2\times \sup_{\bx_i\in \Theta}\frac{1}{n}\sum_{j = 1}^n\left\|\frac{\partial\widetilde{\bg}_{ij}}{\partial\bx_i\transpose}(\bx_i)\right\|_2\\
&\quad\quad + \sup_{\bx_i\in B(\rho_n^{1/2}\bx_{0i}, \delta_n)}\left\|\frac{\partial}{\partial\bx_i\transpose}\widehat{\blambda}_i(\bW\bx_i)\right\|_2\times \sup_{\bx_i\in B(\rho_n^{1/2}\bx_{0i}, \delta_n)}\left\|\frac{1}{n}\sum_{j = 1}^n\widetilde{\bg}_{ij}(\bW\bx_i)\right\|_2.
\end{align*}
By Lemma \ref{lemma:Local_ULLN} and the proof of Lemma \ref{lemma:lambda_gradient_convergence_ETEL} (also see the proof of Lemma \ref{lemma:convergence_lambda_ETEL}), 
\[
\sup_{\bx_i\in B(\rho_n^{1/2}\bx_{0i}, \delta_n)}\left\|\frac{1}{n}\sum_{j = 1}^n\widetilde{\bg}_{ij}(\bW\bx_i)\right\|_2\lesssim \rho_n^{1/2}\delta_n\quad\mbox{w.h.p.}.
\]
Then by Lemma \ref{lemma:Sample_moments_g}, Lemma \ref{lemma:EP_convergence_ETEL}, and Lemma \ref{lemma:lambda_gradient_convergence_ETEL}, we conclude that
\begin{align*}
\sup_{\bx_i\in B(\rho_n^{1/2}\bx_{0i}, \delta_n)}\left\|\frac{1}{n}\sum_{j = 1}^nnp_{ij}(\bW\bx_i)\left\{\frac{\partial(\widehat{\blambda}_i\transpose\widetilde{\bg}_{ij})}{\partial\bx_i}(\bW\bx_i)\right\}\right\|_2\lesssim M_n\sqrt{\frac{(\log n)^{2\xi + 2}}{n\rho_n}}\quad\mbox{w.h.p.}.
\end{align*}
\end{proof}

\begin{lemma}\label{lemma:ETEL_Hessian_term_I}
Suppose Assumptions \ref{assumption:signal_plus_noise}, \ref{assumption:regularity_condition}, and \ref{assumption:weight_functions} hold. Further assume that Assumption \ref{assumption:signal_plus_noise} (vi) is strengthened to Assumption \ref{assumption:signal_plus_noise} (vi) (b). Let $\widehat{\blambda}_i(\bx_i)$ be the Lagrange multiplier given by \eqref{eqn:ETEL_probabilities_dual}
% , $\eps_n = \sqrt{(\log n)/n}$, 
and $\delta_n = M_n\sqrt{(\log n)^{2\xi+1}/(n\rho_n)}$, where $M_n = \log\log n$. Then
\begin{align*}
&\bW\transpose\frac{1}{n}\sum_{j = 1}^nnp_{ij}(\bx_i)\left\{\frac{\partial(\widehat{\blambda}_i\transpose\widetilde{\bg}_{ij})}{\partial\bx_i}(\bx_i)\right\}\left\{\frac{\partial(\widehat{\blambda}_i\transpose\widetilde{\bg}_{ij})}{\partial\bx_i}(\bx_i)\right\}\transpose\bW\\
&\quad = \bG_{in}(\rho_n^{1/2}\bx_{0i})\transpose\bOmega_{in}(\rho_n^{1/2}\bx_{0i})^{-1}\bG_{in}(\rho_n^{1/2}\bx_{0i}) + \bR_{2in}^{(\mathrm{ET})}(\bx_i),
\end{align*}
where 
\[
\sup_{\bx_i\in B(\rho_n^{1/2}\bx_{0i}, \delta_n)}\|\bR_{2in}^{(\mathrm{ET})}(\bW\bx_i)\|_2\lesssim M_n\sqrt{\frac{(\log n)^{2\xi + 2}}{n\rho_n}}\quad\mbox{w.h.p.}.
\]
\end{lemma}

\begin{proof}[\bf Proof of Lemma \ref{lemma:ETEL_Hessian_term_I}]
Denote $\bG_{0in}\overset{\Delta}{=}\bG_{in}(\rho_n^{1/2}\bx_{0i})$ and $\bOmega_{0in}\overset{\Delta}{=}\bOmega_{in}(\rho_n^{1/2}\bx_{0i})$. By Lemma \ref{lemma:covariance_bound_ETEL}, the proof of the second assertion in Lemma \ref{lemma:lambda_gradient_convergence_ETEL}, and the third assertion of Lemma \ref{lemma:lambda_gradient_convergence_ETEL},
\begin{align*}
&\frac{1}{n}\sum_{j = 1}^n\bW\transpose\widetilde{\bg}_{ij}(\bW\bx_i)\widetilde{\bg}_{ij}(\bW\bx_i)\transpose\bW = \bOmega_{0in} + \bR_{in}^{(\bOmega)}(\bW\bx_i),\\
&\bW\transpose\frac{\partial\widehat{\blambda}_i}{\partial\bx_i\transpose}(\bW\bx_i)\bW = -\bOmega_{0in}^{-1}\bG_{0in} + \bR_{in}^{(\blambda)}(\bW\bx_i),
\end{align*}
where $\|\bR_{in}^{(\bOmega)}(\bW\bx_i)\|_2 = O\{\rho_n^{1/2}M_n\sqrt{(\log n)^{2\xi + 1}/n}\}$ and $\|\bR_{in}^{(\blambda)}(\bW\bx_i)\|_2 = O\{\rho_n^{-1}M_n\sqrt{(\log n)^{2\xi + 2}/n}\}$ w.h.p. uniformly in $\bx_i\in B(\rho_n^{1/2}\bx_{0i}, \delta_n)$. 
It follows that
\begin{align*}
&\bW\transpose\left\{\frac{\partial\widehat{\blambda}_i\transpose}{\partial\bx_i}(\bW\bx_i)\right\}\left\{\frac{1}{n}\sum_{j = 1}^n\widetilde{\bg}_{ij}(\bW\bx_i)\widetilde{\bg}_{ij}(\bW\bx_i)\transpose\right\}\left\{\frac{\partial\widehat{\blambda}_i}{\partial\bx_i\transpose}(\bW\bx_i)\right\}\bW\\
&\quad = \{-\bOmega_{0in}^{-1}\bG_{0in} + \bR_{in}^{(\blambda)}(\bW\bx_i)\}\transpose
\{\bOmega_{0in} + \bR_{in}^{(\bOmega)}(\bW\bx_i)\}
\{-\bOmega_{0in}^{-1}\bG_{0in} + \bR_{in}^{(\blambda)}(\bW\bx_i)\}.
\end{align*}
Denote 
\[
\bR_{3in}^{(\mathrm{ET})}(\bW\bx_i) = \bW\transpose\left\{\frac{\partial\widehat{\blambda}_i\transpose}{\partial\bx_i}(\bW\bx_i)\right\}\left\{\frac{1}{n}\sum_{j = 1}^n\widetilde{\bg}_{ij}(\bW\bx_i)\widetilde{\bg}_{ij}(\bW\bx_i)\transpose\right\}\left\{\frac{\partial\widehat{\blambda}_i}{\partial\bx_i\transpose}(\bW\bx_i)\right\}\bW
- \bG_{0in}\transpose\bOmega_{0in}^{-1}\bG_{0in}.
\]
Then
\begin{align*}
\sup_{\bx_i\in B(\rho_n^{1/2}\bx_{0i},\delta_n)}\left\|\bR_{3in}^{(\mathrm{ET})}(\bW\bx_i)
\right\|_2
&\leq \sup_{\bx_i\in B(\rho_n^{1/2}\bx_{0i},\delta_n)}\|\bR_{in}^{(\bOmega)}(\bW\bx_i)\|_2\|-\bOmega_{0in}^{-1}\bG_{0in} + \bR_{in}^{(\blambda)}(\bW\bx_i)\|_2^2\\
&\quad + 2\sup_{\bx_i\in B(\rho_n^{1/2}\bx_{0i},\delta_n)}\|\bR_{in}^{(\blambda)}(\bW\bx_i)\|_2\|\bOmega_{0in}\|_2\|\bOmega_{0in}^{-1}\bG_{0in}\|_2\\
&\quad +  \sup_{\bx_i\in B(\rho_n^{1/2}\bx_{0i},\delta_n)}\|\bOmega_{0in}\|_2\|\bR_{in}^{(\blambda)}(\bW\bx_i)\|_2^2\lesssim M_n\sqrt{\frac{(\log n)^{2\xi + 2}}{n\rho_n}}\quad\mbox{w.h.p.}.
\end{align*}
Now we suppress the argument and compute
\begin{align*}
&\frac{1}{n}\sum_{j = 1}^nnp_{ij}\left\{\frac{\partial(\widehat{\blambda}_i\transpose\widetilde{\bg}_{ij})}{\partial\bx_i}\right\}\left\{\frac{\partial(\widehat{\blambda}_i\transpose\widetilde{\bg}_{ij})}{\partial\bx_i}\right\}\transpose\\
&\quad = \frac{1}{n}\sum_{j = 1}^n\left\{\frac{\partial(\widehat{\blambda}_i\transpose\widetilde{\bg}_{ij})}{\partial\bx_i}\right\}\left\{\frac{\partial(\widehat{\blambda}_i\transpose\widetilde{\bg}_{ij})}{\partial\bx_i}\right\}\transpose - \frac{1}{n}\sum_{j = 1}^n(1 - np_{ij})\left\{\frac{\partial(\widehat{\blambda}_i\transpose\widetilde{\bg}_{ij})}{\partial\bx_i}\right\}\left\{\frac{\partial(\widehat{\blambda}_i\transpose\widetilde{\bg}_{ij})}{\partial\bx_i}\right\}\transpose\\
&\quad = \left(\frac{\partial\widehat{\blambda}_i\transpose}{\partial\bx_i}\right)\left(\frac{1}{n}\sum_{j = 1}^n\widetilde{\bg}_{ij}\widetilde{\bg}_{ij}\transpose\right)\left(\frac{\partial\widehat{\blambda}_i}{\partial\bx_i\transpose}\right) + \left(\frac{\partial\widehat{\blambda}_i\transpose}{\partial\bx_i}\right)\left(\frac{1}{n}\sum_{j = 1}^n\widetilde{\bg}_{ij}\widehat{\blambda}_i\transpose\frac{\partial\widetilde{\bg}_{ij}}{\partial\bx_i\transpose}\right)\\
&\quad\quad + \left(\frac{1}{n}\sum_{j = 1}^n\frac{\partial\widetilde{\bg}_{ij}\transpose}{\partial\bx_i}\widehat{\blambda}_i\widetilde{\bg}_{ij}\transpose\right)\left(\frac{\partial\widehat{\blambda}_i}{\partial\bx_i\transpose}\right)
 + \frac{1}{n}\sum_{j = 1}^n\frac{\partial\widetilde{\bg}_{ij}\transpose}{\partial\bx_i}\widehat{\blambda}_i\widehat{\blambda}_i\transpose\frac{\partial\widetilde{\bg}_{ij}}{\partial\bx_i\transpose}
\\
&\quad\quad + \frac{1}{n}\sum_{j = 1}^n(1 - np_{ij})\left\{\frac{\partial(\widehat{\blambda}_i\transpose\widetilde{\bg}_{ij})}{\partial\bx_i}\right\}\left\{\frac{\partial(\widehat{\blambda}_i\transpose\widetilde{\bg}_{ij})}{\partial\bx_i}\right\}\transpose
\end{align*}
It follows that
\begin{align*}
&\left\|\frac{1}{n}\sum_{j = 1}^nnp_{ij}\left\{\frac{\partial(\widehat{\blambda}_i\transpose\widetilde{\bg}_{ij})}{\partial\bx_i}\right\}\left\{\frac{\partial(\widehat{\blambda}_i\transpose\widetilde{\bg}_{ij})}{\partial\bx_i}\right\}\transpose - \frac{1}{n}\left(\frac{\partial\widehat{\blambda}_i}{\partial\bx_i\transpose}\right)\left(\frac{1}{n}\sum_{j = 1}^n\widetilde{\bg}_{ij}\widetilde{\bg}_{ij}\transpose\right)\left(\frac{\partial\widehat{\blambda}_i\transpose}{\partial\bx_i}\right)\right\|_2\\
&\quad\leq \max_{j\in [n]}|1 - np_{ij}|\frac{1}{n}\sum_{j = 1}^n\left\|\frac{\partial(\widehat{\blambda}_i\transpose\widetilde{\bg}_{ij})}{\partial\bx_i}\right\|_2^2
 + 2\left\|\frac{\partial\widehat{\blambda}_i}{\partial\bx_i\transpose}\right\|_2\|\widehat{\blambda}_i\|_2\left(
\frac{1}{n}\sum_{j = 1}^n\|\widetilde{\bg}_{ij}\|_2^2
 \right)^{1/2}\left(
\frac{1}{n}\sum_{j = 1}^n\left\|\frac{\partial\widetilde{\bg}_{ij}}{\partial\bx_i\transpose}\right\|_2^2
 \right)^{1/2}\\
&\quad\quad + \|\widehat{\blambda}_i\|_2^2\frac{1}{n}\sum_{j = 1}^n\left\|\frac{\partial\widetilde{\bg}_{ij}}{\partial\bx_i\transpose}\right\|_2^2.
\end{align*}
By Lemma \ref{lemma:Sample_moments_g}, Lemma \ref{lemma:EP_convergence_ETEL}, and Lemma \ref{lemma:lambda_gradient_convergence_ETEL}, we obtain that
\begin{align*}
&\frac{1}{n}\sum_{j = 1}^nnp_{ij}(\bW\bx_i)\left\{\frac{\partial(\widehat{\blambda}_i\transpose\widetilde{\bg}_{ij})}{\partial\bx_i}(\bW\bx_i)\right\}\left\{\frac{\partial(\widehat{\blambda}_i\transpose\widetilde{\bg}_{ij})}{\partial\bx_i}(\bW\bx_i)\right\}\transpose\\
&\quad = \frac{1}{n}\left\{\frac{\partial\widehat{\blambda}_i}{\partial\bx_i\transpose}(\bW\bx_i)\right\}\left\{\frac{1}{n}\sum_{j = 1}^n\widetilde{\bg}_{ij}(\bW\bx_i)\widetilde{\bg}_{ij}(\bW\bx_i)\transpose\right\}\left\{\frac{\partial\widehat{\blambda}_i\transpose}{\partial\bx_i}(\bW\bx_i)\right\} + \bR_{4in}^{(\mathrm{ET})}(\bW\bx_i),
\end{align*}
where $\sup_{\bx_i\in B(\rho_n^{1/2}\bx_{0i}, \delta_n)}\|\bR_{4in}^{(\mathrm{ET})}(\bW\bx_i)\|_2\lesssim M_n\sqrt{{(\log n)^{2\xi + 2}}/{(n\rho_n)}}$ w.h.p.. The proof is them completed by noting that $\|\bR_{2in}^{(ET)}(\bW\bx_i)\|_2\leq \|\bR_{3in}^{(\mathrm{ET})}(\bW\bx_i)\|_2 + \|\bR_{4in}^{(\mathrm{ET})}(\bW\bx_i)\|_2$. 
\end{proof}

\begin{proof}[\bf Proof of Proposition \ref{prop:Criterion_satisfies_assumption} (c)]
Let $\eps_n = (\log n)^{(\xi - 1)/4}/\sqrt{n}$ and $\delta_n = M_n\sqrt{(\log n)^{2\xi + 1} / (n\rho_n)}$, where $M_n = \log\log n$. 
By definition of the ETEL criterion function \eqref{eqn:ETEL}, we have
\begin{align*}
\frac{\partial\ell_{in}}{\partial\bx_i}(\bx_i)
& = \sum_{j = 1}^n\{1 - np_{ij}(\bx_i)\}\frac{\partial(\widehat{\blambda}_i\transpose\widetilde{\bg}_{ij})}{\partial\bx_i}(\bx_i)
\\
% \bW\transpose
\frac{1}{n}\frac{\partial^2\ell_{in}(\bW\bx_i)}{\partial\bx_i\partial\bx_i\transpose}
% \bW
& = \frac{1}{n}\sum_{j = 1}^n\{1 - np_{ij}(\bW\bx_i)\}
% \bW\transpose
\frac{\partial^2(\widehat{\blambda}_i\transpose\widetilde{\bg}_{ij})}{\partial\bx_i\partial\bx_i\transpose}(\bW\bx_i)
% \bW
\\
&\quad - 
% \bW\transpose
\frac{1}{n}\sum_{j = 1}^nnp_{ij}(\bW\bx_i)\left\{\frac{\partial(\widehat{\blambda}_i\transpose\widetilde{\bg}_{ij})}{\partial\bx_i}(\bW\bx_i)\right\}\left\{\frac{\partial(\widehat{\blambda}_i\transpose\widetilde{\bg}_{ij})}{\partial\bx_i}(\bW\bx_i)\right\}\transpose
% \bW
\\
&\quad + 
% \bW\transpose
\left\{\frac{1}{n}\sum_{j = 1}^nnp_{ij}(\bW\bx_i) \frac{\partial(\widehat{\blambda}_i\transpose\widetilde{\bg}_{ij})}{\partial\bx_i}(\bW\bx_i) \right\}\left\{\frac{1}{n}\sum_{j = 1}^nnp_{ij}(\bW\bx_i) \frac{\partial(\widehat{\blambda}_i\transpose\widetilde{\bg}_{ij})}{\partial\bx_i}(\bW\bx_i) \right\}\transpose
% \bW
.
\end{align*}
By the eighth assertion of Lemma \ref{lemma:lambda_gradient_convergence_ETEL} and Lemma \ref{lemma:EP_convergence_ETEL}, 
\begin{align*}
\sup_{\bx_i\in B(\rho_n^{1/2}\bx_{0i}, \delta_n)}\left\|
\frac{1}{n}\sum_{j = 1}^n\{1 - np_{ij}(\bW\bx_i)\}
\frac{\partial^2(\widehat{\blambda}_i\transpose\widetilde{\bg}_{ij})}{\partial\bx_i\partial\bx_i\transpose}(\bW\bx_i)
\right\|_2\lesssim M_n\sqrt{\frac{(\log n)^{2\xi + 2}}{n\rho_n}}\quad\mbox{w.p.a.1.}.
\end{align*}
Then by Lemma \ref{lemma:ETEL_Hessian_term_I} and Lemma \ref{lemma:ETEL_Hessian_term_III}, 
\begin{align*}
\sup_{\bx_i\in B(\rho_n^{1/2}\bx_{0i}, \delta_n)}\left\|\bW\transpose \frac{1}{n}\frac{\partial^2\ell_{in}(\bW\bx_i)}{\partial\bx_i\partial\bx_i\transpose}
\bW + \bG_{0in}\transpose\bOmega_{0in}^{-1}\bG_{0in}\right\|_2\lesssim 
M_n\sqrt{\frac{(\log n)^{2\xi + 2}}{n\rho_n}} = o\left(\frac{1}{n\eps_n^2}\right)\quad\mbox{w.p.a.1.},
\end{align*}
where $\bG_{0in}\overset{\Delta}{=}\bG_{in}(\rho_n^{1/2}\bx_{0i})$ and $\bOmega_{0in} = \bOmega_{in}(\rho_n^{1/2}\bx_{0i})$. 
Since $\eps_n\leq \delta_n/3$ and the eigenvalues of $\bG_{0in}\transpose\bOmega_{0in}^{-1}\bG_{0in}$ are bounded away from $0$, this completes the proof of \eqref{eqn:Hessian_A1} and \eqref{eqn:Hessian_A2} in Assumption \ref{assumption:criterion_function} simultaneously. 

\vspace*{2ex}\noindent
It is now sufficient to establish \eqref{eqn:identifiability} in Assumption \ref{assumption:criterion_function}. The argument here is a modification of the proof of Lemma 1 in \cite{tang2021statistical}. We first claim that $p_{ij}(\widehat{\bx}_i) = 1/n$ for all $j\in [n]$. The reasoning is similar to the proof of Proposition 1 in \cite{10.1093/biomet/asaa028}. By \eqref{eqn:ETEL_probabilities_dual}, the empirical probabilities $\{p_{ij}(\widehat{\bx}_i)\}_{j = 1}^n$ can be viewed as the solution to the constrained optimization problem \eqref{eqn:ETEL_probabilities_dual} with $\bx_i$ evaluated at $\widehat{\bx}_i$. The relaxed problem
\begin{align*}
\max_{p_{i1},\ldots,p_{in}}&\sum_{j = 1}^n(-p_{ij}\log p_{ij})\\
\mbox{subject to }&\sum_{j = 1}^np_{ij} = 1,\quad p_{ij}\geq 0,\quad j\in [n]
\end{align*}
is uniquely solved at $p_{i1} = \ldots = p_{in} = 1/n$. Since $\widehat{\bx}_i$ satisfies $\sum_{j = 1}^n(1/n)\widetilde{\bg}_{ij}(\widehat{\bx}_i) = \zero_d$, we then see that the solution $p_{i1} = \ldots = p_{in} = 1/n$ also satisfies the additional constraint that $\sum_{j = 1}^np_{ij}\widetilde{\bg}_{ij}(\widehat{\bx}_i) = \zero_d$. This implies that $p_{ij}(\widehat{\bx}_i) = 1/n$, $j\in [n]$ solves \eqref{eqn:ETEL_probabilities_dual} with $\bx_i = \widehat{\bx}_i$. By definition of the ETEL criterion function \eqref{eqn:ETEL}, it follows immediately that $\ell_{in}(\widehat{\bx}_i) = -n\log{n}$. We now focus on $\ell_{in}(\bW\bx_i)$ outside $B(\rho_n^{1/2}\bx_{0i}, \delta_n)$. Let $\bp_i(\bx_i) = [p_{i1}(\bx_i),\ldots,p_{in}(\bx_i)]\transpose$ and $\bp_i^{(-n)}(\bx_i) = [p_{i1}(\bx_i),\ldots,p_{i(n - 1)}(\bx_i)]\transpose$. We consider two cases:
\begin{itemize}
	\item[$\blacksquare$] \textbf{Case I: $\|\bp_i(\bx_i)\|_\infty > 2(1 + \alpha)d(\log n)/n$. }This implies that there exists some index $k\in [n]$, such that $p_{ik}(\bx_i)\geq 2(1 + \alpha)d(\log n)/n$, and by the constraint that $\sum_{j = 1}^np_{ij}(\bx_i) = 1$, we see that $\sum_{j\neq k}p_{ij}(\bx_i)\leq 1 - 2(1 + \alpha)d(\log n)/n$. By the algorithmic-geometric inequality and the fact that $\log p_{ik}(\bx_i)\leq 0$, 
	\begin{align*}
	\sum_{j = 1}^n\log p_{ij}(\bx_i)&\leq \log\left\{\prod_{j\neq k}p_{ij}(\bx_i)\right\}
	\leq (n - 1)\log\left\{\frac{1}{n - 1}\sum_{j\neq k}p_{ij}(\bx_i)\right\}\\
	&\leq (n - 1)\log\left\{\frac{1}{n - 1}\left(1 - \frac{2(1 + \alpha)d\log n}{n}\right)\right\}.
	\end{align*}
	Therefore, by the basic inequality $\log(1 + x)\leq x$ for any $x > -1$, 
	\begin{align*}
	\sum_{j = 1}^n\log p_{ij}(\bx_i) + n\log n
	&\leq \log n + (n - 1)\log\left\{\frac{n}{n - 1}\left(1 - \frac{2(1 + \alpha)d\log n}{n}\right)\right\}\\
	&\leq \log n + (n - 1)\log\left\{1 - \frac{2(1 + \alpha)d\log n}{n - 1}\right\}\\
	&\leq \log n - 2(1 + \alpha)d\log n\leq -(1 + \alpha)d\log n.
	\end{align*}
  Namely,
  \[
  \inf_{\|\bW\transpose\bx_i - \rho_n^{1/2}\bx_{0i}\| > \delta_n, \|\bp_i(\bx_i)\|_\infty > 2(1 + \alpha)d(\log n)/n}\left\{\ell_{in}(\widehat{\bx}_i) - \ell_{in}(\bx_i)\right\}\geq (1 + \alpha)d\log n.
  \]

	\item[$\blacksquare$] \textbf{Case II: $\|\bp_i(\bx_i)\|_\infty \leq 2(1 + \alpha)d(\log n)/n$. }By Assumption \ref{assumption:regularity_condition}, 
	\[
	\inf_{\bx_i\notin B(\rho_n^{1/2}\bx_{0i}, \delta_n)}\left\|\frac{1}{n}\sum_{j = 1}^n\expect_0\{\bg_{ij}(\bx_i)\}\right\|\geq \rho_n^{1/2}\delta_0\delta_n,
	\]
	where $\delta_0 > 0$ is a constant. Then by Lemma \ref{lemma:ULLN}, 
	\[
	\inf_{\bx_i\notin B(\rho_n^{1/2}\bx_{0i}, \delta_n)}\left\|\frac{1}{n}\sum_{j = 1}^n\widetilde{\bg}_{ij}(\bW\bx_i)\right\|_2^2\geq \frac{1}{4}\rho_n\delta_0^2\delta_n^2\geq \frac{\delta_0^2 M_n^2(\log n)^{2\xi + 1}}{4n}\quad\mbox{w.h.p.}.
	\]
	By the definition of the empirical probabilities $\{p_{ij}(\bx_i)\}_{j = 1}^n$, we have $\sum_{j = 1}^np_{ij}(\bx_i)\widetilde{\bg}_{ij}(\bx_i) = \zero_d$. Namely
	\begin{align*}
	\inf_{\bW\transpose\bx_i\notin B(\rho_n^{1/2}\bx_{0i}, \delta_n)}\left\|\frac{1}{n}\sum_{j = 1}^n\widetilde{\bg}_{ij}(\bx_i)\right\|_2^2
	& = \inf_{\bW\transpose\bx_i\notin B(\rho_n^{1/2}\bx_{0i}, \delta_n)}\left\|\sum_{j = 1}^n\left\{p_{ij}(\bx_i) - \frac{1}{n}\right\}\widetilde{\bg}_{ij}(\bx_i)\right\|_2^2.
	\end{align*}
	It follows from Cauchy-Schwarz inequality and Lemma \ref{lemma:Sample_moments_g} that
	\begin{align*}
	&\inf_{\bW\transpose\bx_i\notin B(\rho_n^{1/2}\bx_{0i}, \delta_n)}\left\|\sum_{j = 1}^n\left\{p_{ij}(\bx_i) - \frac{1}{n}\right\}\widetilde{\bg}_{ij}(\bx_i)\right\|_2^2\\
	&\quad\leq \inf_{\bW\transpose\bx_i\notin B(\rho_n^{1/2}\bx_{0i}, \delta_n)}\sum_{j = 1}^n\left\{p_{ij}(\bx_i) - \frac{1}{n}\right\}^2\sup_{\bx_i\in\Theta}\sum_{j = 1}^n\left\|\widetilde{\bg}_{ij}(\bx_i)\right\|_2^2\\
	&\quad\lesssim \inf_{\bW\transpose\bx_i\notin B(\rho_n^{1/2}\bx_{0i}, \delta_n)}\sum_{j = 1}^n\left\{p_{ij}(\bx_i) - \frac{1}{n}\right\}^2n\rho_n \quad\mbox{w.h.p.},
	\end{align*}
	implying that
	\[
	\inf_{\bW\transpose\bx_i\notin B(\rho_n^{1/2}\bx_{0i}, \delta_n)}\sum_{j = 1}^n\left\{p_{ij}(\bx_i) - \frac{1}{n}\right\}^2\geq \frac{cM_n^2(\log n)^{2\xi + 1}}{n^2\rho_n}\quad\mbox{w.p.a.1}, 
	\]
  where $c > 0$ is some constant. Since $\|\bp_i(\bx_i)\|_\infty\leq 2(1 + \alpha)d(\log n)/n$, it follows that $\{(p_{in}(\bx_i) - 1/n\}^2\leq 9(1 + \alpha)^2d^2(\log n)^2/n^2$, implying that
  \[
  \inf_{\|\bW\transpose\bx_i - \rho_n^{1/2}\bx_{0i}\| > \delta_n, \|\bp_i(\bx_i)\|_\infty\leq 2(1 + \alpha)d(\log n)/n}\sum_{j = 1}^{n - 1}\left\{p_{ij}(\bx_i) - \frac{1}{n}\right\}^2\geq \frac{cM_n^2(\log n)^{2\xi + 1}}{4n^2\rho_n}\quad\mbox{w.p.a.1}.
  \]
  Denote the function $q(\bp_i^{(-n)}) = \sum_{j = 1}^{n - 1}\log p_{ij} + \log(1 - \sum_{j = 1}^{n - 1}p_{ij})$, where $\bp_i^{(-n)} = [p_{i1},\ldots,p_{i(n - 1)}]\transpose$. By the definition of the ETEL criterion function \eqref{eqn:ETEL} and the result that $\bp_{i}^{(-n)}(\widehat{\bx}_i) = \one_{n - 1}/n$, 
  \[
  \ell_{in}(\widehat{\bx}_i) - \ell_{in}(\bx_i) = q(\one_{n - 1}/n) - q(\bp_i^{(-n)}(\bx_i)).
  \]
  Denote $[\bx]_j$ the $j$th coordinate of a vector $\bx$. 
  The gradient and Hessian of the $q$ function can be obtained directly:
  \begin{align*}
  \left[\frac{\partial q}{\partial\bp_i^{(-n)}}(\bp_i^{(-n)})\right]_j & = \frac{1}{p_{ij}} - \frac{1}{1 - \sum_{j = 1}^{n - 1}p_{ij}},\\
  \frac{\partial^2 q}{\partial\bp_i^{(-n)}\partial\bp_i^{(-n)\mathrm{T}}}(\bp_i^{(-n)}) & = - \mathrm{diag}\left\{\frac{1}{p_{i1}^2},\ldots,\frac{1}{p_{i(n - 1)}^2}\right\} - \frac{1}{(1 - \sum_{j = 1}^{n - 1}p_{ij})^2}\one_{n - 1}\one_{n - 1}\transpose.
  \end{align*}
  It follows that
  \begin{align*}
  \frac{\partial q}{\partial\bp_i^{(-n)}}\left(\frac{1}{n}\one_{n - 1}\right) & = \zero_{n - 1},\quad
  \lambda_{\min}\left\{-\frac{\partial^2 q}{\partial\bp_i^{(-n)}\partial\bp_i^{(-n)\mathrm{T}}}(\bp_i^{(-n)})\right\}\geq \frac{1}{\|\bp_i\|_\infty}.
  \end{align*}
  By Taylor's theorem, there exists some $\theta\in [0, 1]$, such that $\bar{\bp}_i^{(-n)}(\bx_i) = \theta\one_{n - 1}/n + (1 - \theta)\bp_i^{(-n)}(\bx_i)$, and 
  \begin{align*}
  \ell_{in}(\widehat{\bx}_i) - \ell_{in}(\bx_i)& = q\left(\frac{1}{n}\one_{n - 1}\right) - q(\bp_i^{(-n)}(\bx_i))\\
  & = -\frac{\partial q}{\partial\bp_i^{(-n)}}\left(\frac{1}{n}\one_{n - 1}\right)\left\{\bp_i^{(-n)}(\bx_i) - \one_{n - 1}/n\right\}\\
  &\quad - \frac{1}{2}\left\{\bp_i^{(-n)}(\bx_i) - \one_{n - 1}/n\right\}\transpose
  \frac{\partial^2 q}{\partial\bp_i^{(-n)}\partial\bp_i^{(-n)\mathrm{T}}}(\bar{\bp}_i^{(-n)})
   \left\{\bp_i^{(-n)}(\bx_i) - \one_{n - 1}/n\right\}\\
  &\geq \frac{\left\|\bp_i^{(-n)}(\bx_i) - \one_{n - 1}/n\right\|_2^2}{2\{\theta/n + (1 - \theta)\|{\bp}_i(\bx_i)\|_\infty\}^2}\geq \frac{1}{2}\left(\frac{n}{2(1 + \alpha)d\log n}\right)^2\sum_{j = 1}^{n - 1}\left\{p_{ij}(\bx_i) - \frac{1}{n}\right\}^2.
  \end{align*}
  Therefore, 
  \begin{align*}
  &\inf_{\|\bW\transpose\bx_i - \rho_n^{1/2}\bx_{0i}\| > \delta_n, \|\bp_i(\bx_i)\|_\infty\leq 2(1 + \alpha)d(\log n)/n}\left\{\ell_{in}(\widehat{\bx}_i) - \ell_{in}(\bx_i)\right\}\\
  &\quad \geq \frac{1}{2}\left(\frac{n}{2(1 + \alpha)d\log n}\right)^2\left[\inf_{\|\bW\transpose\bx_i - \rho_n^{1/2}\bx_{0i}\| > \delta_n, \|\bp_i(\bx_i)\|_\infty\leq 2(1 + \alpha)d(\log n)/n}\sum_{j = 1}^{n - 1}\left\{p_{ij}(\bx_i) - \frac{1}{n}\right\}^2\right]\\
  &\quad \geq \frac{1}{32}\left\{\frac{n^2}{(1 + \alpha)^2d^2(\log n)^2}\right\}\left\{\frac{cM_n^2(\log n)^{2\xi + 1}}{n^2\rho_n}\right\}\geq(1 + \alpha)d\log n\quad\mbox{w.p.a.1}.
  \end{align*}
  The proof is thus completed.

\end{itemize}
\end{proof}
% subsection technical_results_for_the_exponentially_tilted_empirical_likelihood (end)

\section{Computational Details}
\label{sec:computational_details}

\subsection{Detailed Metropolis-Hastings Algorithm} % (fold)
\label{sub:MH_algorithm}

This subsection provides the detailed Metropolis-Hastings algorithm for computing the generalized posterior distribution defined in \eqref{eqn:generalized_posterior} in Section \ref{sub:generalized_bayesian_estimation_with_moment_conditions} of the manuscript. The algorithm applies to a generic criterion function $\ell_{in}(\cdot)$ for $\bx_i$, including the M-criterion \eqref{eqn:M_estimation}, the GMM criterion \eqref{eqn:GMM}, and the ETEL criterion \eqref{eqn:ETEL}. See Algorithm \ref{alg:MCMC} below for details. 

\begin{algorithm}[htbp] %Beginning of the algo
  \renewcommand{\algorithmicrequire}{\textbf{Input:}}
  \renewcommand{\algorithmicensure}{\textbf{Output:} }
  \caption{Metropolis-Hastings algorithm for the generalized Bayesian estimation} %Title of the algo
  \label{alg:MCMC} %Assign a label to the algo
  \begin{algorithmic}[1] %means that everyline is accompanied with a number
    % \Require\\
    \State{\textbf{Input:}
      Data matrix $\bA = [A_{ij}]_{n\times n}$, rank $d$, proposal distribution $q(\cdot; \cdot)$, Number of burn-in iterations $B$, Number of post-burn-in iterations $T$}
    % \State{\textbf{Set:}
    %   Tuning parameters $a_0,b_0 > 0$, $\eps\in (0, 1/2]$.}
    \State{Compute the truncated eigen-decomposition $(\bU_\bA, \bS_\bA)$ of the $\bA$: $\bA\bU_\bA = \bU_\bA\bS_\bA$, where $\bU_\bA\in\mathbb{O}(n, d)$, $\bS_\bA = \mathrm{diag}(\widehat{\lambda}_1,\ldots,\widehat{\lambda}_d)$, and  $|\widehat\lambda_1|\geq|\widehat\lambda_2|\geq\ldots\geq|\widehat\lambda_n|$.
    % \[
    % \bA = \sum_{i = 1}^n\widehat{\lambda}_i\widehat{\bu}_i\widehat{\bu}_i\transpose,
    % \]
    % where $|\widehat\lambda_1|\geq|\widehat\lambda_2|\geq\ldots\geq|\widehat\lambda_n|$, and $\widehat{\bu}_i\transpose\widehat{\bu}_j = \mathbbm{1}(i = j)$ for all $i,j\in[n]$. 
    Compute the spectral embedding
        % \[
        $\widetilde\bX  = \bU_\bA\bS_\bA^{1/2}$ 
        % [\widehat{\mathbf{u}}_1, \ldots, \widehat{\mathbf{u}}_d] \cdot \mathrm{diag}(|\widehat{\lambda}_1|^{1/2}, \ldots, |\widehat{\lambda}_d|^{1/2})$,
        % \]
        and write $\widetilde\bX = [\widetilde\bx_1,\ldots,\widetilde\bx_n]\transpose\in\mathbb{R}^{n\times d}$. } 
    \State{For $i = 1,2,\ldots,n$

    Initialize ${\bx}_i^{(1)} = \widetilde{\bx}_i$.

    For $t = 1,2,\ldots,B + T$

    \quad\quad Generate $\bx_i^*\sim q(\bx_i; \bx_i^{(t - 1)})$. 

    \quad\quad Generate $u\sim \mathrm{Unif}(0, 1)$ independent of $\bx_i^*$. 

    \quad\quad Compute the logarithmic Metropolis-Hastings ratio
    \[
    \log\alpha = \log\frac{\pi(\bx_i^*)}{\pi(\bx_i^{(t - 1)})} + \ell_{in}(\bx_i^*) - \ell_{in}(\bx_i^{(t - 1)}) - \log \frac{q(\bx_i^*\mid\bx_i^{(t - 1)})}{q(\bx_i^{(t - 1)}\mid \bx_i^*)}.
    \]

    \quad\quad If $\log u\leq \log\alpha$ then

    \quad\quad\quad $\bx_i^{(t)} = \bx_i^*$

    \quad\quad else

    \quad\quad\quad $\bx_i^{(t)} = \bx_i^{(t - 1)}$

    End For
    
    \noindent
    End For}
    \State{\textbf{Output: }The MCMC samples $(\bX^{(t)})_{t = B + 1}^{B + T}$, where $\bX^{(t)} = [\bx_1^{(t)},\ldots,\bx_n^{(t)}]\transpose$. }
  \end{algorithmic}
\end{algorithm}

\subsection{MCMC Convergence diagnostics}
\label{sub:mcmc_diagnostics}

\begin{figure}[htbp]
  \centerline{\includegraphics[width=1\textwidth]{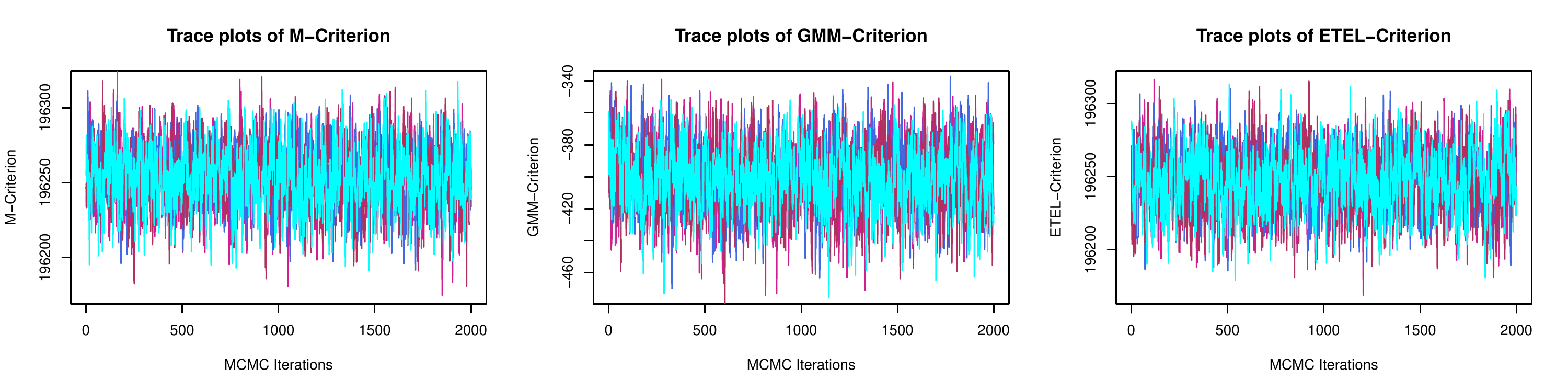}}
  \caption{Trace plots of the three criterion functions across the post-burn-in MCMC samples for the synthetic example under scenario I. Four different colors are used to highlight the trace plots of four different MCMC chains with different initializations. }
  \label{fig:graph_simulation_traceplots}
\end{figure}

\begin{figure}[htbp]
  \centerline{\includegraphics[width=1\textwidth]{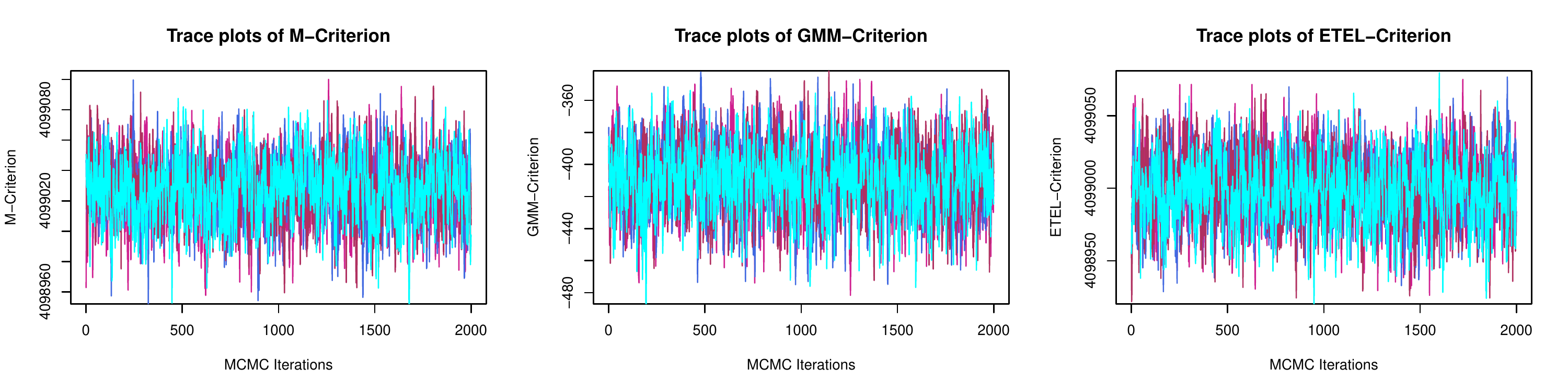}}
  \caption{Trace plots of the three criterion functions across the post-burn-in MCMC samples for the synthetic example under scenario II. Four different colors are used to highlight the trace plots of four different MCMC chains with different initializations. }
  \label{fig:SNMC_simulation_traceplots}
\end{figure}

\begin{figure}[htbp]
  \centerline{\includegraphics[width=1\textwidth]{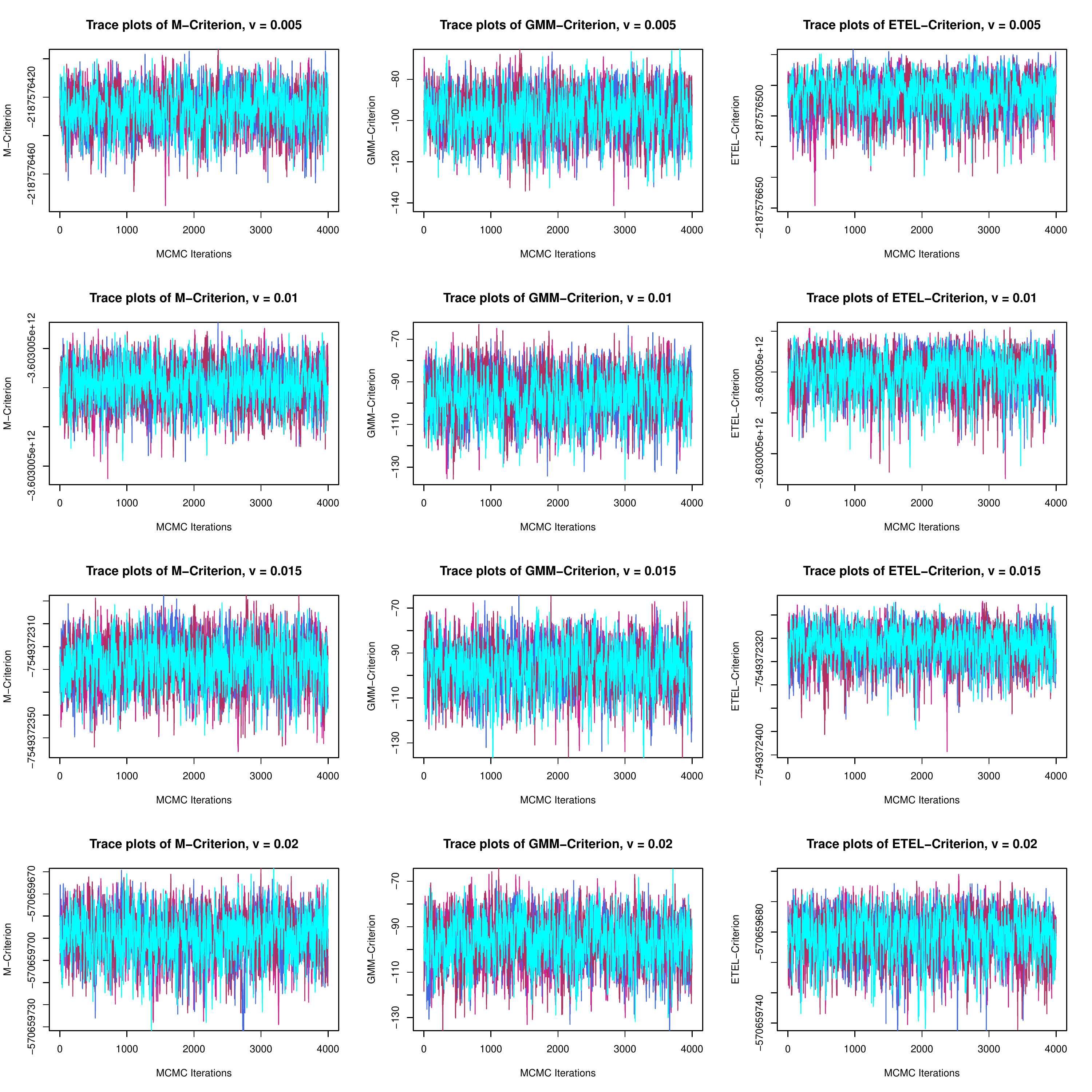}}
  \caption{Trace plots of the three criterion functions across the post-burn-in MCMC samples for ENZYMES 118 network data with different $v$. Four different colors are used to highlight the trace plots of four different MCMC chains with different initializations. }
  \label{fig:ENZYMES118_convergence}
\end{figure}

\begin{figure}[htbp]
  \centerline{\includegraphics[width=1\textwidth]{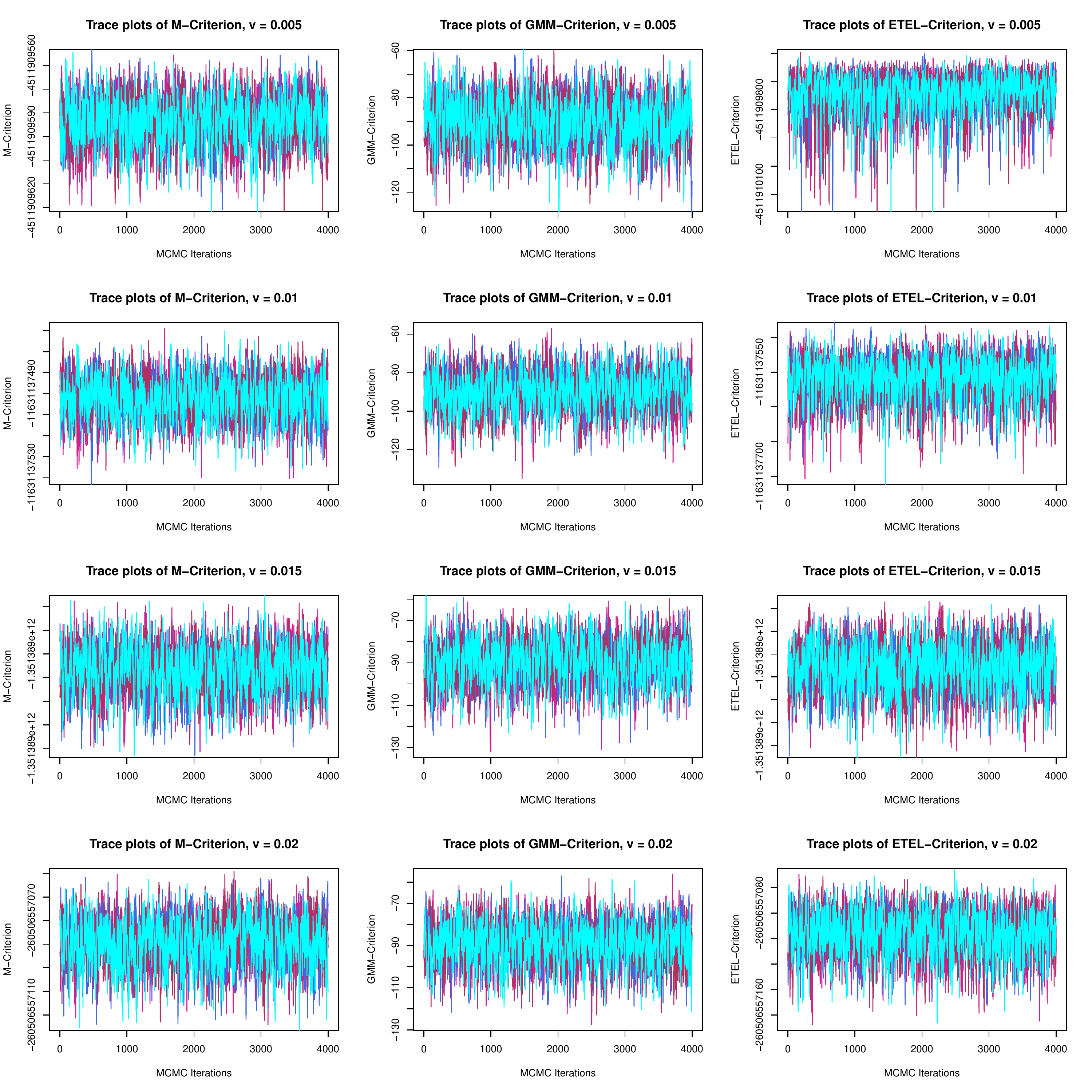}}
  \caption{Trace plots of the three criterion functions across the post-burn-in MCMC samples for ENZYMES 123 network data with different $v$. Four different colors are used to highlight the trace plots of four different MCMC chains with different initializations. }
  \label{fig:ENZYMES123_convergence}
\end{figure}

\begin{figure}[htbp]
  \centerline{\includegraphics[width=1\textwidth]{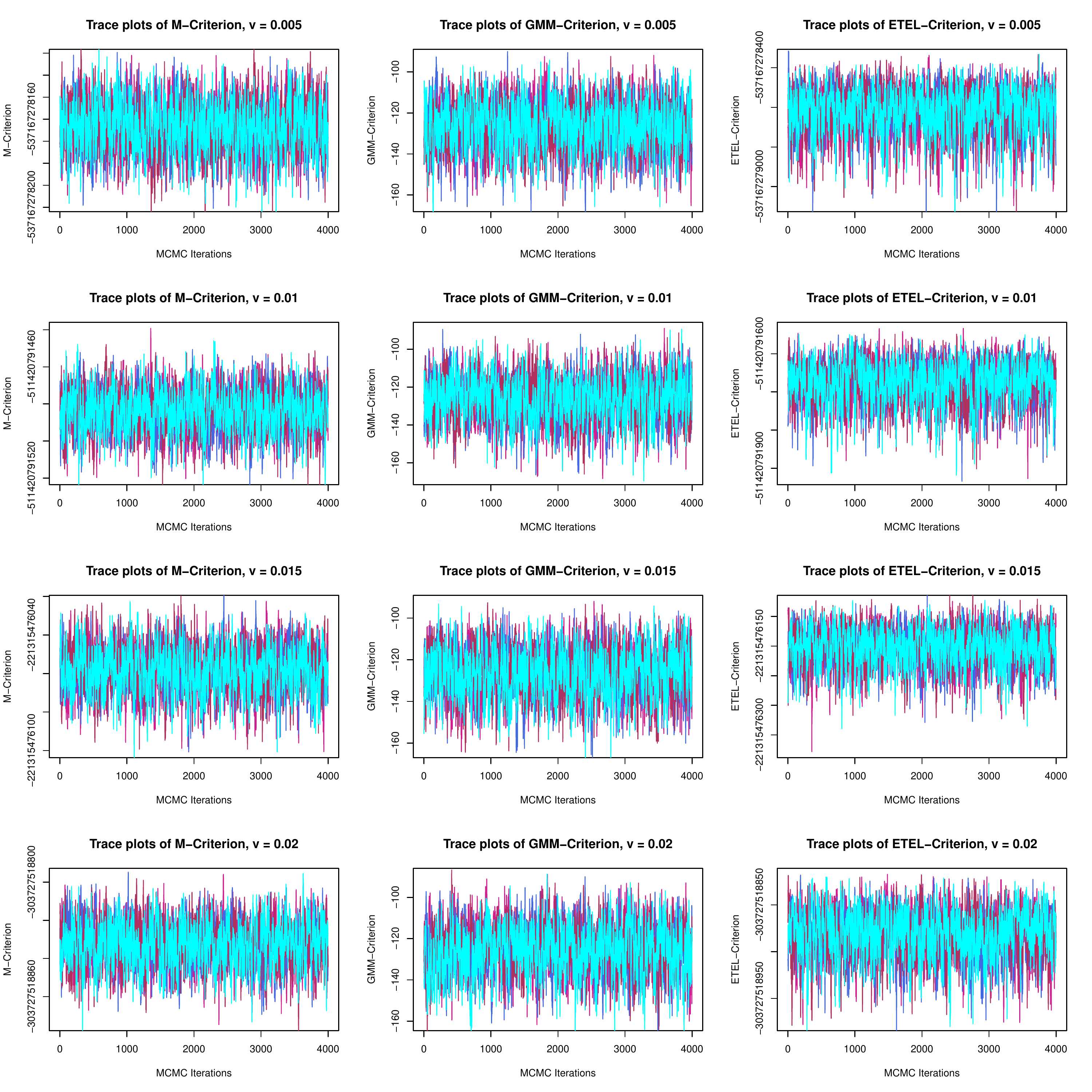}}
  \caption{Trace plots of the three criterion functions across the post-burn-in MCMC samples for ENZYMES 296 network data with different $v$. Four different colors are used to highlight the trace plots of four different MCMC chains with different initializations. }
  \label{fig:ENZYMES296_convergence}
\end{figure}

\begin{figure}[htbp]
  \centerline{\includegraphics[width=1\textwidth]{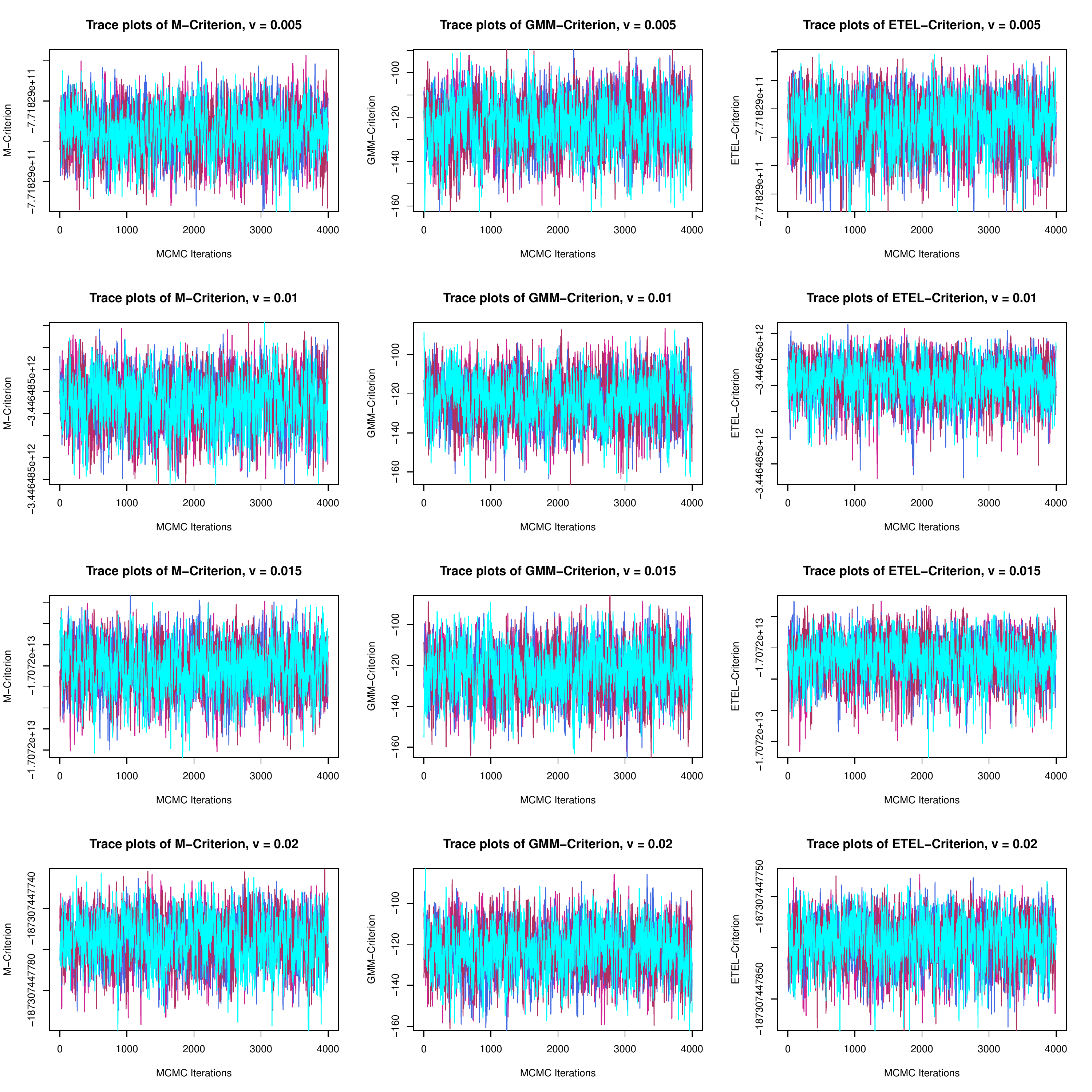}}
  \caption{Trace plots of the three criterion functions across the post-burn-in MCMC samples for ENZYMES 297 network data with different $v$. Four different colors are used to highlight the trace plots of four different MCMC chains with different initializations. }
  \label{fig:ENZYMES297_convergence}
\end{figure}

In this subsection, we provides the convergence diagnostics for the Metropolis-Hastings samplers implemented in Section \ref{sec:numerical_examples} of the manuscript. For each dataset (including the synthetic datasets and the real-world ENZYMES network datasets), the Markov chain Monte Carlo (MCMC) sampler is implemented with $1000$ burn-in iterations and $2000$ post-bur-in MCMC samples. To assess the convergence of the Markov chains, we adopt the trace plots and the Gelman-Rubin convergence diagnostics with $4$ parallel chains for each MCMC implementation. The trace plots of the MCMC implementations are provided in Figures \ref{fig:graph_simulation_traceplots}, \ref{fig:SNMC_simulation_traceplots}, \ref{fig:ENZYMES118_convergence}, \ref{fig:ENZYMES123_convergence}, \ref{fig:ENZYMES296_convergence}, \ref{fig:ENZYMES297_convergence}, showing that the Markov chains mix well in all cases. The summary statistics of the Gelman-Rubin diagnostics are provided in Tables \ref{table:Gelman_convergence_diag_simulation}, \ref{table:Gelman_convergence_diag_ENZYMES118}, \ref{table:Gelman_convergence_diag_ENZYMES123}, \ref{table:Gelman_convergence_diag_ENZYMES296}, \ref{table:Gelman_convergence_diag_ENZYMES297}. In particular, the point estimates of the potential scale reduction factors are close to $1$, and the upper limits of the $95\%$ confidence intervals are no greater than $1.1$ in all circumstances. These convergence diagnostics summaries show no signs of non-convergence of the Markov chains in the involved MCMC implementations. 

\begin{table}[htbp]
  \centering
  \caption{Gelman-Rubin convergence diagnostics for the synthetic example with the point estimates and the upper $95\%$ confidence limits of the potential scale reduction factor given by the Gelman-Rubin convergence diagnostics implemented in the \texttt{coda} package. }
  \begin{tabular}{c | c c c | c c c }
    \hline\hline
    & \multicolumn{3}{c|}{Scenario I} & \multicolumn{3}{c}{Scenario II} \\
    \hline
    Criterion & M & GMM & ETEL & M & GMM & ETEL\\
    \hline
    Point est. & 1.05 & 1.05  &  1.03 & 1.07 & 1.08 & 1.05 \\
    Upper CI   & 1.00 & 1.01  &  1.00 & 1.01 & 1.01 & 1.02 \\
    \hline\hline
  \end{tabular}%
  \label{table:Gelman_convergence_diag_simulation}
\end{table}%

\begin{table}[htbp]
  \centering
  \caption{Gelman-Rubin convergence diagnostics for the ENZYMES 118 network with the point estimates and the upper $95\%$ confidence limits of the potential scale reduction factor given by the Gelman-Rubin convergence diagnostics implemented in the \texttt{coda} package. }
  \begin{tabular}{c | c  c  c| c  c c | c  c  c| c  c c | c c c }
    \hline\hline
    $v$& \multicolumn{3}{c|}{0.005} & \multicolumn{3}{c|}{0.010} & \multicolumn{3}{c|}{0.015} & \multicolumn{3}{c}{0.020} \\
    \hline
    Criterion & M & GMM & ETEL & M & GMM & ETEL & M & GMM & ETEL & M & GMM & ETEL\\
    \hline
    Point est. & 1.05 & 1.02 & 1.01 & 1.04 & 1.04 & 1.01 & 1.03 & 1.06 & 1.03 & 1.05 & 1.03 & 1.02\\
    Upper CI   & 1.09 & 1.04 & 1.02 & 1.07 & 1.07 & 1.02 & 1.05 & 1.10 & 1.04 & 1.09 & 1.05 & 1.03\\
    \hline\hline
  \end{tabular}%
  \label{table:Gelman_convergence_diag_ENZYMES118}
\end{table}%

\begin{table}[htbp]
  \centering
  \caption{Gelman-Rubin convergence diagnostics for the ENZYMES 123 network with the point estimates and the upper $95\%$ confidence limits of the potential scale reduction factor given by the Gelman-Rubin convergence diagnostics implemented in the \texttt{coda} package. }
  \begin{tabular}{c | c  c  c| c  c c | c  c  c| c  c c | c c c }
    \hline\hline
    $v$& \multicolumn{3}{c|}{0.005} & \multicolumn{3}{c|}{0.010} & \multicolumn{3}{c|}{0.015} & \multicolumn{3}{c}{0.020} \\
    \hline
    Criterion & M & GMM & ETEL & M & GMM & ETEL & M & GMM & ETEL & M & GMM & ETEL\\
    \hline
    Point est. & 1.02 & 1.02 & 1.01 & 1.05 & 1.04 & 1.01 & 1.03 & 1.03 & 1.01 & 1.02 & 1.05 & 1.02\\
    Upper CI   & 1.03 & 1.03 & 1.02 & 1.08 & 1.07 & 1.02 & 1.06 & 1.05 & 1.02 & 1.03 & 1.08 & 1.03\\
    \hline\hline
  \end{tabular}%
  \label{table:Gelman_convergence_diag_ENZYMES123}
\end{table}%

\begin{table}[htbp]
  \centering
  \caption{Gelman-Rubin convergence diagnostics for the ENZYMES 296 network with the point estimates and the upper $95\%$ confidence limits of the potential scale reduction factor given by the Gelman-Rubin convergence diagnostics implemented in the \texttt{coda} package. }
  \begin{tabular}{c | c  c  c| c  c c | c  c  c| c  c c | c c c }
    \hline\hline
    $v$& \multicolumn{3}{c|}{0.005} & \multicolumn{3}{c|}{0.010} & \multicolumn{3}{c|}{0.015} & \multicolumn{3}{c}{0.020} \\
    \hline
    Criterion & M & GMM & ETEL & M & GMM & ETEL & M & GMM & ETEL & M & GMM & ETEL\\
    \hline
    Point est. & 1.03 & 1.03 & 1.02 & 1.05 & 1.03 & 1.01 & 1.04 & 1.02 & 1.03 & 1.02 & 1.04 & 1.02\\
    Upper CI   & 1.04 & 1.06 & 1.03 & 1.08 & 1.06 & 1.02 & 1.07 & 1.03 & 1.05 & 1.03 & 1.06 & 1.03\\
    \hline\hline
  \end{tabular}%
  \label{table:Gelman_convergence_diag_ENZYMES296}
\end{table}%

\begin{table}[htbp]
  \centering
  \caption{Gelman-Rubin convergence diagnostics for the ENZYMES 297 network with the point estimates and the upper $95\%$ confidence limits of the potential scale reduction factor given by the Gelman-Rubin convergence diagnostics implemented in the \texttt{coda} package. }
  \begin{tabular}{c | c  c  c| c  c c | c  c  c| c  c c | c c c }
    \hline\hline
    $v$& \multicolumn{3}{c|}{0.005} & \multicolumn{3}{c|}{0.010} & \multicolumn{3}{c|}{0.015} & \multicolumn{3}{c}{0.020} \\
    \hline
    Criterion & M & GMM & ETEL & M & GMM & ETEL & M & GMM & ETEL & M & GMM & ETEL\\
    \hline
    Point est. & 1.04 & 1.05 & 1.02 & 1.03 & 1.02 & 1.02 & 1.02 & 1.04 & 1.02 & 1.02 & 1.02 & 1.02\\
    Upper CI   & 1.06 & 1.08 & 1.04 & 1.05 & 1.04 & 1.03 & 1.04 & 1.07 & 1.04 & 1.04 & 1.03 & 1.03\\
    \hline\hline
  \end{tabular}%
  \label{table:Gelman_convergence_diag_ENZYMES297}
\end{table}%

%% if your bibliography is in bibtex format, uncomment commands:

\clearpage
\bibliographystyle{apalike} % Style BST file
\bibliography{reference1, reference2}       % Bibliography file (usually '*.bib')

\end{document}